# The Physicalization of Metamathematics and Its Implications for the Foundations of Mathematics

Stephen Wolfram*

*Both metamathematics and physics are posited to emerge from samplings by observers of the unique ruliad structure that corresponds to the entangled limit of all possible computations. The possibility of higher-level mathematics accessible to humans is posited to be the analog for mathematical observers of the perception of physical space for physical observers. A physicalized analysis is given of the bulk limit of traditional axiomatic approaches to the foundations of mathematics, together with explicit empirical metamathematics of some examples of formalized mathematics. General physicalized laws of mathematics are discussed, associated with concepts such as metamathematical motion, inevitable dualities, proof topology and metamathematical singularities. It is argued that mathematics as currently practiced can be viewed as derived from the ruliad in a direct Platonic fashion analogous to our experience of the physical world, and that axiomatic formulation, while often convenient, does not capture the ultimate character of mathematics. Among the implications of this view is that only certain collections of axioms may be consistent with inevitable features of human mathematical observers. A discussion is included of historical and philosophical connections, as well as of foundational implications for the future of mathematics.*

## 1 | Mathematics and Physics Have the Same Foundations

One of the many surprising (and to me, unexpected) implications of our Physics Project is its suggestion of a very deep correspondence between the foundations of physics and mathematics. We might have imagined that physics would have certain laws, and mathematics would have certain theories, and that while they might be historically related, there wouldn't be any fundamental formal correspondence between them.

But what our Physics Project suggests is that underneath everything we physically experience there is a single very general abstract structure—that we call the ruliad—and that our physical laws arise in an inexorable way from the particular samples we take of this structure. We can think of the ruliad as the entangled limit of all possible computations—or in effect a representation of all possible formal processes. And this then leads us to the idea that perhaps the ruliad might underlie not only physics but also mathematics—and that everything in mathematics, like everything in physics, might just be the result of sampling the ruliad.

Of course, mathematics as it's normally practiced doesn't look the same as physics. But the idea is that they can both be seen as views of the same underlying structure. What makes them different is that physical and mathematical observers sample this structure in somewhat different ways. But since in the end both kinds of observers are associated with human experience they inevitably have certain core characteristics in common. And the result is that there should be "fundamental laws of mathematics" that in some sense mirror the perceived laws of physics that we derive from our physical observation of the ruliad.

So what might those fundamental laws of mathematics be like? And how might they inform our conception of the foundations of mathematics, and our view of what mathematics really is?

The most obvious manifestation of the mathematics that we humans have developed over the course of many centuries is the few million mathematical theorems that have been published in the literature of mathematics. But what can be said in generality about this thing we call mathematics? Is there some notion of what mathematics is like "in bulk"? And what might we be able to say, for example, about the structure of mathematics in the limit of infinite future development?

When we do physics, the traditional approach has been to start from our basic sensory experience of the physical world, and of concepts like space, time and motion—and then to try to formalize our descriptions of these things, and build on these formalizations. And in its early development—for example by Euclid—mathematics took the same basic approach. But beginning a little more than a century ago there emerged the idea that one could build mathematics purely from formal axioms, without necessarily any reference to what is accessible to sensory experience.

And in a way our Physics Project begins from a similar place. Because at the outset it just considers purely abstract structures and abstract rules—typically described in terms of hypergraph rewriting—and then tries to deduce their consequences. Many of these consequences are incredibly complicated, and full of computational irreducibility. But the remarkable discovery is that when sampled by observers with certain general characteristics that make them like us, the behavior that emerges must generically have regularities that we can recognize, and in fact must follow exactly known core laws of physics.

And already this begins to suggest a new perspective to apply to the foundations of mathematics. But there's another piece, and that's the idea of the ruliad. We might have supposed that our universe is based on some particular chosen underlying rule, like an axiom system we might choose in mathematics. But the concept of the ruliad is in effect to represent the entangled result of "running all possible rules". And the key point is then that it turns out that an "observer like us" sampling the ruliad must perceive behavior that corresponds to known laws of physics. In other words, without "making any choice" it's inevitable—given what we're like as observers—that our "experience of the ruliad" will show fundamental laws of physics.

But now we can make a bridge to mathematics. Because in embodying all possible computational processes the ruliad also necessarily embodies the consequences of all possible axiom systems. As humans doing physics we're effectively taking a certain sampling of the ruliad. And we realize that as humans doing mathematics we're also doing essentially the same kind of thing.

But will we see "general laws of mathematics" in the same kind of way that we see "general laws of physics"? It depends on what we're like as "mathematical observers". In physics, there turn out to be general laws—and concepts like space and motion—that we humans can assimilate. And in the abstract it might not be that anything similar would be true in mathematics. But it seems as if the thing mathematicians typically call mathematics is something for which it is—and where (usually in the end leveraging our experience of physics) it's possible to successfully carve out a sampling of the ruliad that's again one we humans can assimilate.

When we think about physics we have the idea that there's an actual physical reality that exists—and that we experience physics within this. But in the formal axiomatic view of mathematics, things are different. There's no obvious "underlying reality" there; instead there's just a certain choice we make of axiom system. But now, with the concept of the ruliad, the story is different. Because now we have the idea that "deep underneath" both physics and mathematics there's the same thing: the ruliad. And that means that insofar as physics is "grounded in reality", so also must mathematics be.

When most working mathematicians do mathematics it seems to be typical for them to reason as if the constructs they're dealing with (whether they be numbers or sets or whatever) are "real things". But usually there's a concept that in principle one could "drill down" and formalize everything in terms of some axiom system. And indeed if one wants to get a global view of mathematics and its structure as it is today, it seems as if the best approach is to work from the formalization that's been done with axiom systems.

In starting from the ruliad and the ideas of our Physics Project we're in effect positing a certain "theory of mathematics". And to validate this theory we need to study the "phenomena of mathematics". And, yes, we could do this in effect by directly "reading the whole literature of mathematics". But it's more efficient to start from what's in a sense the "current prevailing underlying theory of mathematics" and to begin by building on the methods of formalized mathematics and axiom systems.

Over the past century a certain amount of metamathematics has been done by looking at the general properties of these methods. But most often when the methods are systematically used today, it's to set up some particular mathematical derivation, normally with the aid of a computer. But here what we want to do is think about what happens if the methods are used

"in bulk". Underneath there may be all sorts of specific detailed formal derivations being done. But somehow what emerges from this is something higher level, something "more human"—and ultimately something that corresponds to our experience of pure mathematics.

How might this work? We can get an idea from an analogy in physics. Imagine we have a gas. Underneath, it consists of zillions of molecules bouncing around in detailed and complicated patterns. But most of our "human" experience of the gas is at a much more coarse-grained level—where we perceive not the detailed motions of individual molecules, but instead continuum fluid mechanics.

And so it is, I think, with mathematics. All those detailed formal derivations—for example of the kind automated theorem proving might do—are like molecular dynamics. But most of our "human experience of mathematics"—where we talk about concepts like integers or morphisms—is like fluid dynamics. The molecular dynamics is what builds up the fluid, but for most questions of "human interest" it's possible to "reason at the fluid dynamics level", without dropping down to molecular dynamics.

It's certainly not obvious that this would be possible. It could be that one might start off describing things at a "fluid dynamics" level—say in the case of an actual fluid talking about the motion of vortices—but that everything would quickly get "shredded", and that there'd soon be nothing like a vortex to be seen, only elaborate patterns of detailed microscopic molecular motions. And similarly in mathematics one might imagine that one would be able to prove theorems in terms of things like real numbers but actually find that everything gets "shredded" to the point where one has to start talking about elaborate issues of mathematical logic and different possible axiomatic foundations.

But in physics we effectively have the Second Law of thermodynamics—which we now understand in terms of computational irreducibility—that tells us that there's a robust sense in which the microscopic details are systematically "washed out" so that things like fluid dynamics "work". Just sometimes—like in studying Brownian motion, or hypersonic flow—the molecular dynamics level still "shines through". But for most "human purposes" we can describe fluids just using ordinary fluid dynamics.

So what's the analog of this in mathematics? Presumably it's that there's some kind of "general law of mathematics" that explains why one can so often do mathematics "purely in the large". Just like in fluid mechanics there can be "corner-case" questions that probe down to the "molecular scale"—and indeed that's where we can expect to see things like undecidability, as a rough analog of situations where we end up tracing the potentially infinite paths of single molecules rather than just looking at "overall fluid effects". But somehow in most cases there's some much stronger phenomenon at work—that effectively aggregates low-level details to allow the kind of "bulk description" that ends up being the essence of what we normally in practice call mathematics.

But is such a phenomenon something formally inevitable, or does it somehow depend on us humans "being in the loop"? In the case of the Second Law it's crucial that we only get to track coarse-grained features of a gas—as we humans with our current technology typically do. Because if instead we watched and decoded what every individual molecule does, we wouldn't end up identifying anything like the usual bulk "Second-Law" behavior. In other words, the emergence of the Second Law is in effect a direct consequence of the fact that it's us humans—with our limitations on measurement and computation—who are observing the gas.

So is something similar happening with mathematics? At the underlying "molecular level" there's a lot going on. But the way we humans think about things, we're effectively taking just particular kinds of samples. And those samples turn out to give us "general laws of mathematics" that give us our usual experience of "human-level mathematics".

To ultimately ground this we have to go down to the fully abstract level of the ruliad, but we'll already see many core effects by looking at mathematics essentially just at a traditional "axiomatic level", albeit "in bulk".

The full story—and the full correspondence between physics and mathematics—requires in a sense "going below" the level at which we have recognizable formal axiomatic mathematical structures; it requires going to a level at which we're just talking about making everything out of completely abstract elements, which in physics we might interpret as "atoms of space" and in mathematics as some kind of "symbolic raw material" below variables and operators and everything else familiar in traditional axiomatic mathematics.

The deep correspondence we're describing between physics and mathematics might make one wonder to what extent the methods we use in physics can be applied to mathematics, and vice versa. In axiomatic mathematics the emphasis tends to be on looking at particular theorems and seeing how they can be knitted together with proofs. And one could certainly imagine an analogous "axiomatic physics" in which one does particular experiments, then sees how they can "deductively" be knitted together. But our impression that there's an "actual reality" to physics makes us seek broader laws. And the correspondence between physics and mathematics implied by the ruliad now suggests that we should be doing this in mathematics as well.

What will we find? Some of it in essence just confirms impressions that working pure mathematicians already have. But it provides a definite framework for understanding these impressions and for seeing what their limits may be. It also lets us address questions like why undecidability is so comparatively rare in practical pure mathematics, and why it is so common to discover remarkable correspondences between apparently quite different areas of mathematics. And beyond that, it suggests a host of new questions and approaches both to mathematics and metamathematics—that help frame the foundations of the remarkable intellectual edifice that we call mathematics.

# 2 | The Underlying Structure of Mathematics and Physics

If we "drill down" to what we've called above the "molecular level" of mathematics, what will we find there? There are many technical details (some of which we'll discuss later) about the historical conventions of mathematics and its presentation. But in broad outline we can think of there as being a kind of "gas" of "mathematical statements"—like 1+1=2 or *x+y=y+x*—represented in some specified symbolic language. (And, yes, Wolfram Language provides a well-developed example of what that language can be like.)

But how does the "gas of statements" behave? The essential point is that new statements are derived from existing ones by "interactions" that implement laws of inference (like that *q* can be derived from the statement *p* and the statement "*p* implies *q*"). And if we trace the paths by which one statement can be derived from others, these correspond to proofs. And the whole graph of all these derivations is then a representation of the possible historical development of mathematics—with slices through this graph corresponding to the sets of statements reached at a given stage.

By talking about things like a "gas of statements" we're making this sound a bit like physics. But while in physics a gas consists of actual, physical molecules, in mathematics our statements are just abstract things. But this is where the discoveries of our Physics Project start to be important. Because in our project we're "drilling down" beneath for example the usual notions of space and time to an "ultimate machine code" for the physical universe. And we can think of that ultimate machine code as operating on things that are in effect just abstract constructs—very much like in mathematics.

In particular, we imagine that space and everything in it is made up of a giant network (hypergraph) of "atoms of space"—with each "atom of space" just being an abstract element that has certain relations with other elements. The evolution of the universe in time then corresponds to the application of computational rules that (much like laws of inference) take abstract relations and yield new relations—thereby progressively updating the network that represents space and everything in it.

But while the individual rules may be very simple, the whole detailed pattern of behavior to which they lead is normally very complicated—and typically shows computational irreducibility, so that there's no way to systematically find its outcome except in effect by explicitly tracing each step. But despite all this underlying complexity it turns out—much like in the case of an ordinary gas—that at a coarse-grained level there are much simpler ("bulk") laws of behavior that one can identify. And the remarkable thing is that these turn out to be exactly general relativity and quantum mechanics (which, yes, end up being the same theory when looked at in terms of an appropriate generalization of the notion of space).

But down at the lowest level, is there some specific computational rule that's "running the universe"? I don't think so. Instead, I think that in effect all possible rules are always being applied. And the result is the ruliad: the entangled structure associated with performing all possible computations.

But what then gives us our experience of the universe and of physics? Inevitably we are observers embedded within the ruliad, sampling only certain features of it. But what features we sample are determined by the characteristics of us as observers. And what seem to be critical to have "observers like us" are basically two characteristics. First, that we are computationally bounded. And second, that we somehow persistently maintain our coherence—in the sense that we can consistently identify what constitutes "us" even though the detailed atoms of space involved are continually changing.

But we can think of different "observers like us" as taking different specific samples, corresponding to different reference frames in rulial space, or just different positions in rulial space. These different observers may describe the universe as evolving according to different specific underlying rules. But the crucial point is that the general structure of the ruliad implies that so long as the observers are "like us", it's inevitable that their perception of the universe will be that it follows things like general relativity and quantum mechanics.

It's very much like what happens with a gas of molecules: to an "observer like us" there are the same gas laws and the same laws of fluid dynamics essentially independent of the detailed structure of the individual molecules.

So what does all this mean for mathematics? The crucial and at first surprising point is that the ideas we're describing in physics can in effect immediately be carried over to mathematics. And the key is that the ruliad represents not only all physics, but also all mathematics—and it shows that these are not just related, but in some sense fundamentally the same.

In the traditional formulation of axiomatic mathematics, one talks about deriving results from particular axiom systems—say [Peano Arithmetic](), or [ZFC set theory](), or the axioms of [Euclidean geometry](). But the ruliad in effect represents the entangled consequences not just of specific axiom systems but of all possible axiom systems (as well as all possible laws of inference).

But from this structure that in a sense corresponds to all possible mathematics, how do we pick out any particular mathematics that we're interested in? The answer is that just as we are limited observers of the physical universe, so we are also limited observers of the "mathematical universe".

But what are we like as "mathematical observers"? As I'll argue in more detail later, we inherit our core characteristics from those we exhibit as "physical observers". And that means that when we "do mathematics" we're effectively sampling the ruliad in much the same way as when we "do physics".

We can operate in different rulial reference frames, or at different locations in rulial space, and these will correspond to picking out different underlying "rules of mathematics", or essentially using different axiom systems. But now we can make use of the correspondence with physics to say that we can also expect there to be certain "overall laws of mathematics" that are the result of general features of the ruliad as perceived by observers like us.

And indeed we can expect that in some formal sense these overall laws will have exactly the same structure as those in physics—so that in effect in mathematics we'll have something like the notion of space that we have in physics, as well as formal analogs of things like general relativity and quantum mechanics.

What does this mean? It implies that—just as it's possible to have coherent "higher-level descriptions" in physics that don't just operate down at the level of atoms of space, so also this should be possible in mathematics. And this in a sense is why we can expect to consistently do what I described above as "human-level mathematics", without usually having to drop down to the "molecular level" of specific axiomatic structures (or below).

Say we're talking about the Pythagorean theorem. Given some particular detailed axiom system for mathematics we can imagine using it to build up a precise—if potentially very long and pedantic—representation of the theorem. But let's say we change some detail of our axioms, say associated with the way they talk about sets, or real numbers. We'll almost certainly still be able to build up something we consider to be "the Pythagorean theorem"—even though the details of the representation will be different.

In other words, this thing that we as humans would call "the Pythagorean theorem" is not just a single point in the ruliad, but a whole cloud of points. And now the question is: what happens if we try to derive other results from the Pythagorean theorem? It might be that each particular representation of the theorem—corresponding to each point in the cloud—would lead to quite different results. But it could also be that essentially the whole cloud would coherently lead to the same results.

And the claim from the correspondence with physics is that there should be "general laws of mathematics" that apply to "observers like us" and that ensure that there'll be coherence between all the different specific representations associated with the cloud that we identify as "the Pythagorean theorem".

In physics it could have been that we'd always have to separately say what happens to every atom of space. But we know that there's a coherent higher-level description of space—in which for example we can just imagine that objects can move while somehow maintaining their identity. And we can now expect that it's the same kind of thing in mathematics: that just as there's a coherent notion of space in physics where things can for example move without being "shredded", so also this will happen in mathematics. And this is why it's possible to do "higher-level mathematics" without always dropping down to the lowest level of axiomatic derivations.

It's worth pointing out that even in physical space a concept like "pure motion" in which objects can move while maintaining their identity doesn't always work. For example, close to a spacetime singularity, one can expect to eventually be forced to see through to the discrete structure of space—and for any "object" to inevitably be "shredded". But most of the time it's possible for observers like us to maintain the idea that there are coherent large-scale features whose behavior we can study using "bulk" laws of physics.

And we can expect the same kind of thing to happen with mathematics. Later on, we'll discuss more specific correspondences between phenomena in physics and mathematics—and we'll see the effects of things like general relativity and quantum mechanics in mathematics, or, more precisely, in metamathematics.

But for now, the key point is that we can think of mathematics as somehow being made of exactly the same stuff as physics: they're both just features of the ruliad, as sampled by observers like us. And in what follows we'll see the great power that arises from using this to combine the achievements and intuitions of physics and mathematics—and how this lets us think about new "general laws of mathematics", and view the ultimate foundations of mathematics in a different light.

## 3 | The Metamodeling of Axiomatic Mathematics

Consider all the mathematical statements that have appeared in mathematical books and papers. We can view these in some sense as the "observed phenomena" of (human) mathematics. And if we're going to make a "general theory of mathematics" a first step is to do something like we'd typically do in natural science, and try to "drill down" to find a uniform underlying model—or at least representation—for all of them.

At the outset, it might not be clear what sort of representation could possibly capture all those different mathematical statements. But what's emerged over the past century or so—with particular clarity in Mathematica and the Wolfram Language—is that there is in fact a rather simple and general representation that works remarkably well: a representation in which everything is a symbolic expression.

One can view a symbolic expression such as **f[g[x][y, h[z]], w]** as a hierarchical or tree structure, in which at every level some particular "head" (like **f**) is "applied to" one or more arguments. Often in practice one deals with expressions in which the heads have "known meanings"—as in **Times**[**Plus**[2, 3], 4] in Wolfram Language. And with this kind of setup symbolic expressions are reminiscent of human natural language, with the heads basically corresponding to "known words" in the language.

And presumably it's this familiarity from human natural language that's caused "human natural mathematics" to develop in a way that can so readily be represented by symbolic expressions.

But in typical mathematics there's an important wrinkle. One often wants to make statements not just about particular things but about whole classes of things. And it's common to then just declare that some of the "symbols" (like, say, $x$) that appear in an expression are "variables", while others (like, say, **Plus**) are not. But in our effort to capture the essence of mathematics as uniformly as possible it seems much better to burn the idea of an object representing a whole class of things right into the structure of the symbolic expression.

And indeed this is a core idea in the Wolfram Language, where something like $x$ or **f** is just a "symbol that stands for itself", while $x\_$ is a pattern (named $x$) that can stand for anything. (More precisely, _ on its own is what stands for "anything", and $x\_$—which can also be written $x:\_$—just says that whatever _ stands for in a particular instance will be called $x$.)

Then with this notation an example of a "mathematical statement" might be:

$x\_ \circ y\_ = (y\_ \circ x\_) \circ y\_$

In more explicit form we could write this as **Equal**[f[x\_, y\_], f[f[y\_, x\_],y\_]]—where **Equal** (=) has the "known meaning" of representing equality. But what can we do with this statement? At a "mathematical level" the statement asserts that $x\_ \circ y\_$ and $(y\_ \circ x\_) \circ y\_$ should be considered equivalent. But thinking in terms of symbolic expressions there's now a more explicit, lower-level, "structural" interpretation: that any expression whose structure matches $x\_ \circ y\_$ can equivalently be replaced by $(y\_ \circ x\_) \circ y\_$ (or, in Wolfram Language notation, just $(y \circ x) \circ y$) and vice versa. We can indicate this interpretation using the notation

$$x\_ \circ y\_ \longleftrightarrow (y\_ \circ x\_) \circ y\_$$

which can be viewed as a shorthand for the pair of Wolfram Language rules:

$$x\_ \circ y\_ \to (y \circ x) \circ y, \quad (y\_ \circ x\_) \circ y\_ \to x \circ y$$

OK, so let's say we have the expression $(a \circ b) \circ a$. Now we can just apply the rules defined by our statement. Here's what happens if we do this just once in all possible ways:

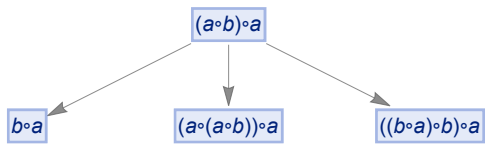

And here we see, for example, that $(a \circ b) \circ a$ can be transformed to $b \circ a$. Continuing this we build up a whole multiway graph. After just one more step we get:

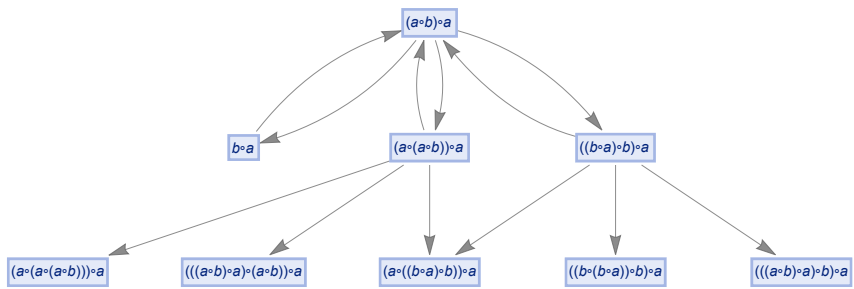

Continuing for a few more steps we then get

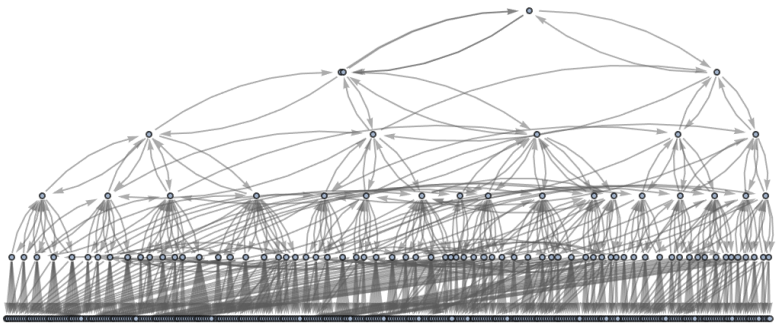

or in a different rendering:

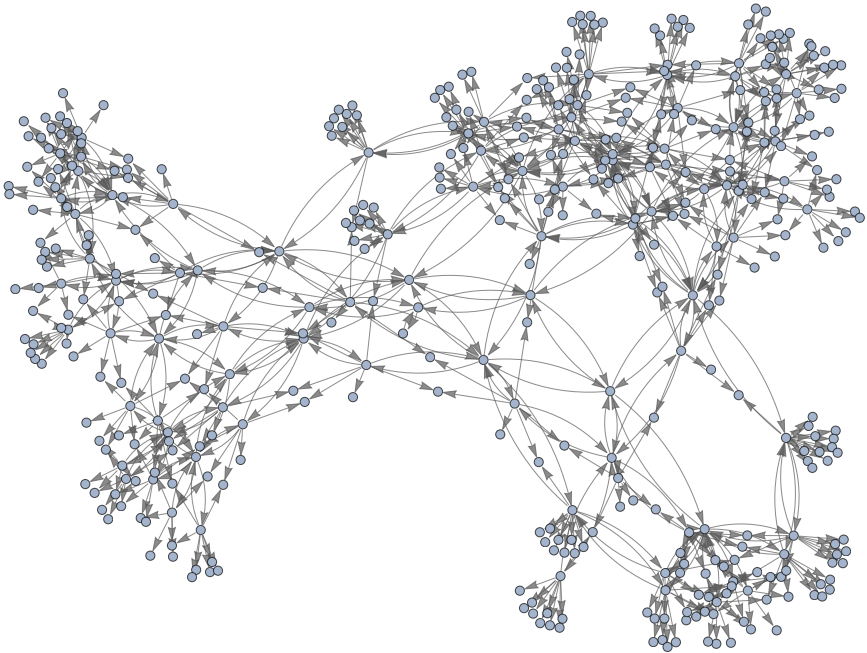

But what does this graph mean? Essentially it gives us a map of equivalences between expressions—with any pair of expressions that are connected being equivalent. So, for example, it turns out that the expressions (*a*∘((*b*∘*a*)∘(*a*∘*b*)))∘*a* and *b*∘*a* are equivalent, and we can "prove this" by exhibiting a path between them in the graph:

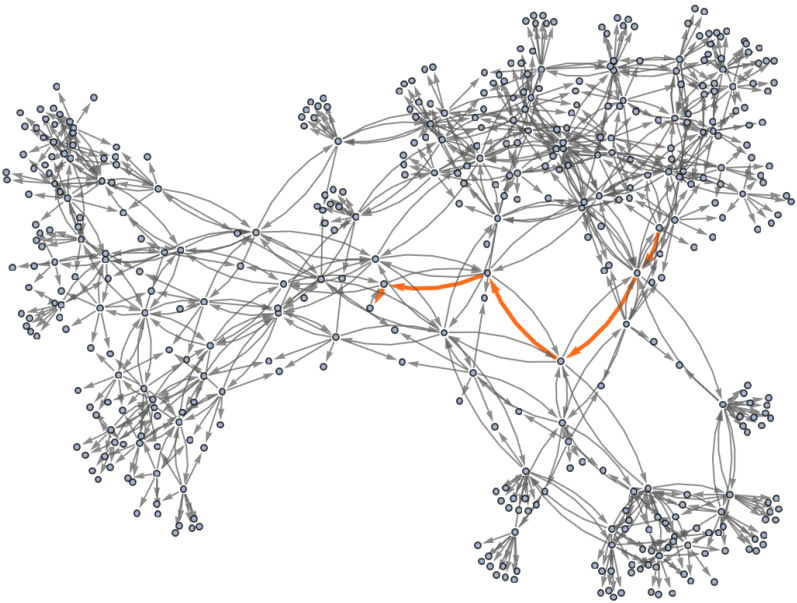

The steps on the path can then be viewed as steps in the proof, where here at each step we've indicated where the transformation in the expression took place:

| |
|---|
| (a ∘ ((b ∘ a) ∘ (a ∘ b))) ∘ a |
| (a ∘ (((a ∘ b) ∘ a) ∘ (a ∘ b))) ∘ a |
| (((a ∘ b) ∘ a) ∘ (a ∘ b)) ∘ a |
| (a ∘ (a ∘ b)) ∘ a |
| (a ∘ b) ∘ a |
| b ∘ a |

In mathematical terms, we can then say that starting from the "axiom" $x \circ y = (y \circ x) \circ y$ we were able to prove a certain equivalence theorem between two expressions. We gave a particular proof. But there are others, for example the "less efficient" 35-step one

| |
|---|
| (a∘((b∘a)∘(a∘b)))∘a |
| ((b∘a)∘(a∘b))∘a |
| (((a∘b)∘(b∘a))∘(a∘b))∘a |
| (((a∘b)∘((a∘b)∘a))∘(a∘b))∘a |
| ((((b∘a)∘b)∘((a∘b)∘a))∘(a∘b))∘a |
| ((((b∘a)∘b)∘((a∘b)∘a))∘((b∘a)∘b))∘a |
| ((((b∘a)∘b)∘(((b∘a)∘b)∘a))∘((b∘a)∘b))∘a |
| ((((b∘a)∘b)∘a)∘((b∘a)∘b))∘a |
| (((((a∘b)∘a)∘b)∘a)∘((b∘a)∘b))∘a |
| (((((a∘b)∘a)∘b)∘a)∘(((a∘b)∘a)∘b))∘a |
| ((((b∘a)∘b)∘a)∘(((a∘b)∘a)∘b))∘a |
| (((a∘b)∘a)∘(((a∘b)∘a)∘b))∘a |
| (((a∘b)∘a)∘((b∘a)∘b))∘a |
| (((a∘(a∘b))∘a)∘((b∘a)∘b))∘a |
| (((a∘((b∘a)∘b))∘a)∘((b∘a)∘b))∘a |
| (((a∘((b∘a)∘b))∘a)∘(a∘b))∘a |
| ((((b∘a)∘b)∘a)∘(a∘b))∘a |
| ((((b∘(b∘a))∘b)∘a)∘(a∘b))∘a |
| ((((b∘(b∘a))∘b)∘a)∘((b∘a)∘b))∘a |
| ((((b∘(b∘a))∘b)∘a)∘((b∘(b∘a))∘b))∘a |
| (a∘((b∘(b∘a))∘b))∘a |
| ((b∘(b∘a))∘b)∘a |
| ((((b∘a)∘b)∘(b∘a))∘b)∘a |
| (((((a∘b)∘a)∘b)∘(b∘a))∘b)∘a |
| (((((a∘b)∘a)∘b)∘((a∘b)∘a))∘b)∘a |
| ((b∘((a∘b)∘a))∘b)∘a |
| (((a∘b)∘a)∘b)∘a |
| ((((b∘a)∘b)∘a)∘b)∘a |
| (((a∘((b∘a)∘b))∘a)∘b)∘a |
| (((a∘(a∘b))∘a)∘b)∘a |
| (a∘(((a∘(a∘b))∘a)∘b))∘a |
| (a∘(((a∘b)∘a)∘b))∘a |
| (a∘((b∘a)∘b))∘a |
| ((b∘a)∘b)∘a |
| (a∘b)∘a |
| b∘a |

corresponding to the path:

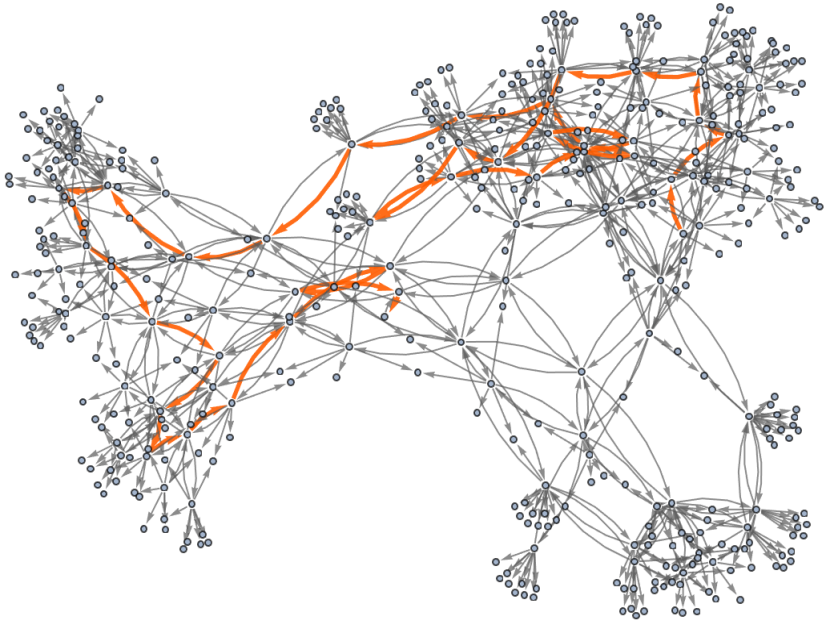

For our later purposes it's worth talking in a little bit more detail here about how the steps in these proofs actually proceed. Consider the expression:

$(a \circ ((b \circ a) \circ (a \circ b))) \circ a$

We can think of this as a tree:

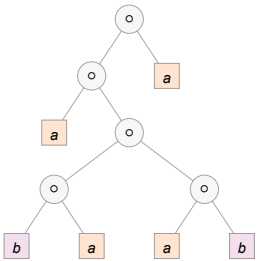

Our axiom can then be represented as:

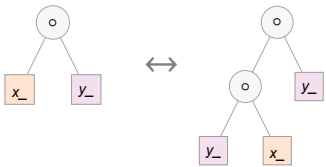

In terms of trees, our first proof becomes

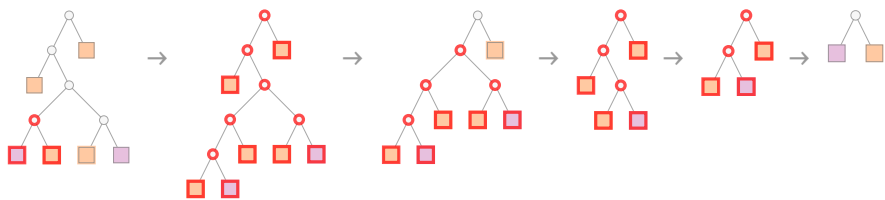

where we're indicating at each step which piece of tree gets "substituted for" using the axiom.

What we've done so far is to generate a multiway graph for a certain number of steps, and then to see if we can find a "proof path" in it for some particular statement. But what if we are given a statement, and asked whether it can be proved within the specified axiom system? In effect this asks whether if we make a sufficiently large multiway graph we can find a path of any length that corresponds to the statement.

If our system was computationally reducible we could expect always to be able to find a finite answer to this question. But in general—with the Principle of Computational Equivalence and the ubiquitous presence of computational irreducibility—it'll be common that there is no fundamentally better way to determine whether a path exists than effectively to try explicitly generating it. If we knew, for example, that the intermediate expressions generated always remained of bounded length, then this would still be a bounded problem. But in general the expressions can grow to any size—with the result that there is no general upper bound on the length of path necessary to prove even a statement about equivalence between small expressions.

For example, for the axiom we are using here, we can look at statements of the form (a∘b)∘a=*expr*. Then this shows how many expressions *expr* of what sizes have shortest proofs of (a∘b)∘a=*expr* with progressively greater lengths:

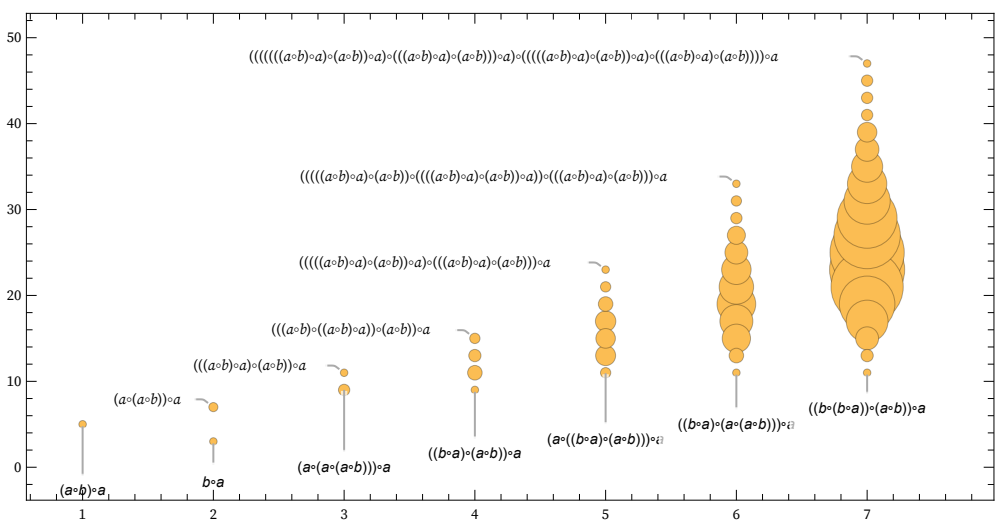

And for example if we look at the statement

$$(a \circ b) \circ a = (a \circ ((b \circ (b \circ a)) \circ (a \circ b))) \circ a$$

its shortest proof is

| |
|---|
| $(a \circ b) \circ a$ |
| $(a \circ (a \circ b)) \circ a$ |
| $(((a \circ b) \circ a) \circ (a \circ b)) \circ a$ |
| $(((a \circ b) \circ ((a \circ b) \circ a)) \circ (a \circ b)) \circ a$ |
| $(a \circ (((a \circ b) \circ ((a \circ b) \circ a)) \circ (a \circ b))) \circ a$ |
| $(a \circ ((((b \circ a) \circ b) \circ ((a \circ b) \circ a)) \circ (a \circ b))) \circ a$ |
| $(a \circ ((((b \circ a) \circ b) \circ (b \circ a)) \circ (a \circ b))) \circ a$ |
| $(a \circ ((b \circ (b \circ a)) \circ (a \circ b))) \circ a$ |

where, as is often the case, there are intermediate expressions that are longer than the final result.

## 4 | Some Simple Examples with Mathematical Interpretations

The multiway graphs in the previous section are in a sense fundamentally metamathematical. Their "raw material" is mathematical statements. But what they represent are the results of operations—like substitution—that are defined at a kind of meta level, that "talks about mathematics" but isn't itself immediately "representable as mathematics". But to help understand this relationship it's useful to look at simple cases where it's possible to make at least some kind of correspondence with familiar mathematical concepts.

Consider for example the axiom

$$x\_ \circ y\_ = y\_ \circ x\_$$

that we can think of as representing commutativity of the binary operator ∘. Now consider using substitution to "apply this axiom", say starting from the expression $(a \circ b) \circ (c \circ d)$. The result is the (finite) multiway graph:

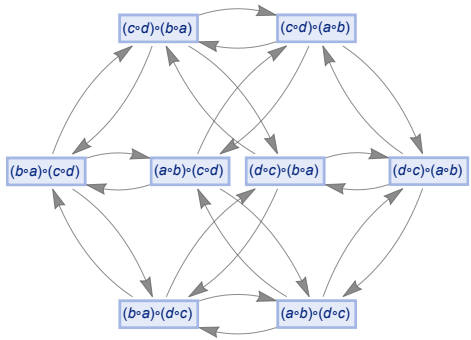

Conflating the pairs of edges going in opposite directions, the resulting graphs starting from any expression involving $s$ ∘'s (and $s+1$ distinct variables) are:

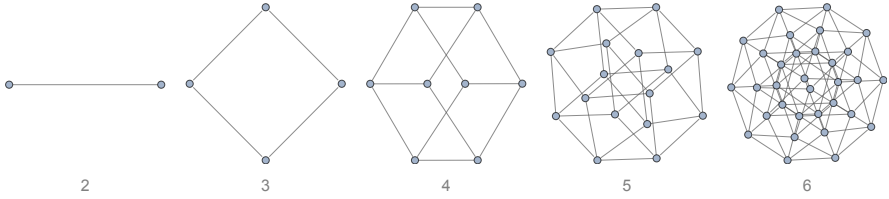

And these are just the Boolean hypercubes, each with $2^s$ nodes.

If instead of commutativity we consider the associativity axiom

$$x\_ \circ (y\_ \circ z\_) = (x\_ \circ y\_) \circ z\_$$

then we get a simple "ring" multiway graph:

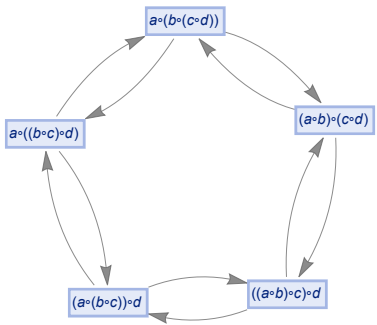

With both associativity and commutativity we get:

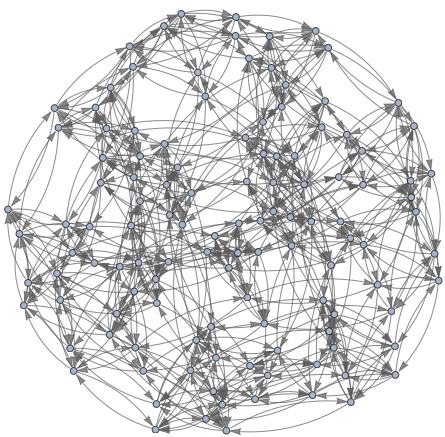

What is the mathematical significance of this object? We can think of our axioms as being the general axioms for a commutative semigroup. And if we build a multiway graph—say starting with (a∘b)∘b—we'll find out what expressions are equivalent to (a∘b)∘b in any commutative semigroup—or, in other words, we'll get a collection of theorems that are "true for any commutative semigroup":

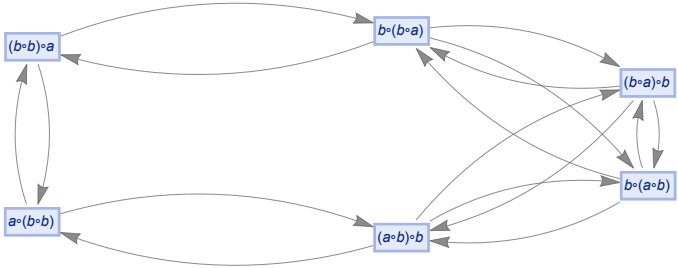

But what if we want to deal with a "specific semigroup" rather than a generic one? We can think of our symbols *a* and *b* as generators of the semigroup, and then we can add relations, as in:

$$x\_ \circ y\_ = y\_ \circ x\_, \; x\_ \circ (y\_ \circ z\_) = (x\_ \circ y\_) \circ z\_, \; a \circ a = b \circ b$$

And the result of this will be that we get more equivalences between expressions:

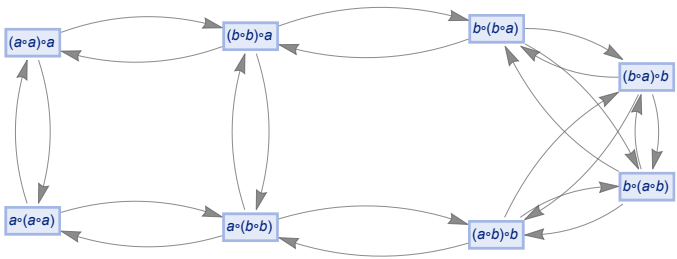

The multiway graph here is still finite, however, giving a finite number of equivalences. But let's say instead that we add the relations:

$$a = b \circ b, \; b = a \circ b$$

Then if we start from *a* we get a multiway graph that begins like

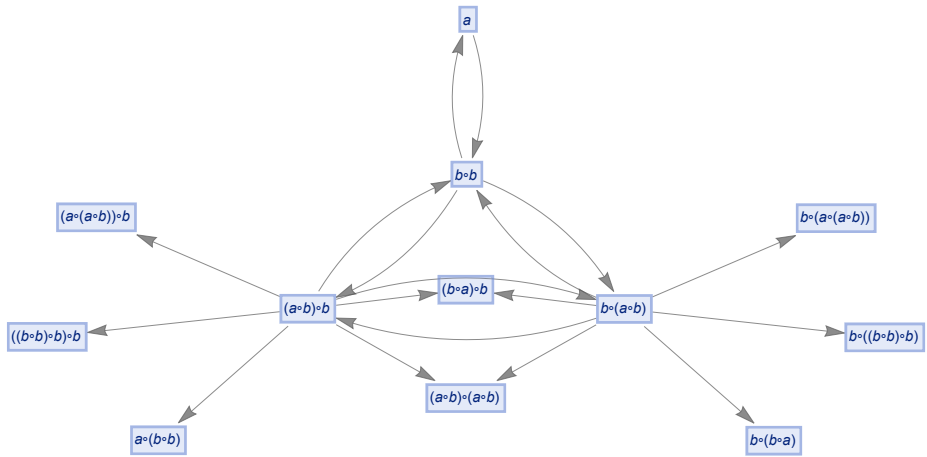

but just keeps growing forever (here shown after 6 steps):

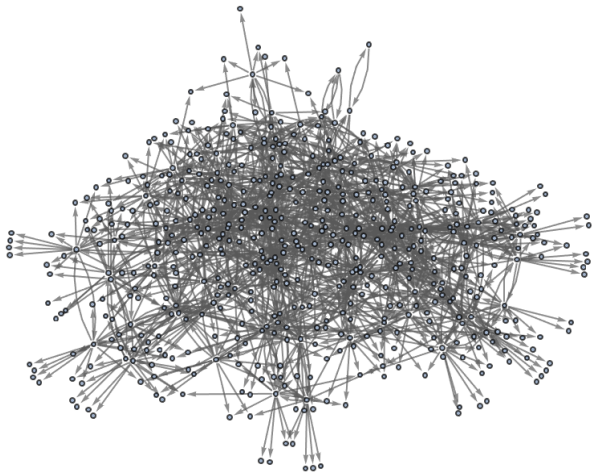

And what this then means is that there are an infinite number of equivalences between expressions. We can think of our basic symbols *a* and *b* as being generators of our semigroup. Then our expressions correspond to "words" in the semigroup formed from these generators. The fact that the multiway graph is infinite then tells us that there are an infinite number of equivalences between words.

But when we think about the semigroup mathematically we're typically not so interested in specific words as in the overall "distinct elements" in the semigroup, or in other words, in those "clusters of words" that don't have equivalences between them. And to find these we can imagine starting with all possible expressions, then building up multiway graphs from them. Many of the graphs grown from different expressions will join up. But what we want to know in the end is how many disconnected graph components are ultimately formed. And each of these will correspond to an element of the semigroup.

As a simple example, let's start from all words of length 2:

a∘a   a∘b   b∘a   b∘b

The multiway graphs formed from each of these after 1 step are:

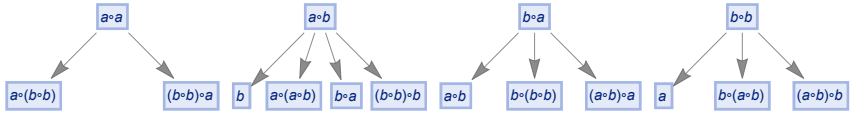

But these graphs in effect "overlap", leaving three disconnected components:

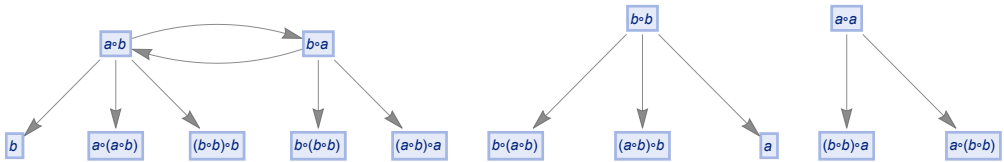

After 2 steps the corresponding result has two components:

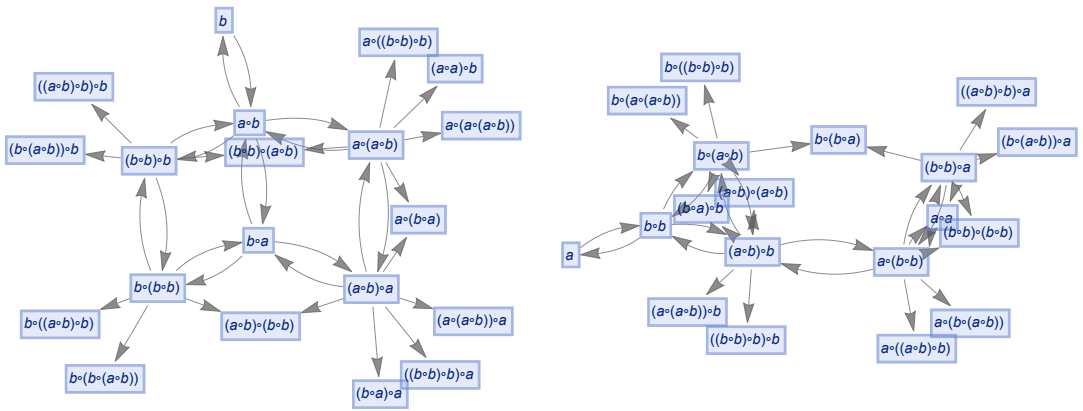

And if we start with longer (or shorter) words, and run for more steps, we'll keep finding the same result: that there are just two disconnected "droplets" that "condense out" of the "gas" of all possible initial words:

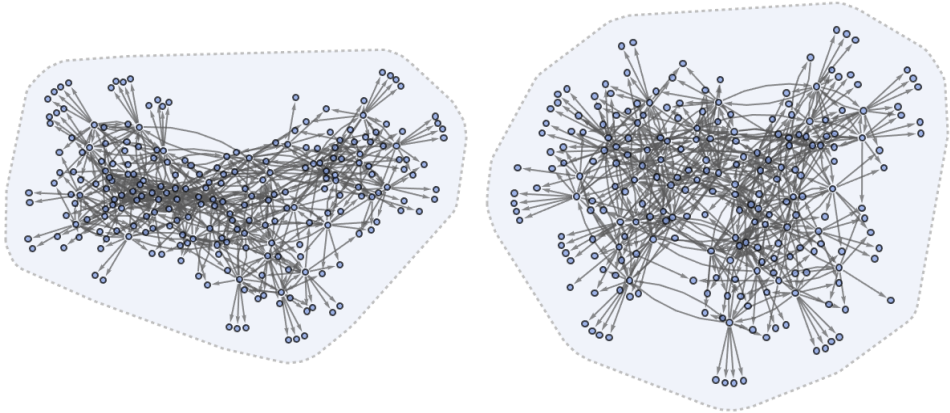

And what this means is that our semigroup ultimately has just two distinct elements—each of which can be represented by any of the different ("equivalent") words in each "droplet". (In this particular case the droplets just contain respectively all words with an odd or even number of *b*'s.)

In the mathematical analysis of semigroups (as well as groups), it's common ask what happens if one forms products of elements. In our setting what this means is in effect that one wants to "combine droplets using ∘". The simplest words in our two droplets are respectively *a* and *b*. And we can use these as "representatives of the droplets". Then we can see how multiplication by *a* and by *b* transforms words from each droplet:

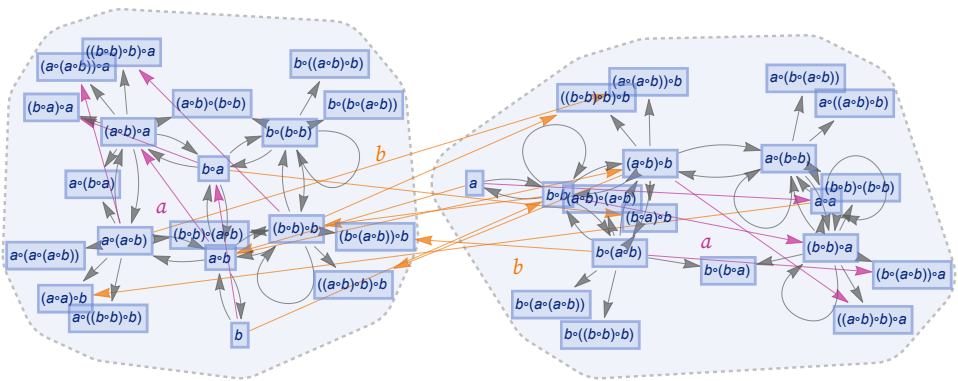

With only finite words the multiplications will sometimes not "have an immediate target" (so they are not indicated here). But in the limit of an infinite number of multiway steps, every multiplication will "have a target" and we'll be able to summarize the effect of multiplication in our semigroup by the graph:

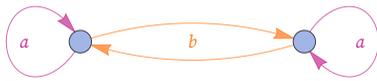

More familiar as mathematical objects than semigroups are groups. And while their axioms are slightly more complicated, the basic setup we've discussed for semigroups also applies to groups. And indeed the graph we've just generated for our semigroup is very much like a standard Cayley graph that we might generate for a group—in which the nodes are elements of the group and the edges define how one gets from one element to another by multiplying by a generator. (One technical detail is that in Cayley graphs identity-element self-loops are normally dropped.)

Consider the group $D_2$ (the "Klein four-group"). In our notation the axioms for this group can be written:

$$x_\_ \circ (y_\_ \circ z_\_) = (x_\_ \circ y_\_) \circ z_\_, \; x_\_ \circ y_\_ = y_\_ \circ x_\_, \; x_\_ = a \circ x_\_, \; a = b \circ b, \; a = c \circ c$$

Given these axioms we do the same construction as for the semigroup above. And what we find is that now four "droplets" emerge, corresponding to the four elements of $D_2$

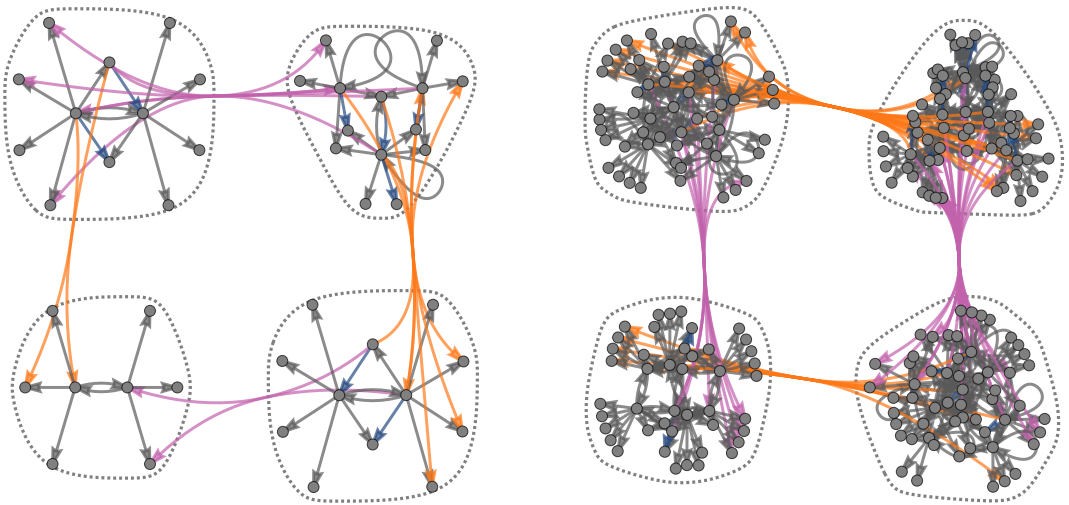

and the pattern of connections between them in the limit yields exactly the Cayley graph for $D_2$:

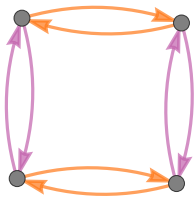

We can view what's happening here as a first example of something we'll return to at length later: the idea of "parsing out" recognizable mathematical concepts (here things like elements of groups) from lower-level "purely metamathematical" structures.

## 5 | Metamathematical Space

In multiway graphs like those we've shown in previous sections we routinely generate very large numbers of "mathematical" expressions. But how are these expressions related to each other? And in some appropriate limit can we think of them all being embedded in some kind of "metamathematical space"?

It turns out that this is the direct analog of what in our Physics Project we call branchial space, and what in that case defines a map of the entanglements between branches of quantum history. In the mathematical case, let's say we have a multiway graph generated using the axiom:

$$x\_ \circ y\_ \leftrightarrow (y\_ \circ x\_) \circ y\_$$

After a few steps starting from $a \circ b$ we have:

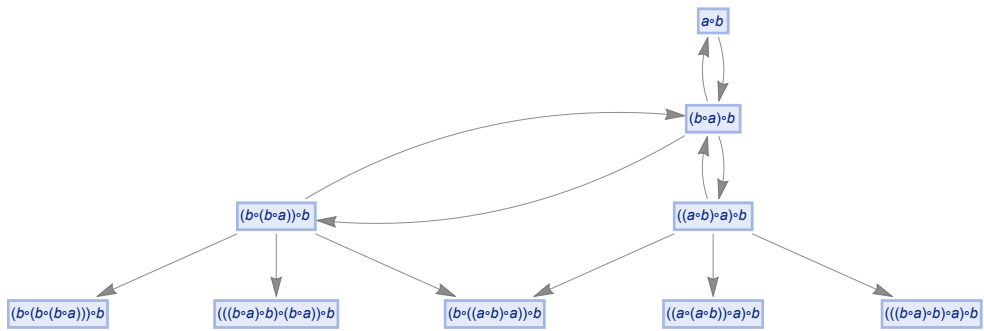

Now—just as in our Physics Project—let's form a branchial graph by looking at the final expressions here and connecting them if they are "entangled" in the sense that they share an ancestor on the previous step:

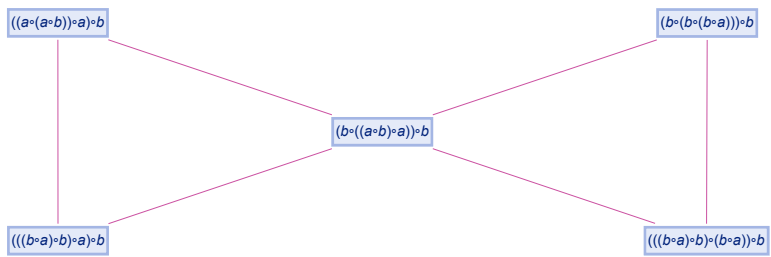

There's some trickiness here associated with loops in the multiway graph (which are the analog of closed timelike curves in physics) and what it means to define different "steps in evolution". But just iterating once more the construction of the multiway graph, we get a branchial graph:

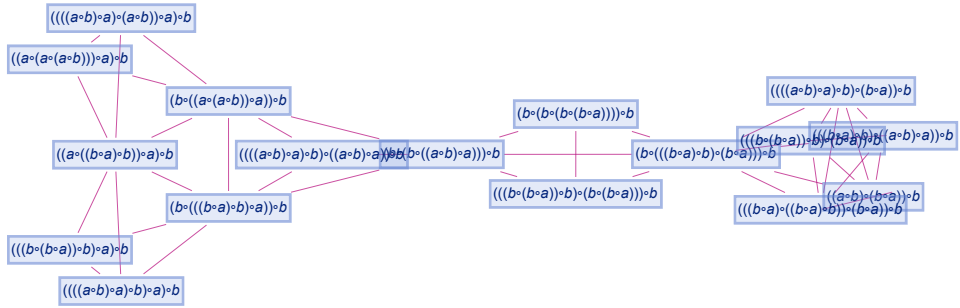

After a couple more iterations the structure of the branchial graph is (with each node sized according to the size of expression it represents):

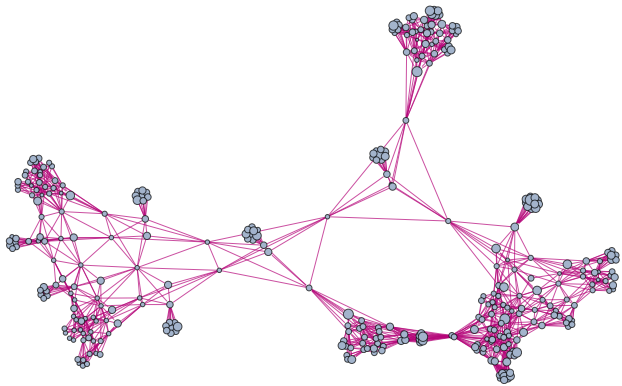

Continuing another iteration, the structure becomes:

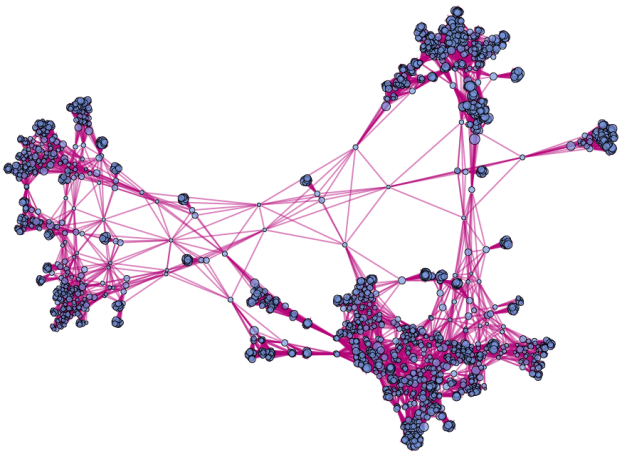

And in essence this structure can indeed be thought of as defining a kind of "metamathematical space" in which the different expressions are embedded. But what is the "geography" of this space? This shows how expressions (drawn as trees) are laid out on a particular branchial graph and we see that there is at least a general clustering of similar trees on the graph—indicating that "similar expressions" tend to be "nearby" in the metamathematical space defined by this axiom system.

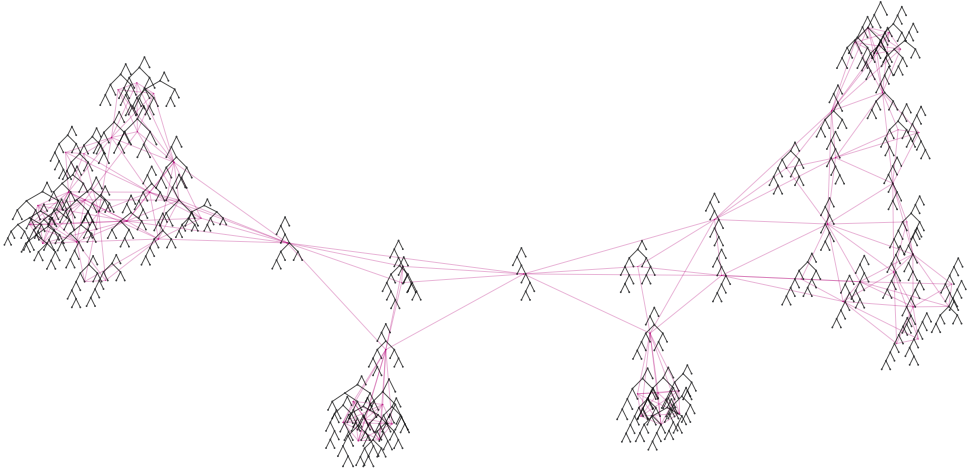

An important feature of branchial graphs is that effects are—essentially by construction—always local in the branchial graph. For example, if one changes an expression at a particular step in the evolution of a multiway system, it can only affect a region of the branchial graph that essentially expands by one edge per step.

One can think of the affected region—in analogy with a light cone in spacetime—as being the "entailment cone" of a particular expression. The edge of the entailment cone in effect expands at a certain "maximum metamathematical speed" in metamathematical (i.e. branchial) space—which one can think of as being measured in units of "expression change per multiway step".

By analogy with physics one can start talking in general about motion in metamathematical space. A particular proof path in the multiway graph will progressively "move around" in the branchial graph that defines metamathematical space. (Yes, there are many subtle issues here, not least the fact that one has to imagine a certain kind of limit being taken so that the structure of the branchial graph is "stable enough" to "just be moving around" in something like a "fixed background space".)

By the way, the shortest proof path in the multiway graph is the analog of a geodesic in spacetime. And later we'll talk about how the "density of activity" in the branchial graph is the analog of energy in physics, and how it can be seen as "deflecting" the path of geodesics, just as gravity does in spacetime.

It's worth mentioning just one further subtlety. Branchial graphs are in effect associated with "transverse slices" of the multiway graph—but there are many consistent ways to make these slices. In physics terms one can think of the foliations that define different choices of sequences of slices as being like "reference frames" in which one is specifying a sequence of "simultaneity surfaces" (here "branchtime hypersurfaces"). The particular branchial graphs we've shown here are ones associated with what in physics might be called the cosmological rest frame in which every node is the result of the same number of updates since the beginning.

## 6 | The Issue of Generated Variables

A rule like

$$x\_ \circ y\_ \leftrightarrow (y\_ \circ x\_) \circ y\_$$

defines transformations for any expressions $x\_$ and $y\_$. So, for example, if we use the rule from left to right on the expression $a \circ (b \circ a)$ the "pattern variable" $x\_$ will be taken to be $a$ while $y\_$ will be taken to be $b \circ a$, and the result of applying the rule will be $((b \circ a) \circ a) \circ (b \circ a)$.

But consider instead the case where our rule is:

$$x\_ \circ y\_ \leftrightarrow (y\_ \circ x\_) \circ z\_$$

Applying this rule (from left to right) to $a \circ (b \circ a)$ we'll now get $((b \circ a) \circ a) \circ z\_$. And applying the rule to $a \circ b$ we'll get $(b \circ a) \circ z\_$. But what should we make of those $z\_$'s? And in particular, are they "the same", or not?

A pattern variable like $z\_$ can stand for any expression. But do two different $z\_$'s have to stand for the same expression? In a rule like $z\_ \circ z\_ \leftrightarrow ...$ we're assuming that, yes, the two $z\_$'s always stand for the same expression. But if the $z\_$'s appear in different rules it's a different story. Because in that case we're dealing with two separate and unconnected $z\_$'s—that can stand for completely different expressions.

To begin seeing how this works, let's start with a very simple example. Consider the (for now, one-way) rule

$$a \to x\_$$

where $a$ is the literal symbol $a$, and $x\_$ is a pattern variable. Applying this to $a \circ a$ we might think we could just write the result as:

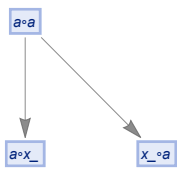

Then if we apply the rule again both branches will give the same expression $x\_\circ x\_$, so there'll be a merge in the multiway graph:

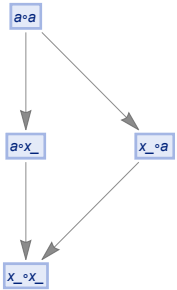

But is this really correct? Well, no. Because really those should be two different $x\_$'s, that could stand for two different expressions. So how can we indicate this? One approach is just to give every "generated" $x\_$ a new name:

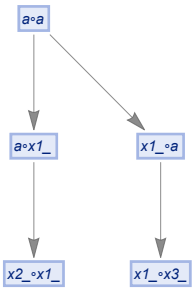

But this result isn't really correct either. Because if we look at the second step we see the two expressions $x1\_\circ x3\_$ and $x2\_\circ x4\_$. But what's really the difference between these? The names $xi$ are arbitrary; the only constraint is that within any given expression they have to be different. But between expressions there's no such constraint. And in fact $x1\_\circ x3\_$ and $x2\_\circ x4\_$ both represent exactly the same class of expressions: any expression of the form $u\_\circ v\_$.

So in fact it's not correct that there are two separate branches of the multiway system producing two separate expressions. Because those two branches produce equivalent expressions, which means they can be merged. And turning both equivalent expressions into the same canonical form we get:

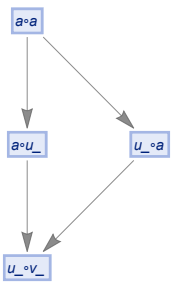

It's important to notice that this isn't the same result as what we got when we assumed that every $x\_$ was the same. Because then our final result was the expression $x\_ \circ x\_$ which can match $a \circ a$ but not $a \circ b$—whereas now the final result is $u\_ \circ v\_$ which can match both $a \circ a$ and $a \circ b$.

This may seem like a subtle issue. But it's critically important in practice. Not least because generated variables are in effect what make up all "truly new stuff" that can be produced. With a rule like $x\_ \circ y\_ \leftrightarrow (y\_ \circ x\_) \circ y\_$ one's essentially just taking whatever one started with, and successively rearranging the pieces of it. But with a rule like $x\_ \circ y\_ \leftrightarrow (y\_ \circ x\_) \circ z\_$ there's something "truly new" generated every time $z\_$ appears.

By the way, the basic issue of "generated variables" isn't something specific to the particular symbolic expression setup we've been using here. For example, there's a direct analog of it in the hypergraph rewriting systems that appear in our Physics Project. But in that case there's a particularly clear interpretation: the analog of "generated variables" are new "atoms of space" produced by the application of rules. And far from being some kind of footnote, these "generated atoms of space" are what make up everything we have in our universe today.

The issue of generated variables—and especially their naming—is the bane of all sorts of formalism for mathematical logic and programming languages. As we'll see later, it's perfectly possible to "go to a lower level" and set things up with no names at all, for example using combinators. But without names, things tend to seem quite alien to us humans—and certainly if we want to understand the correspondence with standard presentations of mathematics it's pretty necessary to have names. So at least for now we'll keep names, and handle the issue of generated variables by uniquifying their names, and canonicalizing every time we have a complete expression.

Let's look at another example to see the importance of how we handle generated variables. Consider the rule:

$x\_ \circ y\_ \leftrightarrow y\_ \circ z\_$

If we start with $a \circ a$ and do no uniquification, we'll get:

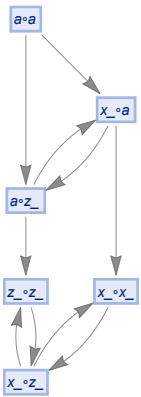

With uniquification, but not canonicalization, we'll get a pure tree:

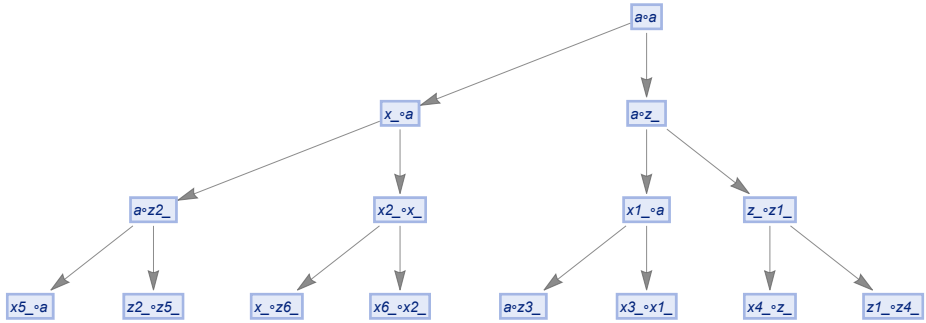

But with canonicalization this is reduced to:

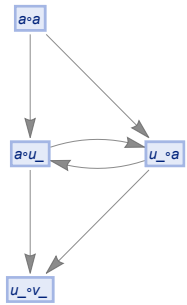

A confusing feature of this particular example is that this same result would have been obtained just by canonicalizing the original "assume-all-$x\_$'s-are-the-same" case.

But things don't always work this way. Consider the rather trivial rule

$x\_ \leftrightarrow y\_$

starting from $a{\circ}x\_$. If we don't do uniquification, and don't do canonicalization, we get:

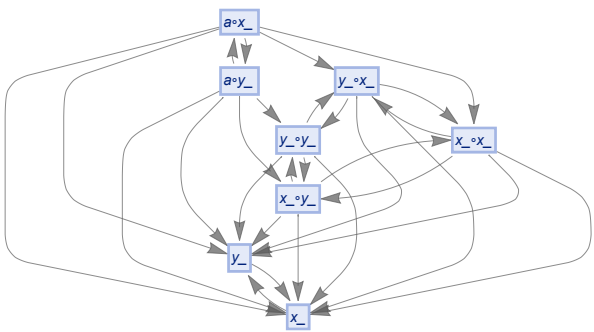

If we do uniquification (but not canonicalization), we get a pure tree:

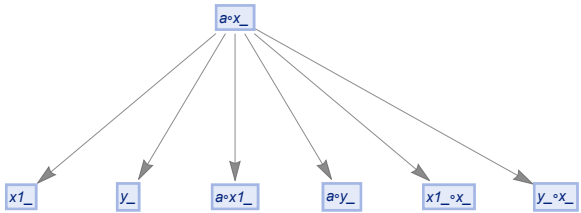

But if we now canonicalize this, we get:

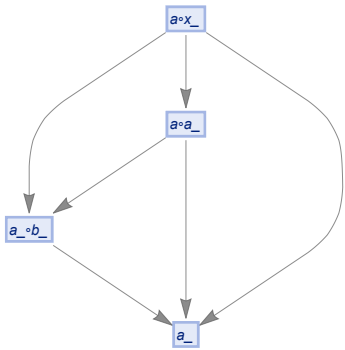

And this is now not the same as what we would get by canonicalizing, without uniquifying:

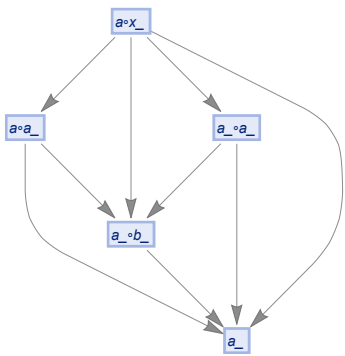

# 7 | Rules Applied to Rules

In what we've done so far, we've always talked about applying rules (like $x\_\circ y\_ \leftrightarrow (y\_\circ x\_)\circ y\_$) to expressions (like $((b\circ a)\circ a)\circ b$ or $x\_\circ(y\_\circ x\_)$). But if everything is a symbolic expression there shouldn't really need to be a distinction between "rules" and "ordinary expressions". They're all just expressions. And so we should as well be able to apply rules to rules as to ordinary expressions.

And indeed the concept of "applying rules to rules" is something that has a familiar analog in standard mathematics. The "two-way rules" we've been using effectively define equivalences—which are very common kinds of statements in mathematics, though in mathematics they're usually written with = rather than with ↔. And indeed, many axioms and many theorems are specified as equivalences—and in equational logic one takes everything to be defined using equivalences. And when one's dealing with theorems (or axioms) specified as equivalences, the basic way one derives new theorems is by applying one theorem to another—or in effect by applying rules to rules.

As a specific example, let's say we have the "axiom":

$$x\_ \circ y\_ \leftrightarrow (y\_ \circ x\_) \circ y\_$$

We can now apply this to the rule

$$a \circ a \leftrightarrow b \circ b$$

to get (where since $u \leftrightarrow v$ is equivalent to $v \leftrightarrow u$ we're sorting each two-way rule that arises)

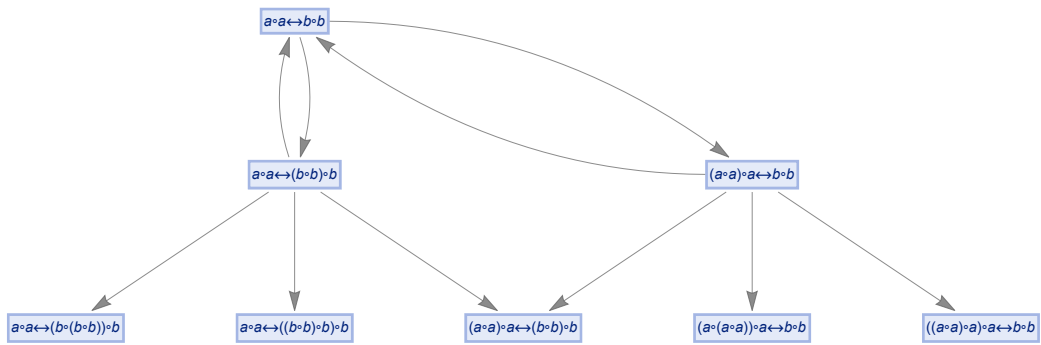

or after a few more steps:

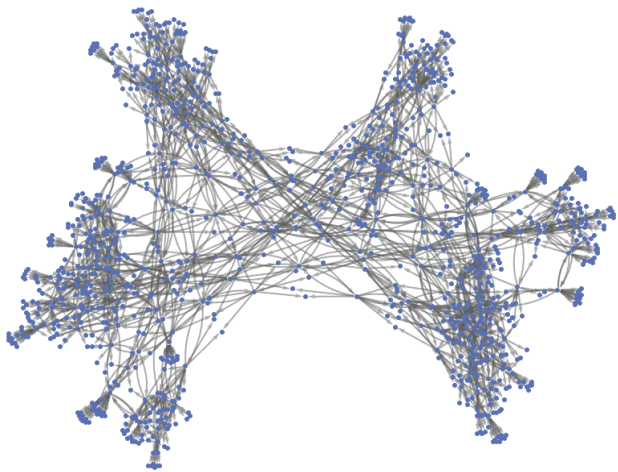

In this example all that's happening is that the substitutions specified by the axiom are getting separately applied to the left- and right-hand sides of each rule that is generated. But if we really take seriously the idea that everything is a symbolic expression, things can get a bit more complicated.

Consider for example the rule:

$x\_ \leftrightarrow x\_ \circ y\_$

If we apply this to

$a \circ a \leftrightarrow b \circ b$

then if $x\_$ "matches any expression" it can match the whole expression $a \circ a \leftrightarrow b \circ b$ giving the result:

$(a \circ a \leftrightarrow b \circ b) \circ y$

Standard mathematics doesn't have an obvious meaning for something like this—although as soon as one "goes metamathematical" it's fine. But in an effort to maintain contact with standard mathematics we'll for now have the "meta rule" that $x\_$ can't match an expression whose top-level operator is $\leftrightarrow$. (As we'll discuss later, including such matches would allow us to do exotic things like encode set theory within arithmetic, which is again something usually considered to be "syntactically prevented" in mathematical logic.)

Another—still more obscure—meta rule we have is that $x\_$ can't "match inside a variable". In Wolfram Language, for example, $a\_$ has the full form **Pattern**[a,**Blank**[]], and one could imagine that $x\_$ could match "internal pieces" of this. But for now, we're going to treat all variables as atomic—even though later on, when we "descend below the level of variables", the story will be different.

When we apply a rule like $x\_ \leftrightarrow (x\_ \circ y\_)$ to $a \circ a \leftrightarrow b \circ b$ we're taking a rule with pattern variables, and doing substitutions with it on a "literal expression" without pattern variables. But it's also perfectly possible to apply pattern rules to pattern rules—and indeed that's what we'll mostly do below. But in this case there's another subtle issue that can arise. Because if our rule generates variables, we can end up with two different kinds of variables with "arbitrary names": generated variables, and pattern variables from the rule we're operating on. And when we canonicalize the names of these variables, we can end up with identical expressions that we need to merge.

Here's what happens if we apply the rule $x\_ \leftrightarrow (x\_ \circ y\_)$ to the literal rule $a \circ b \leftrightarrow a$:

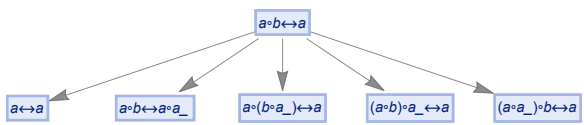

If we apply it to the pattern rule $a\_\circ b\_\leftrightarrow a\_$ but don't do canonicalization, we'll just get the same basic result:

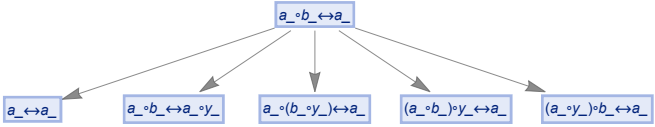

But if we canonicalize we get instead:

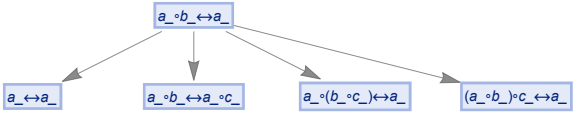

The effect is more dramatic if we go to two steps. When operating on the literal rule we get:

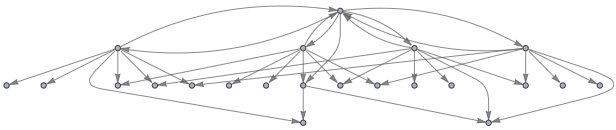

Operating on the pattern rule, but without canonicalization, we get

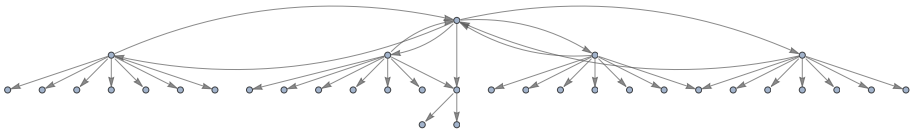

while if we include canonicalization many rules merge and we get:

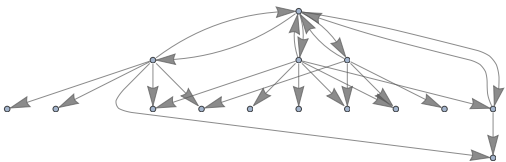

# 8 | Accumulative Evolution

We can think of "ordinary expressions" like $a \circ b$ as being like "data", and rules as being like "code". But when everything is a symbolic expression, it's perfectly possible—as we saw above—to "treat code like data", and in particular to generate rules as output. But this now raises a new possibility. When we "get a rule as output", why not start "using it like code" and applying it to things?

In mathematics we might apply some theorem to prove a lemma, and then we might subsequently use that lemma to prove another theorem—eventually building up a whole "accumulative structure" of lemmas (or theorems) being used to prove other lemmas. In any given proof we can in principle always just keep using the axioms over and over again—but it'll be much more efficient to progressively build a library of more and more lemmas, and use these. And in general we'll build up a richer structure by "accumulating lemmas" than always just going back to the axioms.

In the multiway graphs we've drawn so far, each edge represents the application of a rule, but that rule is always a fixed axiom. To represent accumulative evolution we need a slightly more elaborate structure—and it'll be convenient to use token-event graphs rather than pure multiway graphs.

Every time we apply a rule we can think of this as an event. And with the setup we're describing, that event can be thought of as taking two tokens as input: one the "code rule" and the other the "data rule". The output from the event is then some collection of rules, which can then serve as input (either "code" or "data") to other events.

Let's start with the very simple example of the rule

$$x \leftrightarrow y$$

where for now there are no patterns being used. Starting from this rule, we get the token-event graph (where now we're indicating the initial "axiom" statement using a slightly different color):

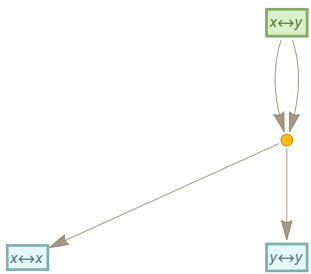

One subtlety here is that the $x \leftrightarrow y$ is applied to itself—so there are two edges going into the event from the node representing the rule. Another subtlety is that there are two different ways the rule can be applied, with the result that there are two output rules generated.

Here's another example, based on the two rules:

$$x \leftrightarrow y, y \leftrightarrow z$$

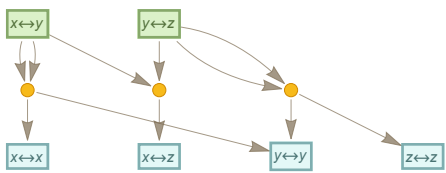

Continuing for another step we get:

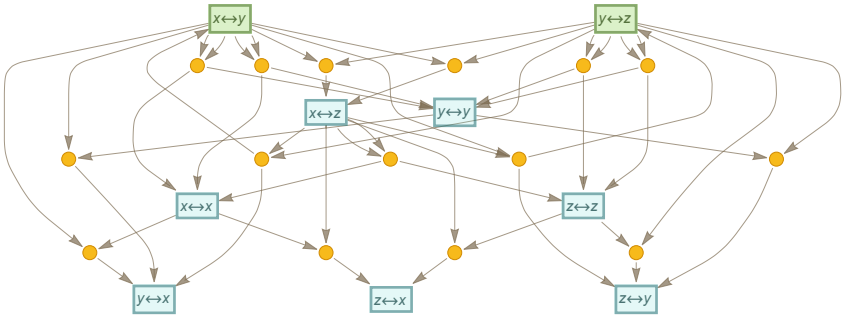

Typically we will want to consider ⟷ as "defining an equivalence", so that $v \leftrightarrow u$ means the same as $u \leftrightarrow v$, and can be conflated with it—yielding in this case:

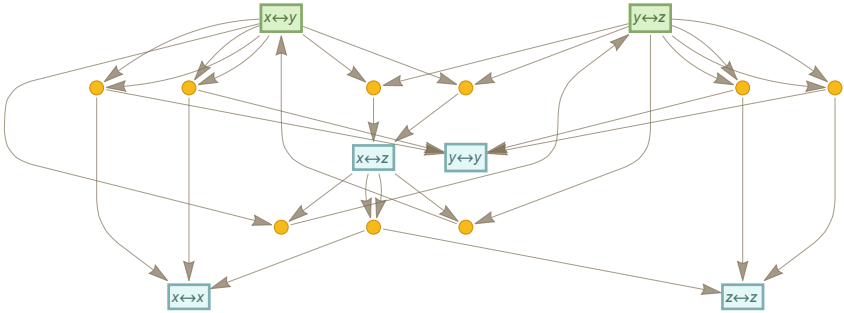

Now let's consider the rule:

$a \circ b \leftrightarrow b$

After one step we get:

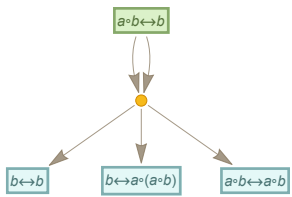

After 2 steps we get:

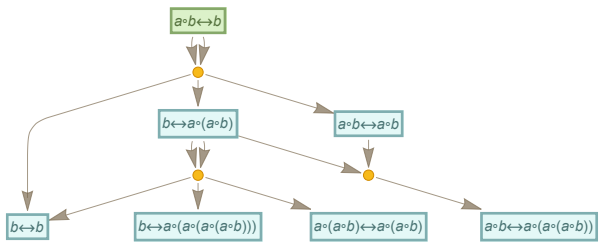

The token-event graphs after 3 and 4 steps in this case are (where now we've deduplicated events):

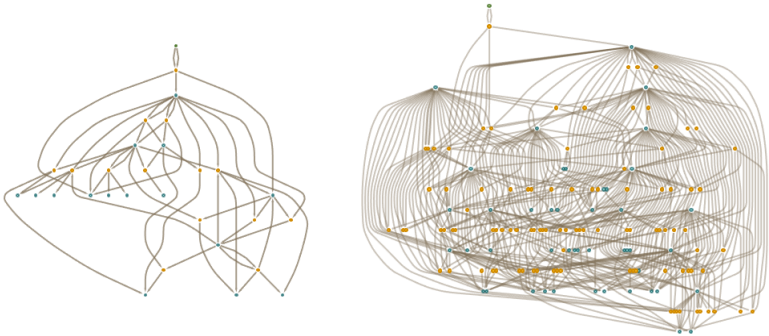

Let's now consider a rule with the same structure, but with pattern variables instead of literal symbols:

$a\_ \circ b\_ \leftrightarrow b\_$

Here's what happens after one step (note that there's canonicalization going on, so $a\_$'s in different rules aren't "the same")

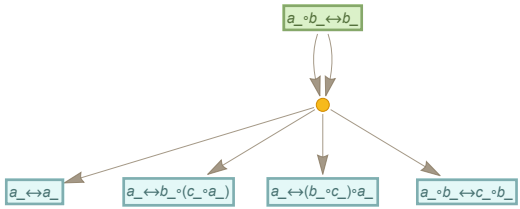

and we see that there are different theorems from the ones we got without patterns. After 2 steps with the pattern rule we get

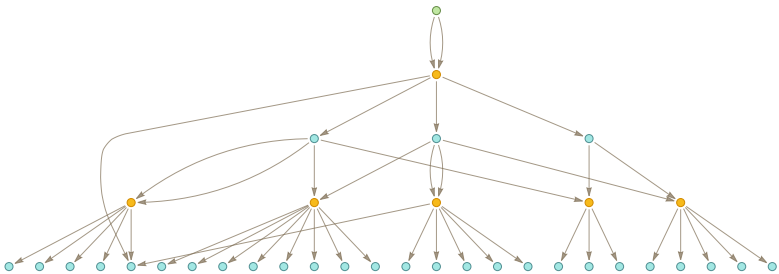

where now the complete set of "theorems that have been derived" is (dropping the _'s for readability)

| | | | | |
|---|---|---|---|---|
| $a \circ b \leftrightarrow b$ | $a \leftrightarrow a$ | $a \leftrightarrow b \circ (c \circ a)$ | $a \leftrightarrow (b \circ c) \circ a$ | $a \circ b \leftrightarrow c \circ b$ |
| $a \leftrightarrow b \circ a$ | $a \leftrightarrow b \circ (c \circ (d \circ (e \circ a)))$ | $a \leftrightarrow b \circ (c \circ ((d \circ e) \circ a))$ | $a \leftrightarrow b \circ ((c \circ d) \circ (e \circ a))$ | $a \leftrightarrow b \circ ((c \circ (d \circ e)) \circ a)$ |
| $a \leftrightarrow b \circ (((c \circ d) \circ e) \circ a)$ | $a \leftrightarrow (b \circ c) \circ (d \circ (e \circ a))$ | $a \leftrightarrow (b \circ c) \circ ((d \circ e) \circ a)$ | $a \leftrightarrow (b \circ (c \circ d)) \circ (e \circ a)$ | $a \leftrightarrow (b \circ (c \circ (d \circ e))) \circ a$ |
| $a \leftrightarrow (b \circ ((c \circ d) \circ e)) \circ a$ | $a \leftrightarrow ((b \circ c) \circ d) \circ (e \circ a)$ | $a \leftrightarrow ((b \circ c) \circ (d \circ e)) \circ a$ | $a \leftrightarrow ((b \circ (c \circ d)) \circ e) \circ a$ | $a \leftrightarrow (((b \circ c) \circ d) \circ e) \circ a$ |
| $a \circ b \leftrightarrow c \circ (d \circ (e \circ b))$ | $a \circ b \leftrightarrow c \circ ((d \circ e) \circ b)$ | $a \circ b \leftrightarrow (c \circ d) \circ (e \circ b)$ | $a \circ b \leftrightarrow (c \circ (d \circ e)) \circ b$ | $a \circ b \leftrightarrow ((c \circ d) \circ e) \circ b$ |
| $a \circ (b \circ c) \leftrightarrow d \circ (e \circ c)$ | $a \circ (b \circ c) \leftrightarrow (d \circ e) \circ c$ | $a \circ (b \circ (c \circ d)) \leftrightarrow e \circ d$ | $a \circ ((b \circ c) \circ d) \leftrightarrow e \circ d$ | $(a \circ b) \circ c \leftrightarrow (d \circ e) \circ c$ |

or as trees:

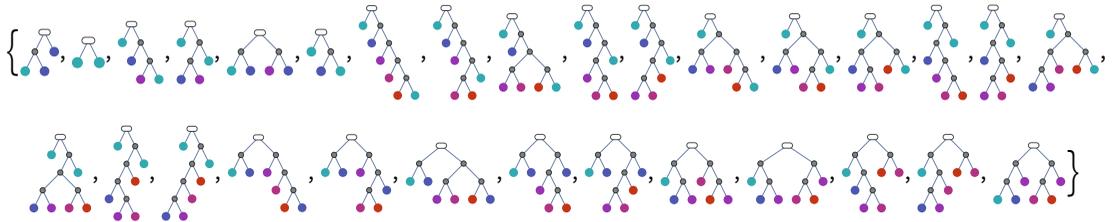

After another step one gets

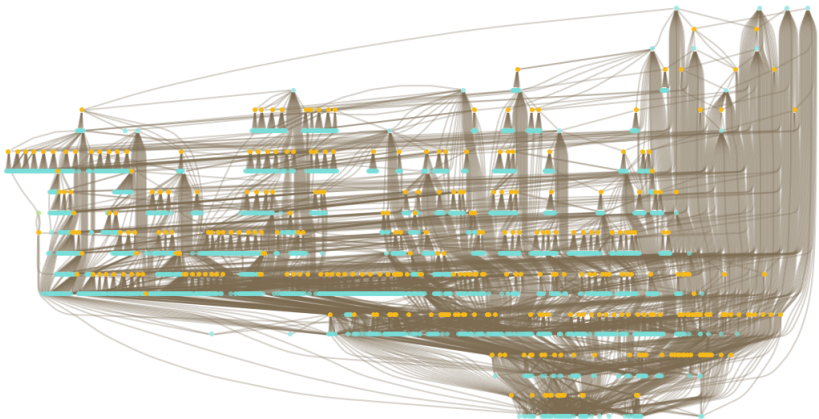

where now there are 2860 "theorems", roughly exponentially distributed across sizes according to

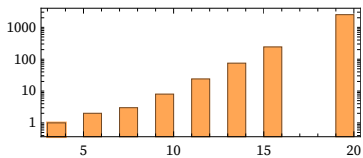

and with a typical "size-19" theorem being:

$$a\_ \circ b\_ \leftrightarrow c\_ \circ \big(d\_ \circ \big(e\_ \circ \big(f\_ \circ \big(g\_ \circ (h\_ \circ (i\_ \circ b\_))\big)\big)\big)\big)$$

In effect we can think of our original rule (or "axiom") as having initiated some kind of "mathematical Big Bang" from which an increasing number of theorems are generated. Early on we described having a "gas" of mathematical theorems that—a little like molecules—can interact and create new theorems. So now we can view our accumulative evolution process as a concrete example of this.

Let's consider the rule from previous sections:

$$x\_ \circ y\_ \leftrightarrow (y\_ \circ x\_) \circ y\_$$

After one step of accumulative evolution according to this rule we get:

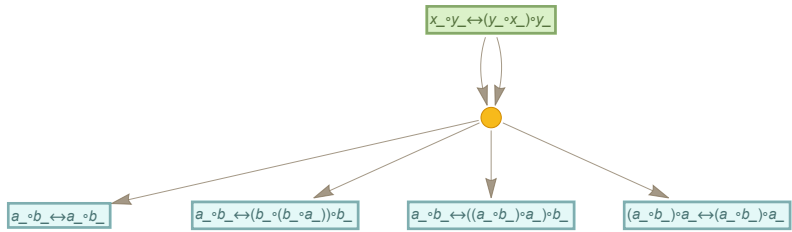

After 2 and 3 steps the results are:

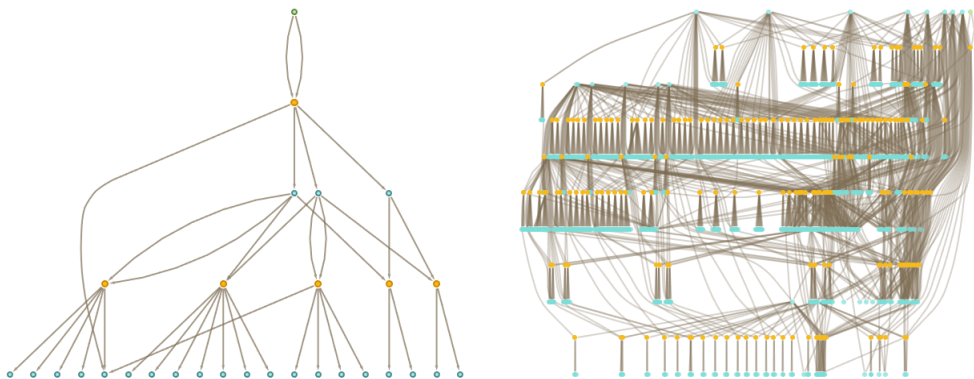

What is the significance of all this complexity? At a basic level, it's just an example of the ubiquitous phenomenon in the computational universe (captured in the Principle of Computational Equivalence) that even systems with very simple rules can generate behavior as complex as anything. But the question is whether—on top of all this complexity—there are simple "coarse-grained" features that we can identify as "higher-level mathematics"; features that we can think of as capturing the "bulk" behavior of the accumulative evolution of axiomatic mathematics.

## 9 | Accumulative String Systems

As we've just seen, the accumulative evolution of even very simple transformation rules for expressions can quickly lead to considerable complexity. And in an effort to understand the essence of what's going on, it's useful to look at the slightly simpler case not of rules for "tree-structured expressions" but instead at rules for strings of characters.

Consider the seemingly trivial case of the rule:

A ↔ B

After one step this gives

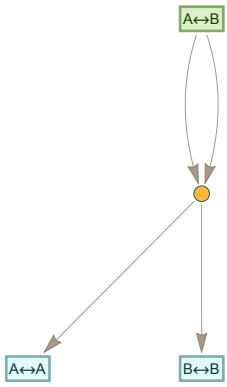

while after 2 steps we get

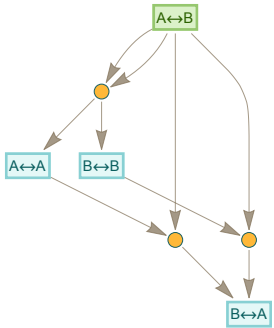

though treating *u*↔*v* as the same as *v*↔*u* this just becomes:

Here's what happens with the rule:

AB ↔ B

After 2 steps we get

and after 3 steps

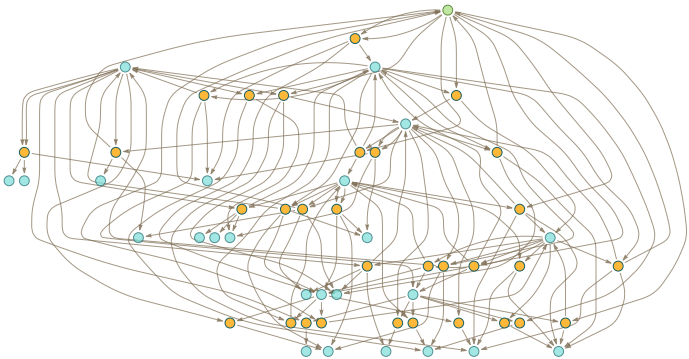

where now there are a total of 25 "theorems", including (unsurprisingly) things like:

AAAAAAB ↔ AAB

It's worth noting that despite the "lexical similarity" of the string rule AB↔B we're now using to the expression rule $a \circ b \leftrightarrow b$ from the previous section, these rules actually work in very different ways. The string rule can apply to characters anywhere within a string, but what it inserts is always of fixed size. The expression rule deals with trees, and only applies to "whole subtrees", but what it inserts can be a tree of any size. (One can align these setups by thinking of strings as expressions in which characters are "bound together" by an associative operator, as in A·B·A·A. But if one explicitly gives associativity axioms these will lead to additional pieces in the token-event graph.)

A rule like $a\_ \circ b\_ \leftrightarrow b\_$ also has the feature of involving patterns. In principle we could include patterns in strings too—both for single characters (as with _) and for sequences of characters (as with __)—but we won't do this here. (We can also consider one-way rules, using → instead of ↔.)

To get a general sense of the kinds of things that happen in accumulative (string) systems, we can consider enumerating all possible distinct two-way string transformation rules. With only a single character A, there are only two distinct cases

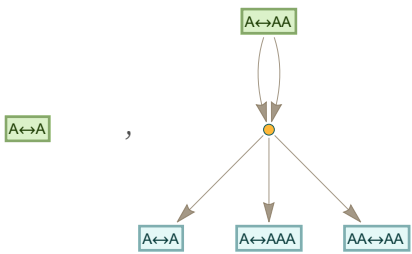

because A↔AA systematically generates all possible $A^n \leftrightarrow A^m$ rules

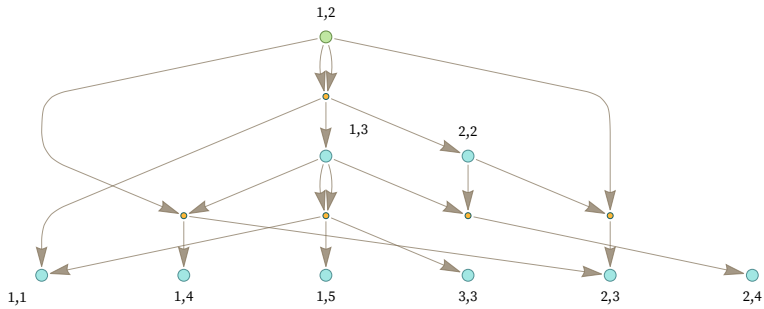

and at *t* steps gives a total number of rules equal to:

$(2^{t-1}+1)^2$

With characters A and B the distinct token-event graphs generated starting from rules with a total of at most 5 characters are:

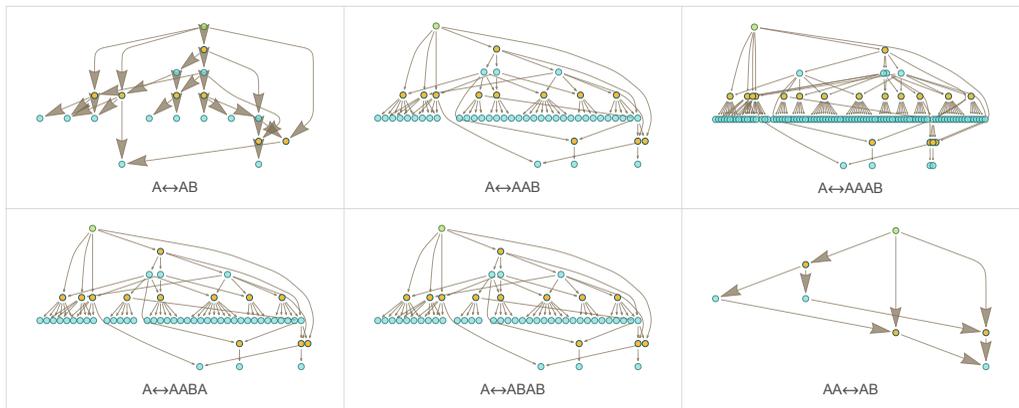

Note that when the strings in the initial rule are the same length, only a rather trivial finite token-event graph is ever generated, as in the case of AA↔AB:

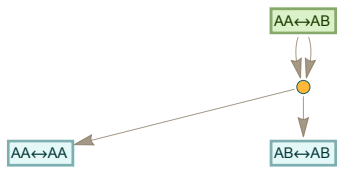

But when the strings are of different lengths, there is always unbounded growth.

## 10 | The Case of Hypergraphs

We've looked at accumulative versions of expression and string rewriting systems. So what about accumulative versions of hypergraph rewriting systems of the kind that appear in our Physics Project?

Consider the very simple hypergraph rule

{{*x*, *x*}} → {{*x*, *x*}, {*x*, *x*}}

or pictorially:

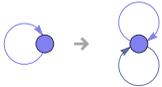

(Note that the nodes that are named 1 here are really like pattern variables, that could be named for example *x_*.)

We can now do accumulative evolution with this rule, at each step combining results that involve equivalent (i.e. isomorphic) hypergraphs:

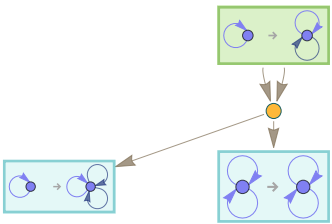

After two steps this gives:

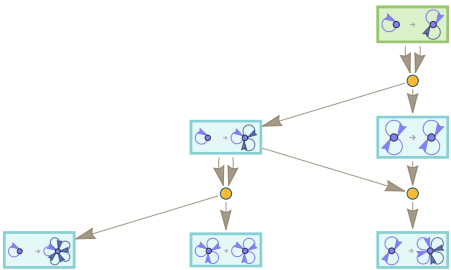

And after 3 steps:

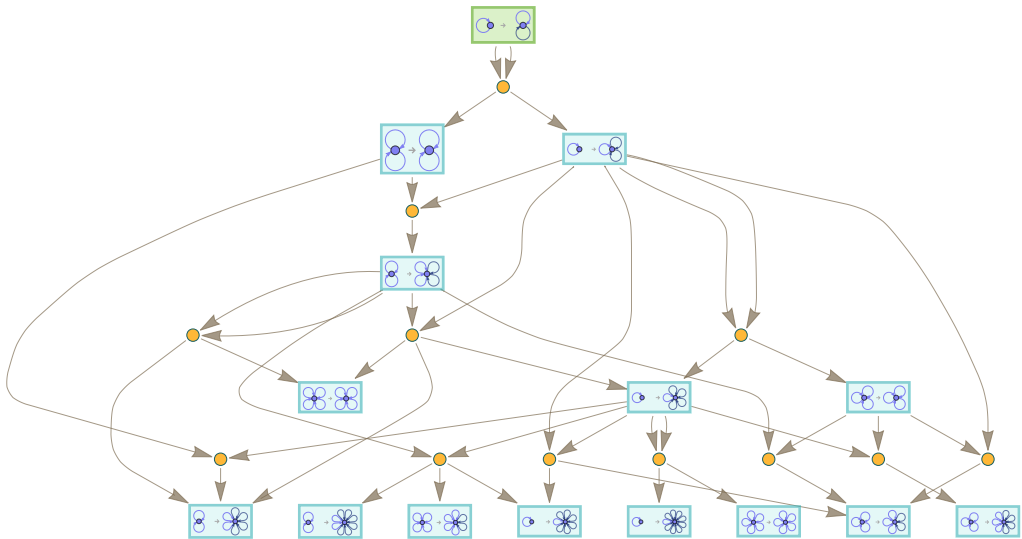

How does all this compare to "ordinary" evolution by hypergraph rewriting? Here's a multiway graph based on applying the same underlying rule repeatedly, starting from an initial condition formed from the rule:

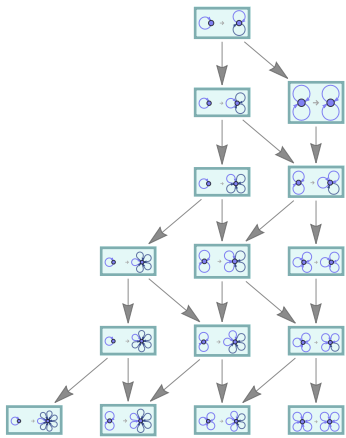

What we see is that the accumulative evolution in effect "shortcuts" the ordinary multiway evolution, essentially by "caching" the result of every piece of every transformation between states (which in this case are rules), and delivering a given state in fewer steps.

In our typical investigation of hypergraph rewriting for our Physics Project we consider one-way transformation rules. Inevitably, though, the ruliad contains rules that go both ways. And here, in an effort to understand the correspondence with our metamodel of mathematics, we can consider two-way hypergraph rewriting rules. An example is the two-way version of the rule above:

$\{\{x, x\}\} \leftrightarrow \{\{x, x\}, \{x, x\}\}$

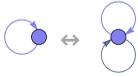

Now the token-event graph becomes

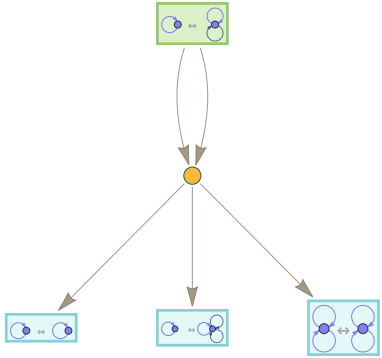

or after 2 steps (where now the transformations from "later states" to "earlier states" have started to fill in):

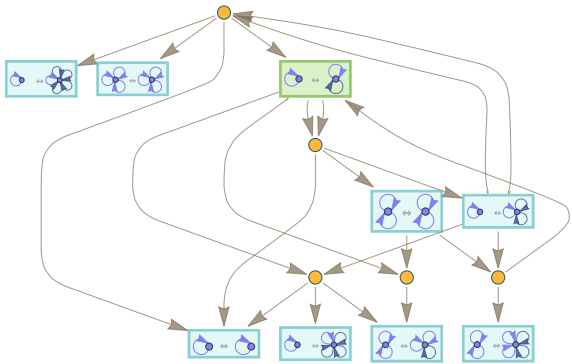

Just like in ordinary hypergraph evolution, the only way to get hypergraphs with additional hyperedges is to start with a rule that involves the addition of new hyperedges—and the same is true for the addition of new elements. Consider the rule:

$\{\{x, y\}, \{x, z\}\} \leftrightarrow \{\{x, y\}, \{x, w\}, \{y, w\}, \{z, w\}\}$

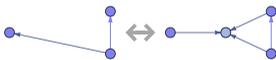

After 1 step this gives

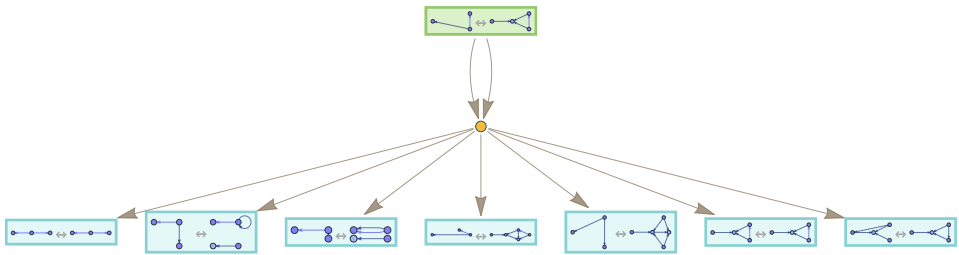

while after 2 steps it gives:

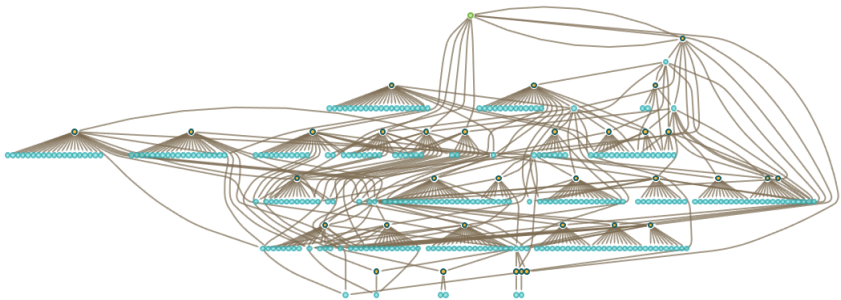

The general appearance of this token-event graph is not much different from what we saw with string rewrite or expression rewrite systems. So what this suggests is that it doesn't matter much whether we're starting from our metamodel of axiomatic mathematics or from any other reasonably rich rewriting system: we'll always get the same kind of "large-scale" token-event graph structure. And this is an example of what we'll use to argue for general laws of metamathematics.

## 11 | Proofs in Accumulative Systems

In an earlier section, we discussed how paths in a multiway graph can represent proofs of "equivalence" between expressions (or the "entailment" of one expression by another). For example, with the rule (or "axiom")

{A → BBB, BB → A}

this shows a path that "proves" that "BA entails AAB":

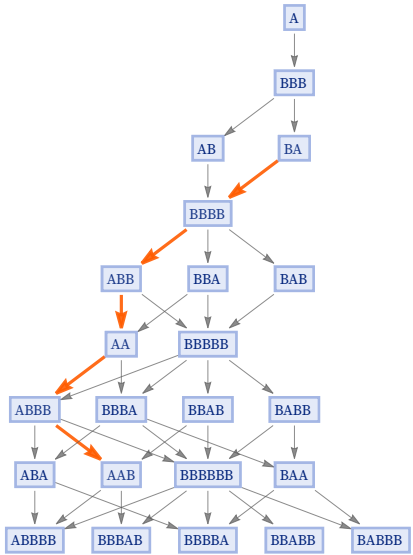

But once we know this, we can imagine adding this result (as what we can think of as a "lemma") to our original rule:

{A → BBB, BB → A, BA → AAB}

And now (the "theorem") "BA entails AAB" takes just one step to prove—and all sorts of other proofs are also shortened:

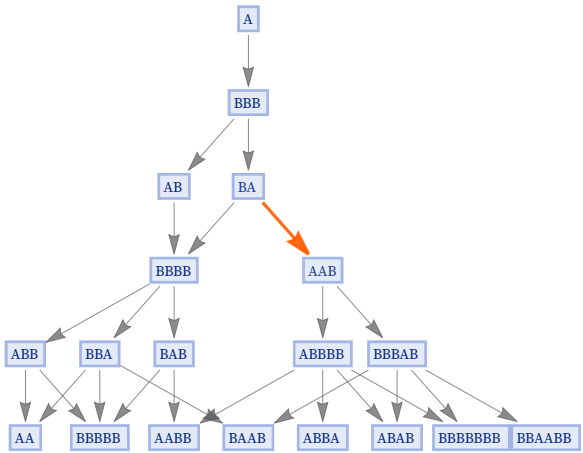

It's perfectly possible to imagine evolving a multiway system with a kind of "caching-based" speed-up mechanism where every new entailment discovered is added to the list of underlying rules. And, by the way, it's also possible to use two-way rules throughout the multiway system:

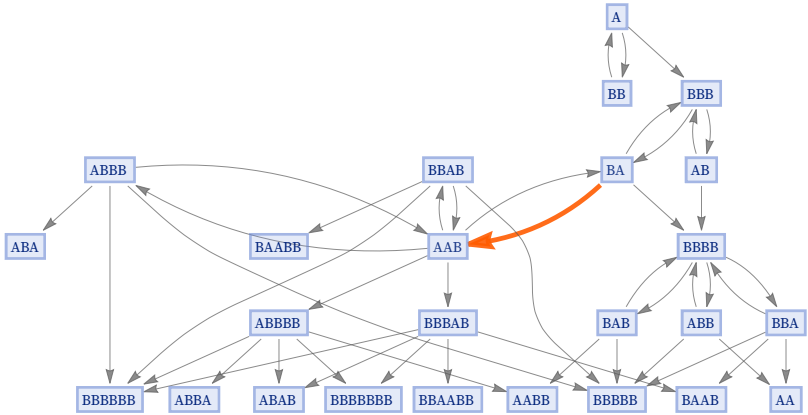

But accumulative systems provide a much more principled way to progressively "add what's discovered". So what do proofs look like in such systems?

Consider the rule:

A ↔ AB

Running it for 2 steps we get the token-event graph:

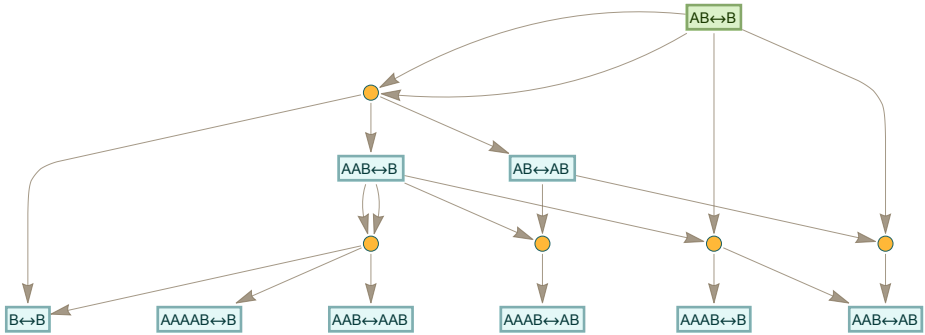

Now let's say we want to prove that the original "axiom" AB↔B implies (or "entails") the "theorem" AAAB↔AB. Here's the subgraph that demonstrates the result:

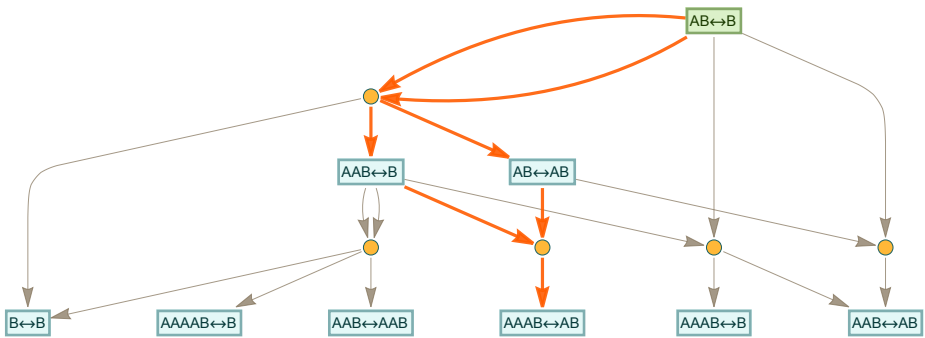

And here it is as a separate "proof graph"

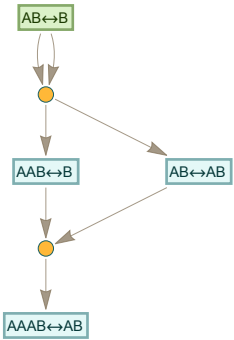

where each event takes two inputs—the "rule to be applied" and the "rule to apply to"—and the output is the derived (i.e. entailed or implied) new rule or rules.

If we run the accumulative system for another step, we get:

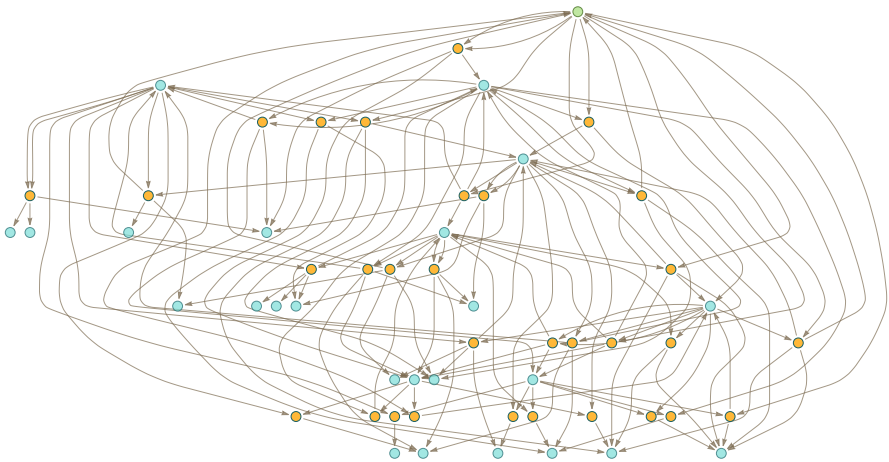

Now there are additional "theorems" that have been generated. An example is:

AAAB ↔ AAB

And now we can find a proof of this theorem:

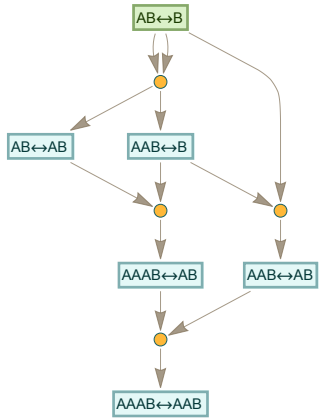

This proof exists as a subgraph of the token-event graph:

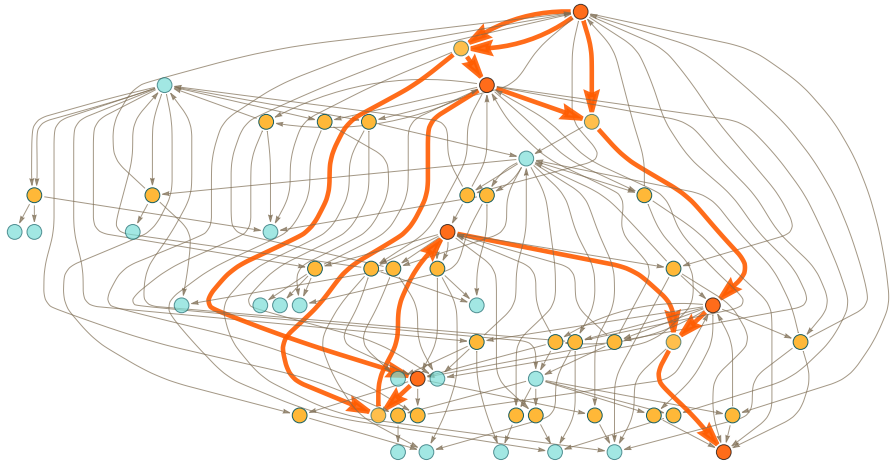

The proof just given has the fewest events—or "proof steps"—that can be used. But altogether there are 50 possible proofs, other examples being:

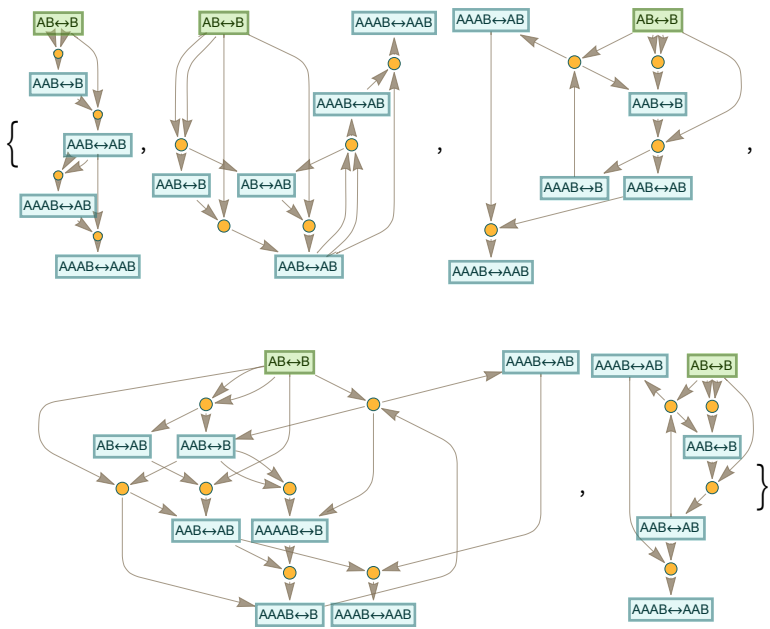

These correspond to the subgraphs:

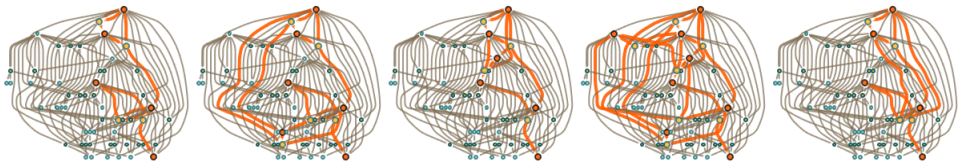

How much has the accumulative character of these token-event graphs contributed to the structure of these proofs? It's perfectly possible to find proofs that never use "intermediate lemmas" but always "go back to the original axiom" at every step. In this case examples are which all in effect require at least one more "sequential event" than our shortest proof using intermediate lemmas.

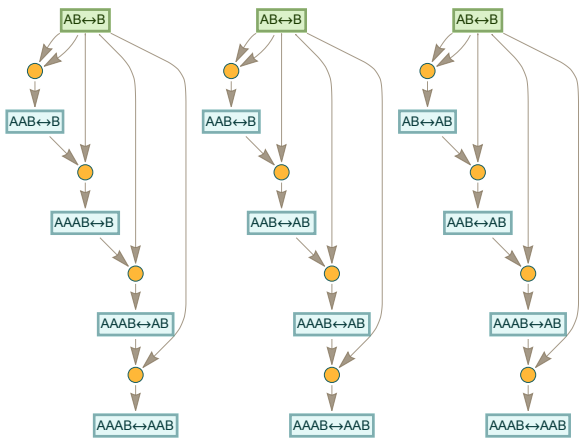

A slightly more dramatic example occurs for the theorem

AAAAAB ↔ AAAAB

where now without intermediate lemmas the shortest proof is

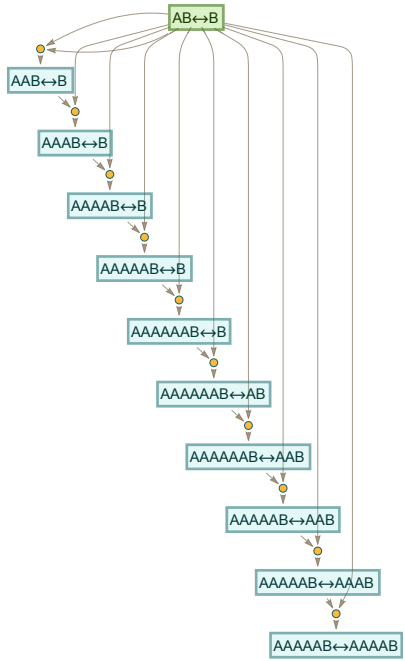

but with intermediate lemmas it becomes:

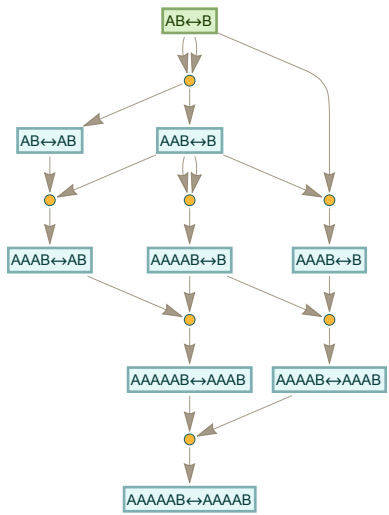

What we've done so far here is to generate a complete token-event graph for a certain number of steps, and then to see if we can find a proof in it for some particular statement. The proof is a subgraph of the "relevant part" of the full token-event graph. Often—in analogy to the simpler case of finding proofs of equivalences between expressions in a multiway graph—we'll call this subgraph a "proof path".

But in addition to just "finding a proof" in a fully constructed token-event graph, we can ask whether, given a statement, we can directly construct a proof for it. As discussed in the context of proofs in ordinary multiway graphs, computational irreducibility implies that in general there's no "shortcut" way to find a proof. In addition, for any statement, there may be no upper bound on the length of proof that will be required (or on the size or number of intermediate "lemmas" that will have to be used). And this, again, is the shadow of undecidability in our systems: that there can be statements whose provability may be arbitrarily difficult to determine.

## 12 | Beyond Substitution: Cosubstitution and Bisubstitution

In making our "metamodel" of mathematics we've been discussing the rewriting of expressions according to rules. But there's a subtle issue that we've so far avoided, that has to do with the fact that the expressions we're rewriting are often themselves patterns that stand for whole classes of expressions. And this turns out to allow for additional kinds of transformations that we'll call cosubstitution and bisubstitution.

Let's talk first about cosubstitution. Imagine we have the expression **f[a]**. The rule $a \to b$ would do a substitution for $a$ to give **f[b]**. But if we have the expression **f[c]** the rule $a \to b$ will do nothing.

Now imagine that we have the expression **f[x_]**. This stands for a whole class of expressions, including **f[a]**, **f[c]**, etc. For most of this class of expressions, the rule $a \to b$ will do nothing. But in the specific case of **f[a]**, it applies, and gives the result **f[b]**.

If our rule is **f[x_] → s** then this will apply as an ordinary substitution to **f[a]**, giving the result $s$. But if the rule is **f[b] → s** this will not apply as an ordinary substitution to **f[a]**. However, it can apply as a cosubstitution to **f[x_]** by picking out the specific case where **x_** stands for $b$, then using the rule to give $s$.

In general, the point is that ordinary substitution specializes patterns that appear in rules—while what one can think of as the "dual operation" of cosubstitution specializes patterns that appear in the expressions to which the rules are being applied. If one thinks of the rule that's being applied as like an operator, and the expression to which the rule is being applied as an operand, then in effect substitution is about making the operator fit the operand, and cosubstitution is about making the operand fit the operator.

It's important to realize that as soon as one's operating on expressions involving patterns, cosubstitution is not something "optional": it's something that one has to include if one is really going to interpret patterns—wherever they occur—as standing for classes of expressions.

When one's operating on a literal expression (without patterns) only substitution is ever possible, as in

f[x_] → g[x]:  f[a] →●→ g[a]

corresponding to this fragment of a token-event graph:

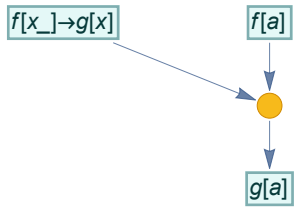

Let's say we have the rule **f[a] → s** (where **f[a]** is a literal expression). Operating on **f[b]** this rule will do nothing. But what if we apply the rule to **f[x_]**? Ordinary substitution still does nothing. But cosubstitution can do something. In fact, there are two different cosubstitutions that can be done in this case:

f[a] → s:  f[x_] →●→ s

f[a] → s:  f[x_] →●→ f[s]

What's going on here? In the first case, **f[x_]** has the "special case" **f[a]**, to which the rule applies ("by cosubstitution")—giving the result *s*. In the second case, however, it's *x_* on its own which has the special case **f[a]**, that gets transformed by the rule to *s*, giving the final cosubstitution result **f[s]**.

There's an additional wrinkle when the same pattern (such as *x_*) appears multiple times:

a → b:  f[x_, x_, x_] →●→ f[b, a, a]

a → b:  f[x_, x_, x_] →●→ f[a, b, a]

a → b:  f[x_, x_, x_] →●→ f[a, a, b]

In all cases, *x_* is matched to *a*. But which of the *x_*'s is actually replaced is different in each case.

Here's a slightly more complicated example:

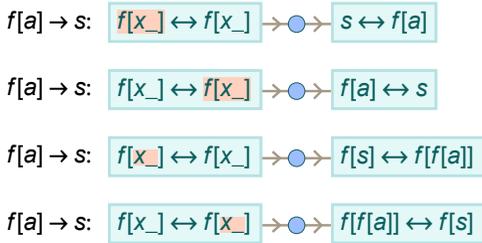

In ordinary substitution, replacements for patterns are in effect always made "locally", with each specific pattern separately being replaced by some expression. But in cosubstitution, a "special case" found for a pattern will get used throughout when the replacement is done.

Let's see how this all works in an accumulative axiomatic system. Consider the very simple rule:

$x\_ \circ y\_ \leftrightarrow y\_ \circ x\_$

One step of substitution gives the token-event graph (where we've canonicalized the names of pattern variables to $a\_$ and $b\_$):

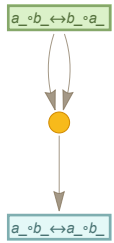

But one step of cosubstitution gives instead:

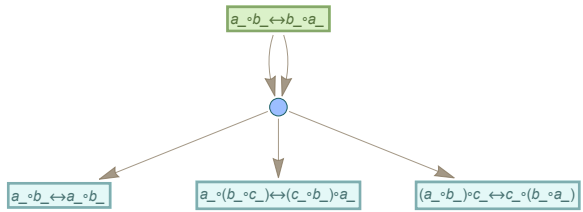

Here are the individual transformations that were made (with the rule at least nominally being applied only in one direction):

$x\_ \circ y\_ \to y \circ x$: $\boxed{a\_ \circ b\_ \leftrightarrow b\_ \circ a\_} \to \circ \to y\_ \circ x\_ \leftrightarrow y\_ \circ x\_$

$x\_ \circ y\_ \to y \circ x$: $a\_ \circ b\_ \leftrightarrow b\_ \circ a\_ \to \circ \to y\_ \circ x\_ \leftrightarrow y\_ \circ x\_$

$x\_ \circ y\_ \to y \circ x$: $a\_ \circ b\_ \leftrightarrow b\_ \circ a\_ \to \circ \to (y\_ \circ x\_) \circ b\_ \leftrightarrow b\_ \circ (x\_ \circ y\_)$

$x\_ \circ y\_ \to y \circ x$: $a\_ \circ b\_ \leftrightarrow b\_ \circ a\_ \to \circ \to a\_ \circ (y\_ \circ x\_) \leftrightarrow (x\_ \circ y\_) \circ a\_$

$x\_ \circ y\_ \to y \circ x$: $a\_ \circ b\_ \leftrightarrow b\_ \circ a\_ \to \circ \to a\_ \circ (x\_ \circ y\_) \leftrightarrow (y\_ \circ x\_) \circ a\_$

$x\_ \circ y\_ \to y \circ x$: $a\_ \circ b\_ \leftrightarrow b\_ \circ a\_ \to \circ \to (x\_ \circ y\_) \circ b\_ \leftrightarrow b\_ \circ (y\_ \circ x\_)$

The token-event graph above is then obtained by canonicalizing variables, and combining identical expressions (though for clarity we don't merge rules of the form $a \leftrightarrow b$ and $b \leftrightarrow a$).

If we go another step with this particular rule using only substitution, there are additional events (i.e. transformations) but no new theorems produced:

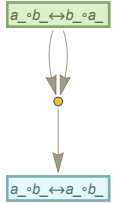

Cosubstitution, however, produces another 27 theorems

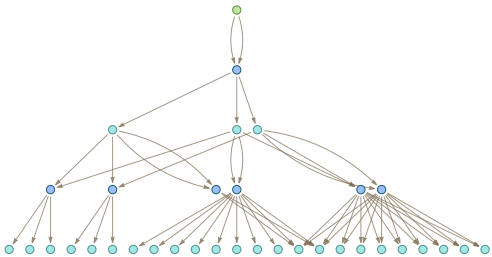

or altogether

| a∘b ↔ b∘a | a∘b ↔ a∘b | a∘(b∘c) ↔ (c∘b)∘a | (a∘b)∘c ↔ c∘(b∘a) |
|---|---|---|---|
| a∘(b∘c) ↔ a∘(b∘c) | a∘(b∘c) ↔ (b∘c)∘a | a∘(b∘(c∘d)) ↔ (b∘(d∘c))∘a | a∘(b∘(c∘d)) ↔ ((c∘d)∘b)∘a |
| a∘(b∘(c∘d)) ↔ ((e∘c)∘b)∘a | a∘(b∘(c∘(d∘e))) ↔ (((e∘d)∘c)∘b)∘a | a∘(b∘((c∘d)∘e)) ↔ ((e∘(d∘c))∘b)∘a | a∘((b∘c)∘d) ↔ (d∘(b∘c))∘a |
| a∘((b∘c)∘d) ↔ (d∘(c∘b))∘a | a∘((b∘c)∘d) ↔ ((c∘b)∘d)∘a | a∘((b∘(c∘d))∘e) ↔ (e∘((d∘c)∘b))∘a | a∘(((b∘c)∘d)∘e) ↔ (e∘(d∘(c∘b)))∘a |
| (a∘b)∘c ↔ c∘(a∘b) | (a∘b)∘c ↔ (a∘b)∘c | (a∘b)∘(c∘d) ↔ (c∘d)∘(b∘a) | (a∘b)∘(c∘d) ↔ (d∘c)∘(a∘b) |
| (a∘b)∘(c∘d) ↔ (d∘c)∘(b∘a) | (a∘b)∘(c∘(d∘e)) ↔ ((e∘d)∘c)∘(b∘a) | (a∘b)∘((c∘d)∘e) ↔ (e∘(d∘c))∘(b∘a) | (a∘(b∘c))∘d ↔ d∘(a∘(c∘b)) |
| (a∘(b∘c))∘d ↔ d∘((b∘c)∘a) | (a∘(b∘c))∘d ↔ d∘((c∘b)∘a) | (a∘(b∘c))∘(d∘e) ↔ (e∘d)∘((c∘b)∘a) | (a∘(b∘(c∘d)))∘e ↔ e∘(((d∘c)∘b)∘a) |
| (a∘((b∘c)∘d))∘e ↔ e∘((d∘(c∘b))∘a) | ((a∘b)∘c)∘d ↔ d∘(c∘(a∘b)) | ((a∘b)∘c)∘d ↔ d∘(c∘(b∘a)) | ((a∘b)∘c)∘d ↔ d∘((b∘a)∘c) |
| ((a∘b)∘c)∘(d∘e) ↔ (e∘d)∘(c∘(b∘a)) | ((a∘(b∘c))∘d)∘e ↔ e∘(d∘((c∘b)∘a)) | (((a∘b)∘c)∘d)∘e ↔ e∘(d∘(c∘(b∘a))) | |

or as trees:

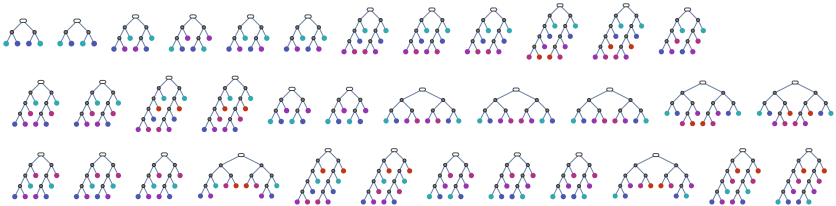

We've now seen examples of both substitution and cosubstitution in action. But in our metamodel for mathematics we're ultimately dealing not with each of these individually, but rather with the "symmetric" concept of bisubstitution, in which both substitution and cosubstitution can be mixed together, and applied even to parts of the same expression.

In the particular case of $a\_\circ b\_ \leftrightarrow b\_\circ a\_$, bisubstitution adds nothing beyond cosubstitution. But often it does. Consider the rule:

$x\_ \circ a \to b$

Here's the result of applying this to three different expressions using substitution, cosubstitution and bisubstitution (where we consider only matches for "whole ∘ expressions", not subparts):

$x\_ \circ b \to c$: {a∘b, x_∘y_, a∘y_} → {c, x_∘y_, a∘y_}

$x\_ \circ b \to c$: {a∘b, x_∘y_, a∘y_} → {a∘b, c, a∘b}

$x\_ \circ b \to c$: {a∘b, x_∘y_, a∘y_} → {c, x_∘y_, a∘y_}

$x\_ \circ b \to c$: {a∘b, x_∘y_, a∘y_} → {a∘b, c, a∘b}

$x\_ \circ b \to c$: {a∘b, x_∘y_, a∘y_} → {a∘b, x_∘b, c}

Cosubstitution very often yields substantially more transformations than substitution—bisubstitution then yielding modestly more than cosubstitution. For example, for the axiom system

$$x\_ \circ y\_ \leftrightarrow (y\_ \circ x\_) \circ y\_$$

the number of theorems derived after 1 and 2 steps is given by:

|  | step 1 | step 2 |
|---|---|---|
| ● substitution | 5 | 24 |
| ● cosubstitution | 14 | 1630 |
| ● bisubstitution | 14 | 1885 |

In some cases there are theorems that can be produced by full bisubstitution, but not—even after any number of steps—by substitution or cosubstitution alone. However, it is also common to find that theorems can in principle be produced by substitution alone, but that this just takes more steps (and sometimes vastly more) than when full bisubstitution is used. (It's worth noting, however, that the notion of "how many steps" it takes to "reach" a given theorem depends on the foliation one chooses to use in the token-event graph.)

The various forms of substitution that we've discussed here represent different ways in which one theorem can entail others. But our overall metamodel of mathematics—based as it is purely on the structure of symbolic expressions and patterns—implies that bisubstitution covers all entailments that are possible.

In the history of metamathematics and mathematical logic, a whole variety of "laws of inference" or "methods of entailment" have been considered. But with the modern view of symbolic expressions and patterns (as used, for example, in the Wolfram Language), bisubstitution emerges as the fundamental form of entailment, with other forms of entailment corresponding to the use of particular types of expressions or the addition of further elements to the pure substitutions we've used here.

It should be noted, however, that when it comes to the ruliad different kinds of entailments correspond merely to different foliations—with the form of entailment that we're using representing just a particularly straightforward case.

The concept of bisubstitution has arisen in the theory of term rewriting, as well as in automated theorem proving (where it is often viewed as a particular "strategy", and called "paramodulation"). In term rewriting, bisubstitution is closely related to the concept of unification—which essentially asks what assignment of values to pattern variables is needed in order to make different subterms of an expression be identical.

# 13 | Some First Metamathematical Phenomenology

Now that we've finished describing the many technical issues involved in constructing our metamodel of mathematics, we can start looking at its consequences. We discussed above how multiway graphs formed from expressions can be used to define a branchial graph that represents a kind of "metamathematical space". We can now use a similar approach to set up a metamathematical space for our full metamodel of the "progressive accumulation" of mathematical statements.

Let's start by ignoring cosubstitution and bisubstitution and considering only the process of substitution—and beginning with the axiom:

$a\_ \circ b\_ \leftrightarrow b\_$

Doing accumulative evolution from this axiom we get the token-event graph

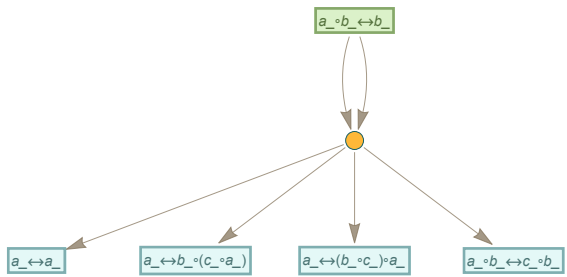

or after 2 steps:

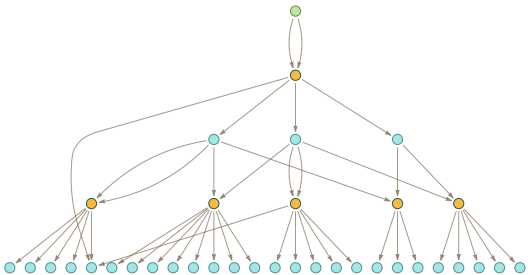

From this we can derive an "effective multiway graph" by directly connecting all input and output tokens involved in each event:

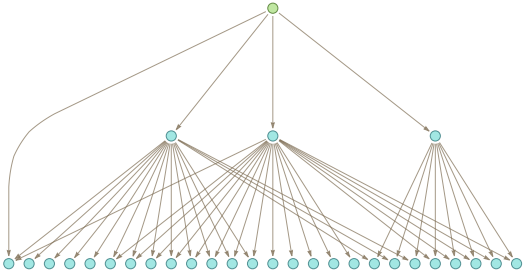

And then we can produce a branchial graph, which in effect yields an approximation to the "metamathematical space" generated by our axiom:

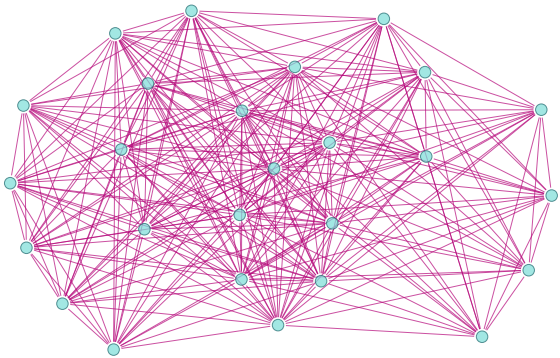

Showing the statements produced in the form of trees we get (with the top node representing ↔):

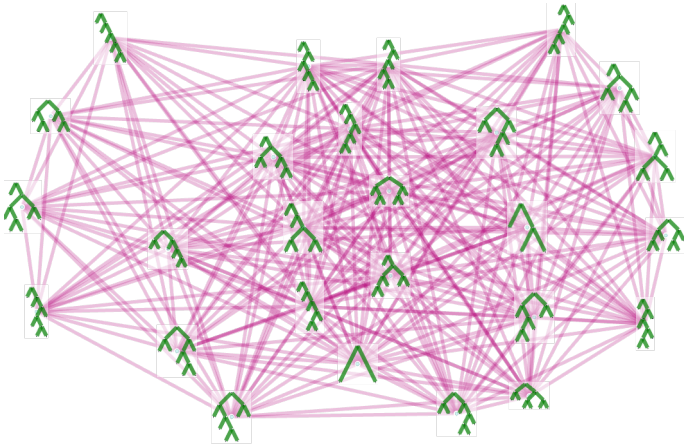

If we do the same thing with full bisubstitution, then even after one step we get a slightly larger token-event graph:

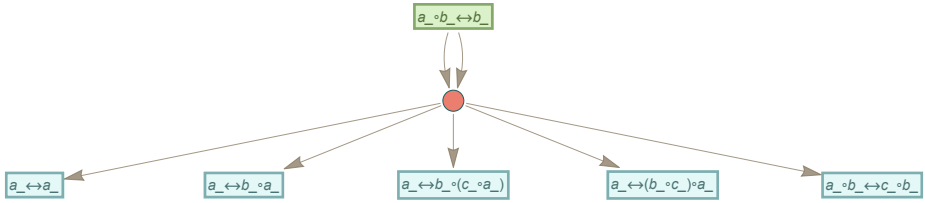

After two steps, we get

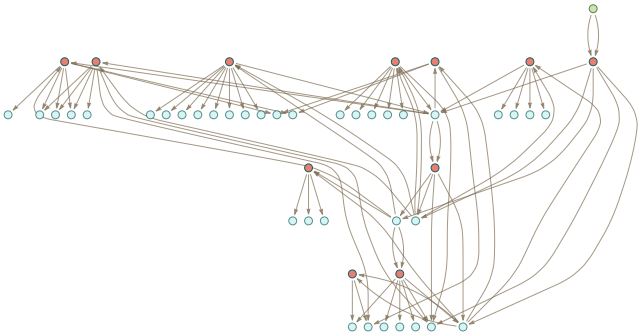

which contains 46 statements, compared to 42 if only substitution is used. The corresponding branchial graph is:

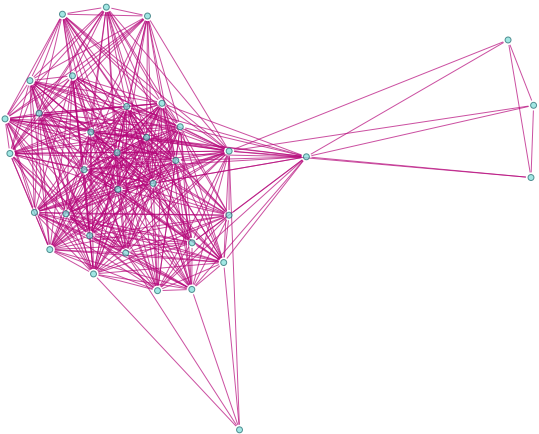

The adjacency matrices for the substitution and bisubstitution cases are then

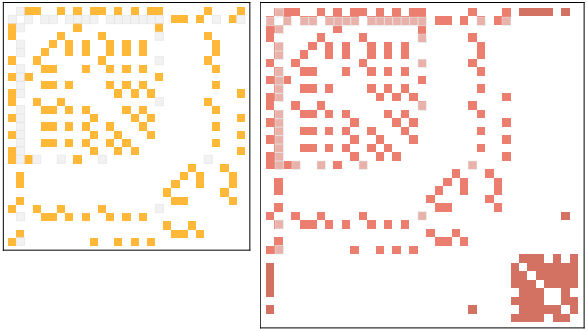

which have 80% and 85% respectively of the number of edges in complete graphs of these sizes.

Branchial graphs are usually quite dense, but they nevertheless do show definite structure. Here are some results after 2 steps:

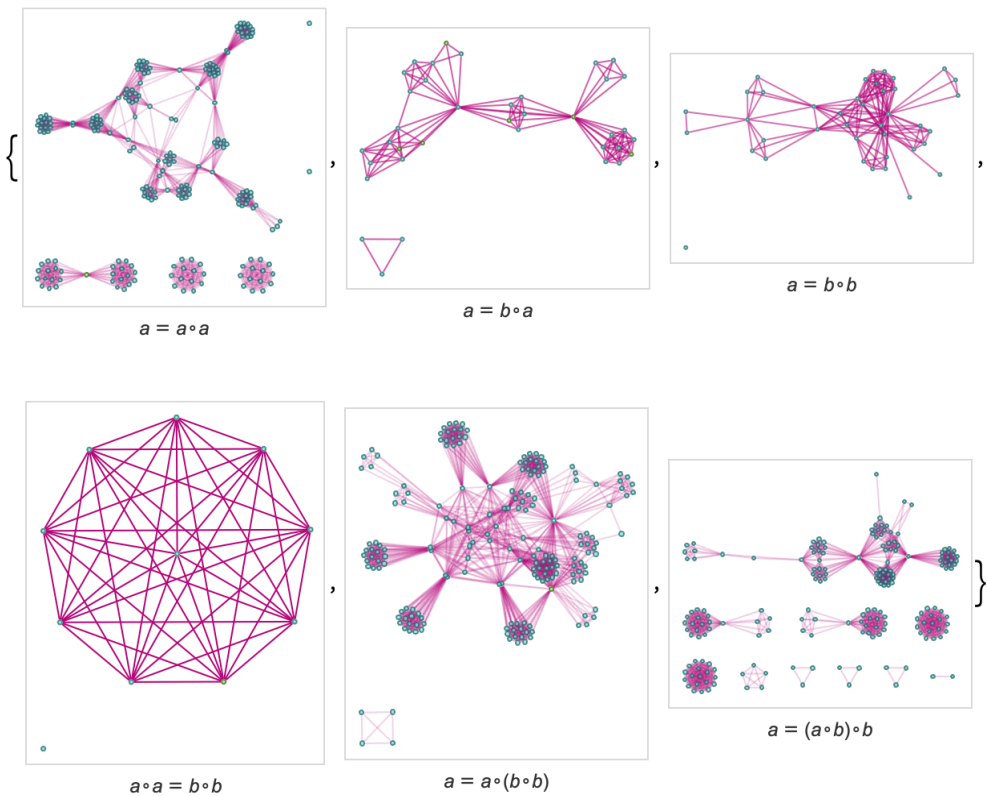

## 14 | Relations to Automated Theorem Proving

We've discussed at some length what happens if we start from axioms and then build up an "entailment cone" of all statements that can be derived from them. But in the actual practice of mathematics people often want to just look at particular target statements, and see if they can be derived (i.e. proved) from the axioms.

But what can we say "in bulk" about this process? The best source of potential examples we have right now come from the practice of automated theorem proving—as for example implemented in the Wolfram Language function **FindEquationalProof**. As a simple example of how this works, consider the axiom

$x\_ \circ y\_ = (y\_ \circ x\_) \circ y\_$

and the theorem:

$(a \circ ((b \circ a) \circ (a \circ b))) \circ a = b \circ a$

Automated theorem proving (based on **FindEquationalProof**) finds the following proof of this theorem:

$(a \circ ((b \circ a) \circ (a \circ b))) \circ a \longrightarrow ((b \circ a) \circ (a \circ b)) \circ a \longrightarrow (((a \circ b) \circ a) \circ (a \circ b)) \circ a \longrightarrow (a \circ (a \circ b)) \circ a \longrightarrow (a \circ b) \circ a \longrightarrow b \circ a$

Needless to say, this isn't the only possible proof. And in this very simple case, we can construct the full entailment cone—and determine that there aren't any shorter proofs, though there are two more of the same length:

$(a \circ ((b \circ a) \circ (a \circ b))) \circ a \longrightarrow (a \circ (((a \circ b) \circ a) \circ (a \circ b))) \circ a \longrightarrow (((a \circ b) \circ a) \circ (a \circ b)) \circ a \longrightarrow (a \circ (a \circ b)) \circ a \longrightarrow (a \circ b) \circ a \longrightarrow b \circ a$

$(a \circ ((b \circ a) \circ (a \circ b))) \circ a \longrightarrow (a \circ (((a \circ b) \circ a) \circ (a \circ b))) \circ a \longrightarrow (a \circ (a \circ (a \circ b))) \circ a \longrightarrow (a \circ (a \circ b)) \circ a \longrightarrow (a \circ b) \circ a \longrightarrow b \circ a$

All three of these proofs can be seen as paths in the entailment cone:

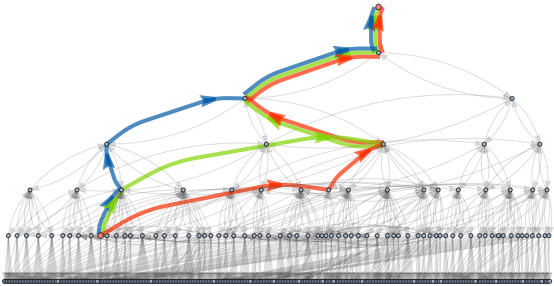

How "complicated" are these proofs? In addition to their lengths, we can for example ask how big the successive intermediate expressions they involve become, where here we are including not only the proofs already shown, but also some longer ones as well:

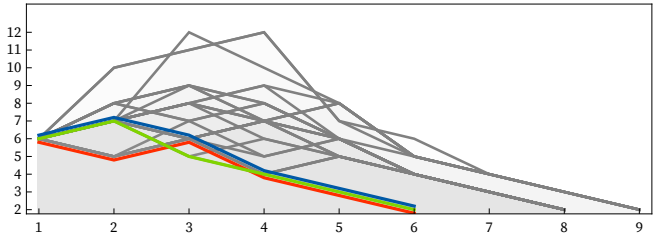

In the setup we're using here, we can find a proof of *lhs* = *rhs* by starting with *lhs*, building up an entailment cone, and seeing whether there's any path in it that reaches *rhs*. In general there's no upper bound on how far one will have to go to find such a path—or how big the intermediate expressions may need to get.

One can imagine all kinds of optimizations, for example where one looks at multistep consequences of the original axioms, and treats these as "lemmas" that we can "add as axioms" to provide new rules that jump multiple steps on a path at a time. Needless to say, there are lots of tradeoffs in doing this. (Is it worth the memory to store the lemmas? Might we "jump" past our target? etc.)

But typical actual automated theorem provers tend to work in a way that is much closer to our accumulative rewriting systems—in which the "raw material" on which one operates is statements rather than expressions.

Once again, we can in principle always construct a whole entailment cone, and then look to see whether a particular statement occurs there. But then to give a proof of that statement it's sufficient to find the subgraph of the entailment cone that leads to that statement. For example, starting with the axiom

$$a \circ b = (b \circ a) \circ b$$

we get the entailment cone (shown here as a token-event graph, and dropping _'s):

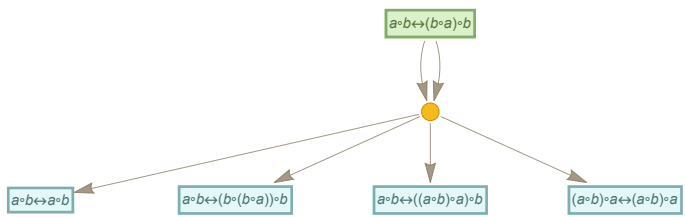

After 2 steps the statement

$$(a \circ b) \circ a = (a \circ (a \circ (a \circ b))) \circ a$$

shows up in this entailment cone

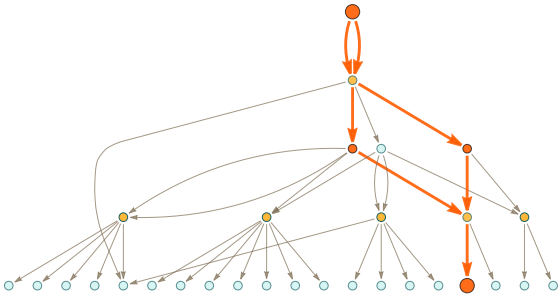

where we're indicating the subgraph that leads from the original axiom to this statement. Extracting this subgraph we get

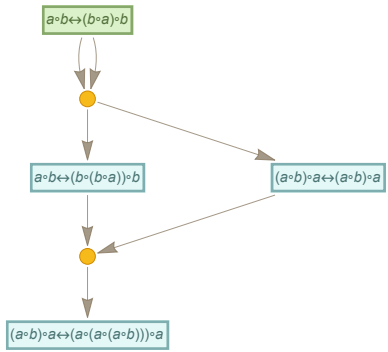

which we can view as a proof of the statement within this axiom system.

But now let's use traditional automated theorem proving (in the form of **FindEquationalProof**) to get a proof of this same statement. Here's what we get:

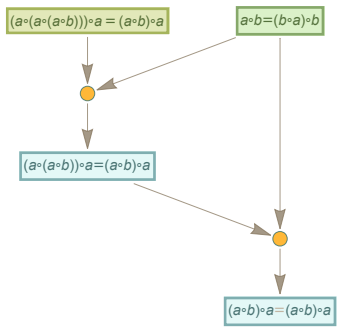

This is again a token-event graph, but its structure is slightly different from the one we "fished out of" the entailment cone. Instead of starting from the axiom and "progressively deriving" our statement we start from both the statement and the axiom and then show that together they lead "merely via substitution" to a statement of the form $x=x$, which we can take as an "obviously derivable tautology".

Sometimes the minimal "direct proof" found from the entailment cone can be considerably simpler than the one found by automated theorem proving. For example, for the statement

$$a \circ b = (((b \circ (b \circ a)) \circ b) \circ (b \circ a)) \circ b$$

the minimal direct proof is

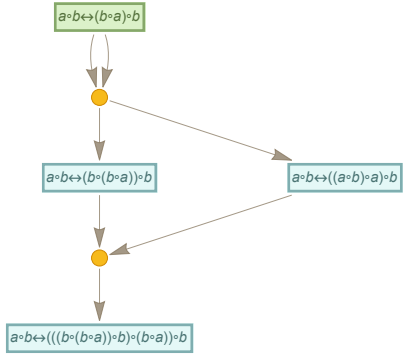

while the one found by **FindEquationalProof** is:

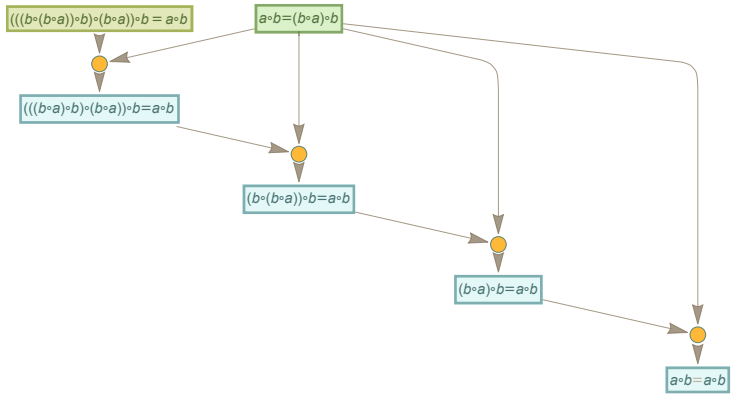

But the great advantage of automated theorem proving is that it can "directedly" search for proofs instead of just "fishing them out of" the entailment cone that contains all possible exhaustively generated proofs. To use automated theorem proving you have to "know where you want to go"—and in particular identify the theorem you want to prove.

Consider the axiom

$$(b \circ a) \circ ((c \circ a) \circ c) = a$$

and the statement:

$$a \circ (a \circ (b \circ b)) = b \circ b$$

This statement doesn't show up in the first few steps of the entailment cone for the axiom, even though millions of other theorems do. But automated theorem proving finds a proof of it—and rearranging the "prove-a-tautology proof" so that we just have to feed in a tautology somewhere in the proof, we get:

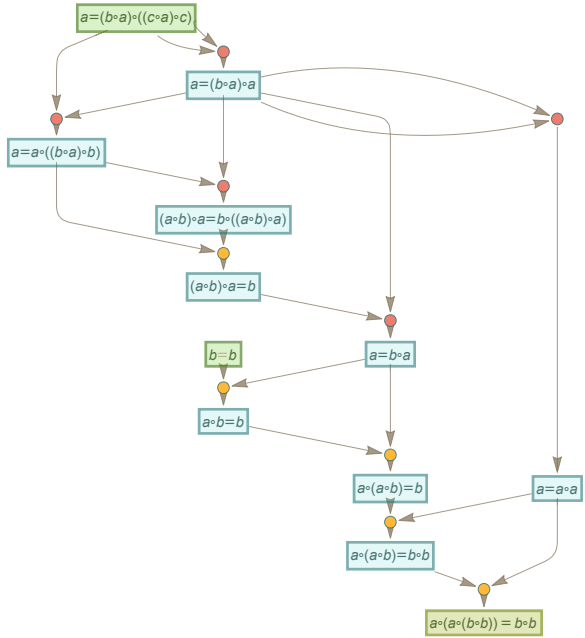

The model-theoretic methods we'll discuss a little later allow one effectively to "guess" theorems that might be derivable from a given axiom system. So, for example, for the axiom system

$$(b \circ a) \circ (c \circ ((c \circ a) \circ c)) = a$$

here's a "guess" at a theorem

$$a \circ ((a \circ b) \circ a) = (b \circ b) \circ (a \circ b)$$

and here's a representation of its proof found by automated theorem proving—where now the length of an intermediate "lemma" is indicated by the size of the corresponding node

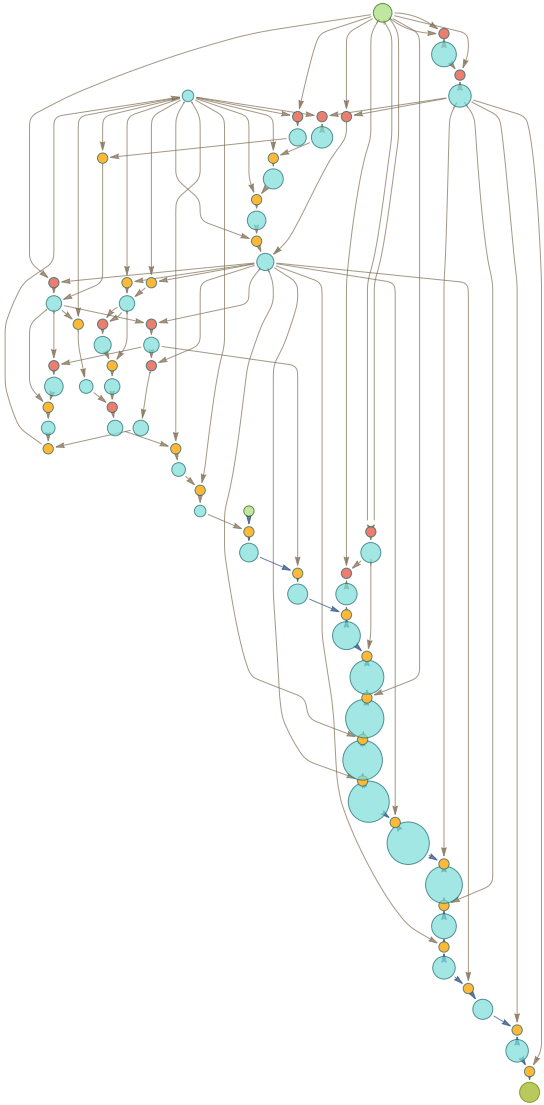

and in this case the longest intermediate lemma is of size 67 and is:

$(b \circ b) \circ (a \circ b) = ((b \circ ((b \circ b) \circ b)) \circ b) \circ$

$((b \circ ((b \circ ((b \circ b) \circ b)) \circ b)) \circ (((b \circ ((b \circ b) \circ b)) \circ b) \circ (((b \circ b) \circ b) \circ ((b \circ ((b \circ b) \circ b)) \circ (((b \circ b) \circ b) \circ (b \circ ((b \circ b) \circ b))))))$

In principle it's possible to rearrange token-event graphs generated by automated theorem proving to have the same structure as the ones we get directly from the entailment cone—with axioms at the beginning and the theorem being proved at the end. But typical strategies for automated theorem proving don't naturally produce such graphs. In principle automated theorem proving could work by directly searching for a "path" that leads to the theorem one's trying to prove. But usually it's much easier instead to have as the "target" a simple tautology.

At least conceptually automated theorem proving must still try to "navigate" through the full token-event graph that makes up the entailment cone. And the main issue in doing this is that there are many places where one does not know "which branch to take". But here there's a crucial—if at first surprising—fact: at least so long as one is using full bisubstitution it ultimately doesn't matter which branch one takes; there'll always be a way to "merge back" to any other branch.

This is a consequence of the fact that the accumulative systems we're using automatically have the property of confluence which says that every branch is accompanied by a subsequent merge. There's an almost trivial way in which this is true by virtue of the fact that for every edge the system also includes the reverse of that edge. But there's a more substantial reason as well: that given any two statements on two different branches, there's always a way to combine them using a bisubstitution to get a single statement.

In our Physics Project, the concept of causal invariance—which effectively generalizes confluence—is an important one, that leads among other things to ideas like relativistic invariance. Later on we'll discuss the idea that "regardless of what order you prove theorems in, you'll always get the same math", and its relationship to causal invariance and to the notion of relativity in metamathematics. But for now the importance of confluence is that it has the potential to simplify automated theorem proving—because in effect it says one can never ultimately "make a wrong turn" in getting to a particular theorem, or, alternatively, that if one keeps going long enough every path one might take will eventually be able to reach every theorem.

And indeed this is exactly how things work in the full entailment cone. But the challenge in automated theorem proving is to generate only a tiny part of the entailment cone, yet still "get to" the theorem we want. And in doing this we have to carefully choose which "branches" we should try to merge using bisubstitution events. In automated theorem proving these bisubstitution events are typically called "critical pair lemmas", and there are a variety of strategies for defining an order in which critical pair lemmas should be tried.

It's worth pointing out that there's absolutely no guarantee that such procedures will find the shortest proof of any given theorem (or in fact that they'll find a proof at all with a given amount of computational effort). One can imagine "higher-order proofs" in which one attempts to transform not just statements of the form *lhs*=*rhs*, but full proofs (say represented as token-event graphs). And one can imagine using such transformations to try to simplify proofs.

A general feature of the proofs we've been showing is that they are accumulative, in the sense they continually introduce lemmas which are then reused. But in principle any proof can be "unrolled" into one that just repeatedly uses the original axioms (and in fact, purely by substitution)—and never introduces other lemmas. The necessary "cut elimination" can effectively be done by always recreating each lemma from the axioms whenever it's needed—a process which can become exponentially complex.

As an example, from the axiom above we can generate the proof

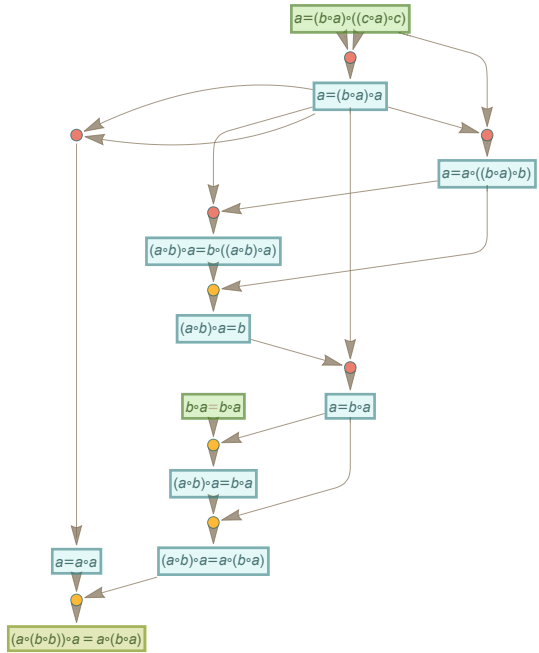

where for example the first lemma at the top is reused in four events. But now by cut elimination we can "unroll" this whole proof into a "straight-line" sequence of substitutions on expressions done just using the original axiom

| |
|---|
| (a∘(b∘b))∘a |
| (a∘(((b∘b)∘((((c∘b)∘c)∘b)∘((c∘b)∘c)))∘b))∘a |
| (a∘(((b∘b)∘b)∘b))∘a |
| (a∘(((b∘b)∘b)∘(((((d∘b)∘d)∘b)∘((d∘b)∘d))))∘a |
| (a∘b)∘a |
| ((((a∘a)∘a)∘((b∘a)∘b))∘b)∘a |
| ((((a∘a)∘((((e∘a)∘e)∘a)∘((e∘a)∘e)))∘((b∘a)∘b))∘b)∘a |
| ((a∘((b∘a)∘b))∘b)∘a |
| (((((e∘a)∘a)∘((b∘a)∘b))∘((b∘a)∘b))∘b)∘a |
| (((((e∘a)∘(((((f∘a)∘f)∘a)∘((f∘a)∘f)))∘((b∘a)∘b))∘((b∘a)∘b))∘b)∘a |
| (((a∘((b∘a)∘b))∘((b∘a)∘b))∘b)∘a |
| (((a∘((b∘a)∘b))∘((((c∘((b∘a)∘b))∘c)∘((b∘a)∘b))∘((c∘((b∘a)∘b))∘c)))∘b)∘a |
| (((b∘a)∘b)∘b)∘a |
| (((b∘a)∘b)∘(((((d∘b)∘d)∘b)∘((d∘b)∘d)))∘a |
| b∘a |
| (((b∘a)∘a)∘(b∘a))∘((((d∘(b∘a))∘d)∘(b∘a))∘((d∘(b∘a))∘d)) |
| (((b∘a)∘a)∘(b∘a))∘(b∘a) |
| ((a∘(((b∘a)∘a)∘(b∘a)))∘((((c∘(((b∘a)∘a)∘(b∘a)))∘c)∘(((b∘a)∘a)∘(b∘a)))∘((c∘(((b∘a)∘a)∘(b∘a)))∘c)))∘(b∘a) |
| ((a∘(((b∘a)∘a)∘(b∘a)))∘(((b∘a)∘a)∘(b∘a)))∘(b∘a) |
| (((((e∘a)∘(((((f∘a)∘f)∘a)∘((f∘a)∘f)))∘(((b∘a)∘a)∘(b∘a)))∘(((b∘a)∘a)∘(b∘a)))∘(b∘a) |
| (((((e∘a)∘a)∘(((b∘a)∘a)∘(b∘a)))∘(((b∘a)∘a)∘(b∘a)))∘(b∘a) |
| (a∘(((b∘a)∘a)∘(b∘a)))∘(b∘a) |
| (((a∘a)∘((((e∘a)∘e)∘a)∘((e∘a)∘e)))∘(((b∘a)∘a)∘(b∘a)))∘(b∘a) |
| (((a∘a)∘a)∘(((b∘a)∘a)∘(b∘a)))∘(b∘a) |
| a∘(b∘a) |

and we see that our final theorem is the statement that the first expression in the sequence is equivalent under the axiom to the last one.

As is fairly evident in this example, a feature of automated theorem proving is that its result tends to be very "non-human". Yes, it can provide incontrovertible evidence that a theorem is valid. But that evidence is typically far away from being any kind of "narrative" suitable for human consumption. In the analogy to molecular dynamics, an automated proof gives detailed "turn-by-turn instructions" that show how a molecule can reach a certain place in a gas. Typical "human-style" mathematics, on the other hand, operates on a higher level, analogous to talking about overall motion in a fluid. And a core part of what's achieved by our physicalization of metamathematics is understanding why it's possible for mathematical observers like us to perceive mathematics as operating at this higher level.

## 15 | Axiom Systems of Present-Day Mathematics

The axiom systems we've been talking about so far were chosen largely for their axiomatic simplicity. But what happens if we consider axiom systems that are used in practice in present-day mathematics?

The simplest common example are the axioms (actually, a single axiom) of semigroup theory, stated in our notation as:

$$a\_ \circ (b\_ \circ c\_) \leftrightarrow (a\_ \circ b\_) \circ c\_$$

Using only substitution, all we ever get after any number of steps is the token-event graph (i.e. "entailment cone"):

But with bisubstitution, even after one step we already get the entailment cone which contains such theorems as:

$a \circ ((b \circ c) \circ d) = (a \circ b) \circ (c \circ d)$

After 2 steps, the entailment cone becomes

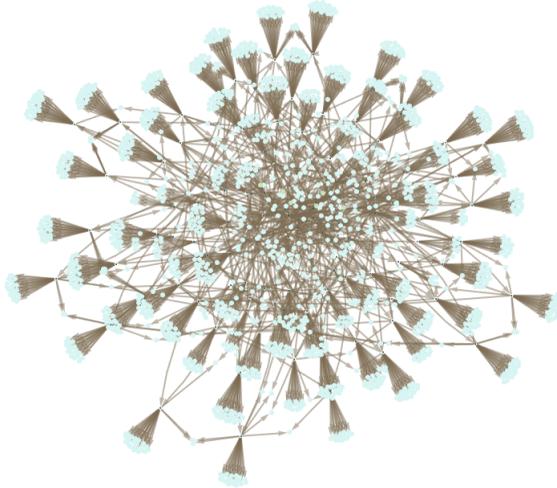

which contains 1617 theorems such as

| |
|---|
| $a \circ (b \circ (c \circ (d \circ e))) \circ f) = (a \circ ((b \circ (c \circ d)) \circ e)) \circ f$ |
| $(a \circ (b \circ c)) \circ (((d \circ e) \circ f) \circ g) = (a \circ (b \circ (c \circ (d \circ (e \circ f))))) \circ g$ |
| $a \circ (b \circ ((((c \circ d) \circ e) \circ f) \circ g)) = ((a \circ b) \circ (c \circ (d \circ (e \circ f)))) \circ g$ |
| $(a \circ (b \circ c)) \circ (d \circ ((e \circ (f \circ g)) \circ (h \circ i))) = (((a \circ b) \circ c) \circ d) \circ ((((e \circ f) \circ g) \circ h) \circ i)$ |
| $((a \circ ((b \circ ((c \circ d) \circ e)) \circ f)) \circ g) \circ h = ((a \circ (((b \circ c) \circ d) \circ e)) \circ f) \circ (g \circ h)$ |

with sizes distributed as follows:

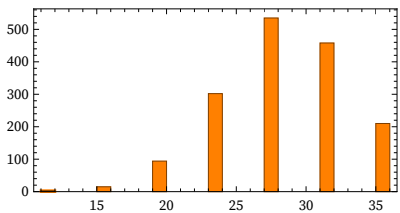

Looking at these theorems we can see that—in fact by construction—they are all just statements of the associativity of ∘. Or, put another way, they state that under this axiom all expression trees that have the same sequence of leaves are equivalent.

What about group theory? The standard axioms can be written

| |
|---|
| $a \circ (b \circ c) = (a \circ b) \circ c$ |
| $a \circ \diamond = a$ |
| $a \circ \overline{a} = \diamond$ |

where ∘ is interpreted as the binary group multiplication operation, overbar as the unary inverse operation, and 1 as the constant identity element (or, equivalently, zero-argument function).

One step of substitution already gives:

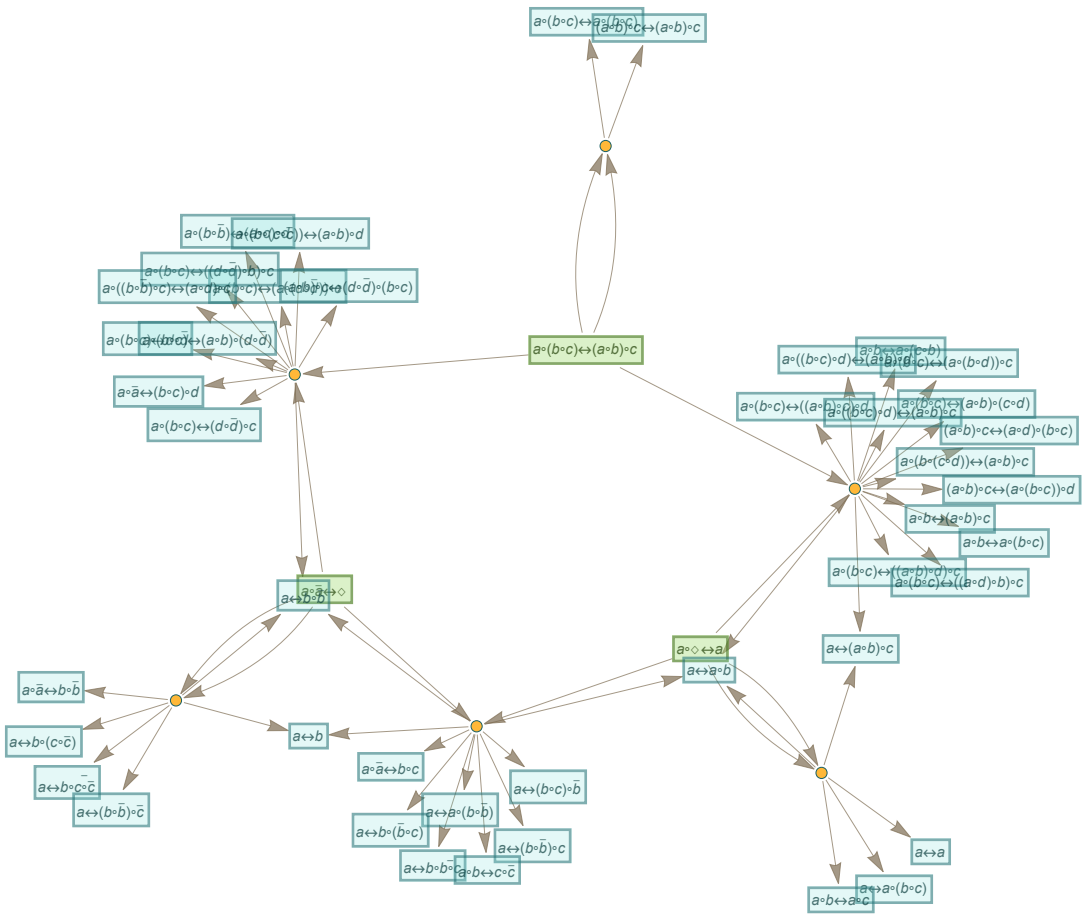

It's notable that in this picture one can already see "different kinds of theorems" ending up in different "metamathematical locations". One also sees some "obvious" tautological "theorems", like *a=a* and 1=1.

If we use full bisubstitution, we get 56 rather than 27 theorems, and many of the theorems are more complicated:

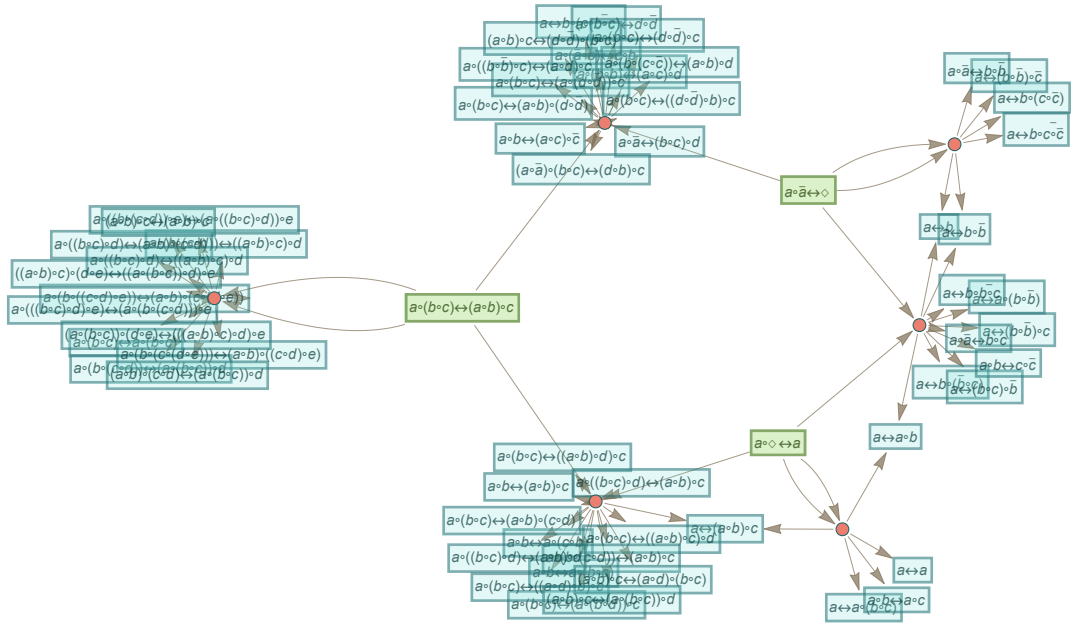

After 2 steps of pure substitution, the entailment cone in this case becomes

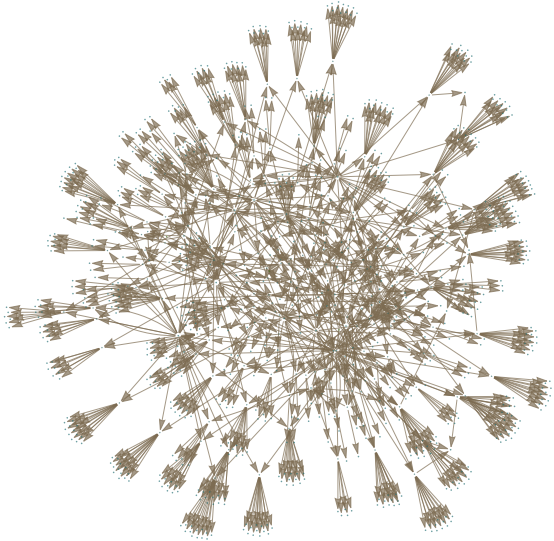

which includes 792 theorems with sizes distributed according to:

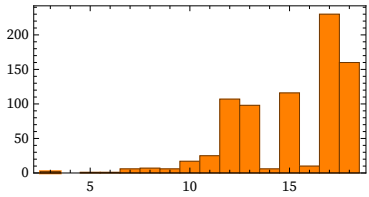

But among all these theorems, do straightforward "textbook theorems" appear, like?

| | |
|---|---|
| inverse of inverse | $\bar{\bar{a}}=a$ |
| left identity | $\diamond \circ a = a$ |
| left inverse | $\bar{a} \circ a = \diamond$ |
| inverse of composite | $\overline{a \circ b} = \bar{b} \circ \bar{a}$ |

The answer is no. It's inevitable that in the end all such theorems must appear in the entailment cone. But it turns out that it takes quite a few steps. And indeed with automated theorem proving we can find "paths" that can be taken to prove these theorems—involving significantly more than two steps:

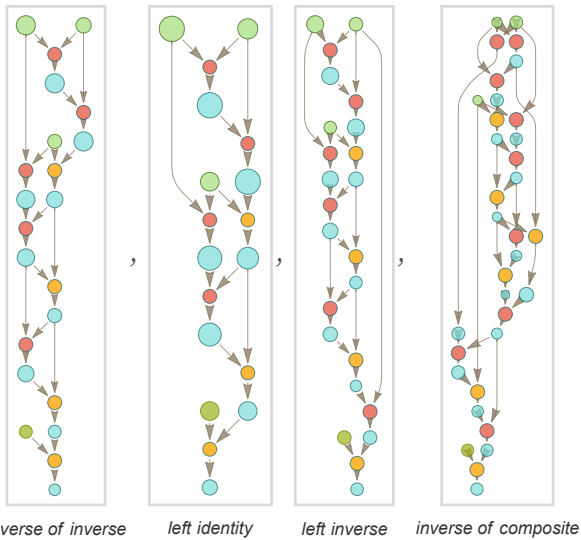

inverse of inverse, left identity, left inverse, inverse of composite

So how about logic, or, more specifically Boolean algebra? A typical textbook axiom system for this (represented in terms of **And** ∧, **Or** ∨ and **Not** $\bar{\Box}$) is:

| | |
|---|---|
| $a \vee b = b \vee a$ | $a \wedge b = b \wedge a$ |
| $a \vee (b \wedge \bar{b}) = a$ | $a \wedge (b \vee \bar{b}) = a$ |
| $a \vee (b \wedge c) = (a \vee b) \wedge (a \vee c)$ | $a \wedge (b \vee c) = (a \wedge b) \vee (a \wedge c)$ |

After one step of substitution from these axioms we get

or in our more usual rendering:

So what happens here with "named textbook theorems" (excluding commutativity and distributivity, which already appear in the particular axioms we're using)?

| | | | |
|---|---|---|---|
| idempotence of And | $a = a \wedge a$ | de Morgan law | $\overline{a \vee b} = \overline{a} \wedge \overline{b}$ |
| idempotence of Or | $a = a \vee a$ | de Morgan law | $\overline{a \wedge b} = \overline{a} \vee \overline{b}$ |
| law of double negation | $a = \overline{\overline{a}}$ | absorption law | $a = a \wedge (a \vee b)$ |
| law of noncontradiction | $\overline{a} \wedge a = \overline{b} \wedge b$ | absorption law | $a = a \vee (a \wedge b)$ |
| law of excluded middle | $\overline{a} \vee a = \overline{b} \vee b$ | associativity of And | $(a \wedge b) \wedge c = a \wedge (b \wedge c)$ |
| | | associativity of Or | $(a \vee b) \vee c = a \vee (b \vee c)$ |

Once again none of these appear in the first step of the entailment cone. But at step 2 with full bisubstitution the idempotence laws show up

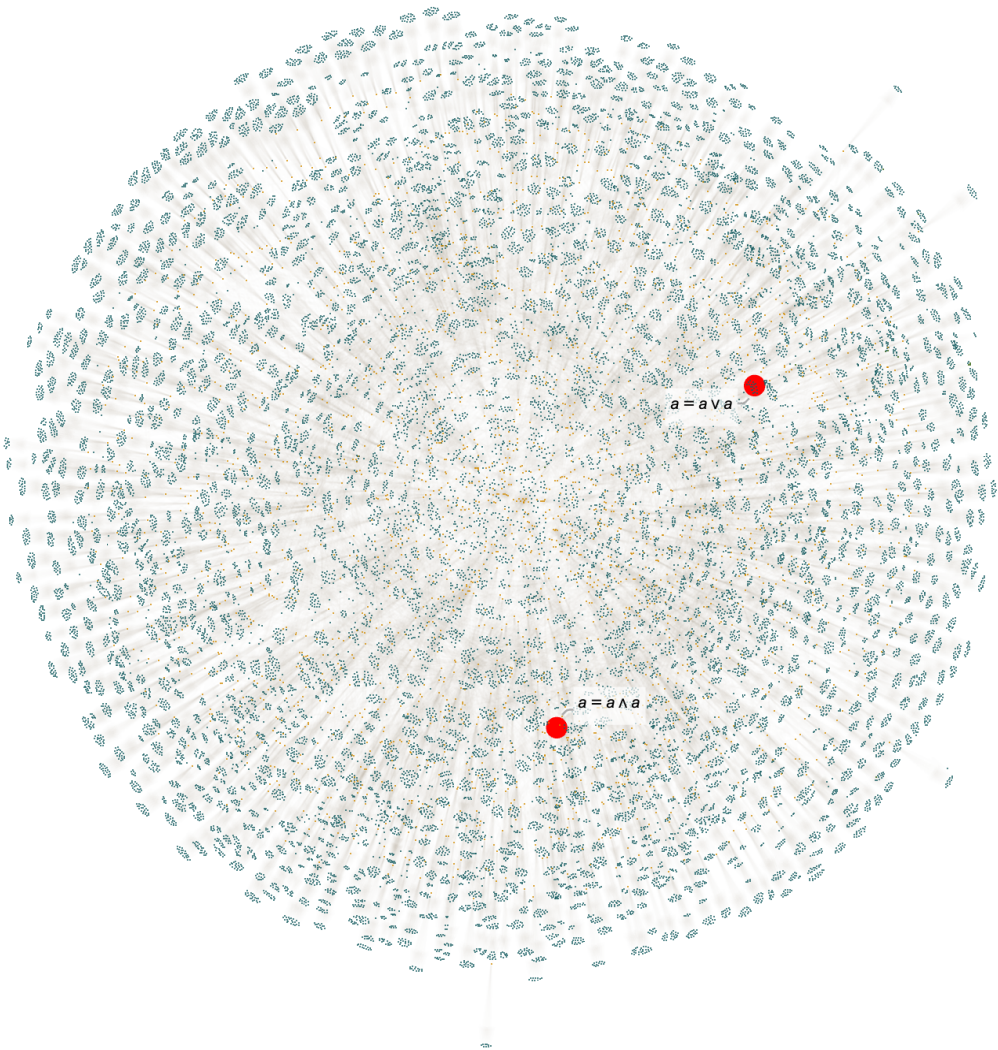

where here we're only operating on theorems with leaf count below 14 (of which there are a total of 27,953).

And if we go to step 3—and use leaf count below 9—we see the law of excluded middle and the law of noncontradiction show up:

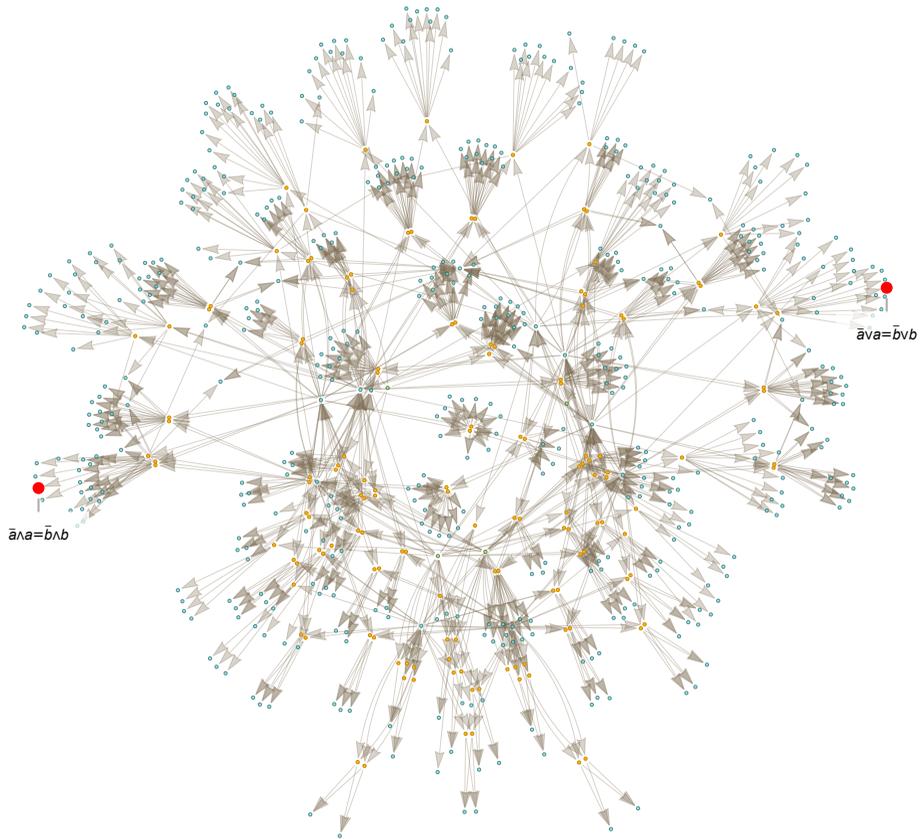

How are these reached? Here's the smallest fragment of token-event graph ("shortest path") within this entailment cone from the axioms to the law of excluded middle:

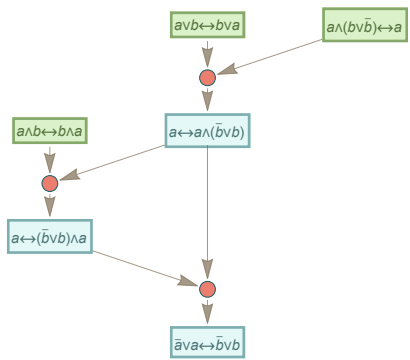

There are actually many possible "paths" (476 in all with our leaf count restriction); the next smallest ones with distinct structures are:

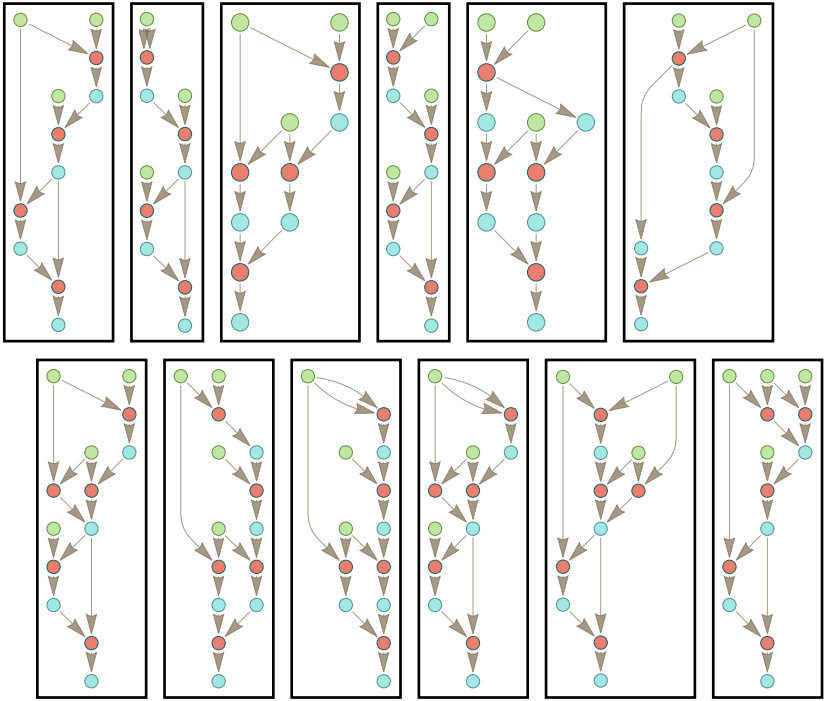

Here's the "path" for this theorem found by automated theorem proving:

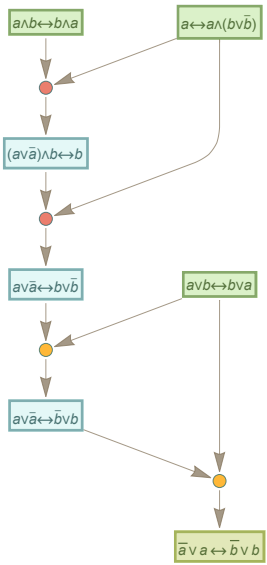

Most of the other "named theorems" involve longer proofs—and so won't show up until much later in the entailment cone:

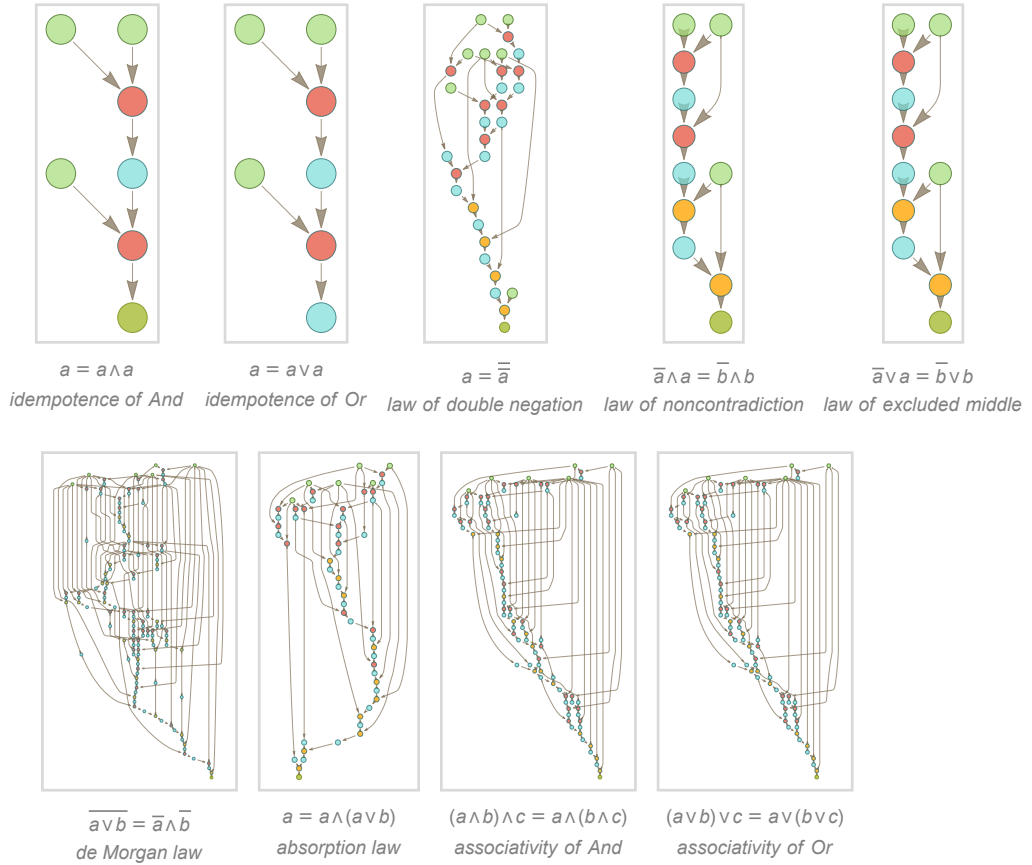

The axiom system we've used for Boolean algebra here is by no means the only possible one. For example, it's stated in terms of **And**, **Or** and **Not**—but one doesn't need all those operators; any Boolean expression (and thus any theorem in Boolean algebra) can also be stated just in terms of the single operator **Nand**.

And in terms of that operator the very simplest axiom system for Boolean algebra contains (as I found in 2000) just one axiom (where here ∘ is now interpreted as **Nand**):

$((b \circ c) \circ a) \circ (b \circ ((b \circ a) \circ b)) = a$

Here's one step of the substitution entailment cone for this axiom:

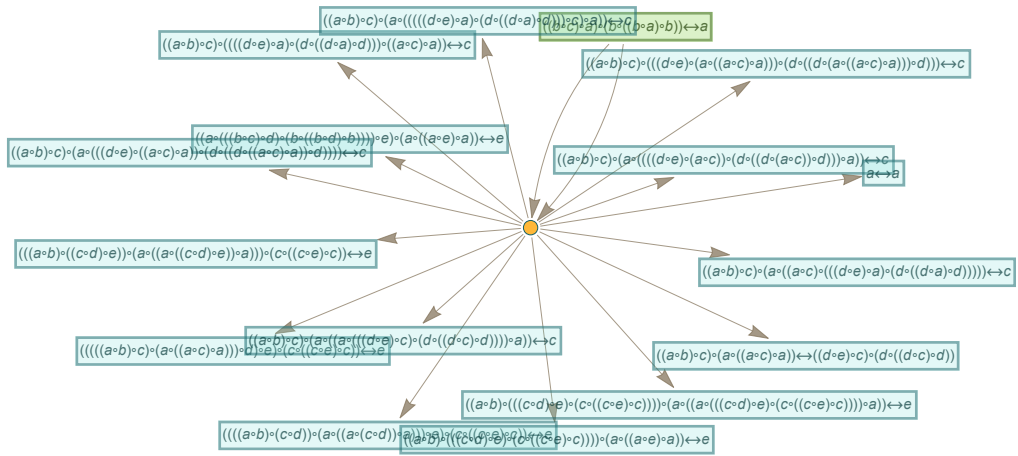

After 2 steps this gives an entailment cone with 5486 theorems

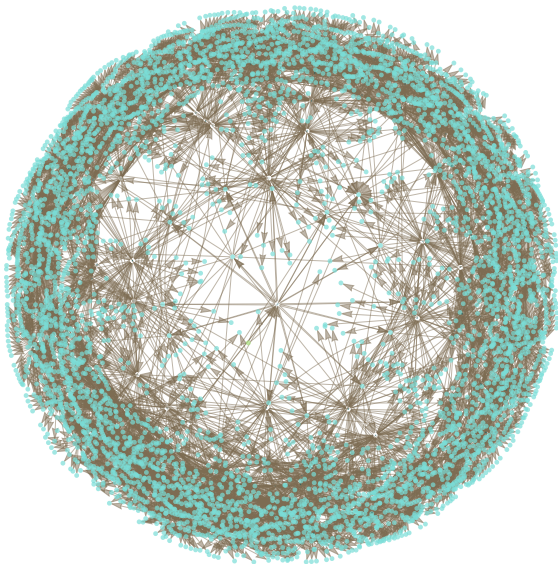

with size distribution:

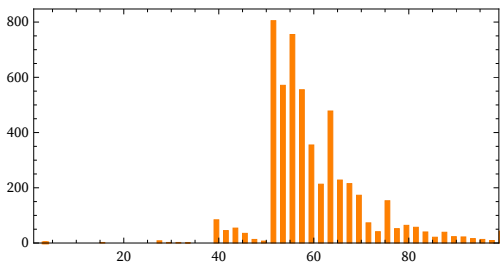

When one's working with **Nand**, it's less clear what one should consider to be "notable theorems". But an obvious one is the commutativity of **Nand**:

$p \circ q = q \circ p$

Here's a proof of this obtained by automated theorem proving (tipped on its side for readability):

Eventually it's inevitable that this theorem must show up in the entailment cone for our axiom system. But based on this proof we would expect it only after something like 102 steps. And with the entailment cone growing exponentially this means that by the time $p \circ q = q \circ p$ shows up, perhaps $10^{100}$ other theorems would have done so—though most vastly more complicated.

We've looked at axioms for group theory and for Boolean algebra. But what about other axiom systems from present-day mathematics? In a sense it's remarkable how few of these there are—and indeed I was able to list essentially all of them in just two pages in *A New Kind of Science*:

The longest axiom system listed here is a precise version of Euclid's original axioms

| |
|---|
| $\forall_{\{x,y,z\}}$ implies[congruent[line[x, y], line[z, z]], congruent[x, y]] |
| $\forall_{\{x,y,z,u,v,w\}}$ implies[<br>    and[congruent[line[x, y], line[z, u]], congruent[line[x, y], line[v, w]]], congruent[line[z, u], line[v, w]]] |
| $\forall_{\{x,y,z\}}$ implies[between[x, y, z], equal[x, y]] |
| $\forall_{\{x,y,z,u,v\}}$ implies[and[between[x, u, z], between[y, v, z]], $\exists_a$ and[between[u, a, y], between[v, a, x]]] |
| $\forall_{\{x,y,z,u,v\}}$ implies[and[and[and[congruent[line[x, u], line[x, v]], congruent[line[y, u], line[y, v]]],<br>      congruent[line[z, u], line[z, v]]], not[equal[u, v]]],<br>    or[or[between[x, y, z], between[y, z, x]], between[z, x, y]]] |
| $\forall_{\{x,y,z,u,v,w\}}$ implies[and[and[and[between[x, y, w], congruent[line[x, y], line[y, w]]],<br>      and[between[x, u, v], congruent[line[x, u], line[u, v]]]],<br>    and[between[y, u, z], congruent[line[y, u], line[z, u]]]], congruent[line[y, z], line[v, w]]] |
| $\forall_{\{x,y,z,a,b,c,u,v\}}$ implies[<br>    and[and[and[and[and[and[not[equal[x, y]], between[x, y, z]], between[a, b, c]], congruent[line[x, y],<br>        line[a, b]]], congruent[line[y, z], line[b, c]]], congruent[line[x, u], line[a, v]]],<br>    congruent[line[y, u], line[b, v]]], congruent[line[z, u], line[c, v]]] |
| $\forall_{\{x,y\}}$ implies[equal[x, y], equal[y, x]] |
| $\forall_{\{x,y,z\}}$ implies[and[equal[x, y], equal[y, z]], equal[x, z]] |
| $\forall_x$ equal[x, x] |
| $\forall_{\{a,b\}}$ and[a, b] = and[b, a] |
| $\forall_{\{a,b\}}$ or[a, b] = or[b, a] |
| $\forall_{\{a,b\}}$ and[a, or[b, not[b]]] = a |
| $\forall_{\{a,b\}}$ or[a, and[b, not[b]]] = a |
| $\forall_{\{a,b,c\}}$ and[a, or[b, c]] = or[and[a, b], and[a, c]] |
| $\forall_{\{a,b,c\}}$ or[a, and[b, c]] = and[or[a, b], or[a, c]] |
| $\forall_{\{a,b\}}$ implies[a, b] = or[not[a], b] |
| $\forall_{\{\alpha,\beta,y,z\}}$ implies[$\exists_x$ implies[and[$\alpha$[y], $\beta$[z]], between[x, y, z]], $\exists_u$ implies[and[$\alpha$[y], $\beta$[z]], between[y, u, z]]] |

where we are listing everything (even logic) in explicit (Wolfram Language) functional form. Given these axioms we should now be able to prove all theorems in Euclidean geometry. As an example (that's already complicated enough) let's take Euclid's very first "proposition" (Book 1, Proposition 1) which states that it's possible "with a ruler and compass" (i.e. with lines and circles) to construct an equilateral triangle based on any line segment—as in:

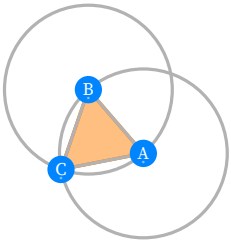

We can write this theorem by saying that given the axioms together with the "setup"

∀_{a,b,c} implies[equal[circle[a, b], circle[a, c]], congruent[line[a, b], line[a, c]]]

it's possible to derive:

and[congruent[line[a, b], line[a, c]], congruent[line[b, a], line[b, c]]]

We can now use automated theorem proving to generate a proof and in this case the proof takes 272 steps. But the fact that it's possible to generate this proof shows that (up to various issues about the "setup conditions") the theorem it proves must eventually "occur naturally" in the entailment cone of the original axioms—though along with an absolutely immense number of other theorems that Euclid didn't "call out" and write down in his books.

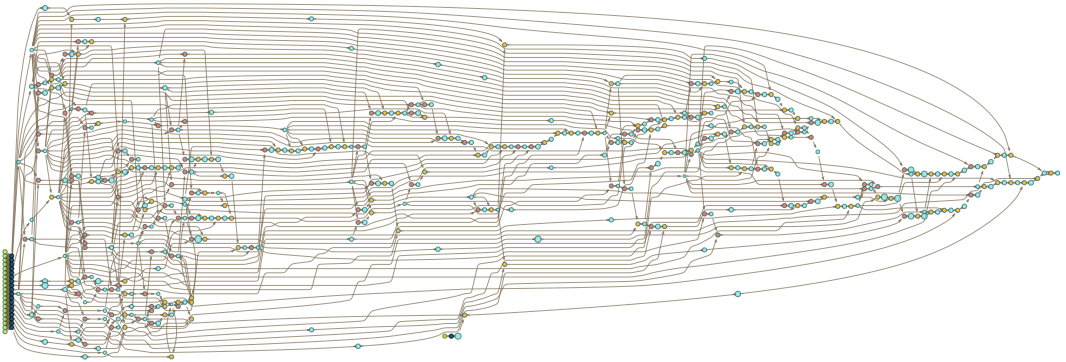

Looking at the collection of axiom systems from *A New Kind of Science* (and a few related ones) for many of them we can just directly start generating entailment cones—here shown after one step, using substitution only:

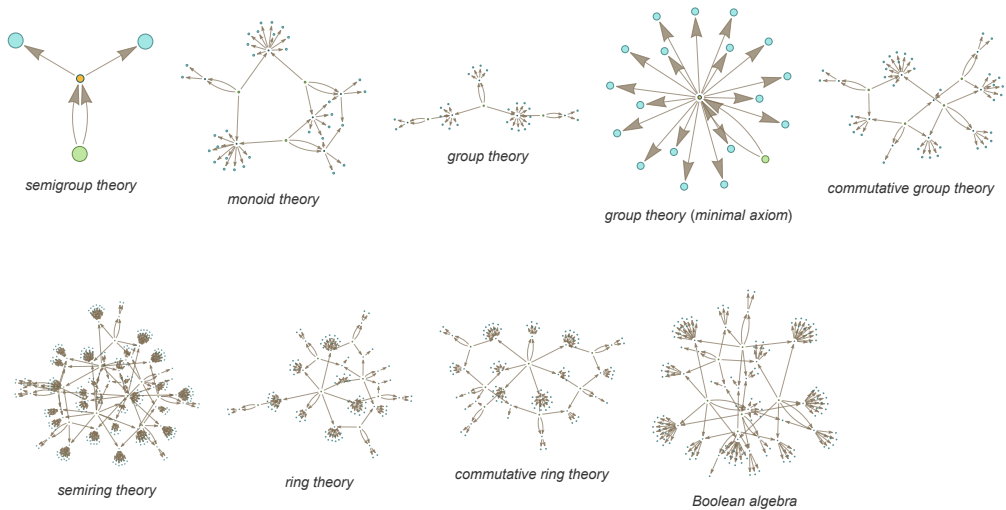

But if we're going to make entailment cones for all axiom systems there are a few other technical wrinkles we have to deal with. The axiom systems shown above are all "straightforwardly equational" in the sense that they in effect state what amount to "algebraic relations" (in the sense of universal algebra) universally valid for all choices of variables. But some axiom systems traditionally used in mathematics also make other kinds of statements. In the traditional formalism and notation of mathematical logic these can look quite complicated and abstruse. But with a metamodel of mathematics like ours it's possible to untangle things to the point where these different kinds of statements can also be handled in a streamlined way.

In standard mathematical notation one might write

$\forall_a \forall_b \ a \circ b = b \circ a$

which we can read as "for all $a$ and $b$, $a \circ b$ equals $b \circ a$"—and which we can interpret in our "metamodel" of mathematics as the (two-way) rule:

$a\_ \circ b\_ \leftrightarrow b\_ \circ a\_$

What this says is just that any time we see an expression that matches the pattern $a\_ \circ b\_$ we can replace it by $b\_ \circ a\_$ (or in Wolfram Language notation just $b \circ a$), and vice versa, so that in effect $a\_ \circ b\_$ can be said to entail $b\_ \circ a\_$.

But what if we have axioms that involve not just universal statements ("for all ...") but also existential statements ("there exists...") ? In a sense we're already dealing with these. Whenever we write $a\_ \circ b\_$—or in explicit functional form, say $o[a\_, b\_]$—we're effectively asserting that there exists some operator $o$ that we can do operations with. It's important to note that once we introduce $o$ (or $\circ$) we imagine that it represents the same thing wherever it appears (in contrast to a pattern variable like $a\_$ that can represent different things in different instances).

Now consider an "explicit existential statement" like

$$\exists_a \, a \circ a = a$$

which we can read as "there exists something *a* for which *a*∘*a* equals *a*". To represent the "something" we just introduce a "constant", or equivalently an expression with head, say, *α*, and zero arguments: *α*[ ]. Now we can write our existential statement as

$$\alpha[\,] \circ \alpha[\,] \leftrightarrow \alpha[\,]$$

or:

$$\circ[\alpha[\,], \, \alpha[\,]] \leftrightarrow \alpha[\,]$$

We can operate on this using rules like $a\_ \circ b\_ \leftrightarrow b\_ \circ a\_$, with *α*[] always "passing through" unchanged—but with its mere presence asserting that "it exists".

A very similar setup works even if we have both universal and existential quantifiers. For example, we can represent

$$\forall_a \, \exists_b \, a \circ b = a$$

as just

$$a\_ \circ \beta[a\_] \leftrightarrow a\_$$

where now there isn't just a single object, say *β*[], that we assert exists; instead there are "lots of different *β*'s", "parametrized" in this case by *a*.

We can apply our standard accumulative bisubstitution process to this statement—and after one step we get:

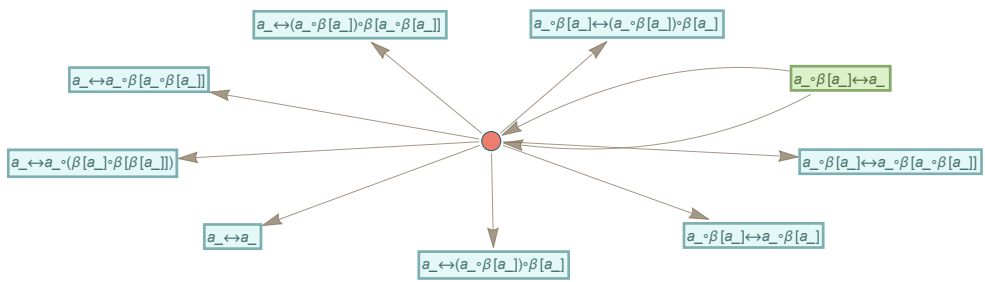

Note that this is a very different result from the one for the "purely universal" statement:

$a\_ \circ b\_ \leftrightarrow a\_$

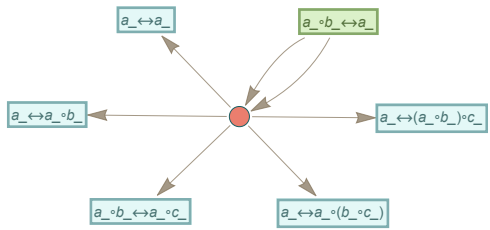

In general, we can "compile" any statement in terms of quantifiers into our metamodel, essentially using the standard technique of Skolemization from mathematical logic. Thus for example

$\forall_a \forall_c \exists_b \, a \circ b = c \circ b$

can be "compiled into"

$a\_ \circ \beta[a\_, c\_] \leftrightarrow c\_ \circ \beta[a\_, c\_]$

while

$\forall_a \exists_b \forall_c \, a \circ b = c$

can be compiled into:

$a\_ \circ \beta[a\_] \leftrightarrow c\_$

If we look at the actual axiom systems used in current mathematics there's one more issue to deal with—which doesn't affect the axioms for logic or group theory, but does show up, for example, in the Peano axioms for arithmetic. And the issue is that in addition to quantifying over "variables", we also need to quantify over "functions". Or formulated differently, we need to set up not just individual axioms, but a whole "axiom schema" that can generate an infinite sequence of "ordinary axioms", one for each possible "function".

In our metamodel of mathematics, we can think of this in terms of "parametrized functions", or in Wolfram Language, just as having functions whose heads are themselves patterns, as in f[n_][a_].

Using this setup we can then "compile" the standard induction axiom of Peano arithmetic

$\forall_f \forall_y (f[0, y] \wedge \forall_x (f[x, y] \Rightarrow f[s[x], y]) \Rightarrow \forall_x f[x, y])$

into the (Wolfram Language) metamodel form

$f\_[0, y\_] \wedge (f\_[x\_, y\_] \rightarrow f\_[s[x\_], y\_]) \rightarrow f\_[z\_, y\_]$

where the "implications" in the original axiom have been converted into one-way rules, so that what the axiom can now be seen to do is to define a transformation for something that is not an "ordinary mathematical-style expression" but rather an expression that is itself a rule.

But the important point is that our whole setup of doing substitutions in symbolic expressions—like Wolfram Language—makes no fundamental distinction between dealing with "ordinary expressions" and with "rules" (in Wolfram Language, for example, $a \to b$ is just **Rule**[**a,b**]). And as a result we can expect to be able to construct token-event graphs, build entailment cones, etc. just as well for axiom systems like Peano arithmetic, as for ones like Boolean algebra and group theory.

The actual number of nodes that appear even in what might seem like simple cases can be huge, but the whole setup makes it clear that exploring an axiom system like this is just another example—that can be uniformly represented with our metamodel of mathematics—of a form of sampling of the ruliad.

## 16 | The Model-Theoretic Perspective

We've so far considered something like

$$x \circ y = (y \circ x) \circ y$$

just as an abstract statement about arbitrary symbolic variables $x$ and $y$, and some abstract operator $\circ$. But can we make a "model" of what $x$, $y$, and $\circ$ could "explicitly be"?

Let's imagine for example that $x$ and $y$ can take 2 possible values, say 0 or 1. (We'll use numbers for notational convenience, though in principle the values could be anything we want.) Now we have to ask what $\circ$ can be in order to have our original statement always hold. It turns out in this case that there are several possibilities, that can be specified by giving possible "multiplication tables" for $\circ$:

| ∘ | 0 | 1 |   | ∘ | 0 | 1 |   | ∘ | 0 | 1 |   | ∘ | 0 | 1 |   | ∘ | 0 | 1 |   | ∘ | 0 | 1 |
|---|---|---|---|---|---|---|---|---|---|---|---|---|---|---|---|---|---|---|---|---|---|---|
| 0 | 0 | 0 |   | 0 | 0 | 0 |   | 0 | 0 | 1 |   | 0 | 0 | 1 |   | 0 | 1 | 0 |   | 0 | 1 | 1 |
| 1 | 0 | 0 |   | 1 | 0 | 1 |   | 1 | 0 | 1 |   | 1 | 1 | 1 |   | 1 | 1 | 0 |   | 1 | 1 | 1 |

(For convenience we'll often refer to such multiplication tables by numbers **FromDigits**[**Flatten**[**m**],**k**], here 0, 1, 5, 7, 10, 15.) Using let's say the second multiplication table we can then "evaluate" both sides of the original statement for all possible choices of $x$ and $y$, and verify that the statement always holds:

| x | y | x∘y | (y∘x)∘y |
|---|---|-----|---------|
| 1 | 1 | 1   | 1       |
| 1 | 0 | 0   | 0       |
| 0 | 1 | 0   | 0       |
| 0 | 0 | 0   | 0       |

If we allow, say, 3 possible values for *x* and *y*, there turn out to be 221 possible forms for ∘. The first few are:

| ∘ | 0 | 1 | 2 |
|---|---|---|---|
| 0 | 0 | 0 | 0 |
| 1 | 0 | 0 | 0 |
| 2 | 0 | 0 | 0 |

| ∘ | 0 | 1 | 2 |
|---|---|---|---|
| 0 | 0 | 0 | 0 |
| 1 | 0 | 0 | 0 |
| 2 | 0 | 0 | 2 |

| ∘ | 0 | 1 | 2 |
|---|---|---|---|
| 0 | 0 | 0 | 0 |
| 1 | 0 | 0 | 2 |
| 2 | 0 | 1 | 0 |

| ∘ | 0 | 1 | 2 |
|---|---|---|---|
| 0 | 0 | 0 | 0 |
| 1 | 0 | 0 | 2 |
| 2 | 0 | 1 | 2 |

| ∘ | 0 | 1 | 2 |
|---|---|---|---|
| 0 | 0 | 0 | 0 |
| 1 | 0 | 0 | 2 |
| 2 | 0 | 2 | 2 |

| ∘ | 0 | 1 | 2 |
|---|---|---|---|
| 0 | 0 | 0 | 0 |
| 1 | 0 | 1 | 0 |
| 2 | 0 | 0 | 0 |

| ∘ | 0 | 1 | 2 |
|---|---|---|---|
| 0 | 0 | 0 | 0 |
| 1 | 0 | 1 | 0 |
| 2 | 0 | 0 | 2 |

| ∘ | 0 | 1 | 2 |
|---|---|---|---|
| 0 | 0 | 0 | 0 |
| 1 | 0 | 1 | 1 |
| 2 | 0 | 1 | 0 |

…

As another example, let's consider the simplest axiom for Boolean algebra (that I discovered in 2000):

$((b \circ c) \circ a) \circ (b \circ ((b \circ a) \circ b)) = a$

Here are the "size-2" models for this

| ∘ | 0 | 1 |
|---|---|---|
| 0 | 1 | 0 |
| 1 | 0 | 0 |

| ∘ | 0 | 1 |
|---|---|---|
| 0 | 1 | 1 |
| 1 | 1 | 0 |

and these, as expected, are the truth tables for **Nand** and **Nor** respectively. (In this particular case, there are no size-3 models, 12 size-4 models, and in general $\frac{2^n!}{n!}$ models of size $2^n$—and no finite models of any other size.)

Looking at this example suggests a way to talk about models for axiom systems. We can think of an axiom system as defining a collection of abstract constraints. But what can we say about objects that might satisfy those constraints? A model is in effect telling us about these objects. Or, put another way, it's telling what "things" the axiom system "describes". And in the case of my axiom for Boolean algebra, those "things" would be Boolean variables, operated on using **Nand** or **Nor**.

As another example, consider the axioms for group theory

$\{a \circ (b \circ c) = (a \circ b) \circ c, \, a \circ \diamond = a, \, a \circ \overline{a} = \diamond\}$

| ∘ | 0 | 1 |
|---|---|---|
| 0 | 0 | 1 |
| 1 | 1 | 0 |

| ∘ | 0 | 1 | 2 |
|---|---|---|---|
| 0 | 0 | 1 | 2 |
| 1 | 1 | 2 | 0 |
| 2 | 2 | 0 | 1 |

| ∘ | 0 | 1 | 2 | 3 |
|---|---|---|---|---|
| 0 | 0 | 1 | 2 | 3 |
| 1 | 1 | 0 | 3 | 2 |
| 2 | 2 | 3 | 0 | 1 |
| 3 | 3 | 2 | 1 | 0 |

| ∘ | 0 | 1 | 2 | 3 |
|---|---|---|---|---|
| 0 | 0 | 1 | 2 | 3 |
| 1 | 1 | 0 | 3 | 2 |
| 2 | 2 | 3 | 1 | 0 |
| 3 | 3 | 2 | 0 | 1 |

| ∘ | 0 | 1 | 2 | 3 |
|---|---|---|---|---|
| 0 | 0 | 1 | 2 | 3 |
| 1 | 1 | 3 | 0 | 2 |
| 2 | 2 | 0 | 3 | 1 |
| 3 | 3 | 2 | 1 | 0 |

| ∘ | 0 | 1 | 2 | 3 |
|---|---|---|---|---|
| 0 | 0 | 1 | 2 | 3 |
| 1 | 1 | 2 | 3 | 0 |
| 2 | 2 | 3 | 0 | 1 |
| 3 | 3 | 0 | 1 | 2 |

Is there a mathematical interpretation of these? Well, yes. They essentially correspond to (representations of) particular finite groups. The original axioms define constraints to be satisfied by any group. These models now correspond to particular groups with specific finite numbers of elements (and in fact specific representations of these groups). And just like in the Boolean algebra case this interpretation now allows us to start saying what the models are "about". The first three, for example, correspond to cyclic groups which can be thought of as being "about" addition of integers mod $k$.

For axiom systems that haven't traditionally been studied in mathematics, there typically won't be any such preexisting identification of what they're "about". But we can still think of models as being a way that a mathematical observer can characterize—or summarize—an axiom system. And in a sense we can see the collection of possible finite models for an axiom system as being a kind of "model signature" for the axiom system.

But let's now consider what models tell us about "theorems" associated with a given axiom system. Take for example the axiom:

$$x = (x \circ y) \circ x$$

Here are the size-2 models for this axiom system:

| ∘ | 0 | 1 |   | ∘ | 0 | 1 |   | ∘ | 0 | 1 |   | ∘ | 0 | 1 |   | ∘ | 0 | 1 |
|---|---|---|---|---|---|---|---|---|---|---|---|---|---|---|---|---|---|---|
| 0 | 0 | 0 |   | 0 | 0 | 1 |   | 0 | 0 | 1 |   | 0 | 1 | 1 |   | 0 | 1 | 1 |
| 1 | 1 | 1 |   | 1 | 0 | 0 |   | 1 | 0 | 1 |   | 1 | 0 | 0 |   | 1 | 0 | 1 |

Let's now pick the last of these models. Then we can take any symbolic expression involving ∘, and say what its values would be for every possible choice of the values of the variables that appear in it:

| a | b | a∘a | a∘b | b∘a | b∘b | (a∘a)∘a | (a∘a)∘b | (a∘b)∘a | (a∘b)∘b | (b∘a)∘a | (b∘a)∘b | (b∘b)∘a | (b∘b)∘b | a∘(a∘a) | a∘(a∘b) |
|---|---|-----|-----|-----|-----|---------|---------|---------|---------|---------|---------|---------|---------|---------|---------|
| 0 | 0 | 1 | 1 | 1 | 1 | 0 | 0 | 0 | 0 | 0 | 0 | 0 | 0 | 1 | 1 |
| 0 | 1 | 1 | 1 | 0 | 1 | 0 | 1 | 0 | 1 | 1 | 1 | 0 | 1 | 1 | 1 |
| 1 | 0 | 1 | 0 | 1 | 1 | 1 | 0 | 1 | 1 | 1 | 0 | 1 | 0 | 1 | 0 |
| 1 | 1 | 1 | 1 | 1 | 1 | 1 | 1 | 1 | 1 | 1 | 1 | 1 | 1 | 1 | 1 |
| 3 | 5 | 15 | 13 | 11 | 15 | 3 | 5 | 3 | 7 | 7 | 5 | 3 | 5 | 15 | 13 |

The last row here gives an "expression code" that summarizes the values of each expression in this particular model. And if two expressions have different codes in the model then this tells us that these expressions cannot be equivalent according to the underlying axiom system.

But if the codes are the same, then it's at least possible that the expressions are equivalent in the underlying axiom system. So as an example, let's take the equivalences associated with pairs of expressions that have code 3 (according to the model we're using):

$\{a, (a \circ a) \circ a, (a \circ b) \circ a, (b \circ b) \circ a\}$

So now let's compare with an actual entailment cone for our underlying axiom system (where to keep the graph of modest size we have dropped expressions involving more than 3 variables):

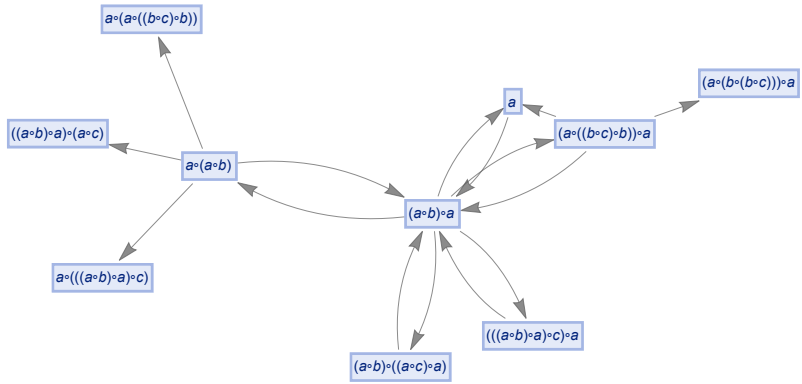

So far this doesn't establish equivalence between any of our code-3 expressions. But if we generate a larger entailment cone (here using a different initial expression) we get

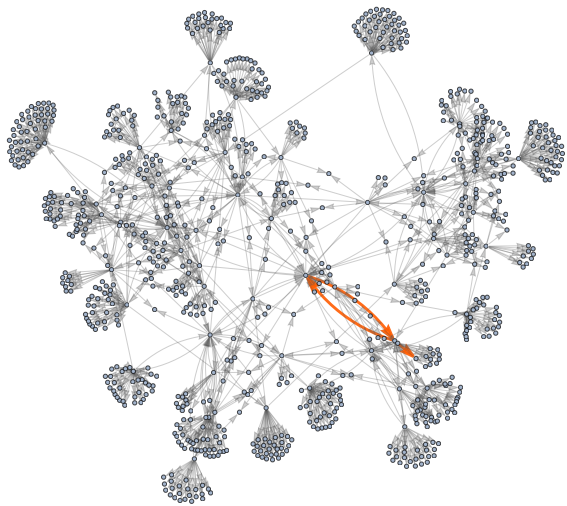

where the path shown corresponds to the statement

$(a \circ a) \circ a = (a \circ b) \circ a$

demonstrating that this is an equivalence that holds in general for the axiom system.

But let's take another statement implied by the model, such as:

$(a \circ b) \circ a = (b \circ b) \circ a$

Yes, it's valid in the model. But it's not something that's generally valid for the underlying axiom system, or could ever be derived from it. And we can see this for example by picking another model for the axiom system, say the second-to-last one in our list above

| ∘ | 0 | 1 |
|---|---|---|
| 0 | 1 | 1 |
| 1 | 0 | 0 |

and finding out that the values for the two expressions here are different in that model:

| a | b | (a∘b)∘a | (b∘b)∘a |
|---|---|---------|---------|
| 0 | 0 | 0 | 0 |
| 0 | 1 | 0 | 1 |
| 1 | 0 | 1 | 0 |
| 1 | 1 | 1 | 1 |
| 3 | 5 | 3 | 5 |

The definitive way to establish that a particular statement follows from a particular axiom system is to find an explicit proof for it, either directly by picking it out as a path in the entailment cone or by using automated theorem proving methods. But models in a sense give one a way to "get an approximate result".

As an example of how this works, consider a collection of possible expressions, with pairs of them joined whenever they can be proved equal in the axiom system we're discussing:

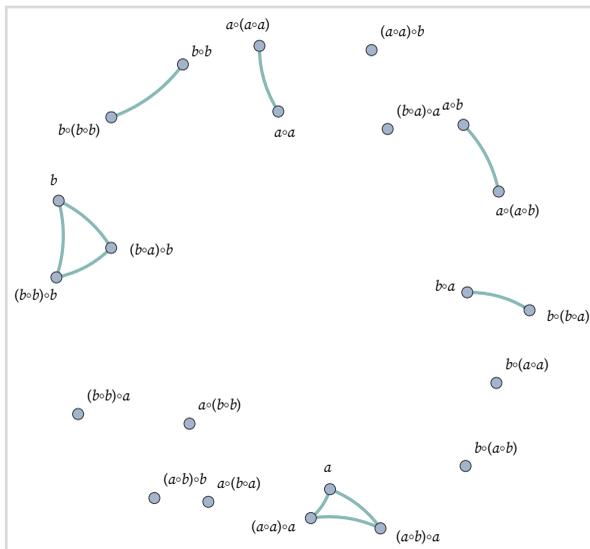

Now let's indicate what codes two models of the axiom system assign to the expressions:

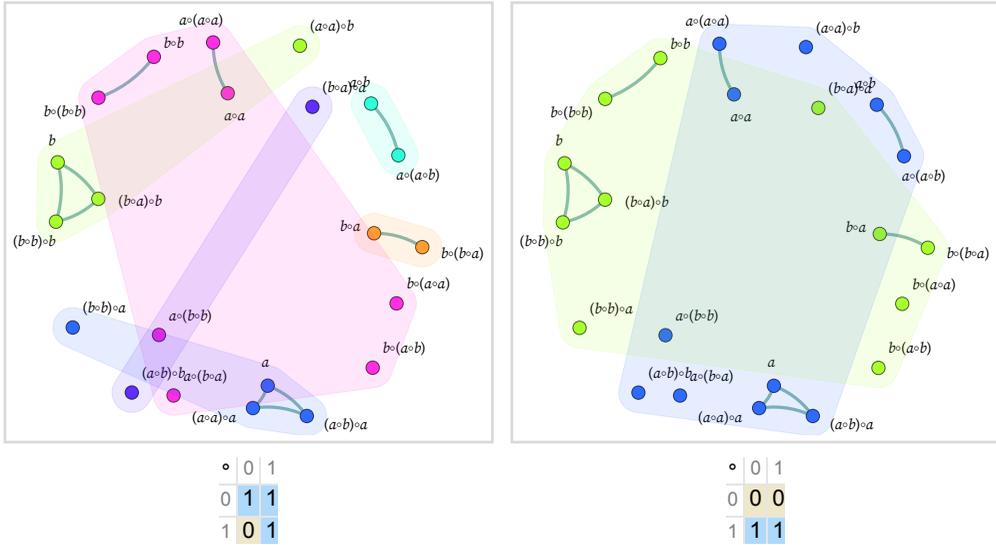

The expressions within each connected graph component are equal according to the underlying axiom system, and in both models they are always assigned the same codes. But sometimes the models "overshoot", assigning the same codes to expressions not in the same connected component—and therefore not equal according to the underlying axiom system.

The models we've shown so far are ones that are valid for the underlying axiom system. If we use a model that isn't valid we'll find that even expressions in the same connected component of the graph (and therefore equal according to the underlying axiom system) will be assigned different codes (note the graphs have been rearranged to allow expressions with the same code to be drawn in the same "patch"):

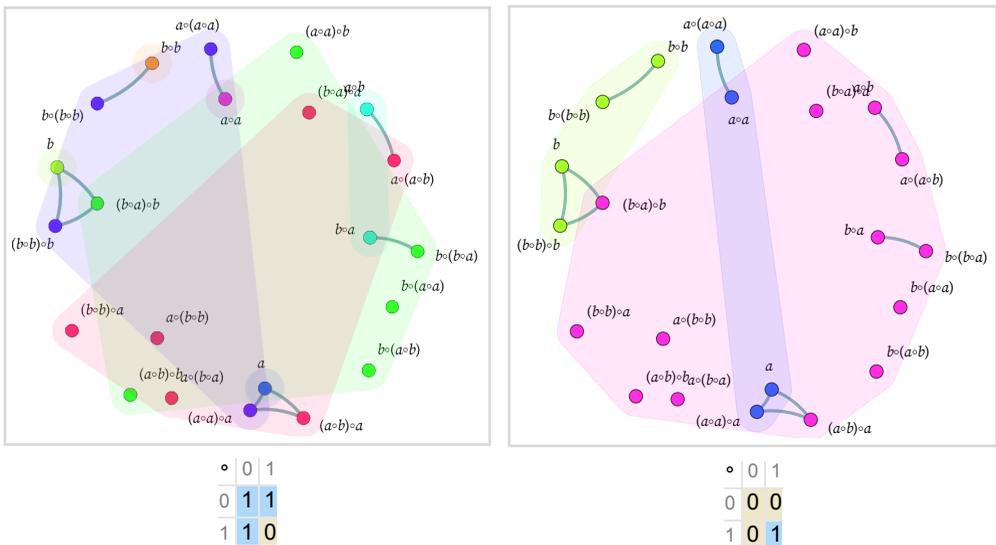

We can think of our graph of equivalences between expressions as corresponding to a slice through an entailment graph—and essentially being "laid out in metamathematical space", like a branchial graph, or what we'll later call an "entailment fabric". And what we see is that when we have a valid model different codes yield different patches that in effect cover metamathematical space in a way that respects the equivalences implied by the underlying axiom system.

But now let's see what happens if we make an entailment cone, tagging each node with the code corresponding to the expression it represents, first for a valid model, and then for non-valid ones:

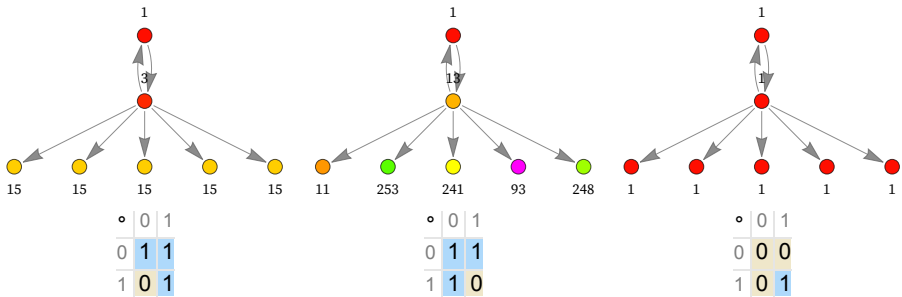

With the valid model, the whole entailment cone is tagged with the same code (and here also same color). But for the non-valid models, different "patches" in the entailment cone are tagged with different codes.

Let's say we're trying to see if two expressions are equal according to the underlying axiom system. The definitive way to tell this is to find a "proof path" from one expression to the other. But as an "approximation" we can just "evaluate" these two expressions according to a model, and see if the resulting codes are the same. Even if it's a valid model, though, this can only definitively tell us that two expressions aren't equal; it can't confirm that they are. In principle we can refine things by checking in multiple models—particularly ones with more elements. But without essentially pre-checking all possible equalities we can't in general be sure that this will give us the complete story.

Of course, generating explicit proofs from the underlying axiom system can also be hard—because in general the proof can be arbitrarily long. And in a sense there is a tradeoff. Given a particular equivalence to check we can either search for a path in the entailment graph, often effectively having to try many possibilities. Or we can "do the work up front" by finding a model or collection of models that we know will correctly tell us whether the equivalence is correct.

Later we'll see how these choices relate to how mathematical observers can "parse" the structure of metamathematical space. In effect observers can either explicitly try to trace out "proof paths" formed from sequences of abstract symbolic expressions—or they can "globally predetermine" what expressions "mean" by identifying some overall model.

In general there may be many possible choices of models—and what we'll see is that these different choices are essentially analogous to different choices of reference frames in physics.

One feature of our discussion of models so far is that we've always been talking about making models for axioms, and then applying these models to expressions. But in the accumulative systems we've discussed above (and that seem like closer metamodels of actual mathematics), we're only ever talking about "statements"—with "axioms" just being statements we happen to start with. So how do models work in such a context?

Here's the beginning of the token-event graph starting with

$x \circ y \leftrightarrow (y \circ x) \circ y$

produced using one step of entailment by substitution:

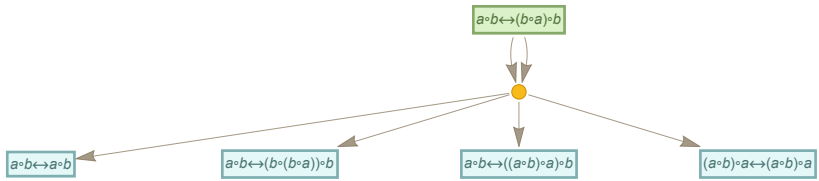

For each of the statements given here, there are certain size-2 models (indicated here by their multiplication tables) that are valid—or in some cases all models are valid:

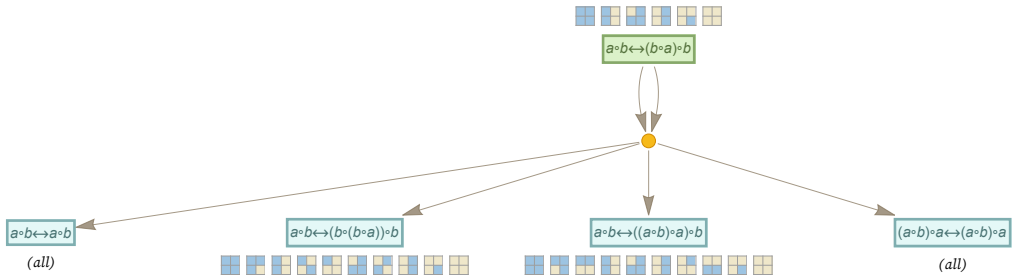

We can summarize this by indicating in a 4×4 grid which of the 16 possible size-2 models are consistent with each statement generated so far in the entailment cone:

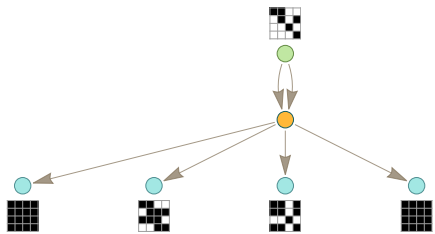

Continuing one more step we get:

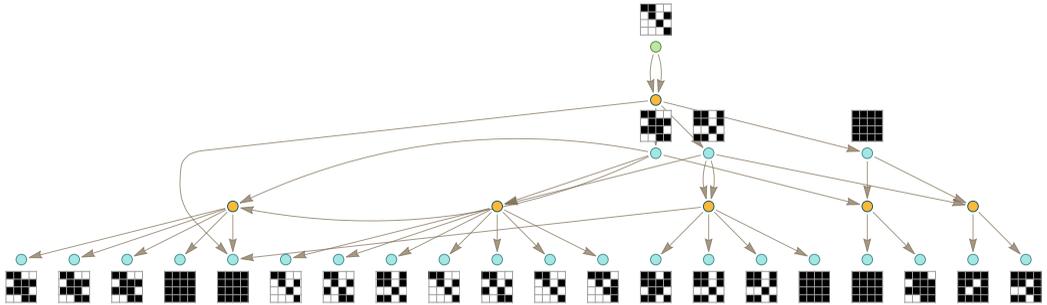

It's often the case that statements generated on successive steps in the entailment cone in essence just "accumulate more models". But—as we can see from the right-hand edge of this graph—it's not always the case—and sometimes a model valid for one statement is no longer valid for a statement it entails. (And the same is true if we use full bisubstitution rather than just substitution.)

Everything we've discussed about models so far here has to do with expressions. But there can also be models for other kinds of structures. For strings it's possible to use something like the same setup, though it doesn't work quite so well. One can think of transforming the string

ABABBAA

into

(((((A∘B)∘A)∘B)∘B)∘A)∘A

and then trying to find appropriate "multiplication tables" for ∘, but here operating on the specific elements A and B, not on a collection of elements defined by the model.

Defining models for a hypergraph rewriting system is more challenging, if interesting. One can think of the expressions we've used as corresponding to trees—which can be "evaluated" as soon as definite "operators" associated with the model are filled in at each node. If we try to do the same thing with graphs (or hypergraphs) we'll immediately be thrust into issues of the order in which we scan the graph.

At a more general level, we can think of a "model" as being a way that an observer tries to summarize things. And we can imagine many ways to do this, with differing degrees of fidelity, but always with the feature that if the summaries of two things are different, then those two things can't be transformed into each other by whatever underlying process is being used.

Put another way, a model defines some kind of invariant for the underlying transformations in a system. The raw material for computing this invariant may be operators at nodes, or may be things like overall graph properties (like cycle counts).

## 17 | Axiom Systems in the Wild

We've talked about what happens with specific, sample axiom systems, as well as with various axiom systems that have arisen in present-day mathematics. But what about "axiom systems in the wild"—say just obtained by random sampling, or by systematic enumeration? In effect, each possible axiom system can be thought of as "defining a possible field of mathematics"—just in most cases not one that's actually been studied in the history of human mathematics. But the ruliad certainly contains all such axiom systems. And in the style of *A New Kind of Science* we can do ruliology to explore them.

As an example, let's look at axiom systems with just one axiom, one binary operator and one or two variables. Here are the smallest few:

| $a = b$ | $a = a \circ a$ | $a = a \circ b$ | $a = b \circ a$ | $a = b \circ b$ | $a \circ a = a \circ b$ | $a \circ a = b \circ a$ |
|---|---|---|---|---|---|---|
| $a \circ a = b \circ b$ | $a \circ b = b \circ a$ | $a \circ b = b \circ b$ | $a = a \circ (a \circ a)$ | $a = a \circ (a \circ b)$ | $a = a \circ (b \circ a)$ | $a = a \circ (b \circ b)$ |
| $a = b \circ (a \circ a)$ | $a = b \circ (a \circ b)$ | $a = b \circ (b \circ a)$ | $a = b \circ (b \circ b)$ | $a = (a \circ a) \circ a$ | $a = (a \circ a) \circ b$ | $a = (a \circ b) \circ a$ |
| $a = (a \circ b) \circ b$ | $a = (b \circ a) \circ a$ | $a = (b \circ a) \circ b$ | $a = (b \circ b) \circ a$ | $a = (b \circ b) \circ b$ | $a \circ a = a \circ (a \circ a)$ | $a \circ a = a \circ (a \circ b)$ |
| $a \circ a = a \circ (b \circ a)$ | $a \circ a = a \circ (b \circ b)$ | $a \circ a = b \circ (a \circ a)$ | $a \circ a = b \circ (a \circ b)$ | $a \circ a = b \circ (b \circ a)$ | $a \circ a = b \circ (b \circ b)$ | $a \circ a = (a \circ a) \circ a$ |

For each of these axiom systems, we can then ask what theorems they imply. And for example we can enumerate theorems—just as we have enumerated axiom systems—then use automated theorem proving to determine which theorems are implied by which axiom systems. This shows the result, with possible axiom systems going down the page, possible theorems going across, and a particular square being filled in (darker for longer proofs) if a given theorem can be proved from a given axiom system:

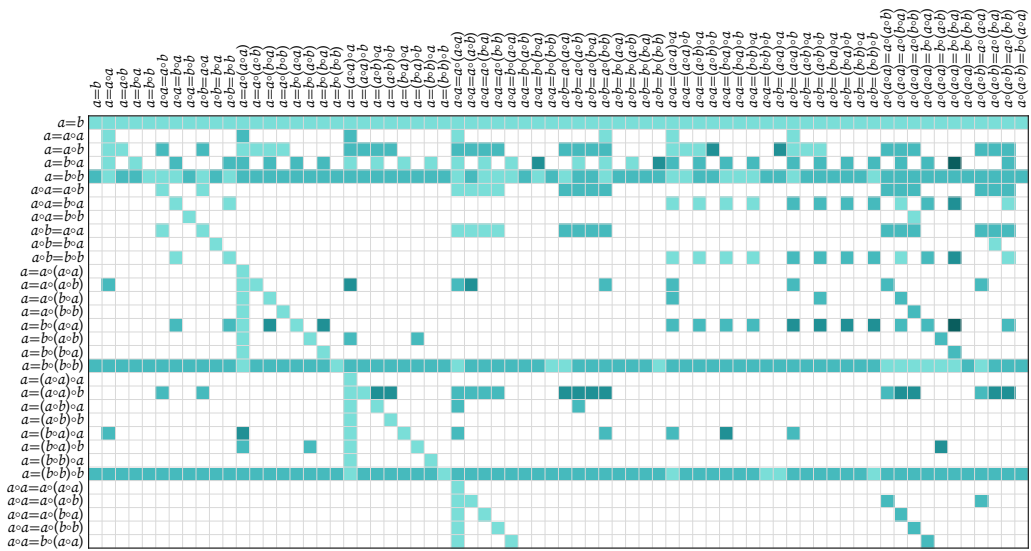

The diagonal on the left is axioms "proving themselves". The lines across are for axiom systems like *a=b* that basically say that any two expressions are equal—so that any theorem that is stated can be proved from the axiom system.

But what if we look at the whole entailment cone for each of these axiom systems? Here are a few examples of the first two steps:

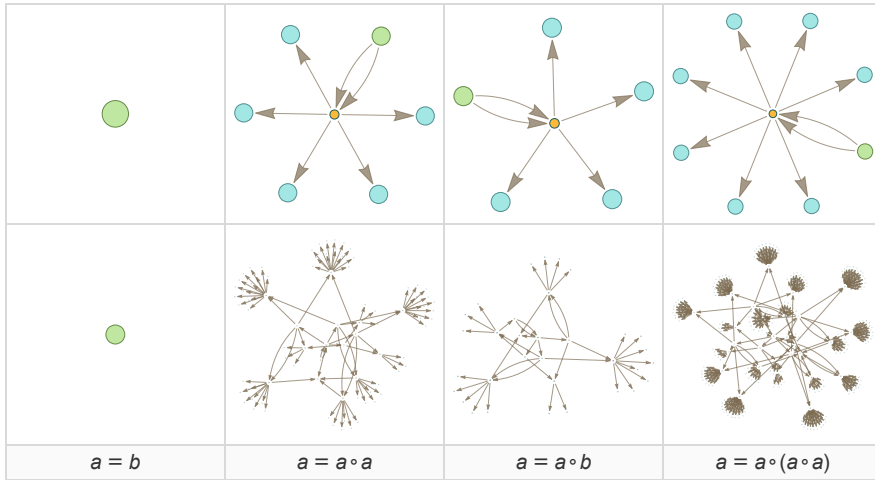

With our method of accumulative evolution the axiom *a=b* doesn't on its own generate a growing entailment cone (though if combined with any axiom containing ∘ it does, and so does *a=b∘c* on its own). But in all the other cases shown the entailment cone grows rapidly (typically at least exponentially)—in effect quickly establishing many theorems. Most of those theorems, however, are "not small"—and for example after 2 steps here are the distributions of their sizes:

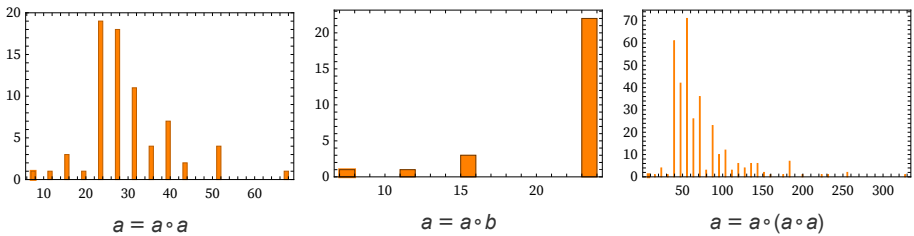

So let's say we generate only one step in the entailment cone. This is the pattern of "small theorems" we establish:

And here is the corresponding result after two steps:

Superimposing this on our original array of theorems we get:

[Figure: grid showing theorem array with highlighted entries]

In other words, there are many small theorems that we can establish "if we look for them", but which won't "naturally be generated" quickly in the entailment cone (though eventually it's inevitable that they will be generated). (Later we'll see how this relates to the concept of "entailment fabrics" and the "knitting together of pieces of mathematics".)

In the previous section we discussed the concept of models for axiom systems. So what models do typical "axiom systems from the wild" have? The number of possible models of a given size varies greatly for different axiom systems:

| | *size 2* | *size 3* |
|---|---|---|
| $a=b$ | [grid] | [grid] |
| $a=a \circ a$ | [grid] | [grid] …(729) |
| $a=a \circ b$ | [grid] | [grid] |
| $a=a \circ (a \circ a)$ | [grid] …(9) | [grid] …(3375) |

|  | size 2 | size 3 | size 4 |
|---|---|---|---|
| $a=b$ | 2 | 3 | 4 |
| $a=a\circ a$ | 4 | 729 | 16 777 216 |
| $a=a\circ b$ | 1 | 1 | 1 |
| $a=b\circ a$ | 1 | 1 | 1 |
| $a=b\circ b$ | 0 | 0 | 0 |
| $a\circ a=a\circ b$ | 4 | 27 | 256 |
| $a\circ a=b\circ a$ | 4 | 27 | 256 |
| $a\circ a=b\circ b$ | 8 | 2187 | 67 108 864 |

,

|  | size 2 | size 3 | size 4 |
|---|---|---|---|
| $a\circ b=a\circ a$ | 4 | 27 | 256 |
| $a\circ b=b\circ a$ | 8 | 729 | 1 048 576 |
| $a\circ b=b\circ b$ | 4 | 27 | 256 |
| $a=a\circ(a\circ a)$ | 9 | 3375 | 157 351 936 |
| $a=a\circ(a\circ b)$ | 1 | 27 | 10 000 |
| $a=a\circ(b\circ a)$ | 5 | 136 | 46 121 |
| $a=a\circ(b\circ b)$ | 5 | 298 | 1 147 649 |
| $a=b\circ(a\circ a)$ | 2 | 4 | 10 |

,

|  | size 2 | size 3 | size 4 |
|---|---|---|---|
| $a=b\circ(a\circ b)$ | 2 | 3 | 18 |
| $a=b\circ(b\circ a)$ | 4 | 64 | 10 000 |
| $a=b\circ(b\circ b)$ | 0 | 0 | 0 |
| $a=(a\circ a)\circ a$ | 9 | 3375 | 157 351 936 |
| $a=(a\circ a)\circ b$ | 2 | 4 | 10 |
| $a=(a\circ b)\circ a$ | 5 | 136 | 46 121 |
| $a=(a\circ b)\circ b$ | 4 | 64 | 10 000 |
| $a=(b\circ a)\circ a$ | 1 | 27 | 10 000 |

,

|  | size 2 | size 3 | size 4 |
|---|---|---|---|
| $a=(b\circ a)\circ b$ | 2 | 3 | 18 |
| $a=(b\circ b)\circ a$ | 5 | 298 | 1 147 649 |
| $a=(b\circ b)\circ b$ | 0 | 0 | 0 |
| $a\circ a=a\circ(a\circ a)$ | 9 | 3375 | 157 351 936 |
| $a\circ a=a\circ(a\circ b)$ | 4 | 343 | 614 656 |
| $a\circ a=a\circ(b\circ a)$ | 7 | 476 | 506 698 |
| $a\circ a=a\circ(b\circ b)$ | 6 | 411 | 1 401 880 |
| $a\circ a=b\circ(a\circ a)$ | 5 | 298 | 1 147 649 |

But for each model we can ask what theorems it implies are valid. And for example combining all models of size 2 yields the following "predictions" for what theorems are valid (with the actual theorems indicated by dots):

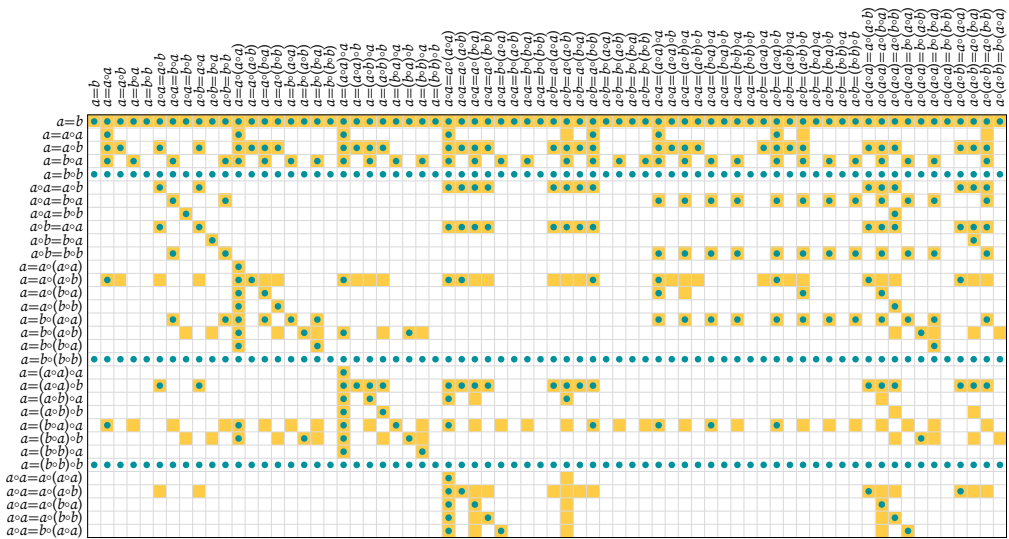

Using instead models of size 3 gives "more accurate predictions":

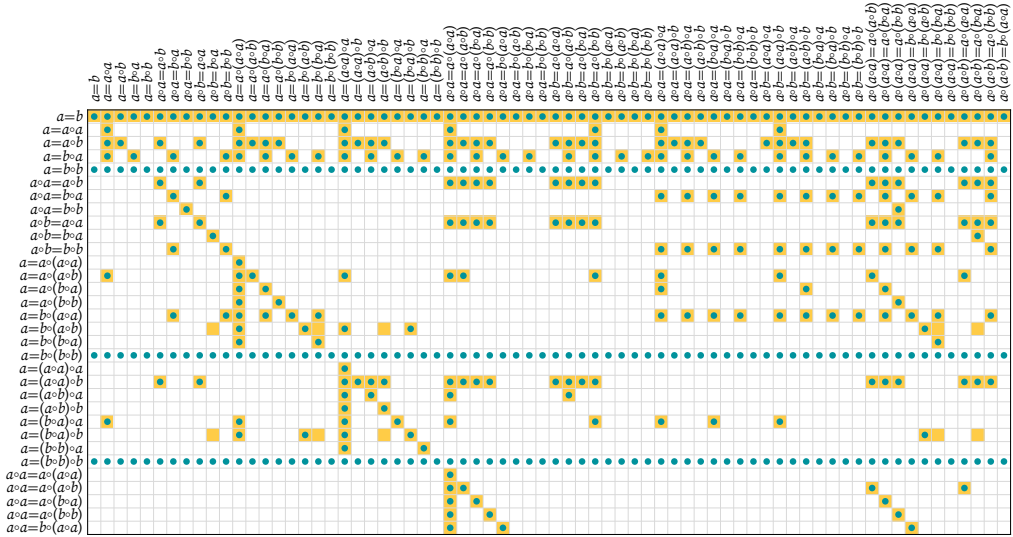

As expected, looking at a fixed number of steps in the entailment cone "underestimates" the number of valid theorems, while looking at finite models overestimates it.

So how does our analysis for "axiom systems from the wild" compare with what we'd get if we considered axiom systems that have been explicitly studied in traditional human mathematics? Here are some examples of "known" axiom systems that involve just a single binary operator

| | | |
|---:|:---|:---|
| (1) | Boolean algebra | $(a \circ a) \circ (a \circ a) = a,\ a \circ (b \circ (b \circ b)) = a \circ a,\ (a \circ (b \circ c)) \circ (a \circ (b \circ c)) = ((b \circ b) \circ a) \circ ((c \circ c) \circ a)$ |
| (2) | Boolean algebra | $(a \circ a) \circ (a \circ b) = a,\ a \circ (a \circ b) = a \circ (b \circ b),\ a \circ (a \circ (b \circ c)) = b \circ (b \circ (a \circ c))$ |
| (3) | Boolean algebra | $a \circ b = b \circ a,\ (a \circ b) \circ (a \circ (b \circ c)) = a$ |
| (4) | Boolean algebra | $((b \circ c) \circ a) \circ (b \circ ((b \circ a) \circ b)) = a$ |
| (5) | Boolean algebra | $(a \circ ((c \circ a) \circ a)) \circ (c \circ (b \circ a)) = c$ |
| (6) | equivalential calculus | $(a \circ b) \circ c = a \circ (b \circ c),\ a \circ b = b \circ a,\ (b \circ b) \circ a = a$ |
| (7) | implicational calculus | $(a \circ b) \circ a = a,\ a \circ (b \circ c) = b \circ (a \circ c),\ (a \circ b) \circ b = (b \circ a) \circ a$ |
| (8) | junctional calculus | $(a \circ b) \circ c = a \circ (b \circ c),\ a \circ b = b \circ a,\ a \circ a = a$ |
| (9) | semigroups | $a \circ (b \circ c) = (a \circ b) \circ c$ |
| (10) | groups | $a \circ ((((a \circ a) \circ b) \circ c) \circ (((a \circ a) \circ a) \circ c)) = b$ |
| (11) | Abelian groups | $a \circ (b \circ (c \circ (a \circ b))) = c$ |
| (12) | central groupoids | $(a \circ b) \circ (b \circ c) = a$ |

and here's the distribution of theorems they give:

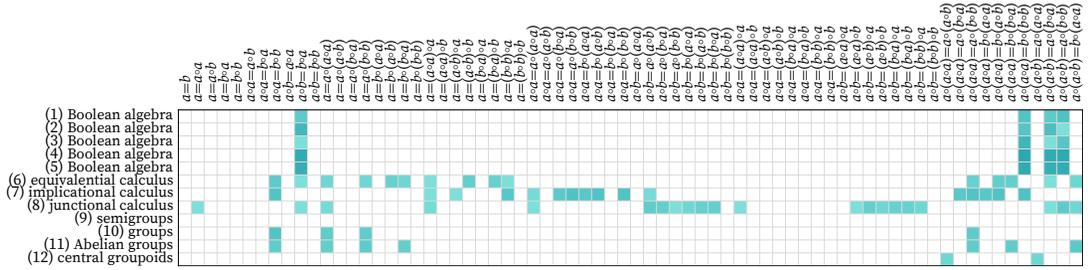

As must be the case, all the axiom systems for Boolean algebra yield the same theorems. But axiom systems for "different mathematical theories" yield different collections of theorems.

What happens if we look at entailments from these axiom systems? Eventually all theorems must show up somewhere in the entailment cone of a given axiom system. But here are the results after one step of entailment:

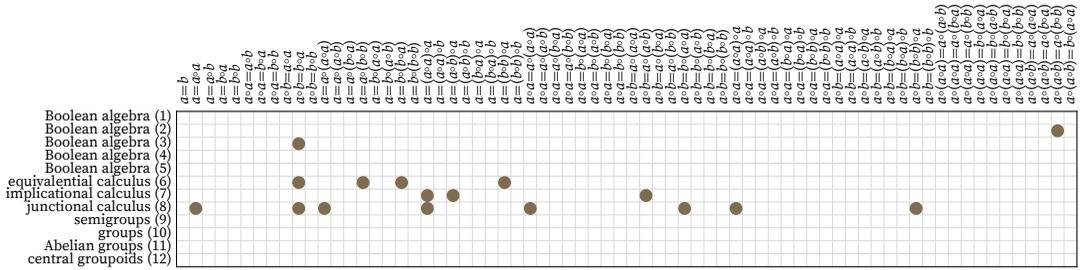

Some theorems have already been generated, but many have not:

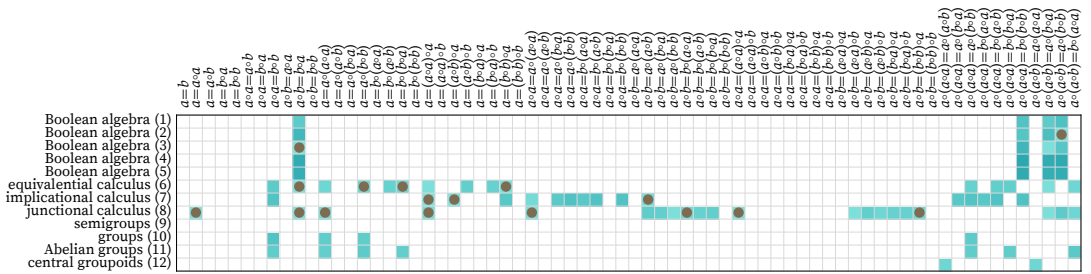

Just as we did above, we can try to "predict" theorems by constructing models. Here's what happens if we ask what theorems hold for all valid models of size 2:

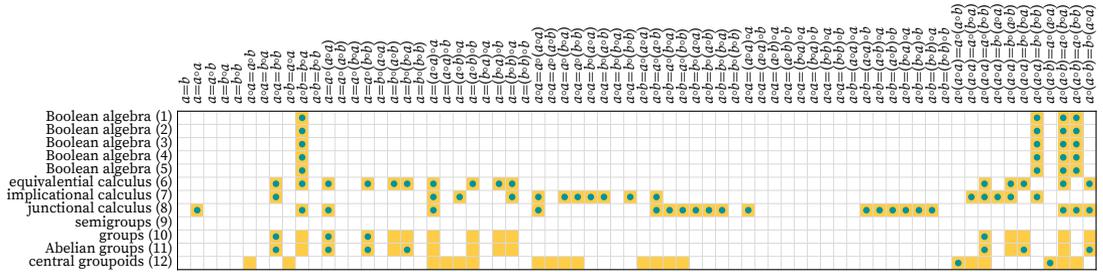

For several of the axiom systems, the models "perfectly predict" at least the theorems we show here. And for Boolean algebra, for example, this isn't surprising: the models just correspond to identifying ∘ as **Nand** or **Nor**, and to say this gives a complete description of Boolean algebra. But in the case of groups, "size-2 models" just capture particular groups that happen to be of size 2, and for these particular groups there are special, extra theorems that aren't true for groups in general.

If we look at models specifically of size 3 there aren't any examples for Boolean algebra so we don't predict any theorems. But for group theory, for example, we start to get a slightly more accurate picture of what theorems hold in general:

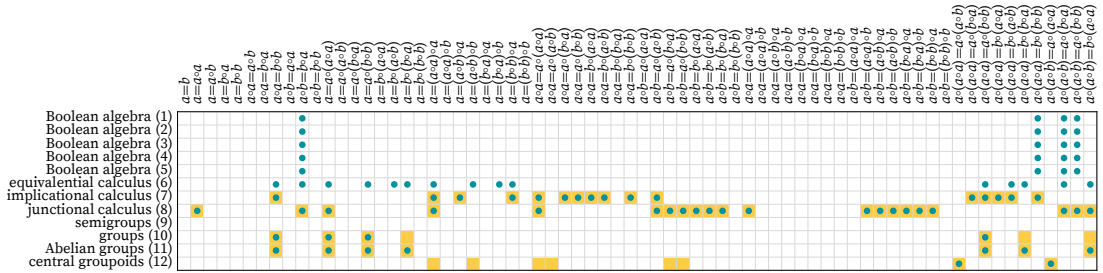

Based on what we've seen here, is there something "obviously special" about the axiom systems that have traditionally been used in human mathematics? There are cases like Boolean algebra where the axioms in effect constrain things so much that we can reasonably say that they're "talking about definite things" (like **Nand** and **Nor**). But there are plenty of other cases, like group theory, where the axioms provide much weaker constraints, and for example allow an infinite number of possible specific groups. But both situations occur among axiom systems "from the wild". And in the end what we're doing here doesn't seem to reveal anything "obviously special" (say in the statistics of models or theorems) about "human" axiom systems.

And what this means is that we can expect that conclusions we draw from looking at the "general case of all axiom systems"—as captured in general by the ruliad—can be expected to hold in particular for the specific axiom systems and mathematical theories that human mathematics has studied.

## 18 | The Topology of Proof Space

In the typical practice of pure mathematics the main objective is to establish theorems. Yes, one wants to know that a theorem has a proof (and perhaps the proof will be helpful in understanding the theorem), but the main focus is on theorems and not on proofs. In our effort to "go underneath" mathematics, however, we want to study not only what theorems there are, but also the process by which the theorems are reached. We can view it as an important simplifying assumption of typical mathematical observers that all that matters is theorems—and that different proofs aren't relevant. But to explore the underlying structure of metamathematics, we need to unpack this—and in effect look directly at the structure of proof space.

Let's consider a simple system based on strings. Say we have the rewrite rule {A→BBB, BB→A} and we want to establish the theorem A→ABA. To do this we have to find some path from A to ABA in the multiway system (or, effectively, in the entailment cone for this axiom system):

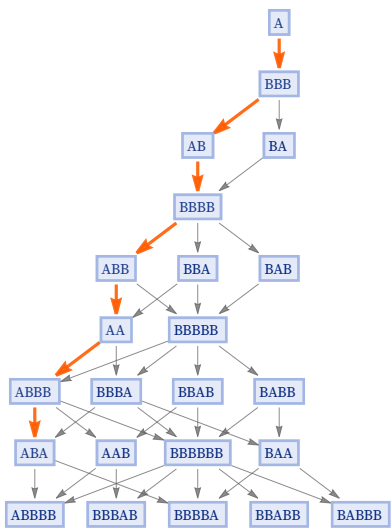

But this isn't the only possible path, and thus the only possible proof. In this particular case, there are 20 distinct paths, each corresponding to at least a slightly different proof:

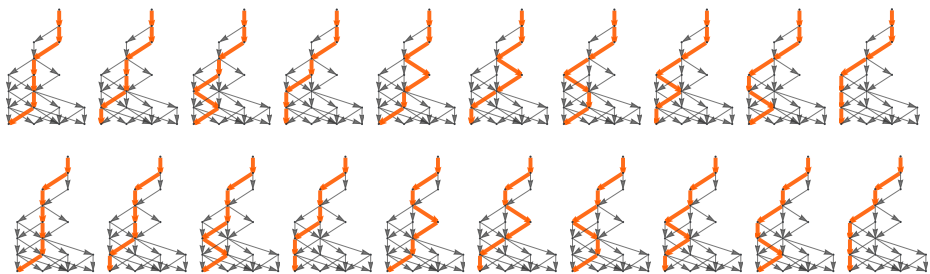

But one feature here is that all these different proofs can in a sense be "smoothly deformed" into each other, in this case by progressively changing just one step at a time. So this means that in effect there is no nontrivial topology to proof space in this case—and "distinctly inequivalent" collections of proofs:

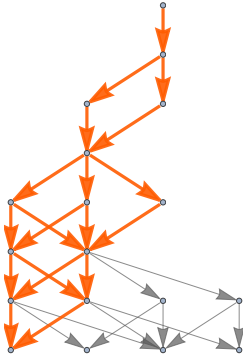

But consider instead the rule {A→AA,A→BAAB}. With this "axiom system" there are 15 possible proofs for the theorem A→BAABAABAAB:

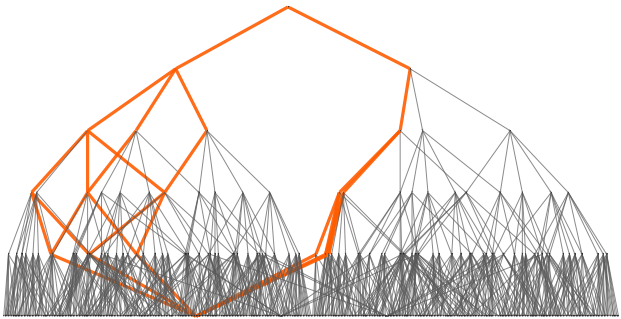

Pulling out just the proofs we get:

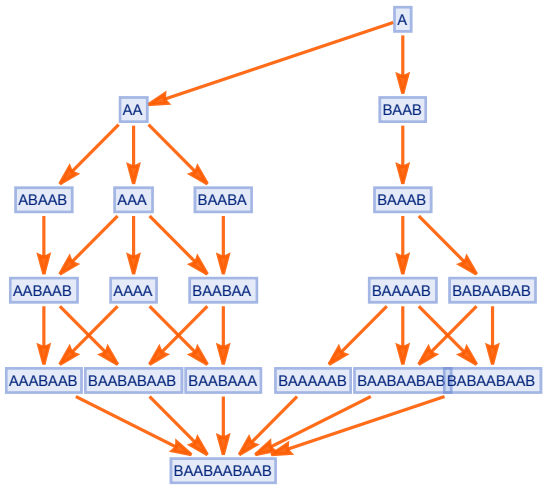

And we see that in a sense there's a "hole" in proof space here—so that there are two distinctly different kinds of proofs that can be done.

One place it's common to see a similar phenomenon is in games and puzzles. Consider for example the Towers of Hanoi puzzle. We can set up a multiway system for the possible moves that can be made. Starting from all disks on the left peg, we get after 1 step:

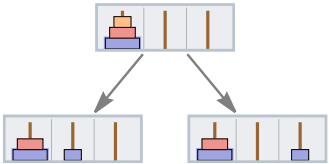

After 2 steps we have:

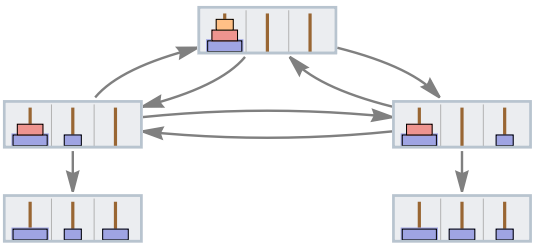

And after 8 steps (in this case) we have the whole "game graph":

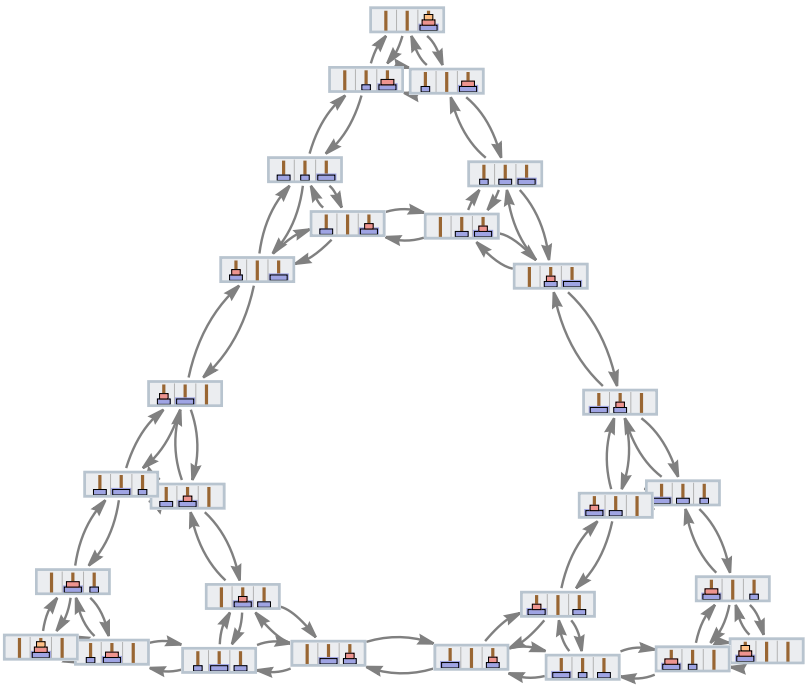

The corresponding result for 4 disks is:

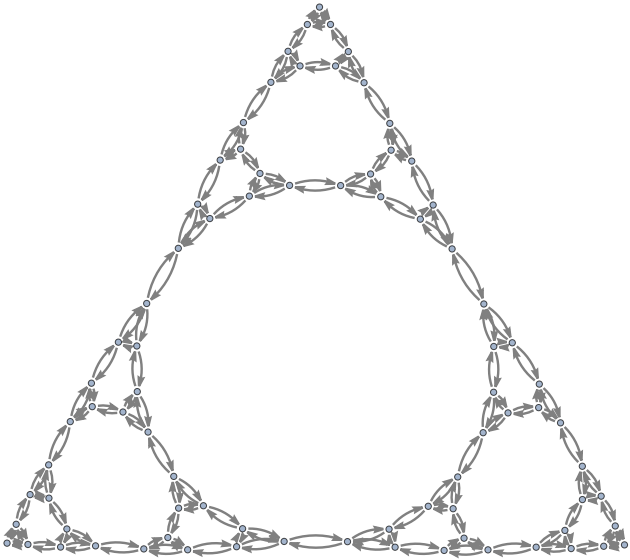

And in each case we see the phenomenon of nontrivial topology. What fundamentally causes this? In a sense it reflects the possibility for distinctly different strategies that lead to the same result. Here, for example, different sides of the "main loop" correspond to the "foundational choice" of whether to move the biggest disk first to the left or to the right. And the same basic thing happens with 4 disks on 4 pegs, though the overall structure is more complicated there:

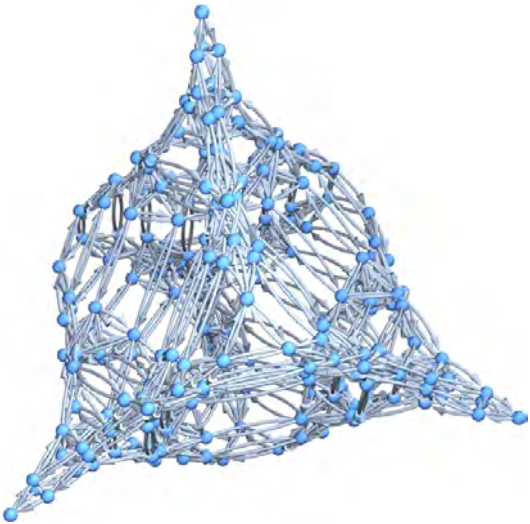

If two paths diverge in a multiway system it could be that it will never be possible for them to merge again. But whenever the system has the property of confluence, it's guaranteed that eventually the paths will merge. And, as it turns out, our accumulative evolution setup guarantees that (at least ignoring generation of new variables) confluence will always be achieved. But the issue is how quickly. If branches always merge after just one step, then in a sense there'll always be topologically trivial proof space. But if the merging can take awhile (and in a continuum limit, arbitrarily long) then there'll in effect be nontrivial topology.

And one consequence of the nontrivial topology we're discussing here is that it leads to disconnection in branchial space. Here are the branchial graphs for the first 3 steps in our original 3-disk 3-peg case:

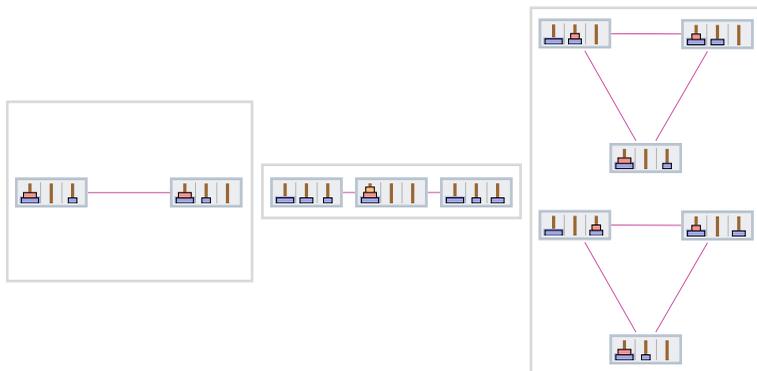

For the first two steps, the branchial graphs stay connected; but on the third step there's disconnection. For the 4-disk 4-peg case the sequence of branchial graphs begins:

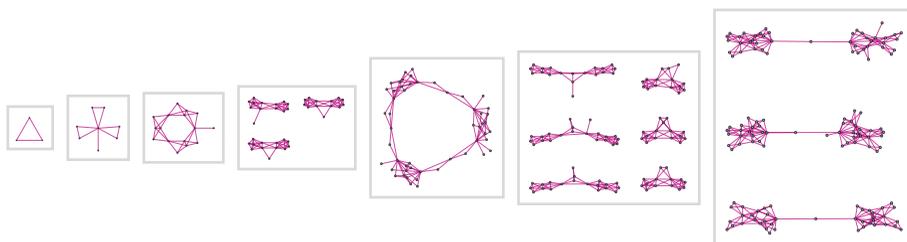

At the beginning (and also the end) there's a single component, that we might think of as a coherent region of metamathematical space. But in the middle it breaks into multiple disconnected components—in effect reflecting the emergence of multiple distinct regions of metamathematical space with something like event horizons temporarily existing between them.

How should we interpret this? First and foremost, it's something that reveals that there's structure "below" the "fluid dynamics" level of mathematics; it's something that depends on the discrete "axiomatic infrastructure" of metamathematics. And from the point of view of our Physics Project, we can think of it as a kind of metamathematical analog of a "quantum effect".

In our Physics Project we imagine different paths in the multiway system to correspond to different possible quantum histories. The observer is in effect spread over multiple paths, which they coarse grain or conflate together. An "observable quantum effect" occurs when there are paths that can be followed by the system, but that are somehow "too far apart" to be immediately coarse-grained together by the observer.

Put another way, there is "noticeable quantum interference" when the different paths corresponding to different histories that are "simultaneously happening" are "far enough apart" to be distinguished by the observer. "Destructive interference" is presumably associated with paths that are so far apart that to conflate them would effectively require conflating essentially every possible path. (And our later discussion of the relationship between falsity and the "principle of explosion" then suggests a connection between destructive interference in physics and falsity in mathematics.)

In essence what determines the extent of "quantum effects" is then our "size" as observers in branchial space relative to the size of features in branchial space such as the "topological holes" we've been discussing. In the metamathematical case, the "size" of us as observers is in effect related to our ability (or choice) to distinguish slight differences in axiomatic formulations of things. And what we're saying here is that when there is nontrivial topology in proof space, there is an intrinsic dynamics in metamathematical entailment that leads to the development of distinctions at some scale—though whether these become "visible" to us as mathematical observers depends on how "strong a metamathematical microscope" we choose to use relative to the scale of the "topological holes".

## 19 | Time, Timelessness and Entailment Fabrics

A fundamental feature of our metamodel of mathematics is the idea that a given set of mathematical statements can entail others. But in this picture what does "mathematical progress" look like?

In analogy with physics one might imagine it would be like the evolution of the universe through time. One would start from some limited set of axioms and then—in a kind of "mathematical Big Bang"—these would lead to a progressively larger entailment cone containing more and more statements of mathematics. And in analogy with physics, one could imagine that the process of following chains of successive entailments in the entailment cone would correspond to the passage of time.

But realistically this isn't how most of the actual history of human mathematics has proceeded. Because people—and even their computers—basically never try to extend mathematics by axiomatically deriving all possible valid mathematical statements. Instead, they come up with particular mathematical statements that for one reason or another they think are valid and interesting, then try to prove these.

Sometimes the proof may be difficult, and may involve a long chain of entailments. Occasionally—especially if automated theorem proving is used—the entailments may approximate a geodesic path all the way from the axioms. But the practical experience of human mathematics tends to be much more about identifying "nearby statements" and then trying to "fit them together" to deduce the statement one's interested in.

And in general human mathematics seems to progress not so much through the progressive "time evolution" of an entailment graph as through the assembly of what one might call an "entailment fabric" in which different statements are being knitted together by entailments.

In physics, the analog of the entailment graph is basically the causal graph which builds up over time to define the content of a light cone (or, more accurately, an entanglement cone). The analog of the entailment fabric is basically the (more-or-less) instantaneous state of space (or, more accurately, branchial space).

In our Physics Project we typically take our lowest-level structure to be a hypergraph—and informally we often say that this hypergraph "represents the structure of space". But really we should be deducing the "structure of space" by taking a particular time slice from the "dynamic evolution" represented by the causal graph—and for example we should think of two "atoms of space" as "being connected" in the "instantaneous state of space" if there's a causal connection between them defined within the slice of the causal graph that occurs within the time slice we're considering. In other words, the "structure of space" is knitted together by the causal connections represented by the causal graph. (In traditional physics, we might say that space can be "mapped out" by looking at overlaps between lots of little light cones.)

Let's look at how this works out in our metamathematical setting, using string rewrites to simplify things. If we start from the axiom A↔AA this is the beginning of the entailment cone it generates:

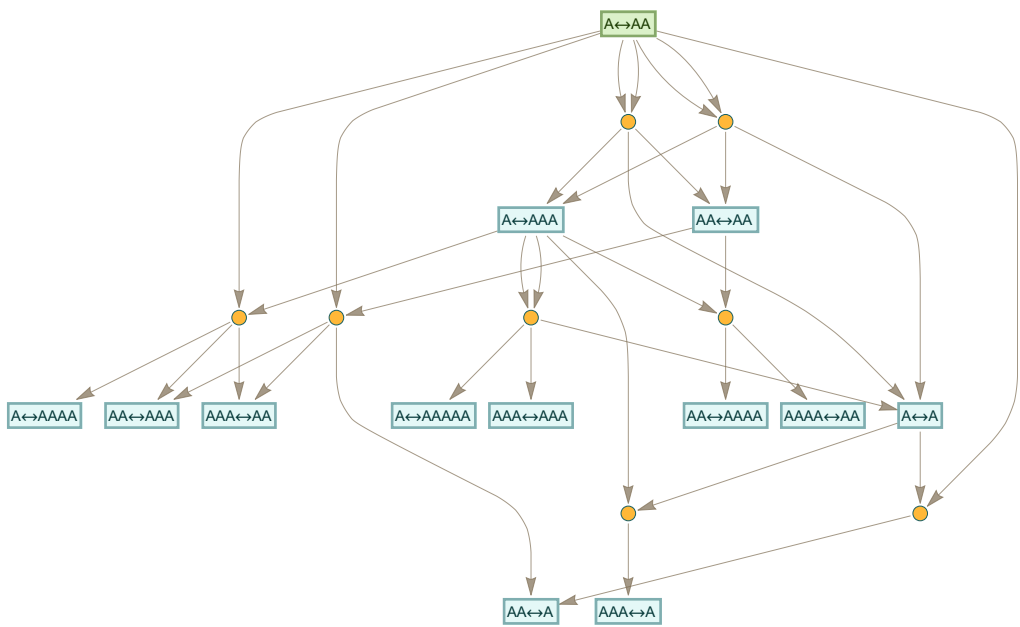

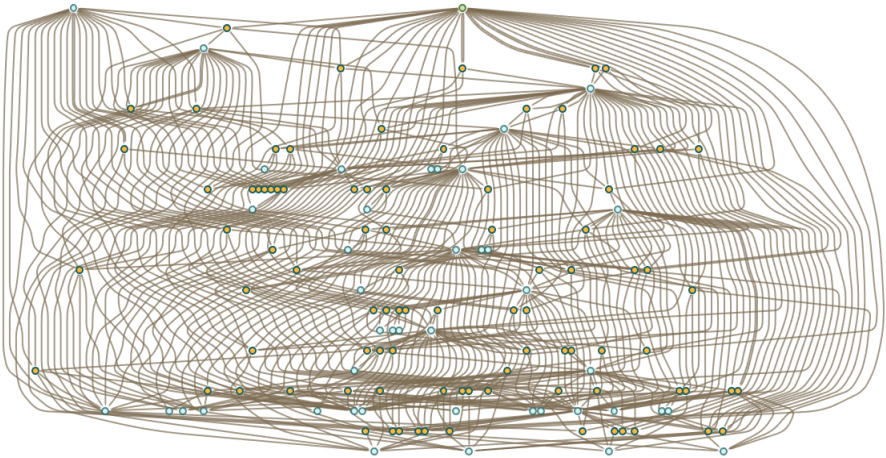

But instead of starting with one axiom and building up a progressively larger entailment cone, let's start with multiple statements, and from each one generate a small entailment cone, say applying each rule at most twice. Here are entailment cones started from several different statements:

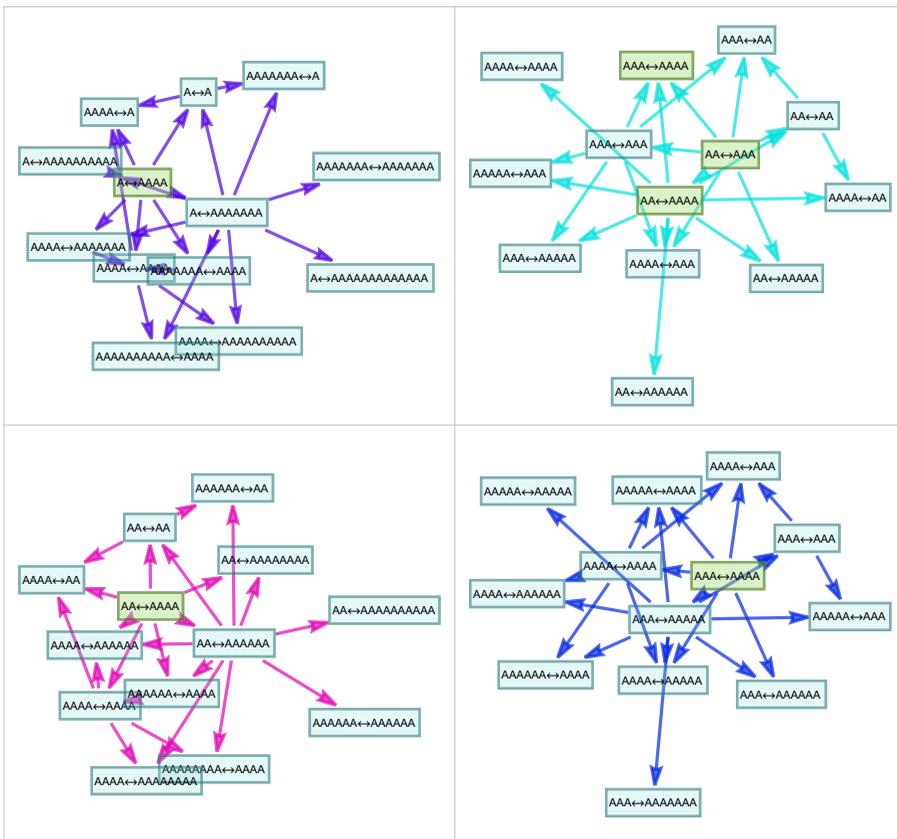

But the crucial point is that these entailment cones overlap—so we can knit them together into an "entailment fabric":

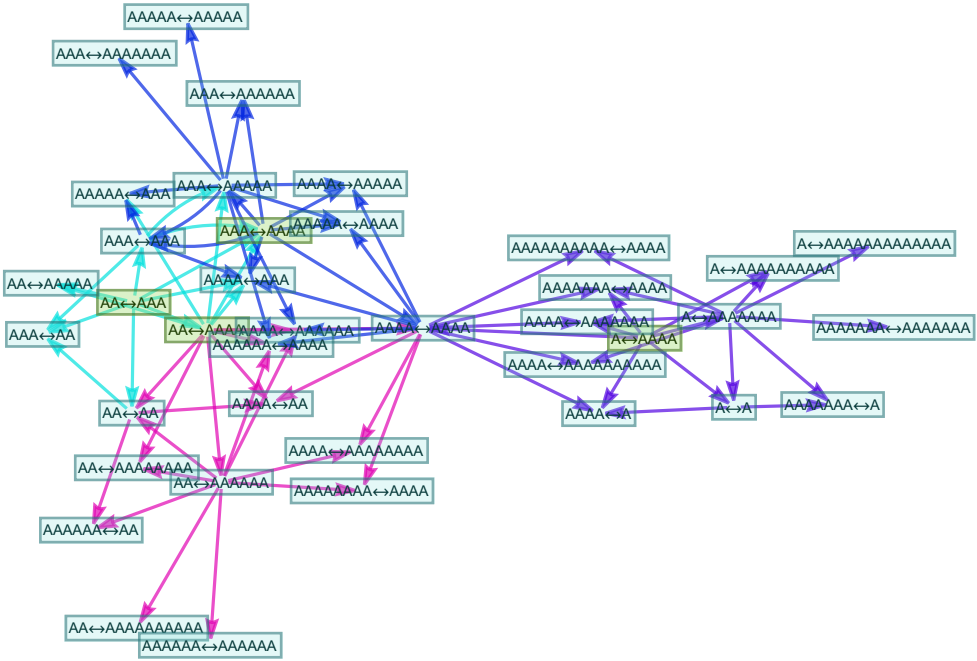

Or with more pieces and another step of entailment:

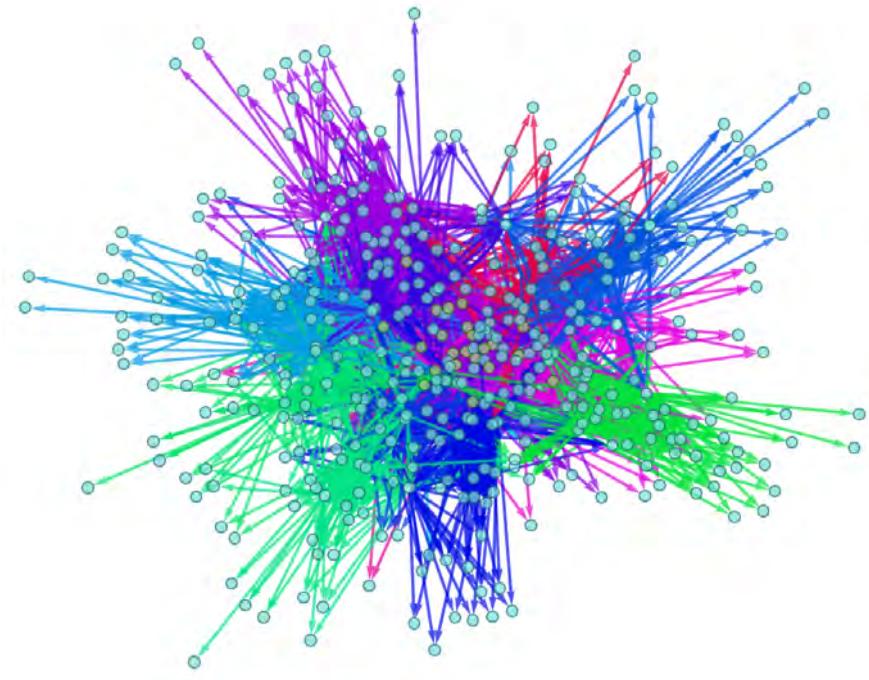

And in a sense this is a "timeless" way to imagine building up mathematics—and metamathematical space. Yes, this structure can in principle be viewed as part of the branchial graph obtained from a slice of an entailment graph (and technically this will be a useful way to think about it). But a different view—closer to the practice of human mathematics—is that it's a "fabric" formed by fitting together many different mathematical statements. It's not something where one's tracking the overall passage of time, and seeing causal connections between things—as one might in "running a program". Rather, it's something where one's fitting pieces together in order to satisfy constraints—as one might in creating a tiling.

Underneath everything is the ruliad. And entailment cones and entailment fabrics can be thought of just as different samplings or slicings of the ruliad. The ruliad is ultimately the entangled limit of all possible computations. But one can think of it as being built up by starting from all possible rules and initial conditions, then running them for an infinite number of steps. An entailment cone is essentially a "slice" of this structure where one's looking at the "time evolution" from a particular rule and initial condition. An entailment fabric is an "orthogonal" slice, looking "at a particular time" across different rules and initial conditions. (And, by the way, rules and initial conditions are essentially equivalent, particularly in an accumulative system.)

One can think of these different slices of the ruliad as being what different kinds of observers will perceive within the ruliad. Entailment cones are essentially what observers who persist through time but are localized in rulial space will perceive. Entailment fabrics are what observers who ignore time but explore more of rulial space will perceive.

Elsewhere I've argued that a crucial part of what makes us perceive the laws of physics we do is that we are observers who consider ourselves to be persistent through time. But now we're seeing that in the way human mathematics is typically done, the "mathematical observer" will be of a different character. And whereas for a physical observer what's crucial is causality through time, for a mathematical observer (at least one who's doing mathematics the way it's usually done) what seems to be crucial is some kind of consistency or coherence across metamathematical space.

In physics it's far from obvious that a persistent observer would be possible. It could be that with all those detailed computationally irreducible processes happening down at the level of atoms of space there might be nothing in the universe that one could consider consistent through time. But the point is that there are certain "coarse-grained" attributes of the behavior that are consistent through time. And it is by concentrating on these that we end up describing things in terms of the laws of physics we know.

There's something very analogous going on in mathematics. The detailed branchial structure of metamathematical space is complicated, and presumably full of computational irreducibility. But once again there are "coarse-grained" attributes that have a certain consistency and coherence across it. And it is on these that we concentrate as human "mathematical observers". And it is in terms of these that we end up being able to do "human-level mathematics"—in effect operating at a "fluid dynamics" level rather than a "molecular dynamics" one.

The possibility of "doing physics in the ruliad" depends crucially on the fact that as physical observers we assume that we have certain persistence and coherence through time. The possibility of "doing mathematics (the way it's usually done) in the ruliad" depends crucially on the fact that as "mathematical observers" we assume that the mathematical statements we consider will have a certain coherence and consistency—or, in effect, that it's possible for us to maintain and grow a coherent body of mathematical knowledge, even as we try to include all sorts of new mathematical statements.

## 20 | The Notion of Truth

Logic was originally conceived as a way to characterize human arguments—in which the concept of "truth" has always seemed quite central. And when logic was applied to the foundations of mathematics, "truth" was also usually assumed to be quite central. But the way we've modeled mathematics here has been much more about what statements can be derived (or entailed) than about any kind of abstract notion of what statements can be "tagged as true". In other words, we've been more concerned with "structurally deriving" that "1+1=2" than in saying that "1+1=2 is true".

But what is the relation between this kind of "constructive derivation" and the logical notion of truth? We might just say that "if we can construct a statement then we should consider it true". And if we're starting from axioms, then in a sense we'll never have an "absolute notion of truth"—because whatever we derive is only "as true as the axioms we started from".

One issue that can come up is that our axioms might be inconsistent—in the sense that from them we can derive two obviously inconsistent statements. But to get further in discussing things like this we really need not only to have a notion of truth, but also a notion of falsity.

In traditional logic it has tended to be assumed that truth and falsity are very much "the same kind of thing"—like 1 and 0. But one feature of our view of mathematics here is that actually truth and falsity seem to have a rather different character. And perhaps this is not surprising—because in a sense if there's one true statement about something there are typically an infinite number of false statements about it. So, for example, the single statement 1+1=2 is true, but the infinite collection of statements 1+1=$n$ for any other $n$ are all false.

There is another aspect to this, discussed since at least the Middle Ages, often under the name of the "principle of explosion": that as soon as one assumes any statement that is false, one can logically derive absolutely any statement at all. In other words, introducing a single "false axiom" will start an explosion that will eventually "blow up everything".

So within our model of mathematics we might say that things are "true" if they can be derived, and are "false" if they lead to an "explosion". But let's say we're given some statement. How can we tell if it's true or false? One thing we can do to find out if it's true is to construct an entailment cone from our axioms and see if the statement appears anywhere in it. Of course, given computational irreducibility there's in general no upper bound on how far we'll need to go to determine this. But now to find out if a statement is false we can imagine introducing the statement as an additional axiom, and then seeing if the entailment cone that's now produced contains an explosion—though once again there'll in general be no upper bound on how far we'll have to go to guarantee that we have a "genuine explosion" on our hands.

So is there any alternative procedure? Potentially the answer is yes: we can just try to see if our statement is somehow equivalent to "true" or "false". But in our model of mathematics where we're just talking about transformations on symbolic expressions, there's no immediate built-in notion of "true" and "false". To talk about these we have to add something. And for example what we can do is to say that "true" is equivalent to what seems like an "obvious tautology" such as $x=x$, or in our computational notation, $x\_ \leftrightarrow x\_$ , while "false" is equivalent to something "obviously explosive", like $x\_ \leftrightarrow y\_$ (or in our particular setup something more like $x\_ \leftrightarrow x\_ \circ y\_$).

But even though something like "Can we find a way to reach $x\_ \leftrightarrow x\_$ from a given statement?" seems like a much more practical question for an actual theorem-proving system than "Can we fish our statement out of a whole entailment cone?", it runs into many of the same issues—in particular that there's no upper limit on the length of path that might be needed.

Soon we'll return to the question of how all this relates to our interpretation of mathematics as a slice of the ruliad—and to the concept of the entailment fabric perceived by a mathematical observer. But to further set the context for what we're doing let's explore how what we've discussed so far relates to things like Gödel's theorem, and to phenomena like incompleteness.

From the setup of basic logic we might assume that we could consider any statement to be either true or false. Or, more precisely, we might think that given a particular axiom system, we should be able to determine whether any statement that can be syntactically constructed with the primitives of that axiom system is true or false. We could explore this by asking whether every statement is either derivable or leads to an explosion—or can be proved equivalent to an "obvious tautology" or to an "obvious explosion".

But as a simple "approximation" to this, let's consider a string rewriting system in which we define a "local negation operation". In particular, let's assume that given a statement like AB↔BBA the "negation" of this statement just exchanges A and B, in this case yielding BA↔AAB.

Now let's ask what statements are generated from a given axiom system. Say we start with AB↔B. After one step of possible substitutions we get

{AAB ↔ B, AB ↔ AB, AB ↔ B, B ↔ B}

while after 2 steps we get:

{AAAAB ↔ B, AAAB ↔ AB, AAAB ↔ B, AAB ↔ AAB, AAB ↔ AB, AAB ↔ B, AB ↔ AB, AB ↔ B, B ↔ B}

And in our setup we're effectively asserting that these are "true" statements. But now let's "negate" the statements, by exchanging A and B. And if we do this, we'll see that there's never a statement where both it and its negation occur. In other words, there's no obvious inconsistency being generated within this axiom system.

But if we consider instead the axiom AB↔BA then this gives:

{AB ↔ AB, AB ↔ BA, BA ↔ BA}

And since this includes both AB↔AB and its "negation" BA↔BA, by our criteria we must consider this axiom system to be inconsistent.

In addition to inconsistency, we can also ask about incompleteness. For all possible statements, does the axiom system eventually generate either the statement or its negation? Or, in other words, can we always decide from the axiom system whether any given statement is true or false?

With our simple assumption about negation, questions of inconsistency and incompleteness become at least in principle very simple to explore. Starting from a given axiom system, we generate its entailment cone, then we ask within this cone what fraction of possible statements, say of a given length, occur.

If the answer is more than 50% we know there's inconsistency, while if the answer is less than 50% that's evidence of incompleteness. So what happens with different possible axiom systems?

Here are some results from *A New Kind of Science*, in each case showing both what amounts to the raw entailment cone (or, in this case, multiway system evolution from "true"), and the number of statements of a given length reached after progressively more steps:

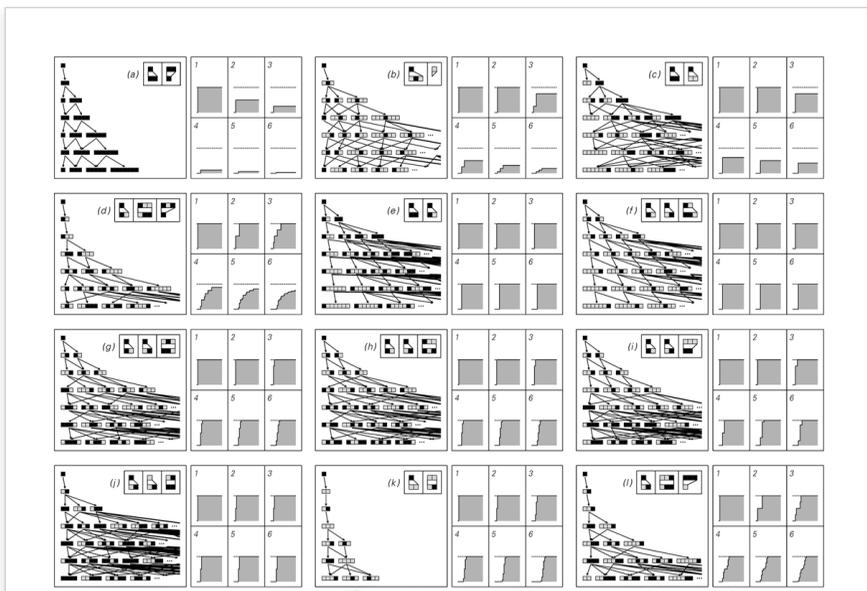

At some level this is all rather straightforward. But from the pictures above we can already get a sense that there's a problem. For most axiom systems the fraction of statements reached of a given length changes as we increase the number of steps in the entailment cone. Sometimes it's straightforward to see what fraction will be achieved even after an infinite number of steps. But often it's not.

And in general we'll run into computational irreducibility—so that in effect the only way to determine whether some particular statement is generated is just to go to ever more steps in the entailment cone and see what happens. In other words, there's no guaranteed-finite way to decide what the ultimate fraction will be—and thus whether or not any given axiom system is inconsistent, or incomplete, or neither.

For some axiom systems it may be possible to tell. But for some axiom systems it's not, in effect because we don't in general know how far we'll have to go to determine whether a given statement is true or not.

A certain amount of additional technical detail is required to reach the standard versions of Gödel's incompleteness theorems. (Note that these theorems were originally stated specifically for the Peano axioms for arithmetic, but the Principle of Computational Equivalence suggests that they're in some sense much more general, and even ubiquitous.) But the important point here is that given an axiom system there may be statements that either can or cannot be reached—but there's no upper bound on the length of path that might be needed to reach them even if one can.

OK, so let's come back to talking about the notion of truth in the context of the ruliad. We've discussed axiom systems that might show inconsistency, or incompleteness—and the difficulty of determining if they do. But the ruliad in a sense contains all possible axiom systems—and generates all possible statements.

So how then can we ever expect to identify which statements are "true" and which are not? When we talked about particular axiom systems, we said that any statement that is generated can be considered true (at least with respect to that axiom system). But in the ruliad every statement is generated. So what criterion can we use to determine which we should consider "true"?

The key idea is any computationally bounded observer (like us) can perceive only a tiny slice of the ruliad. And it's a perfectly meaningful question to ask whether a particular statement occurs within that perceived slice.

One way of picking a "slice" is just to start from a given axiom system, and develop its entailment cone. And with such a slice, the criterion for the truth of a statement is exactly what we discussed above: does the statement occur in the entailment cone?

But how do typical "mathematical observers" actually sample the ruliad? As we discussed in the previous section, it seems to be much more by forming an entailment fabric than by developing a whole entailment cone. And in a sense progress in mathematics can be seen as a process of adding pieces to an entailment fabric: pulling in one mathematical statement after another, and checking that they fit into the fabric.

So what happens if one tries to add a statement that "isn't true"? The basic answer is that it produces an "explosion" in which the entailment fabric can grow to encompass essentially any statement. From the point of view of underlying rules—or the ruliad—there's really nothing wrong with this. But the issue is that it's incompatible with an "observer like us"—or with any realistic idealization of a mathematician.

Our view of a mathematical observer is essentially an entity that accumulates mathematical statements into an entailment fabric. But we assume that the observer is computationally bounded, so in a sense they can only work with a limited collection of statements. So if there's an explosion in an entailment fabric that means the fabric will expand beyond what a mathematical observer can coherently handle. Or, put another way, the only kind of entailment fabrics that a mathematical observer can reasonably consider are ones that "contain no explosions". And in such fabrics, it's reasonable to take the generation or entailment of a statement as a signal that the statement can be considered true.

The ruliad is in a sense a unique and absolute thing. And we might have imagined that it would lead us to a unique and absolute definition of truth in mathematics. But what we've seen is that that's not the case. And instead our notion of truth is something based on how we sample the ruliad as mathematical observers. But now we must explore what this means about what mathematics as we perceive it can be like.

## 21 | What Can Human Mathematics Be Like?

The ruliad in a sense contains all structurally possible mathematics—including all mathematical statements, all axiom systems and everything that follows from them. But mathematics as we humans conceive of it is never the whole ruliad; instead it is always just some tiny part that we as mathematical observers sample.

We might imagine, however, that this would mean that there is in a sense a complete arbitrariness to our mathematics—because in a sense we could just pick any part of the ruliad we want. Yes, we might want to start from a specific axiom system. But we might imagine that that axiom system could be chosen arbitrarily, with no further constraint. And that the mathematics we study can therefore be thought of as an essentially arbitrary choice, determined by its detailed history, and perhaps by cognitive or other features of humans.

But there is a crucial additional issue. When we "sample our mathematics" from the ruliad we do it as mathematical observers and ultimately as humans. And it turns out that even very general features of us as mathematical observers turn out to put strong constraints on what we can sample, and how.

When we discussed physics, we said that the central features of observers are their computational boundedness and their assumption of their own persistence through time. In mathematics, observers are again computationally bounded. But now it is not persistence through time that they assume, but rather a certain coherence of accumulated knowledge.

We can think of a mathematical observer as progressively expanding the entailment fabric that they consider to "represent mathematics". And the question is what they can add to that entailment fabric while still "remaining coherent" as observers. In the previous section, for example, we argued that if the observer adds a statement that can be considered "logically false" then this will lead to an "explosion" in the entailment fabric.

Such a statement is certainly present in the ruliad. But if the observer were to add it, then they wouldn't be able to maintain their coherence—because, whimsically put, their mind would necessarily explode.

In thinking about axiomatic mathematics it's been standard to say that any axiom system that's "reasonable to use" should at least be consistent (even though, yes, for a given axiom system it's in general ultimately undecidable whether this is the case). And certainly consistency is one criterion that we now see is necessary for a "mathematical observer like us". But one can expect that it's not the only criterion.

In other words, although it's perfectly possible to write down any axiom system, and even start generating its entailment cone, only some axiom systems may be compatible with "mathematical observers like us".

And so, for example, something like the Continuum Hypothesis—which is known to be independent of the "established axioms" of set theory—may well have the feature that, say, it has to be assumed to be true in order to get a metamathematical structure compatible with mathematical observers like us.

In the case of physics, we know that the general characteristics of observers lead to certain key perceived features and laws of physics. In statistical mechanics, we're dealing with "coarse-grained observers" who don't trace and decode the paths of individual molecules, and therefore perceive the Second Law of thermodynamics, fluid dynamics, etc. And in our Physics Project we're also dealing with coarse-grained observers who don't track all the details of the atoms of space, but instead perceive space as something coherent and effectively continuous.

And it seems as if in metamathematics there's something very similar going on. As we began to discuss in the very first section above, mathematical observers tend to "coarse grain" metamathematical space. In operational terms, one way they do this is by talking about something like the Pythagorean theorem without always going down to the detailed level of axioms, and for example saying just how real numbers should be defined. And something related is that they tend to concentrate more on mathematical statements and theorems than on their proofs. Later we'll see how in the context of the ruliad there's an even deeper level to which one can go. But the point here is that in actually doing mathematics one tends to operate at the "human scale" of talking about mathematical concepts rather than the "molecular-scale details" of axioms.

But why does this work? Why is one not continually "dragged down" to the detailed axiomatic level—or below? How come it's possible to reason at what we described above as the "fluid dynamics" level, without always having to go down to the detailed "molecular dynamics" level?

The basic claim is that this works for mathematical observers for essentially the same reason as the perception of space works for physical observers. With the "coarse-graining" characteristics of the observer, it's inevitable that the slice of the ruliad they sample will have the kind of coherence that allows them to operate at a higher level. In other words, mathematics can be done "at a human level" for the same basic reason that we have a "human-level experience" of space in physics.

The fact that it works this way depends both on necessary features of the ruliad—and in general of multicomputation—as well as on characteristics of us as observers.

Needless to say, there are "corner cases" where what we've described starts to break down. In physics, for example, the "human-level experience" of space breaks down near spacetime singularities. And in mathematics, there are cases where for example undecidability forces one to take a lower-level, more axiomatic and ultimately more metamathematical view.

But the point is that there are large regions of physical space—and metamathematical space—where these kinds of issues don't come up, and where our assumptions about physical—and mathematical—observers can be maintained. And this is what ultimately allows us to have the "human-scale" views of physics and mathematics that we do.

## 22 | Going below Axiomatic Mathematics

In the traditional view of the foundations of mathematics one imagines that axioms—say stated in terms of symbolic expressions—are in some sense the lowest level of mathematics. But thinking in terms of the ruliad suggests that in fact there is a still-lower "ur level"—a kind of analog of machine code in which everything, including axioms, is broken down into ultimate "raw computation".

Take an axiom like $x \circ y = (y \circ x) \circ y$, or, in more precise computational language:

**x_ ∘ y_ ↔ (y_ ∘ x_) ∘ y_**

Compared to everything we're used to seeing in mathematics this looks simple. But actually it's already got a lot in it. For example, it assumes the notion of a binary operator, which it's in effect naming "∘". And for example it also assumes the notion of variables, and has two distinct pattern variables that are in effect "tagged" with the names $x$ and $y$.

So how can we define what this axiom ultimately "means"? Somehow we have to go from its essentially textual symbolic representation to a piece of actual computation. And, yes, the particular representation we've used here can immediately be interpreted as computation in the Wolfram Language. But the ultimate computational concept we're dealing with is more general than that. And in particular it can exist in any universal computational system.

Different universal computational systems (say particular languages or CPUs or Turing machines) may have different ways to represent computations. But ultimately any computation can be represented in any of them—with the differences in representation being like different "coordinatizations of computation".

And however we represent computations there is one thing we can say for sure: all possible computations are somewhere in the ruliad. Different representations of computations correspond in effect to different coordinatizations of the ruliad. But all computations are ultimately there.

For our Physics Project it's been convenient use a "parametrization of computation" that can be thought of as being based on rewriting of hypergraphs. The elements in these hypergraphs are ultimately purely abstract, but we tend to talk about them as "atoms of space" to indicate the beginnings of our interpretation.

It's perfectly possible to use hypergraph rewriting as the "substrate" for representing axiom systems stated in terms of symbolic expressions. But it's a bit more convenient (though ultimately equivalent) to instead use systems based on expression rewriting—or in effect tree rewriting.

At the outset, one might imagine that different axiom systems would somehow have to be represented by "different rules" in the ruliad. But as one might expect from the phenomenon of universal computation, it's actually perfectly possible to think of different axiom systems as just being specified by different "data" operated on by a single set of rules. There are many rules and structures that we could use. But one set that has the benefit of a century of history are S, K combinators.

The basic concept is to represent everything in terms of "combinator expressions" containing just the two objects S and K. (It's also possible to have just one fundamental object, and indeed S alone may be enough.)

It's worth saying at the outset that when we go this "far down" things get pretty non-human and obscure. Setting things up in terms of axioms may already seem pedantic and low level. But going to a substrate below axioms—that we can think of as getting us to raw "atoms of existence"—will lead us to a whole other level of obscurity and complexity. But if we're going to understand how mathematics can emerge from the ruliad this is where we have to go. And combinators provide us with a more-or-less-concrete example.

Here's an example of a small combinator expression

**S[S[K[S]][S[K[S[K[S]]]][S[K[K]]]][S[S[K[S]][S[K[K]][S]]][K[K]]]**

which corresponds to the "expression tree":

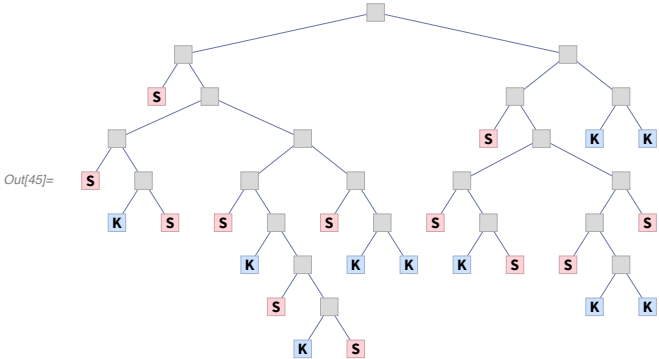

*Out[45]=*

We can write the combinator expression without explicit "function application" [...] by using a (left) application operator •

**S•(S•(K•S)•(S•(K•(S•(K•S)))•(S•(K•K))))•(S•(S•(K•S)•(S•(K•K)•S))•(K•K))**

and it's always unambiguous to omit this operator, yielding the compact representation:

**S(S(KS)(S(K(S(KS)))(S(KK))))(S(S(KS)(S(KK)S))(KK))**

By mapping S, K and the application operator to codewords it's possible to represent this as a simple binary sequence:

1100110010100110010110010100100101011100110010100110010101010010101

But what does our combinator expression mean? The basic combinators are defined to have the rules:

{**S**[x_][y_][z_] → x[z][y[z]], **K**[x_][y_] → x}

These rules on their own don't do anything to our combinator expression. But if we form the expression

**S**[**S**[**K**[**S**]][**S**[**K**[**S**[**K**[**S**]]]][**S**[**K**[**K**]]]]][**S**[**S**[**K**[**S**]][**S**[**K**[**K**]][**S**]]][**K**[**K**]]][c][x][y]

which we can write as

**S(S(KS)(S(K(S(KS)))(S(KK))))(S(S(KS)(S(KK)S))(KK))** $c\,x\,y$

then repeated application of the rules gives:

**S**[**S**[**K**[**S**]]][**S**[**K**[**S**[**K**[**S**]]]][**S**[**K**[**K**]]]][**S**[**S**[**K**[**S**]][**S**[**K**[**K**]][**S**]]][**K**[**K**]]][c][x][y]
**S**[**K**[**S**]][**S**[**K**[**S**[**K**[**S**]]]][**S**[**K**[**K**]]]][c][**S**[**S**[**K**[**S**]][**S**[**K**[**K**]][**S**]]][**K**[**K**]][c]][x][y]
**K**[**S**][c][**S**[**K**[**S**[**K**[**S**]]]][**S**[**K**[**K**]]][c]][**S**[**S**[**K**[**S**]][**S**[**K**[**K**]][**S**]]][**K**[**K**]][c]][x][y]
**S**[**S**[**K**[**S**]]][**S**[**K**[**K**]]][c]][**S**[**S**[**K**[**S**]][**S**[**K**[**K**]][**S**]]][**K**[**K**]][c]][x][y]
**S**[**K**[**S**[**K**[**S**]]]][**S**[**K**[**K**]]][c][x][**S**[**S**[**K**[**S**]][**S**[**K**[**K**]][**S**]]][**K**[**K**]][c][x][y]
**K**[**S**[**K**[**S**]]][c][**S**[**K**[**K**]][c]][x][**S**[**S**[**K**[**S**]][**S**[**K**[**K**]][**S**]]][**K**[**K**]][c][x][y]
**S**[**K**[**S**]][**S**[**K**[**K**]][c]][x][**S**[**S**[**K**[**S**]][**S**[**K**[**K**]][**S**]]][**K**[**K**]][c][x][y]
**K**[**S**][x][**S**[**K**[**K**]][c][x]][**S**[**S**[**K**[**S**]][**S**[**K**[**K**]][**S**]]][**K**[**K**]][c][x][y]
**S**[**S**[**K**[**K**]][c][x]][**S**[**S**[**K**[**S**]][**S**[**K**[**K**]][**S**]]][**K**[**K**]][c][x][y]
**S**[**K**[**K**]][c][x][y][**S**[**S**[**K**[**S**]][**S**[**K**[**K**]][**S**]]][**K**[**K**]][c][x][y]
**K**[**K**][x][c[x]][y][**S**[**S**[**K**[**S**]][**S**[**K**[**K**]][**S**]]][**K**[**K**]][c][x][y]
**K**[c[x]][y][**S**[**S**[**K**[**S**]][**S**[**K**[**K**]][**S**]]][**K**[**K**]][c][x][y]
c[x][**S**[**S**[**K**[**S**]][**S**[**K**[**K**]][**S**]]][**K**[**K**]][c][x][y]]
c[x][**S**[**K**[**S**]][**S**[**K**[**K**]][**S**]][c][**K**[**K**][c]][x][y]]
c[x][**K**[**S**][c][**S**[**K**[**K**]][**S**][c]][**K**[**K**][c]][x][y]]
c[x][**S**[**S**[**K**[**K**]][**S**][c]][**K**[**K**][c]][x][y]]
c[x][**S**[**K**[**K**]][**S**][c][x][**K**[**K**][c][x]][y]]
c[x][**K**[**K**][c][**S**[c]][x][**K**[**K**][c][x]][y]]
c[x][**K**[**S**[c]][x][**K**[**K**][c][x]][y]]
c[x][**S**[c][**K**[**K**][c][x]][y]]
c[x][c[y][**K**[**K**][c][x][y]]]
c[x][c[y][**K**[x][y]]]
c[x][c[y][x]]

We can think of this as "feeding" *c*, *x* and *y* into our combinator expression, then using the "plumbing" defined by the combinator expression to assemble a particular expression in terms of *c*, *x* and *y*.

But what does this expression now mean? Well, that depends on what we think *c*, *x* and *y* mean. We might notice that *c* always appears in the configuration **c**[_][_]. And this means we can interpret it as a binary operator, which we could write in infix form as ∘ so that our expression becomes:

*x* ∘ (*y* ∘ *x*)

And, yes, this is all incredibly low level. But we need to go even further. Right now we're feeding in names like *c*, *x* and *y*. But in the end we want to represent absolutely everything purely in terms of S and K. So we need to get rid of the "human-readable names" and just replace them with "lumps" of S, K combinators that—like the names—get "carried around" when the combinator rules are applied.

We can think about our ultimate expressions in terms of S and K as being like machine code. "One level up" we have assembly language, with the same basic operations, but explicit names. And the idea is that things like axioms—and the laws of inference that apply to them—can be "compiled down" to this assembly language.

But ultimately we can always go further, to the very lowest-level "machine code", in which only S and K ever appear. Within the ruliad as "coordinatized" by S, K combinators, there's an infinite collection of possible combinator expressions. But how do we find ones that "represent something recognizably mathematical"?

As an example let's consider a possible way in which S, K can represent integers, and arithmetic on integers. The basic idea is that an integer *n* can be input as the combinator expression

**Nest[S[S[K[S]][K]], S[K], n]**

which for *n* = 5 gives:

**S[S[K[S]][K]][S[S[K[S]][K]][S[S[K[S]][K]][S[S[K[S]][K]][S[K]]]]]**

But if we now apply this to **[S][K]** what we get reduces to

**S[S[S[S[K]]]]**

which contains 4 S's.

But with this representation of integers it's possible to find combinator expressions that represent arithmetic operations. For example, here's a representation of an addition operator:

**S[K[S]][S[K[S[K[S]]]][S[K[K]]]]**

At the "assembly language" level we might call this **plus**, and apply it to integers *i* and *j* using:

**plus[i][j]**

But at the "pure machine code" level 1+2 can be represented simply by

**S[K[S]][S[K[S[K[S]]]][S[K[K]]]][S[S[K[S]][K]][S[K]]][S[S[K[S]][K]][S[S[K[S]][K]][S[K]]]]**

which when applied to **[S][K]** reduces to the "output representation" of 3:

**S[S[S[K]]]**

As a slightly more elaborate example

**S[K[S[S[K][K]]]][K]**

represents the operation of raising to a power. Then $2^3$ becomes:

**S[K[S[S[K][K]]]][K][S[S[K[S]][K]][S[S[K[S]][K]][S[K]]]][S[S[K[S]][K]][S[S[K[S]][K]][S[S[K[S]][K]][S[K]]]]]**

Applying this to [S][K] repeated application of the combinator rules gives

```
[dense combinator reduction trace]
```

eventually yielding the output representation of 8:

**S[S[S[S[S[S[S[K]]]]]]]**

We could go on and construct any other arithmetic or computational operation we want, all just in terms of the "universal combinators" S and K.

But how should we think about this in terms of our conception of mathematics? Basically what we're seeing is that in the "raw machine code" of S, K combinators it's possible to "find" a representation for something we consider to be a piece of mathematics.

Earlier we talked about starting from structures like axiom systems and then "compiling them down" to raw machine code. But what about just "finding mathematics" in a sense "naturally occurring" in "raw machine code"? We can think of the ruliad as containing "all possible machine code". And somewhere in that machine code must be all the conceivable "structures of mathematics". But the question is: in the wildness of the raw ruliad, what structures can we as mathematical observers successfully pick out?

The situation is quite directly analogous to what happens at multiple levels in physics. Consider for example a fluid full of molecules bouncing around. As we've discussed several times, observers like us usually aren't sensitive to the detailed dynamics of the molecules. But we can still successfully pick out large-scale structures—like overall fluid motions, vortices, etc. And—much like in mathematics—we can talk about physics just at this higher level.

In our Physics Project all this becomes much more extreme. For example, we imagine that space and everything in it is just a giant network of atoms of space. And now within this network we imagine that there are "repeated patterns"—that correspond to things like electrons and quarks and black holes.

In a sense it is the big achievement of natural science to have managed to find these regularities so that we can describe things in terms of them, without always having to go down to the level of atoms of space. But the fact that these are the kinds of regularities we have found is also a statement about us as physical observers.

And the point is that even at the level of the raw ruliad our characteristics as physical observers will inevitably lead us to such regularities. The fact that we are computationally bounded and assume ourselves to have a certain persistence will lead us to consider things that are localized and persistent—that in physics we identify for example as particles.

And it's very much the same thing in mathematics. As mathematical observers we're interested in picking out from the raw ruliad "repeated patterns" that are somehow robust. But now instead of identifying them as particles, we'll identify them as mathematical constructs and definitions. In other words, just as a repeated pattern in the ruliad might in physics be interpreted as an electron, in mathematics a repeated pattern in the ruliad might be interpreted as an integer.

We might think of physics as something "emergent" from the structure of the ruliad, and now we're thinking of mathematics the same way. And of course not only is the "underlying stuff" of the ruliad the same in both cases, but also in both cases it's "observers like us" that are sampling and perceiving things.

There are lots of analogies to the process we're describing of "fishing constructs out of the raw ruliad". As one example, consider the evolution of a ("class 4") cellular automaton in which localized structures emerge:

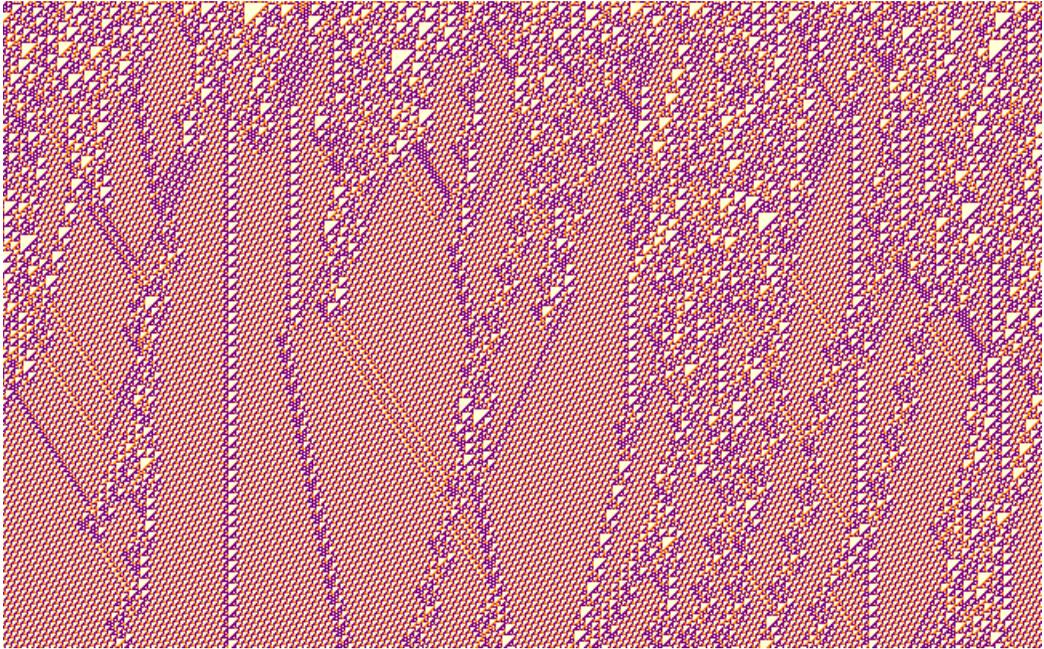

Underneath, just as throughout the ruliad, there's lots of detailed computation going on, with rules repeatedly getting applied to each cell. But out of all this underlying computation we can identify a certain set of persistent structures—which we can use to make a "higher-level description" that may capture the aspects of the behavior that we care about.

Given an "ocean" of S, K combinator expressions, how might we set about "finding mathematics" in them? One straightforward approach is just to identify certain "mathematical properties" we want, and then go searching for S, K combinator expressions that satisfy these.

For example, if we want to "search for (propositional) logic" we first need to pick combinator expressions to symbolically represent "true" and "false". There are many pairs of expressions that will work. As one example, let's pick:

| True | **K** |
|---|---|
| False | **S[K]** |

Now we can just search for combinator expressions which, when applied to all possible pairs of "true" and "false" give truth tables corresponding to particular logical functions. And if we do this, here are examples of the smallest combinator expressions we find:

| And | **S[S][K]** |
|---|---|
| Or | **S[S[S]][S[S[K]]** |
| Implies | **S[S][K[K[K]]]** |
| Nand | **S[S[K[S[S][K[K[K]]]]]][S]** |
| Equal | **S[S][K[S[S[S]][K[K[K]]]][S]]** |

Here's how we can then reproduce the truth table for And:

| And[True][True] | And[True][False] | And[False][True] | And[False][False] |
|---|---|---|---|
| S[S][K][K][K] | S[S][K][K][S[K]] | S[S][K][S[K]][K] | S[S][K][S[K]][S[K]] |
| S[K][K[K]][K] | S[K][K[K]][S[K]] | S[S[K]][K[S[K]]][K] | S[S[K]][K[S[K]]][S[K]] |
| K[K][K[K]][K] | K[S[K]][K[K][S[K]]] | S[K[K]][K[S[K]][K] | S[K[S[K]][K[S[K]][S[K]] |
| K | S[K] | K[K[S[K]][K]][K[K[S[K]][K]]] | K[K[S[K]][S[K]]][S[K][K[S[K]][S[K]]]] |
| | | K[S[K]][K] | K[S[K]][S[K]] |
| | | S[K] | S[K] |

If we just started picking combinator expressions at random, then most of them wouldn't be "interpretable" in terms of this representation of logic. But if we ran across for example

**S[S][K[K[K]]][S[S][K][S[S[S]]][S[S[K]][p][q]][p]][S[S[K[S[S[S]][K[K[K]]]]]][S][p][q]]**

we could recognize in it the combinators for And, Or, etc. that we identified above, and in effect "disassemble" it to give:

Implies[And[Or[p][q]][p]][Nand[p][q]]

It's worth noting, though, that even with the choices we made above for "true" and "false", there's not just a single possible combinator, say for And. Here are a few possibilities:

| S[S][K] | S[S][S[S][K]] | S[S][K[K]][S] | S[S][K[K][S]] | S[S][K[K][K]] | S[K[S][S]][K] | S[K[S][K]][K] |
|---|---|---|---|---|---|---|
| K[S[S][S]][K] | K[S[S]][S[K] | K[S[K]][S] | K[S[S[K]][K] | K[S[S]][K][K] | K[S][K][S][K] | S[S][S[S[K]][K]] |
| S[S[S[K]][S]][K] | S[S][S[K][S[K]] | S[S][S[K][K[K]] | S[S[K][S[S]][K] | S[S[K]][S[S]][K] | S[S[K]][S[S]][K] | S[S[K][S[S]][K]] |
| S[S][K[K][S[S]]] | S[S[K[K]][S]][K] | S[S][K[K[S[K]]]] | S[S][K[K[S[K]]] | S[S][K[K[K][S]]] | S[S][K[K][K[S]]] | S[K[S][S[S][K]] |
| S[K[S][S[S]]][K] | S[K][S[S[S]]][K] | S[K][S[S[S]][K] | S[K[S][K[S]][K] | S[K[S][K[S]][K] | S[K[S][K[K]][K] | S[K[K][S][S][K] |

And there's also nothing unique about the choices for "true" and "false". With the alternative choices

| True | K[K] |
|---|---|
| False | K |

here are the smallest combinator expressions for a few logical functions:

| And | S[S][S[S[K]]] |
|---|---|
| Or | S[S[S[S]]][K] |
| Implies | S[S][S][S[S][K]][S] |
| Nand | S[S[S]][S][S[S][K]][S]] |
| Equal | S[S[S[S[S]]][S[S]]][S[S[K]]]] |

So what can we say in general about the "interpretability" of an arbitrary combinator expression? Obviously any combinator expression does what it does at the level of raw combinators. But the question is whether it can be given a "higher-level"—and potentially "mathematical"—interpretation.

And in a sense this is directly an issue of what a mathematical observer "perceives" in it. Does it contain some kind of robust structure—say a kind of analog for mathematics of a particle in physics?

Axiom systems can be viewed as a particular way to "summarize" certain "raw machine code" in the ruliad. But from the point of a "raw coordinatization of the ruliad" like combinators there doesn't seem to be anything immediately special about them. At least for us humans, however, they do seem to be an obvious "waypoint". Because by distinguishing operators and variables, establishing arities for operators and introducing names for things, they reflect the kind of structure that's familiar from human language.

But now that we think of the ruliad as what's "underneath" both mathematics and physics there's a different path that's suggested. With the axiomatic approach we're effectively trying to leverage human language as a way of summarizing what's going on. But an alternative is to leverage our direct experience of the physical world, and our perception and intuition about things like space. And as we'll discuss later, this is likely in many ways a better "metamodel" of the way pure mathematics is actually practiced by us humans.

In some sense, this goes straight from the "raw machine code" of the ruliad to "human-level mathematics", sidestepping the axiomatic level. But given how much "reductionist" work has already been done in mathematics to represent its results in axiomatic form, there is definitely still great value in seeing how the whole axiomatic setup can be "fished out" of the "raw ruliad".

And there's certainly no lack of complicated technical issues in doing this. As one example, how should one deal with "generated variables"? If one "coordinatizes" the ruliad in terms of something like hypergraph rewriting this is fairly straightforward: it just involves creating new elements or hypergraph nodes (which in physics would be interpreted as atoms of space). But for something like S, K combinators it's a bit more subtle. In the examples we've given above, we have combinators that, when "run", eventually reach a fixed point. But to deal with generated variables we probably also need combinators that never reach fixed points, making it considerably more complicated to identify correspondences with definite symbolic expressions.

Another issue involves rules of entailment, or, in effect, the metalogic of an axiom system. In the full axiomatic setup we want to do things like create token-event graphs, where each event corresponds to an entailment. But what rule of entailment should be used? The underlying rules for S, K combinators, for example, define a particular choice—though they can be used to emulate others. But the ruliad in a sense contains all choices. And, once again, it's up to the observer to "fish out" of the raw ruliad a particular "slice"—which captures not only the axiom system but also the rules of entailment used.

It may be worth mentioning a slightly different existing "reductionist" approach to mathematics: the idea of describing things in terms of types. A type is in effect an equivalence class that characterizes, say, all integers, or all functions from tuples of reals to truth values. But in our terms we can interpret a type as a kind of "template" for our underlying "machine code": we can say that some piece of machine code represents something of a particular type if the machine code matches a particular pattern of some kind. And the issue is then whether that pattern is somehow robust "like a particle" in the raw ruliad.

An important part of what made our Physics Project possible is the idea of going "underneath" space and time and other traditional concepts of physics. And in a sense what we're doing here is something very similar, though for mathematics. We want to go "underneath" concepts like functions and variables, and even the very idea of symbolic expressions. In our Physics Project a convenient "parametrization" of what's "underneath" is a hypergraph made up of elements that we often refer to as "atoms of space". In mathematics we've discussed using combinators as our "parametrization" of what's "underneath".

But what are these "made of"? We can think of them as corresponding to raw elements of metamathematics, or raw elements of computation. But in the end, they're "made of" whatever the ruliad is "made of". And perhaps the best description of the elements of the ruliad is that they are "atoms of existence"—the smallest units of anything, from which everything, in mathematics and physics and elsewhere, must be made.

The atoms of existence aren't bits or points or anything like that. They're something fundamentally lower level that's come into focus only with our Physics Project, and particularly with the identification of the ruliad. And for our purposes here I'll call such atoms of existence "emes" (pronounced "eemes", like phonemes etc.).

Everything in the ruliad is made of emes. The atoms of space in our Physics Project are emes. The nodes in our combinator trees are emes. An eme is a deeply abstract thing. And in a sense all it has is an identity. Every eme is distinct. We could give it a name if we wanted to, but it doesn't intrinsically have one. And in the end the structure of everything is built up simply from relations between emes.

## 23 | The Physicalized Laws of Mathematics

The concept of the ruliad suggests there is a deep connection between the foundations of mathematics and physics. And now that we have discussed how some of the familiar formalism of mathematics can "fit into" the ruliad, we are ready to use the "bridge" provided by the ruliad to start exploring how to apply some of the successes and intuitions of physics to mathematics.

A foundational part of our everyday experience of physics is our perception that we live in continuous space. But our Physics Project implies that at sufficiently small scales space is actually made of discrete elements—and it is only because of the coarse-grained way in which we experience it that we perceive it as continuous.

In mathematics—unlike physics—we've long thought of the foundations as being based on things like symbolic expressions that have a fundamentally discrete structure. Normally, though, the elements of those expressions are, for example, given human-recognizable names (like 2 or **Plus**). But what we saw in the previous section is that these recognizable forms can be thought of as existing in an "anonymous" lower-level substrate made of what we can call atoms of existence or emes.

But the crucial point is that this substrate is directly based on the ruliad. And its structure is identical between the foundations of mathematics and physics. In mathematics the emes aggregate up to give us our universe of mathematical statements. In physics they aggregate up to give us our physical universe.

But now the commonality of underlying "substrate" makes us realize that we should be able to take our experience of physics, and apply it to mathematics. So what is the analog in mathematics of our perception of the continuity of space in physics? We've discussed the idea that we can think of mathematical statements as being laid out in a metamathematical space—or, more specifically, in what we've called an entailment fabric. We initially talked about "coordinatizing" this using axioms, but in the previous section we saw how to go "below axioms" to the level of "pure emes".

When we do mathematics, though, we're sampling this on a much higher level. And just like as physical observers we coarse grain the emes (that we usually call "atoms of space") that make up physical space, so too as "mathematical observers" we coarse grain the emes that make up metamathematical space.

Foundational approaches to mathematics—particularly over the past century or so—have almost always been based on axioms and on their fundamentally discrete symbolic structure. But by going to a lower level and seeing the correspondence with physics we are led to consider what we might think of as a higher-level "experience" of mathematics—operating not at the "molecular dynamics" level of specific axioms and entailments, but rather at what one might call the "fluid dynamics" level of larger-scale concepts.

At the outset one might not have any reason to think that this higher-level approach could consistently be applied. But this is the first big place where ideas from physics can be used. If both physics and mathematics are based on the ruliad, and if our general characteristics as observers apply in both physics and mathematics, then we can expect that similar features will emerge. And in particular, we can expect that our everyday perception of physical space as continuous will carry over to mathematics, or, more accurately, to metamathematical space.

The picture is that we as mathematical observers have a certain "size" in metamathematical space. We identify concepts—like integers or the Pythagorean theorem—as "regions" in the space of possible configurations of emes (and ultimately of slices of the ruliad). At an axiomatic level we might think of ways to capture what a typical mathematician might consider "the same concept" with slightly different formalism (say, different large cardinal axioms or different models of real numbers). But when we get down to the level of emes there'll be vastly more freedom in how we capture a given concept—so that we're in effect using a whole region of "emic space" to do so.

But now the question is what happens if we try to make use of the concept defined by this "region"? Will the "points in the region" behave coherently, or will everything be "shredded", with different specific representations in terms of emes leading to different conclusions?

The expectation is that in most cases it will work much like physical space, and that what we as observers perceive will be quite independent of the detailed underlying behavior at the level of emes. Which is why we can expect to do "higher-level mathematics", without always having to descend to the level of emes, or even axioms.

And this we can consider as the first great "physicalized law of mathematics": that coherent higher-level mathematics is possible for us for the same reason that physical space seems coherent to observers like us.

We've discussed several times before the analogy to the Second Law of thermodynamics— and the way it makes possible a higher-level description of things like fluids for "observers like us". There are certainly cases where the higher-level description breaks down. Some of them may involve specific probes of molecular structure (like Brownian motion). Others may be slightly more "unwitting" (like hypersonic flow).

In our Physics Project we're very interested in where similar breakdowns might occur— because they'd allow us to "see below" the traditional continuum description of space. Potential targets involve various extreme or singular configurations of spacetime, where in effect the "coherent observer" gets "shredded", because different atoms of space "within the observer" do different things.

In mathematics, this kind of "shredding" of the observer will tend to be manifest in the need to "drop below" higher-level mathematical concepts, and go down to a very detailed axiomatic, metamathematical or even eme level—where computational irreducibility and phenomena like undecidability are rampant.

It's worth emphasizing that from the point of view of pure axiomatic mathematics it's not at all obvious that higher-level mathematics should be possible. It could be that there'd be no choice but to work through every axiomatic detail to have any chance of making conclusions in mathematics.

But the point is that we now know there could be exactly the same issue in physics. Because our Physics Project implies that at the lowest level our universe is effectively made of emes that have all sorts of complicated—and computationally irreducible—behavior. Yet we know that we don't have to trace through all the details of this to make conclusions about what will happen in the universe—at least at the level we normally perceive it.

In other words, the fact that we can successfully have a "high-level view" of what happens in physics is something that fundamentally has the same origin as the fact that we can successfully have a high-level view of what happens in mathematics. Both are just features of how observers like us sample the ruliad that underlies both physics and mathematics.

# 24 | Uniformity and Motion in Metamathematical Space

We've discussed how the basic concept of space as we experience it in physics leads us to our first great physicalized law of mathematics—and how this provides for the very possibility of higher-level mathematics. But this is just the beginning of what we can learn from thinking about the correspondences between physical and metamathematical space implied by their common origin in the structure of the ruliad.

A key idea is to think of a limit of mathematics in which one is dealing with so many mathematical statements that one can treat them "in bulk"—as forming something we could consider a continuous metamathematical space. But what might this space be like?

Our experience of physical space is that at our scale and with our means of perception it seems to us for the most part quite simple and uniform. And this is deeply connected to the concept that pure motion is possible in physical space—or, in other words, that it's possible for things to move around in physical space without fundamentally changing their character.

Looked at from the point of view of the atoms of space it's not at all obvious that this should be possible. After all, whenever we move we'll almost inevitably be made up of different atoms of space. But it's fundamental to our character as observers that the features we end up perceiving are ones that have a certain persistence—so that we can imagine that we, and objects around us, can just "move unchanged", at least with respect to those aspects of the objects that we perceive. And this is why, for example, we can discuss laws of mechanics without having to "drop down" to the level of the atoms of space.

So what's the analog of all this in metamathematical space? At the present stage of our physical universe, we seem to be able to experience physical space as having features like being basically three-dimensional. Metamathematical space probably doesn't have such familiar mathematical characterizations. But it seems very likely (and we'll see some evidence of this from empirical metamathematics below) that at the very least we'll perceive metamathematical space as having a certain uniformity or homogeneity.

In our Physics Project we imagine that we can think of physical space as beginning "at the Big Bang" with what amounts to some small collection of atoms of space, but then growing to the vast number of atoms in our current universe through the repeated application of particular rules. But with a small set of rules being applied a vast number of times, it seems almost inevitable that some kind of uniformity must result.

But then the same kind of thing can be expected in metamathematics. In axiomatic mathematics one imagines the mathematical analog of the Big Bang: everything starts from a small collection of axioms, and then expands to a huge number of mathematical statements through repeated application of laws of inference. And from this picture (which gets a bit more elaborate when one considers emes and the full ruliad) one can expect that at least after it's "developed for a while" metamathematical space, like physical space, will have a certain uniformity.

The idea that physical space is somehow uniform is something we take very much for granted, not least because that's our lifelong experience. But the analog of this idea for metamathematical space is something we don't have immediate everyday intuition about—and that in fact may at first seem surprising or even bizarre. But actually what it implies is something that increasingly rings true from modern experience in pure mathematics. Because by saying that metamathematical space is in a sense uniform, we're saying that different parts of it somehow seem similar—or in other words that there's parallelism between what we see in different areas of mathematics, even if they're not "nearby" in terms of entailments.

But this is exactly what, for example, the success of [category theory](category theory) implies. Because it shows us that even in completely different areas of mathematics it makes sense to set up the same basic structures of objects, morphisms and so on. As such, though, category theory defines only the barest outlines of mathematical structure. But what our concept of perceived uniformity in metamathematical space suggests is that there should in fact be closer correspondences between different areas of mathematics.

We can view this as another fundamental "physicalized law of mathematics": that different areas of mathematics should ultimately have structures that are in some deep sense "perceived the same" by mathematical observers. For several centuries we've known there's a certain correspondence between, for example, geometry and algebra. But it's been a major achievement of recent mathematics to identify more and more such correspondences or "dualities".

Often the existence of these has seemed remarkable, and surprising. But what our view of metamathematics here suggests is that this is actually a general physicalized law of mathematics—and that in the end essentially all different areas of mathematics must share a deep structure, at least in some appropriate "bulk metamathematical limit" when enough statements are considered.

But it's one thing to say that two places in metamathematical space are "similar"; it's another to say that "motion between them" is possible. Once again we can make an analogy with physical space. We're used to the idea that we can move around in space, maintaining our identity and structure. But this in a sense requires that we can maintain some kind of continuity of existence on our path between two positions.

In principle it could have been that we would have to be "atomized" at one end, then "reconstituted" at the other end. But our actual experience is that we perceive ourselves to continually exist all the way along the path. In a sense this is just an assumption about how things work that physical observers like us make; but what's nontrivial is that the underlying structure of the ruliad implies that this will always be consistent.

And so we expect it will be in metamathematics. Like a physical observer, the way a mathematical observer operates, it'll be possible to "move" from one area of mathematics to another "at a high level", without being "atomized" along the way. Or, in other words, that a mathematical observer will be able to make correspondences between different areas of mathematics without having to go down to the level of emes to do so.

It's worth realizing that as soon as there's a way of representing mathematics in computational terms the concept of universal computation (and, more tightly, the Principle of Computational Equivalence) implies that at some level there must always be a way to translate between any two mathematical theories, or any two areas of mathematics. But the question is whether it's possible to do this in "high-level mathematical terms" or only at the level of the underlying "computational substrate". And what we're saying is that there's a general physicalized law of mathematics that implies that higher-level translation should be possible.

Thinking about mathematics at a traditional axiomatic level can sometimes obscure this, however. For example, in axiomatic terms we usually think of Peano arithmetic as not being as powerful as ZFC set theory (for example, it lacks transfinite induction)—and so nothing like "dual" to it. But Peano arithmetic can perfectly well support universal computation, so inevitably a "formal emulator" for ZFC set theory can be built in it. But the issue is that to do this essentially requires going down to the "atomic" level and operating not in terms of mathematical constructs but instead directly in terms of "metamathematical" symbolic structure (and, for example, explicitly emulating things like equality predicates).

But the issue, it seems, is that if we think at the traditional axiomatic level, we're not dealing with a "mathematical observer like us". In the analogy we've used above, we're operating at the "molecular dynamics" level, not at the human-scale "fluid dynamics" level. And so we see all sorts of details and issues that ultimately won't be relevant in typical approaches to actually doing pure mathematics.

It's somewhat ironic that our physicalized approach shows this by going below the axiomatic level—to the level of emes and the raw ruliad. But in a sense it's only at this level that there's the uniformity and coherence to conveniently construct a general picture that can encompass observers like us.

Much as with ordinary matter we can say that "everything is made of atoms", we're now saying that everything is "made of computation" (and its structure and behavior is ultimately described by the ruliad). But the crucial idea that emerged from our Physics Project—and that is at the core of what I'm calling the multicomputational paradigm—is that when we ask what observers perceive there is a whole additional level of inexorable structure. And this is what makes it possible to do both human-scale physics and higher-level mathematics—and for there to be what amounts to "pure motion", whether in physical or metamathematical space.

There's another way to think about this, that we alluded to earlier. A key feature of an observer is to have a coherent identity. In physics, that involves having a consistent thread of experience in time. In mathematics, it involves bringing together a consistent view of "what's true" in the space of mathematical statements.

In both cases the observer will in effect involve many separate underlying elements (ultimately, emes). But in order to maintain the observer's view of having a coherent identity, the observer must somehow conflate all these elements, effectively treating them as "the same". In physics, this means "coarse-graining" across physical or branchial (or, in fact, rulial) space. In mathematics, this means "coarse-graining" across metamathematical space—or in effect treating different mathematical statements as "the same".

In practice, there are several ways this happens. First of all, one tends to be more concerned about mathematical results than their proofs, so two statements that have the same form can be considered the same even if the proofs (or other processes) that generated them are different (and indeed this is something we have routinely done in constructing entailment cones here). But there's more. One can also imagine that any statements that entail each other can be considered "the same".

In a simple case, this means that if $a=b$ and $b=c$ then one can always assume $a=c$. But there's a much more general version of this embodied in the univalence axiom of homotopy type theory—that in our terms can be interpreted as saying that mathematical observers consider equivalent things the same.

There's another way that mathematical observers conflate different statements—that's in many ways more important, but less formal. As we mentioned above, when mathematicians talk, say, about the Pythagorean theorem, they typically think they have a definite concept in mind. But at the axiomatic level—and even more so at the level of emes—there are a huge number of different "metamathematical configurations" that are all "considered the same" by the typical working mathematician, or by our "mathematical observer". (At the level of axioms, there might be different axiom systems for real numbers; at the level of emes there might be different ways of representing concepts like addition or equality.)

In a sense we can think of mathematical observers as having a certain "extent" in metamathematical space. And much like human-scale physical observers see only the aggregate effects of huge numbers of atoms of space, so also mathematical observers see only the "aggregate effects" of huge numbers of emes of metamathematical space.

But now the key question is whether a "whole mathematical observer" can "move in metamathematical space" as a single "rigid" entity, or whether it will inevitably be distorted—or shredded—by the structure of metamathematical space. In the next section we'll discuss the analog of gravity—and curvature—in metamathematical space. But our physicalized approach tends to suggest that in "most" of metamathematical space, a typical mathematical observer will be able to "move around freely", implying that there will indeed be paths or "bridges" between different areas of mathematics, that involve only higher-level mathematical constructs, and don't require dropping down to the level of emes and the raw ruliad.

## 25 | Gravitational and Relativistic Effects in Metamathematics

If metamathematical space is like physical space, does that mean that it has analogs of gravity, and relativity? The answer seems to be "yes"—and these provide our next examples of physicalized laws of mathematics.

In the end, we're going to be able to talk about at least gravity in a largely "static" way, referring mostly to the "instantaneous state of metamathematics", captured as an entailment fabric. But in leveraging ideas from physics, it's important to start off formulating things in terms of the analog of time for metamathematics—which is entailment.

As we've discussed above, the entailment cone is the direct analog of the light cone in physics. Starting with some mathematical statement (or, more accurately, some event that transforms it) the forward entailment cone contains all statements (or, more accurately, events) that follow from it. Any possible "instantaneous state of metamathematics" then corresponds to a "transverse slice" through this entailment cone—with the slice in effect being laid out in metamathematical space.

An individual entailment of one statement by another corresponds to a path in the entailment cone, and this path (or, more accurately for accumulative evolution, subgraph) can be thought of as a proof of one statement given another. And in these terms the shortest proof can be thought of as a geodesic in the entailment cone. (In practical mathematics, it's very unlikely one will find—or care about—the strictly shortest proof. But even having a "fairly short proof" will be enough to give the general conclusions we'll discuss here.)

Given a path in the entailment cone, we can imagine projecting it onto a transverse slice, i.e. onto an entailment fabric. Being able to consistently do this depends on having a certain uniformity in the entailment cone, and in the sequence of "metamathematical hypersurfaces" that are defined by whatever "metamathematical reference frame" we're using. But assuming, for example, that underlying computational irreducibility successfully generates a kind of "statistical uniformity" that cannot be "decoded" by the observer, we can expect to have meaningful paths—and geodesics—on entailment fabrics.

But what these geodesics are like then depends on the emergent geometry of entailment fabrics. In physics, the limiting geometry of the analog of this for physical space is presumably a fairly simple 3D manifold. For branchial space, it's more complicated, probably for example being "exponential dimensional". And for metamathematics, the limiting geometry is also undoubtedly more complicated—and almost certainly exponential dimensional.

We've argued that we expect metamathematical space to have a certain perceived uniformity. But what will affect this, and therefore potentially modify the local geometry of the space? The basic answer is exactly the same as in our Physics Project. If there's "more activity" somewhere in an entailment fabric, this will in effect lead to "more local connections", and thus effective "positive local curvature" in the emergent geometry of the network. Needless to say, exactly what "more activity" means is somewhat subtle, especially given that the fabric in which one is looking for this is itself defining the ambient geometry, measures of "area", etc.

In our Physics Project we make things more precise by associating "activity" with energy density, and saying that energy effectively corresponds to the flux of causal edges through spacelike hypersurfaces. So this suggests that we think about an analog of energy in metamathematics: essentially defining it to be the density of update events in the entailment fabric. Or, put another way, energy in metamathematics depends on the "density of proofs" going through a region of metamathematical space, i.e. involving particular "nearby" mathematical statements.

There are lots of caveats, subtleties and details. But the notion that "activity AKA energy" leads to increasing curvature in an emergent geometry is a general feature of the whole multicomputational paradigm that the ruliad captures. And in fact we expect a quantitative relationship between energy density (or, strictly, energy-momentum) and induced curvature of the "transversal space"—that corresponds exactly to Einstein's equations in general relativity. It'll be more difficult to see this in the metamathematical case because metamathematical space is geometrically more complicated—and less familiar—than physical space.

But even at a qualitative level, it seems very helpful to think in terms of physics and spacetime analogies. The basic phenomenon is that geodesics are deflected by the presence of "energy", in effect being "attracted to it". And this is why we can think of regions of higher energy (or energy-momentum/mass)—in physics and in metamathematics—as "generating gravity", and deflecting geodesics towards them. (Needless to say, in metamathematics, as in physics, the vast majority of overall activity is just devoted to knitting together the structure of space, and when gravity is produced, it's from slightly increased activity in a particular region.)

(In our Physics Project, a key result is that the same kind of dependence of "spatial" structure on energy happens not only in physical space, but also in branchial space—where there's a direct analog of general relativity that basically yields the path integral of quantum mechanics.)

What does this mean in metamathematics? Qualitatively, the implication is that "proofs will tend to go through where there's a higher density of proofs". Or, in an analogy, if you want to drive from one place to another, it'll be more efficient if you can do at least part of your journey on a freeway.

One question to ask about metamathematical space is whether one can always get from any place to any other. In other words, starting from one area of mathematics, can one somehow derive all others? A key issue here is whether the area one starts from is computation universal. Propositional logic is not, for example. So if one starts from it, one is essentially trapped, and cannot reach other areas.

But results in mathematical logic have established that most traditional areas of axiomatic mathematics are in fact computation universal (and the Principle of Computational Equivalence suggests that this will be ubiquitous). And given computation universality there will at least be some "proof path". (In a sense this is a reflection of the fact that the ruliad is unique, so everything is connected in "the same ruliad".)

But a big question is whether the "proof path" is "big enough" to be appropriate for a "mathematical observer like us". Can we expect to get from one part of metamathematical space to another without the observer being "shredded"? Will we be able to start from any of a whole collection of places in metamathematical space that are considered "indistinguishably nearby" to a mathematical observer and have all of them "move together" to reach our destination? Or will different specific starting points follow quite different paths—preventing us from having a high-level ("fluid dynamics") description of what's going on, and instead forcing us to drop down to the "molecular dynamics" level?

In practical pure mathematics, this tends to be an issue of whether there is an "elegant proof using high-level concepts", or whether one has to drop down to a very detailed level that's more like low-level computer code, or the output of an automated theorem proving system. And indeed there's a very visceral sense of "shredding" in cases where one's confronted with a proof that consists of page after page of "machine-like details".

But there's another point here as well. If one looks at an individual proof path, it can be computationally irreducible to find out where the path goes, and the question of whether it ever reaches a particular destination can be undecidable. But in most of the current practice of pure mathematics, one's interested in "higher-level conclusions", that are "visible" to a mathematical observer who doesn't resolve individual proof paths.

Later we'll discuss the dichotomy between explorations of computational systems that routinely run into undecidability—and the typical experience of pure mathematics, where undecidability is rarely encountered in practice. But the basic point is that what a typical mathematical observer sees is at the "fluid dynamics level", where the potentially circuitous path of some individual molecule is not relevant.

Of course, by asking specific questions—about metamathematics, or, say, about very specific equations—it's still perfectly possible to force tracing of individual "low-level" proof paths. But this isn't what's typical in current pure mathematical practice. And in a sense we can see this as an extension of our first physicalized law of mathematics: not only is higher-level mathematics possible, but it's ubiquitously so, with the result that, at least in terms of the questions a mathematical observer would readily formulate, phenomena like undecidability are not generically seen.

But even though undecidability may not be directly visible to a mathematical observer, its underlying presence is still crucial in coherently "knitting together" metamathematical space. Because without undecidability, we won't have computation universality and computational irreducibility. But—just like in our Physics Project—computational irreducibility is crucial in producing the low-level apparent randomness that is needed to support any kind of "continuum limit" that allows us to think of large collections of what are ultimately discrete emes as building up some kind of coherent geometrical space.

And when undecidability is not present, one will typically not end up with anything like this kind of coherent space. An extreme example occurs in rewrite systems that eventually terminate—in the sense that they reach a "fixed-point" (or "normal form") state where no more transformations can be applied.

In our Physics Project, this kind of termination can be interpreted as a spacelike singularity at which "time stops" (as at the center of a non-rotating black hole). But in general decidability is associated with "limits on how far paths can go"—just like the limits on causal paths associated with event horizons in physics.

There are many details to work out, but the qualitative picture can be developed further. In physics, the singularity theorems imply that in essence the eventual formation of spacetime singularities is inevitable. And there should be a direct analog in our context that implies the eventual formation of "metamathematical singularities". In qualitative terms, we can expect that the presence of proof density (which is the analog of energy) will "pull in" more proofs until eventually there are so many proofs that one has decidability and a "proof event horizon" is formed.

In a sense this implies that the long-term future of mathematics is strangely similar to the long-term future of our physical universe. In our physical universe, we expect that while the expansion of space may continue, many parts of the universe will form black holes and essentially be "closed off". (At least ignoring expansion in branchial space, and quantum effects in general.)

The analog of this in mathematics is that while there can be continued overall expansion in metamathematical space, more and more parts of it will "burn out" because they've become decidable. In other words, as more work and more proofs get done in a particular area, that area will eventually be "finished"—and there will be no more "open-ended" questions associated with it.

In physics there's sometimes discussion of white holes, which are imagined to effectively be time-reversed black holes, spewing out all possible material that could be captured in a black hole. In metamathematics, a white hole is like a statement that is false and therefore "leads to an explosion". The presence of such an object in metamathematical space will in effect cause observers to be shredded—making it inconsistent with the coherent construction of higher-level mathematics.

We've talked at some length about the "gravitational" structure of metamathematical space. But what about seemingly simpler things like special relativity? In physics, there's a notion of basic, flat spacetime, for which it's easy to construct families of reference frames, and in which parallel trajectories stay parallel. In metamathematics, the analog is presumably metamathematical space in which "parallel proof geodesics" remain "parallel"—so that in effect one can continue "making progress in mathematics" by just "keeping on doing what you've been doing".

And somehow relativistic invariance is associated with the idea that there are many ways to do math, but in the end they're all able to reach the same conclusions. Ultimately this is something one expects as a consequence of fundamental features of the ruliad—and the inevitability of causal invariance in it resulting from the Principle of Computational Equivalence. It's also something that might seem quite familiar from practical mathematics and, say, from the ability to do derivations using different methods—like from either geometry or algebra—and yet still end up with the same conclusions.

So if there's an analog of relativistic invariance, what about analogs of phenomena like time dilation? In our Physics Project time dilation has a rather direct interpretation. To "progress in time" takes a certain amount of computational work. But motion in effect also takes a certain amount of computational work—in essence to continually recreate versions of something in different places. But from the ruliad on up there is ultimately only a certain amount of computational work that can be done—and if computational work is being "used up" on motion, there is less available to devote to progress in time, and so time will effectively run more slowly, leading to the experience of time dilation.

So what is the metamathematical analog of this? Presumably it's that when you do derivations in math you can either stay in one area and directly make progress in that area, or you can "base yourself in some other area" and make progress only by continually translating back and forth. But ultimately that translation process will take computational work, and so will slow down your progress—leading to an analog of time dilation.

In physics, the speed of light defines the maximum amount of motion in space that can occur in a certain amount of time. In metamathematics, the analog is that there's a maximum "translation distance" in metamathematical space that can be "bridged" with a certain amount of derivation. In physics we're used to measuring spatial distance in meters—and time in seconds. In metamathematics we don't yet have familiar units in which to measure, say, distance between mathematical concepts—or, for that matter, "amount of derivation" being done. But with the empirical metamathematics we'll discuss in the next section we actually have the beginnings of a way to define such things, and to use what's been achieved in the history of human mathematics to at least imagine "empirically measuring" what we might call "maximum metamathematical speed".

It should be emphasized that we are only at the very beginning of exploring things like the analogs of relativity in metamathematics. One important piece of formal structure that we haven't really discussed here is causal dependence, and causal graphs. We've talked at length about statements entailing other statements. But we haven't talked about questions like which part of which statement is needed for some event to occur that will entail some other statement. And—while there's no fundamental difficulty in doing it—we haven't concerned ourselves with constructing causal graphs to represent causal relationships and causal dependencies between events.

When it comes to physical observers, there is a very direct interpretation of causal graphs that relates to what a physical observer can experience. But for mathematical observers—where the notion of time is less central—it's less clear just what the interpretation of causal graphs should be. But one certainly expects that they will enter in the construction of any general "observer theory" that characterizes "observers like us" across both physics and mathematics.

## 26 | Empirical Metamathematics

We've discussed the overall structure of metamathematical space, and the general kind of sampling that we humans do of it (as "mathematical observers") when we do mathematics. But what can we learn from the specifics of human mathematics, and the actual mathematical statements that humans have published over the centuries?

We might imagine that these statements are just ones that—as "accidents of history"—humans have "happened to find interesting". But there's definitely more to it—and potentially what's there is a rich source of "empirical data" relevant to our physicalized laws of mathematics, and to what amounts to their "experimental validation".

The situation with "human settlements" in metamathematical space is in a sense rather similar to the situation with human settlements in physical space. If we look at where humans have chosen to live and build cities, we'll find a bunch of locations in 3D space. The details of where these are depend on history and many factors. But there's a clear overarching theme, that's in a sense a direct reflection of underlying physics: all the locations lie on the more-or-less spherical surface of the Earth.

It's not so straightforward to see what's going on in the metamathematical case, not least because any notion of coordinatization seems to be much more complicated for metamathematical space than for physical space. But we can still begin by doing "empirical metamathematics" and asking questions about for example what amounts to where in metamathematical space we humans have so far established ourselves. And as a first example, let's consider Boolean algebra.

Even to talk about something called "Boolean algebra" we have to be operating at a level far above the raw ruliad—where we've already implicitly aggregated vast numbers of emes to form notions of, for example, variables and logical operations.

But once we're at this level we can "survey" metamathematical space just by enumerating possible symbolic statements that can be created using the operations we've set up for Boolean algebra (here And ∧, Or ∨ and Not ◌̄):

| $a = b$ | $a = \bar{a}$ | $a = \bar{b}$ | $\bar{a} = \bar{b}$ | $a = a \wedge a$ | $\bar{a} = a \wedge a$ | $a = a \vee a$ | $\bar{a} = a \vee a$ |
|---|---|---|---|---|---|---|---|
| $a \wedge a = a \vee a$ | $a = a \wedge b$ | $\bar{a} = a \wedge b$ | $a \wedge a = a \wedge b$ | $a \vee a = a \wedge b$ | $a = a \vee b$ | $\bar{a} = a \vee b$ | $a \wedge a = a \vee b$ |
| $a \vee a = a \vee b$ | $a \wedge b = a \vee b$ | $a \wedge b = a \wedge c$ | $a \vee b = a \wedge c$ | $a \wedge b = a \vee c$ | $a \vee b = a \vee c$ | $a = b \wedge a$ | $\bar{a} = b \wedge a$ |
| $a \wedge a = b \wedge a$ | $a \vee a = b \wedge a$ | $a \wedge b = b \wedge a$ | $a \vee b = b \wedge a$ | $a = b \vee a$ | $\bar{a} = b \vee a$ | $a \wedge a = b \vee a$ | $a \vee a = b \vee a$ |
| $a \wedge b = b \vee a$ | $a \vee b = b \vee a$ | $a = b \wedge b$ | $\bar{a} = b \wedge b$ | $a \wedge a = b \wedge b$ | $a \vee a = b \wedge b$ | $a \wedge b = b \wedge b$ | $a \vee b = b \wedge b$ |
| $a = b \vee b$ | $\bar{a} = b \vee b$ | $a \wedge a = b \vee b$ | $a \vee a = b \vee b$ | $a \wedge b = b \vee b$ | $a \vee b = b \vee b$ | $a = b \wedge c$ | $\bar{a} = b \wedge c$ |
| $a \wedge a = b \wedge c$ | $a \vee a = b \wedge c$ | $a \wedge b = b \wedge c$ | $a \vee b = b \wedge c$ | $a = b \vee c$ | $\bar{a} = b \vee c$ | $a \wedge a = b \vee c$ | $a \vee a = b \vee c$ |
| $a \wedge b = b \vee c$ | $a \vee b = b \vee c$ | $a \wedge b = c \wedge a$ | $a \vee b = c \wedge a$ | $a \wedge b = c \vee a$ | $a \vee b = c \vee a$ | $a \wedge b = c \wedge b$ | $a \vee b = c \wedge b$ |
| $a \wedge b = c \vee b$ | $a \vee b = c \vee b$ | $a \wedge b = c \wedge c$ | $a \vee b = c \wedge c$ | $a \wedge b = c \vee c$ | $a \vee b = c \vee c$ | $a = \bar{\bar{a}}$ | $\bar{a} = \bar{\bar{a}}$ |
| $a \wedge a = \bar{\bar{a}}$ | $a \vee a = \bar{\bar{a}}$ | $a \wedge b = \bar{\bar{a}}$ | $a \vee b = \bar{\bar{a}}$ | $a = \bar{\bar{b}}$ | $\bar{a} = \bar{\bar{b}}$ | $a \wedge a = \bar{\bar{b}}$ | $a \vee a = \bar{\bar{b}}$ |

But so far these are just raw, structural statements. To connect with actual Boolean algebra we must pick out which of these can be derived from the axioms of Boolean algebra, or, put another way, which of them are in the entailment cone of these axioms:

| $a = b$ | $a = \bar{a}$ | $a = \bar{b}$ | $\bar{a} = \bar{b}$ | $a = a \wedge a$ | $\bar{a} = a \wedge a$ | $a = a \vee a$ | $\bar{a} = a \vee a$ |
|---|---|---|---|---|---|---|---|
| $a \wedge a = a \vee a$ | $a = a \wedge b$ | $\bar{a} = a \wedge b$ | $a \wedge a = a \wedge b$ | $a \vee a = a \wedge b$ | $a = a \vee b$ | $\bar{a} = a \vee b$ | $a \wedge a = a \vee b$ |
| $a \vee a = a \vee b$ | $a \wedge b = a \vee b$ | $a \wedge b = a \wedge c$ | $a \vee b = a \wedge c$ | $a \wedge b = a \vee c$ | $a \vee b = a \vee c$ | $a = b \wedge a$ | $\bar{a} = b \wedge a$ |
| $a \wedge a = b \wedge a$ | $a \vee a = b \wedge a$ | $a \wedge b = b \wedge a$ | $a \vee b = b \wedge a$ | $a = b \vee a$ | $\bar{a} = b \vee a$ | $a \wedge a = b \vee a$ | $a \vee a = b \vee a$ |
| $a \wedge b = b \vee a$ | $a \vee b = b \vee a$ | $a = b \wedge b$ | $\bar{a} = b \wedge b$ | $a \wedge a = b \wedge b$ | $a \vee a = b \wedge b$ | $a \wedge b = b \wedge b$ | $a \vee b = b \wedge b$ |
| $a = b \vee b$ | $\bar{a} = b \vee b$ | $a \wedge a = b \vee b$ | $a \vee a = b \vee b$ | $a \wedge b = b \vee b$ | $a \vee b = b \vee b$ | $a = b \wedge c$ | $\bar{a} = b \wedge c$ |
| $a \wedge a = b \wedge c$ | $a \vee a = b \wedge c$ | $a \wedge b = b \wedge c$ | $a \vee b = b \wedge c$ | $a = b \vee c$ | $\bar{a} = b \vee c$ | $a \wedge a = b \vee c$ | $a \vee a = b \vee c$ |
| $a \wedge b = b \vee c$ | $a \vee b = b \vee c$ | $a \wedge b = c \wedge a$ | $a \vee b = c \wedge a$ | $a \wedge b = c \vee a$ | $a \vee b = c \vee a$ | $a \wedge b = c \wedge b$ | $a \vee b = c \wedge b$ |
| $a \wedge b = c \vee b$ | $a \vee b = c \vee b$ | $a \wedge b = c \wedge c$ | $a \vee b = c \wedge c$ | $a \wedge b = c \vee c$ | $a \vee b = c \vee c$ | $a = \bar{\bar{a}}$ | $\bar{a} = \bar{\bar{a}}$ |
| $a \wedge a = \bar{\bar{a}}$ | $a \vee a = \bar{\bar{a}}$ | $a \wedge b = \bar{\bar{a}}$ | $a \vee b = \bar{\bar{a}}$ | $a = \bar{\bar{b}}$ | $\bar{a} = \bar{\bar{b}}$ | $a \wedge a = \bar{\bar{b}}$ | $a \vee a = \bar{\bar{b}}$ |

Of all possible statements, it's only an exponentially small fraction that turn out to be derivable:

| variables | $a$ | | $a, b$ | | $a, b, c$ | |
|---|---|---|---|---|---|---|
| size | total | theorems | total | theorems | total | theorems |
| 2 | 1 | 0 | 4 | 0 | 4 | 0 |
| 3 | 33 | 8 | 132 | 8 | 164 | 8 |
| 4 | 673 | 164 | 5348 | 404 | 9316 | 404 |
| 5 | 15 009 | 3620 | 234 724 | 17 940 | 597 092 | 22 292 |

But in the case of Boolean algebra, we can readily collect such statements:

| | | | | | |
|---|---|---|---|---|---|
| $a = a \wedge a$ | $a = a \vee a$ | $a \wedge a = a \vee a$ | $a \wedge b = b \wedge a$ | $a \vee b = b \vee a$ | $a = \bar{\bar{a}}$ |
| $a \wedge a = \bar{\bar{a}}$ | $a \vee a = \bar{\bar{a}}$ | $\bar{a} \wedge a = a \wedge \bar{a}$ | $\bar{a} \vee a = a \vee \bar{a}$ | $\bar{a} = \overline{a \wedge a}$ | $\bar{a} = \overline{a \vee a}$ |
| $\overline{a \wedge a} = \overline{a \vee a}$ | $a \wedge \bar{b} = \bar{b} \wedge a$ | $\bar{a} \wedge b = b \wedge \bar{a}$ | $a \vee \bar{b} = \bar{b} \vee a$ | $\bar{a} \vee b = b \vee \bar{a}$ | $\overline{a \wedge b} = \overline{b \wedge a}$ |
| $\overline{a \vee b} = \overline{b \vee a}$ | $\bar{a} \wedge a = \bar{b} \wedge b$ | $a \wedge \bar{a} = \bar{b} \wedge b$ | $\bar{a} \wedge a = b \wedge \bar{b}$ | $a \wedge \bar{a} = b \wedge \bar{b}$ | $\bar{a} \vee a = \bar{b} \vee b$ |
| $a \vee \bar{a} = \bar{b} \vee b$ | $\bar{a} \vee a = b \vee \bar{b}$ | $a \vee \bar{a} = b \vee \bar{b}$ | $\bar{\bar{a}} = \bar{a} \wedge \bar{a}$ | $\overline{a \wedge a} = \bar{a} \wedge \bar{a}$ | $\overline{a \vee a} = \bar{a} \wedge \bar{a}$ |
| $\bar{\bar{a}} = \bar{a} \vee \bar{a}$ | $\overline{a \wedge a} = \bar{a} \vee \bar{a}$ | $\overline{a \vee a} = \bar{a} \vee \bar{a}$ | $\bar{a} \wedge a = \bar{a} \vee \bar{a}$ | $\overline{a \vee b} = \bar{a} \wedge \bar{b}$ | $\overline{a \wedge b} = \bar{a} \vee \bar{b}$ |
| $\overline{a \vee b} = \bar{b} \wedge \bar{a}$ | $\overline{a \wedge b} = \bar{b} \wedge \bar{a}$ | $\overline{a \wedge b} = \bar{b} \vee \bar{a}$ | $\overline{a \vee b} = \bar{b} \vee \bar{a}$ | $a = (a \wedge a) \wedge a$ | $a \wedge a = (a \wedge a) \wedge a$ |
| $a \vee a = (a \wedge a) \wedge a$ | $\bar{\bar{a}} = (a \wedge a) \wedge a$ | $a = (a \vee a) \wedge a$ | $a \wedge a = (a \vee a) \wedge a$ | $a \vee a = (a \vee a) \wedge a$ | $\bar{\bar{a}} = (a \vee a) \wedge a$ |
| $(a \wedge a) \wedge a = (a \vee a) \wedge a$ | $a = a \wedge (a \wedge a)$ | $a \wedge a = a \wedge (a \wedge a)$ | $a \vee a = a \wedge (a \wedge a)$ | $\bar{\bar{a}} = a \wedge (a \wedge a)$ | $(a \wedge a) \wedge a = a \wedge (a \wedge a)$ |
| $(a \vee a) \wedge a = a \wedge (a \wedge a)$ | $a = a \wedge (a \vee a)$ | $a \wedge a = a \wedge (a \vee a)$ | $a \vee a = a \wedge (a \vee a)$ | $\bar{\bar{a}} = a \wedge (a \vee a)$ | $(a \wedge a) \wedge a = a \wedge (a \vee a)$ |

We've typically explored entailment cones by looking at slices consisting of collections of theorems generated after a specified number of proof steps. But here we're making a very different sampling of the entailment cone—looking in effect instead at theorems in order of their structural complexity as symbolic expressions.

In doing this kind of systematic enumeration we're in a sense operating at a "finer level of granularity" than typical human mathematics. Yes, these are all "true theorems". But mostly they're not theorems that a human mathematician would ever write down, or specifically "consider interesting". And for example only a small fraction of them have historically been given names—and are called out in typical logic textbooks:

| $a = a \land a$ | $a = a \lor a$ | $a \land a = a \lor a$ | $a \land b = b \land a$ | $a \lor b = b \lor a$ | $a = \overline{\overline{a}}$ | $a \land a = \overline{\overline{a}}$ |
|---|---|---|---|---|---|---|
| $a \lor a = \overline{\overline{a}}$ | $\overline{\overline{a}} = \overline{a \land a}$ | $\overline{\overline{a}} = \overline{a \lor a}$ | $a = (a \land a) \land a$ | $a = (a \lor a) \land a$ | $a = a \land (a \lor a)$ | $a = a \land (a \lor a)$ |
| $a = (a \land a) \lor a$ | $a = (a \lor a) \lor a$ | $a = a \lor (a \land a)$ | $a = a \lor (a \lor a)$ | $a = a \land (a \lor b)$ | $a = a \lor (a \land b)$ | $a = (a \lor b) \land a$ |
| $a = a \land (b \lor a)$ | $a = (a \land b) \lor a$ | $a = a \lor (b \land a)$ | $a = (b \lor a) \land a$ | $a = (b \land a) \lor a$ | $\overline{\overline{a}} = \overline{a \land a}$ | $\overline{\overline{a}} = \overline{a \lor a}$ |
| $\overline{a} \land a = a \land \overline{a}$ | $\overline{a} \lor a = a \lor \overline{a}$ | $\overline{a \land a} = \overline{a \lor a}$ | $\overline{\overline{a}} = (a \land a) \land a$ | $\overline{\overline{a}} = (a \lor a) \land a$ | $\overline{\overline{a}} = a \land (a \lor a)$ | $\overline{\overline{a}} = a \land (a \lor a)$ |
| $\overline{\overline{a}} = (a \land a) \lor a$ | $\overline{\overline{a}} = (a \lor a) \lor a$ | $\overline{\overline{a}} = a \lor (a \land a)$ | $\overline{\overline{a}} = a \lor (a \lor a)$ | $\overline{\overline{a}} = a \land (a \lor b)$ | $\overline{\overline{a}} = a \lor (a \land b)$ | $\overline{\overline{a}} = (a \lor b) \land a$ |
| $\overline{\overline{a}} = a \land (b \lor a)$ | $\overline{\overline{a}} = (a \land b) \lor a$ | $\overline{\overline{a}} = a \lor (b \land a)$ | $\overline{a} \land a = \overline{b} \land b$ | $a \land \overline{a} = \overline{b} \land b$ | $\overline{a} \land a = b \land \overline{b}$ | $a \land \overline{a} = b \land \overline{b}$ |
| $\overline{a} \lor a = \overline{b} \lor b$ | $a \lor \overline{a} = \overline{b} \lor b$ | $\overline{a} \lor a = b \lor \overline{b}$ | $a \lor \overline{a} = b \lor \overline{b}$ | $\overline{a} = (b \lor a) \land a$ | $\overline{a} = (b \land a) \lor a$ | $a \land \overline{b} = \overline{b} \land a$ |
| $\overline{a} \land b = b \land \overline{a}$ | $a \lor \overline{b} = \overline{b} \lor a$ | $\overline{a} \lor b = b \lor \overline{a}$ | $\overline{a \land b} = \overline{b \land a}$ | $\overline{a \lor b} = \overline{b \lor a}$ | $a \land a = (a \land a) \land a$ | $a \lor a = (a \land a) \land a$ |
| $a \land a = (a \lor a) \land a$ | $a \lor a = (a \lor a) \land a$ | $a \land a = a \land (a \land a)$ | $a \lor a = a \land (a \land a)$ | $a \land a = a \land (a \lor a)$ | $a \lor a = a \land (a \lor a)$ | $a \land a = (a \land a) \lor a$ |
| $a \lor a = (a \land a) \lor a$ | $a \land a = (a \lor a) \lor a$ | $a \lor a = (a \lor a) \lor a$ | $a \land a = a \lor (a \land a)$ | $a \lor a = a \lor (a \land a)$ | $a \land a = a \lor (a \lor a)$ | $a \lor a = a \lor (a \lor a)$ |
| $a = (a \land a) \land (a \land a)$ | $a = (a \land a) \land (a \lor a)$ | $a = (a \lor a) \land (a \land a)$ | $a = (a \lor a) \land (a \lor a)$ | $a = (a \land a) \lor (a \land a)$ | $a = (a \land a) \lor (a \lor a)$ | $a = (a \lor a) \lor (a \land a)$ |
| $a = (a \lor a) \lor (a \lor a)$ | $a \land a = a \land (a \lor b)$ | $a \lor a = a \land (a \lor b)$ | $a \land a = a \lor (a \land b)$ | $a \lor a = a \lor (a \land b)$ | $a = (a \land a) \land (a \lor b)$ | $a = (a \lor a) \land (a \lor b)$ |
| $a = (a \land a) \lor (a \land b)$ | $a = (a \lor a) \lor (a \land b)$ | $a \land a = a \land (b \lor a)$ | $a \lor a = (a \lor b) \land a$ | $a \land a = a \lor (b \land a)$ | $a \lor a = a \lor (b \land a)$ | $a \land a = (a \lor a) \lor a$ |
| $a \lor a = (a \land b) \lor a$ | $a \land a = a \lor (b \land a)$ | $a \lor a = a \lor (b \land a)$ | $a \land a = (b \lor a) \land a$ | $a \lor a = (b \lor a) \land a$ | $a = (a \lor b) \land (a \land a)$ | $a = (a \lor b) \land (a \lor a)$ |
| $a \land a = (b \land a) \land a$ | $a \lor a = (b \land a) \land a$ | $a \land a = (b \land a) \lor a$ | $a \lor a = (b \land a) \lor a$ | $a = (b \land a) \lor (a \land a)$ | $a = (b \land a) \lor (a \lor a)$ | $a = (a \land b) \lor (a \land a)$ |
| $a = (a \land b) \lor (a \lor a)$ | $a = (b \lor a) \land (a \land a)$ | $a = (b \lor a) \land (a \lor a)$ | $a = (b \lor a) \land (a \land a)$ | $a = (b \land a) \lor (a \lor a)$ | $a \land b = (a \land a) \land b$ | $a \land b = (a \lor a) \land b$ |
| $a \land b = a \land (a \land b)$ | $a \lor b = (a \land a) \lor b$ | $a \lor b = (a \lor a) \lor b$ | $a \lor b = (a \lor a) \lor b$ | $a \lor b = a \lor (a \lor b)$ | $a \land b = (a \lor b) \land a$ | $a \land b = a \land (b \land a)$ |
| $a \lor b = a \lor (b \lor a)$ | $a \land b = (a \land b) \land b$ | $a \land b = a \land (b \land b)$ | $a \land b = a \land (b \lor b)$ | $a \lor b = (b \lor a) \lor b$ | $a \lor b = a \lor (b \lor b)$ | $a \lor b = a \lor (b \lor b)$ |
| $a \land b = (b \land a) \land a$ | $a \land b = b \land (a \land a)$ | $a \land b = b \land (a \lor a)$ | $a \lor b = (b \lor a) \lor a$ | $a \lor b = b \lor (a \land a)$ | $a \lor b = b \lor (a \lor a)$ | $a \land b = (b \land a) \land b$ |
| $a \land b = b \land (a \land b)$ | $a \lor b = (b \lor a) \lor b$ | $a \lor b = b \lor (a \lor b)$ | $a \land b = (b \land b) \land a$ | $a \land b = (b \lor b) \land a$ | $a \land b = b \land (b \land a)$ | $a \lor b = (b \lor b) \lor a$ |
| $a \lor b = (b \lor b) \lor a$ | $a \lor b = b \lor (b \lor a)$ | $\overline{a \land a} = \overline{a} \land \overline{a}$ | $\overline{a \lor a} = \overline{a} \land \overline{a}$ | $\overline{a \land a} = \overline{a} \lor \overline{a}$ | $\overline{a \lor a} = \overline{a} \lor \overline{a}$ | $\overline{a \lor b} = \overline{a} \land \overline{b}$ |
| $\overline{a \land b} = \overline{a} \lor \overline{b}$ | $\overline{a \lor b} = \overline{b} \land \overline{a}$ | $\overline{a \land b} = \overline{b} \lor \overline{a}$ | $\overline{a \land a} = \overline{a} \lor \overline{a}$ | $\overline{a \land b} = \overline{b} \land \overline{a}$ | $\overline{a \lor b} = \overline{b} \lor \overline{a}$ | $\overline{a} \land a = (a \land a) \land \overline{a}$ |
| $a \land \overline{a} = (a \land a) \land \overline{a}$ | $\overline{a} \land a = (a \lor a) \land \overline{a}$ | $a \land \overline{a} = (a \lor a) \land \overline{a}$ | $\overline{a} \land a = \overline{a} \land (a \land a)$ | $a \land \overline{a} = \overline{a} \land (a \land a)$ | $\overline{a} \land a = \overline{a} \land (a \lor a)$ | $a \land \overline{a} = \overline{a} \land (a \lor a)$ |
| $\overline{a} \lor a = (a \land a) \lor \overline{a}$ | $a \lor \overline{a} = (a \land a) \lor \overline{a}$ | $\overline{a} \lor a = (a \lor a) \lor \overline{a}$ | $a \lor \overline{a} = (a \lor a) \lor \overline{a}$ | $\overline{a} \lor a = \overline{a} \lor (a \land a)$ | $a \lor \overline{a} = \overline{a} \lor (a \land a)$ | $\overline{a} \lor a = \overline{a} \lor (a \lor a)$ |

⋮ 88 lines

| $a \lor b = (b \land b) \lor (b \lor a)$ | $a \lor b = (b \lor b) \lor (b \lor a)$ | $(a \land b) \land b = (b \land b) \land b$ | $(a \land b) \lor b = (b \land b) \land b$ | $(a \lor b) \land b = (b \lor b) \land b$ | $(a \land b) \lor b = (b \lor b) \land b$ |
|---|---|---|---|---|---|
| $(a \lor b) \land b = b \land (b \land b)$ | $(a \land b) \lor b = b \land (b \land b)$ | $(a \lor b) \land b = b \land (b \lor b)$ | $(a \land b) \lor b = b \land (b \lor b)$ | $(a \lor b) \land b = (b \lor b) \land b$ | $(a \land b) \lor b = (b \lor b) \lor b$ |
| $(a \lor b) \land b = (b \lor b) \lor b$ | $(a \land b) \lor b = (b \lor b) \lor b$ | $(a \lor b) \land b = b \lor (b \land b)$ | $(a \land b) \lor b = b \lor (b \land b)$ | $(a \lor b) \land b = b \lor (b \lor b)$ | $(a \land b) \lor b = b \lor (b \lor b)$ |
| $(a \lor b) \land b = b \land (b \lor c)$ | $(a \land b) \lor b = b \land (b \lor c)$ | $(a \lor b) \land b = b \lor (b \land c)$ | $(a \land b) \lor b = b \lor (b \land c)$ | $(a \lor b) \land b = (b \lor c) \land b$ | $(a \land b) \lor b = (b \lor c) \land b$ |
| $(a \lor b) \land b = b \land (c \lor b)$ | $(a \land b) \lor b = b \land (c \lor b)$ | $(a \lor b) \land b = b \lor (c \land b)$ | $(a \land b) \lor b = b \lor (c \land b)$ | $(a \lor b) \land b = b \lor (c \land b)$ | $(a \land b) \lor b = b \lor (c \land b)$ |
| $a \land b = (b \lor c) \land (a \lor b)$ | $a \lor b = (b \lor c) \land (a \lor b)$ | $a \land b = (b \lor c) \land (b \land a)$ | $a \lor b = (b \lor c) \lor (b \land a)$ | $(a \lor b) \land b = (c \lor a) \land b$ | $(a \land b) \lor b = (c \lor b) \land b$ |
| $(a \lor b) \land b = (c \lor b) \lor b$ | $(a \land b) \lor b = (c \lor b) \lor b$ | $a \land b = (c \lor a) \land (a \land b)$ | $a \lor b = (c \land a) \lor (a \lor b)$ | $a \land b = (c \land a) \lor (b \land a)$ | $a \lor b = (c \land a) \lor (b \lor a)$ |
| $(a \land b) \land c = a \land (b \land c)$ | $(a \lor b) \lor c = a \lor (b \lor c)$ | $(a \land b) \land c = (a \land c) \land b$ | $a \land (b \land c) = (a \land c) \land b$ | $(a \land b) \land c = a \land (c \land b)$ | $a \land (b \land c) = a \land (c \land b)$ |
| $a \land (b \lor c) = a \land (c \lor b)$ | $(a \lor b) \lor c = a \lor (c \lor b)$ | $a \lor (b \lor c) = (a \lor c) \lor b$ | $a \lor (b \lor c) = a \lor (c \lor b)$ | $(a \lor b) \lor c = a \lor (c \lor b)$ | $a \lor (b \lor c) = a \lor (c \lor b)$ |
| $a \land b = (c \lor b) \land (a \land b)$ | $a \lor b = (c \lor b) \lor (a \lor b)$ | $(a \land b) \land c = (b \land a) \land c$ | $a \land (b \land c) = (b \land a) \land c$ | $(a \land b) \land c = (b \land a) \land c$ | $(a \land b) \land c = b \land (a \land c)$ |

⋮ 324 lines

| $a \land (b \land c) = (a \lor a) \land (c \land b)$ | $a \land (b \lor c) = (a \lor a) \land (c \lor b)$ | $a \lor (b \land c) = (a \land a) \lor (c \land b)$ | $(a \lor b) \lor c = (a \land a) \lor (c \lor b)$ |
|---|---|---|---|
| $a \lor (b \lor c) = (a \land a) \lor (c \lor b)$ | $a \lor (b \land c) = (a \lor a) \lor (c \land b)$ | $(a \lor b) \lor c = (a \lor a) \lor (c \lor b)$ | $a \lor (b \lor c) = (a \lor a) \lor (c \lor b)$ |
| $(a \land b) \land c = (a \land b) \land (a \land c)$ | $a \land (b \land c) = (a \land b) \land (a \land c)$ | $a \lor (b \land c) = (a \lor b) \land (a \lor c)$ | $a \land (b \lor c) = (a \land b) \lor (a \land c)$ |
| $(a \lor b) \lor c = (a \lor b) \lor (a \lor c)$ | $a \lor (b \lor c) = (a \lor b) \lor (a \lor c)$ | $(a \land b) \land c = (a \land b) \land (b \land c)$ | $a \land (b \land c) = (a \land b) \land (b \land c)$ |
| $(a \lor b) \lor c = (a \lor b) \lor (b \lor c)$ | $a \lor (b \lor c) = (a \lor b) \lor (b \lor c)$ | $(a \land b) \land c = (a \land b) \land (c \land a)$ | $a \land (b \land c) = (a \land b) \land (c \land a)$ |
| $a \lor (b \land c) = (a \lor b) \land (c \lor a)$ | $a \land (b \lor c) = (a \lor b) \lor (c \lor a)$ | $(a \lor b) \lor c = (a \lor b) \lor (c \lor a)$ | $a \lor (b \lor c) = (a \lor b) \lor (c \lor a)$ |

The reduction from all "structurally possible" theorems to just "ones we consider interesting" can be thought of as a form of coarse graining. And it could well be that this coarse graining would depend on all sorts of accidents of human mathematical history. But at least in the case of Boolean algebra there seems to be a surprisingly simple and "mechanical" procedure that can reproduce it.

Go through all theorems in order of increasing structural complexity, in each case seeing whether a given theorem can be proved from ones earlier in the list:

It turns out that the theorems identified by humans as "interesting" coincide almost exactly with "root theorems" that cannot be proved from earlier theorems in the list. Or, put another way, the "coarse graining" that human mathematicians do seems (at least in this case) to essentially consist of picking out only those theorems that represent "minimal statements" of new information—and eliding away those that involve "extra ornamentation".

But how are these "notable theorems" laid out in metamathematical space? Earlier we saw how the simplest of them can be reached after just a few steps in the entailment cone of a typical textbook axiom system for Boolean algebra. The full entailment cone rapidly gets unmanageably large but we can get a first approximation to it by generating individual proofs (using automated theorem proving) of our notable theorems, and then seeing how these "knit together" through shared intermediate lemmas in a token-event graph:

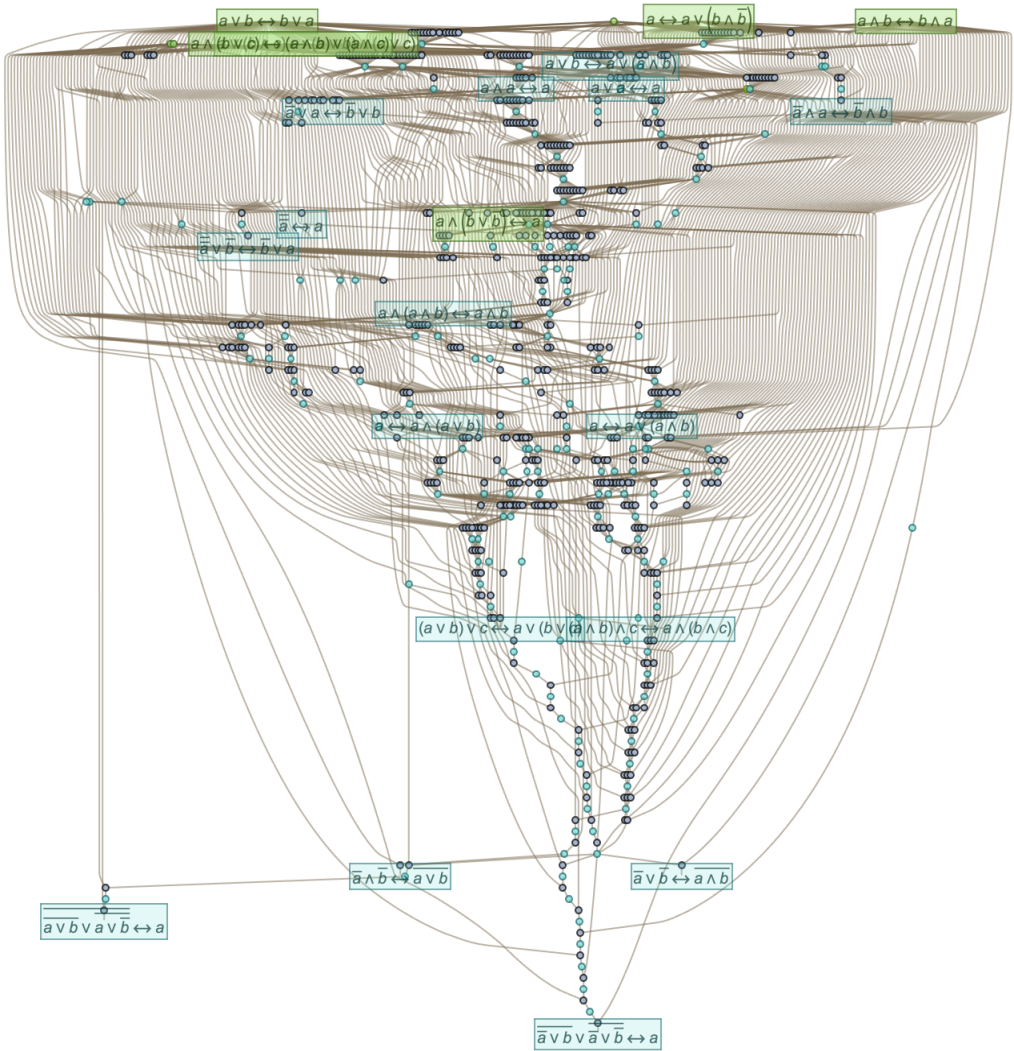

Looking at this picture we see at least a hint that clumps of notable theorems are spread out across the entailment cone, only modestly building on each other—and in effect "staking out separated territories" in the entailment cone. But of the 11 notable theorems shown here, 7 depend on all 6 axioms, while 4 depend only on various different sets of 3 axioms—suggesting at least a certain amount of fundamental interdependence or coherence.

From the token-event graph we can derive a branchial graph that represents a very rough approximation to how the theorems are "laid out in metamathematical space":

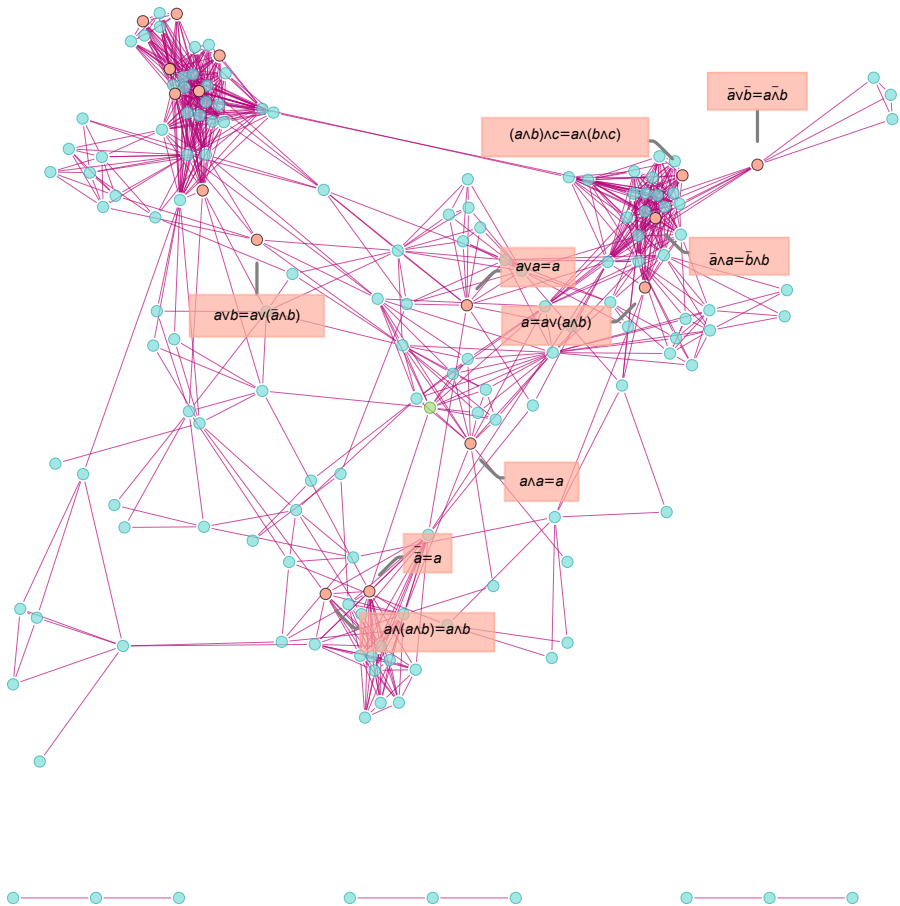

We can get a potentially slightly better approximation by including proofs not just of notable theorems, but of all theorems up to a certain structural complexity. The result shows separation of notable theorems both in the multiway graph

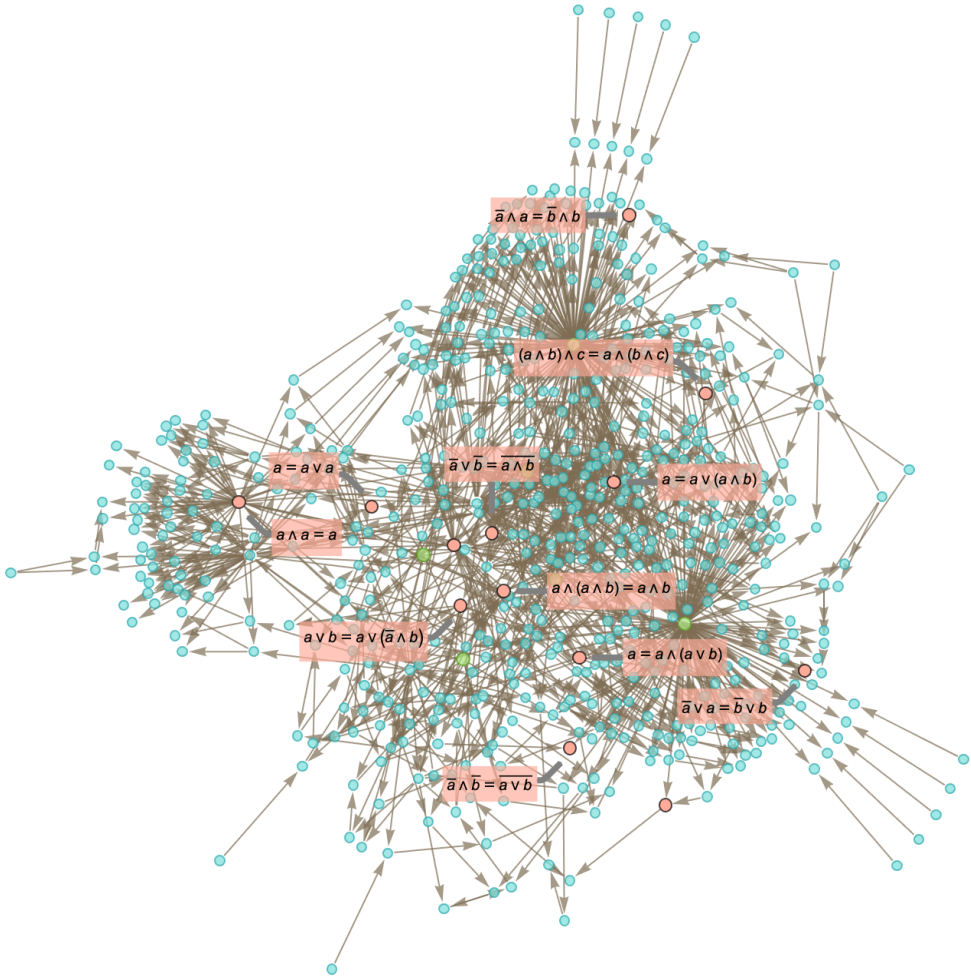

and in the branchial graph:

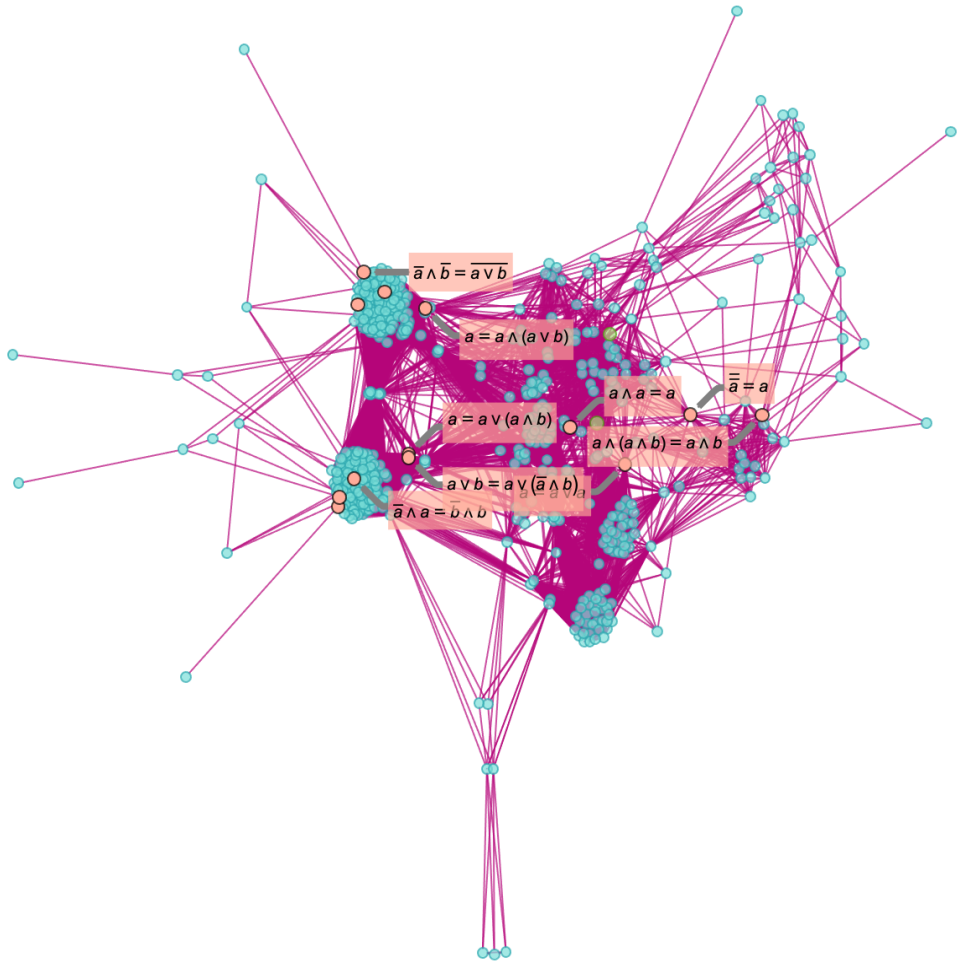

In doing this empirical metamathematics we're including only specific proofs rather than enumerating the whole entailment cone. We're also using only a specific axiom system. And even beyond this, we're using specific operators to write our statements in Boolean algebra.

In a sense each of these choices represents a particular "metamathematical coordinatization"—or particular reference frame or slice that we're sampling in the ruliad.

For example, in what we've done above we've built up statements from **And**, **Or** and **Not**. But we can just as well use any other functionally complete sets of operators, such as the following (here each shown representing a few specific Boolean expressions):

| | Or[p, q] | Nand[p, q] | Xor[p, q] | Implies[p, q] | Or[p, And[q, r]] |
|---|---|---|---|---|---|
| And Or Not<br>∧ ∨ − | p ∨ q | $\overline{p \wedge q}$ | $\overline{p \wedge q} \wedge (p \vee q)$ | $q \vee \overline{p}$ | $p \vee (q \wedge r)$ |
| And Not<br>∧ − | $\overline{\overline{p} \wedge \overline{q} \wedge q}$ | $\overline{p \wedge q}$ | $\overline{\overline{p \wedge q} \wedge \overline{\overline{p} \wedge \overline{q}}}$ | $\overline{p \wedge \overline{q}}$ | $\overline{\overline{p} \wedge \overline{q \wedge r}}$ |
| Or Not<br>∨ − | p ∨ q | $\overline{p} \vee \overline{q}$ | $\overline{\overline{p \vee \overline{q}} \vee \overline{q \vee \overline{p}}}$ | $q \vee \overline{p}$ | $p \vee \overline{\overline{q} \vee \overline{r}}$ |
| Implies Not<br>⇒ − | $\overline{p} \Rightarrow q$ | $p \Rightarrow \overline{q}$ | $(p \Rightarrow q) \Rightarrow \overline{q \Rightarrow p}$ | $p \Rightarrow q$ | $(q \Rightarrow \overline{r}) \Rightarrow p$ |
| Xor Implies<br>⊻ ⇒ | $(p \Rightarrow q) \Rightarrow q$ | $p \Rightarrow p \veebar q$ | $p \veebar q$ | $p \Rightarrow q$ | $(q \Rightarrow (r \Rightarrow p)) \Rightarrow p$ |
| Nand Not<br>∘ − | $\overline{p} \circ \overline{q}$ | $p \circ q$ | $(p \circ \overline{q}) \circ (q \circ \overline{p})$ | $p \circ \overline{q}$ | $(q \circ r) \circ \overline{p}$ |
| Nand<br>∘ | $(p \circ p) \circ (q \circ q)$ | $p \circ q$ | $(p \circ (p \circ q)) \circ (q \circ (p \circ p))$ | $p \circ (p \circ q)$ | $(p \circ p) \circ (q \circ r)$ |
| Nor<br>∘ | $(p \circ q) \circ (p \circ q)$ | $(p \circ (p \circ p)) \circ ((p \circ p) \circ (q \circ q))$ | $(p \circ q) \circ ((p \circ p) \circ (q \circ q))$ | $(p \circ (p \circ p)) \circ (q \circ (p \circ p))$ | $(p \circ q) \circ (p \circ r)$ |

For each set of operators, there are different axiom systems that can be used. And for each axiom system there will be different proofs. Here are a few examples of axiom systems with a few different sets of operators—in each case giving a proof of the law of double negation (which has to be stated differently for different operators):

| operators | axioms | statement | proof | steps |
|---|---|---|---|---|
| And Or Not<br>∧ ∨ ¯ | $a \wedge b = b \wedge a$<br>$a \wedge (b \vee c) = (a \wedge b) \vee (a \wedge c)$<br>$a \vee (b \wedge \overline{b}) = a$<br>$a \wedge (b \vee \overline{b}) = a$<br>$a \vee b = b \vee a$<br>$a \vee (b \wedge c) = (a \vee b) \wedge (a \vee c)$ | $\overline{\overline{p}} = p$ | 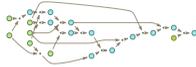 | 19 |
| Or Not<br>∨ ¯ | $a \vee (b \vee c) = (a \vee b) \vee c$<br>$a \vee b = b \vee a$<br>$\overline{\overline{a} \vee b} \vee \overline{\overline{a} \vee \overline{b}} = a$ | $\overline{\overline{p}} = p$ | 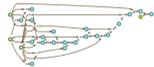 | 24 |
| Or Not<br>∨ ¯ | $\overline{\overline{\overline{c \vee b} \vee a} \vee \overline{\overline{d} \vee d}} \vee (\overline{a} \vee c) = a$ | $\overline{\overline{p}} = p$ | 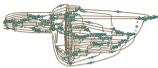 | 116 |
| Nand<br>∘ | $(a \circ a) \circ (a \circ a) = a$<br>$a \circ (b \circ (b \circ b)) = a \circ a$<br>$(a \circ (b \circ c)) \circ (a \circ (b \circ c)) = ((b \circ b) \circ a) \circ ((c \circ c) \circ a)$ | $(p \circ p) \circ (p \circ p) = p$ | 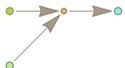 | 3 |
| Nand<br>∘ | $(a \circ a) \circ (a \circ b) = a$<br>$a \circ (a \circ b) = a \circ (b \circ b)$<br>$a \circ (a \circ (b \circ c)) = b \circ (b \circ (a \circ c))$ | $(p \circ p) \circ (p \circ p) = p$ | 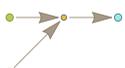 | 3 |
| Nand<br>∘ | $a \circ (b \circ (a \circ c)) = ((c \circ b) \circ b) \circ a$<br>$(a \circ a) \circ (b \circ a) = a$ | $(p \circ p) \circ (p \circ p) = p$ | 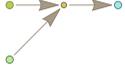 | 3 |
| Nand<br>∘ | $a \circ b = b \circ a$<br>$(a \circ b) \circ (a \circ (b \circ c)) = a$ | $(p \circ p) \circ (p \circ p) = p$ | 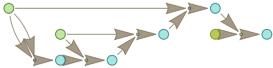 | 8 |
| Nand<br>∘ | $((b \circ c) \circ a) \circ (b \circ ((b \circ a) \circ b)) = a$ | $(p \circ p) \circ (p \circ p) = p$ | 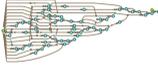 | 54 |
| Nand<br>∘ | $(a \circ ((c \circ a) \circ a)) \circ (c \circ (b \circ a)) = c$ | $(p \circ p) \circ (p \circ p) = p$ | 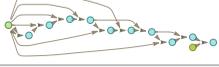 | 12 |

Boolean algebra (or, equivalently, propositional logic) is a somewhat desiccated and thin example of mathematics. So what do we find if we do empirical metamathematics on other areas?

Let's talk first about geometry—for which Euclid's *Elements* provided the very first large-scale historical example of an axiomatic mathematical system. The *Elements* started from 10 axioms (5 "postulates" and 5 "common notions"), then gave 465 theorems.

Each theorem was proved from previous ones, and ultimately from the axioms. Thus, for example, the "proof graph" (or "theorem dependency graph") for Book 1, Proposition 5 (which says that angles at the base of an isosceles triangle are equal) is:

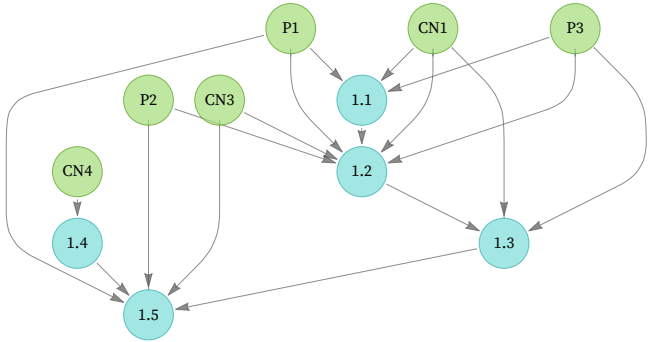

One can think of this as a coarse-grained version of the proof graphs we've used before (which are themselves in turn "slices" of the entailment graph)—in which each node shows how a collection of "input" theorems (or axioms) entails a new theorem.

Here's a slightly more complicated example (Book 1, Proposition 48) that ultimately depends on all 10 of the original axioms:

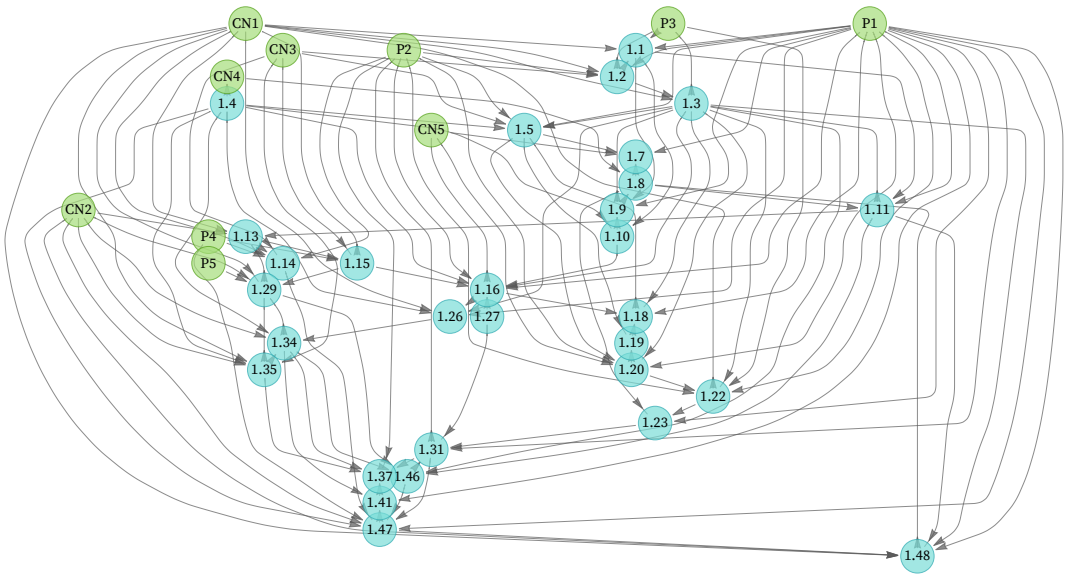

And here's the full graph for all the theorems in Euclid's *Elements*:

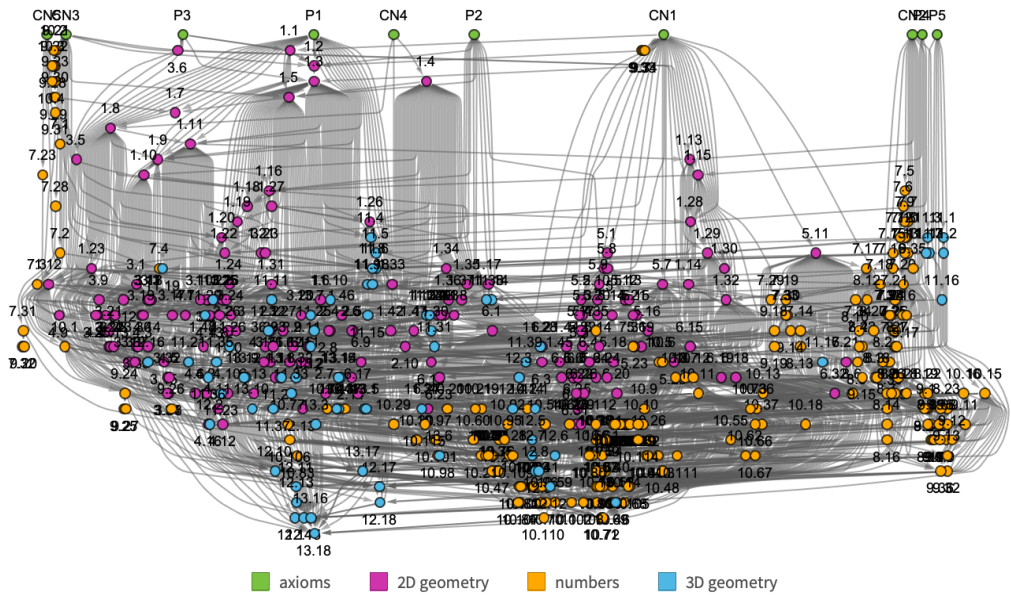

Of the 465 theorems here, 255 (i.e. 55%) depend on all 10 axioms. (For the much smaller number of notable theorems of Boolean algebra above we found that 64% depended on all 6 of our stated axioms.) And the general connectedness of this graph in effect reflects the idea that Euclid's theorems represent a coherent body of connected mathematical knowledge.

The branchial graph gives us an idea of how the theorems are "laid out in metamathematical space":

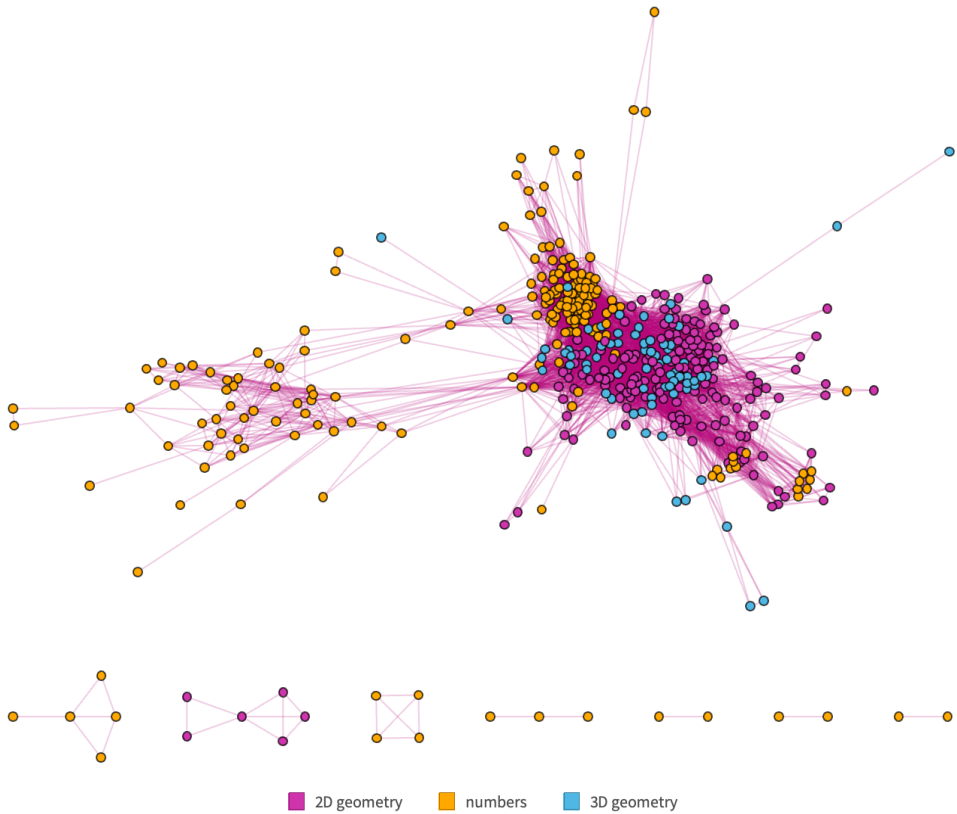

One thing we notice is that theorems about different areas—shown here in different colors—tend to be separated in metamathematical space. And in a sense the seeds of this separation are already evident if we look "textually" at how theorems in different books of Euclid's *Elements* refer to each other:

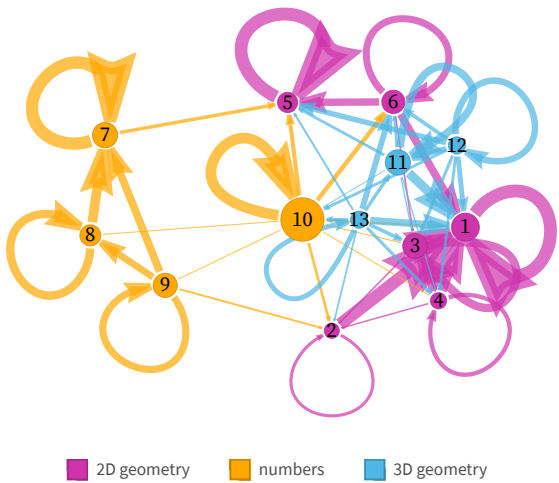

Looking at the overall dependence of one theorem on others in effect shows us a very coarse form of entailment. But can we go to a finer level—as we did above for Boolean algebra? As a first step, we have to have an explicit symbolic representation for our theorems. And beyond that, we have to have a formal axiom system that describes possible transformations between these.

At the level of "whole theorem dependency" we can represent the entailment of Euclid's Book 1, Proposition 1 from axioms as:

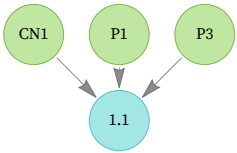

But if we now use the full, formal axiom system for geometry that we discussed in a previous section we can use automated theorem proving to get a full proof of Book 1, Proposition 1:

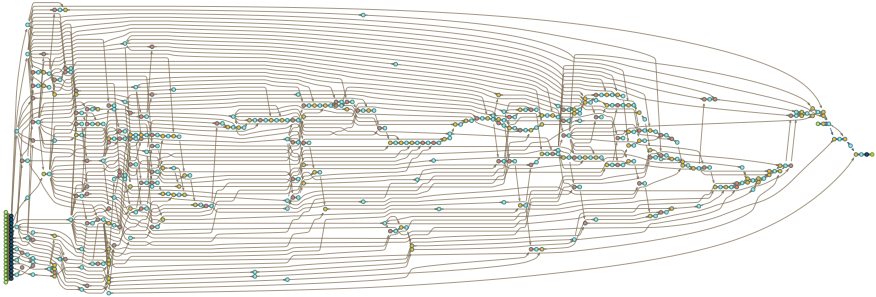

In a sense this is "going inside" the theorem dependency graph to look explicitly at how the dependencies in it work. And in doing this we see that what Euclid might have stated in words in a sentence or two is represented formally in terms of hundreds of detailed intermediate lemmas. (It's also notable that whereas in Euclid's version, the theorem depends only on 3 out of 10 axioms, in the formal version the theorem depends on 18 out of 20 axioms.)

How about for other theorems? Here is the theorem dependency graph from Euclid's *Elements* for the Pythagorean theorem (which Euclid gives as Book 1, Proposition 47):

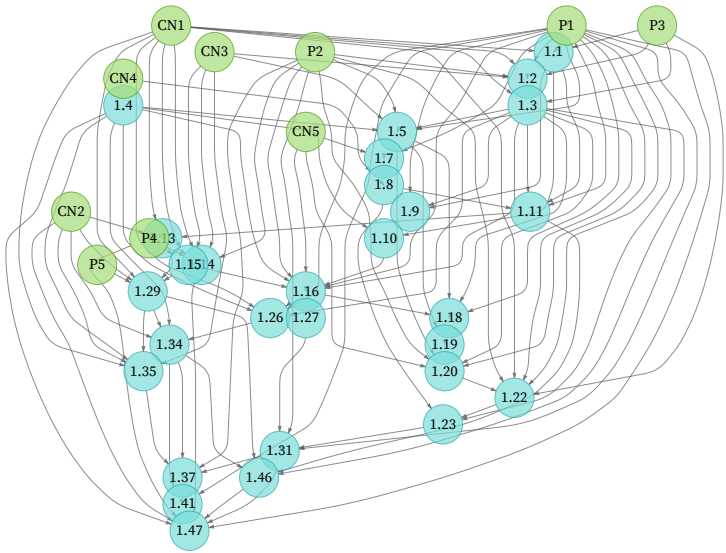

The theorem depends on all 10 axioms, and its stated proof goes through 28 intermediate theorems (i.e. about 6% of all theorems in the *Elements*). In principle we can "unroll" the proof dependency graph to see directly how the theorem can be "built up" just from copies of the original axioms. Doing a first step of unrolling we get:

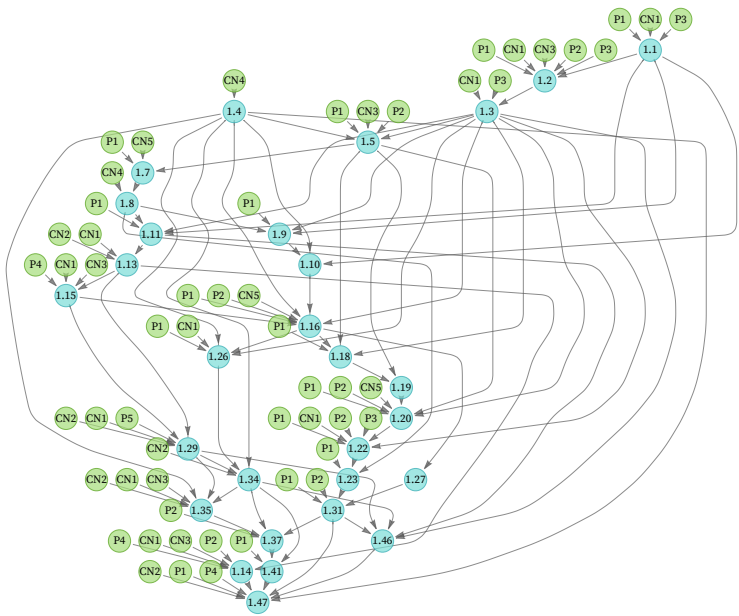

And "flattening everything out" so that we don't use any intermediate lemmas but just go back to the axioms to "re-prove" everything we can derive the theorem from a "proof tree" with the following number of copies of each axiom (and a certain "depth" to reach that axiom):

| axiom  | CN1  | CN2 | CN3 | CN4 | CN5 | P1   | P2  | P3   | P4 | P5 |
|--------|------|-----|-----|-----|-----|------|-----|------|----|----|
| depth  | 2    | 1   | 2   | 2   | 4   | 1    | 2   | 3    | 1  | 3  |
| copies | 1235 | 81  | 310 | 217 | 82  | 1094 | 321 | 1096 | 27 | 10 |

So how about a more detailed and formal proof? We could certainly in principle construct this using the axiom system we discussed above.

But an important general point is that the thing we in practice call "the Pythagorean theorem" can actually be set up in all sorts of different axiom systems. And as an example let's consider setting it up in the main actual axiom system that working mathematicians typically imagine they're (usually implicitly) using, namely ZFC set theory.

Conveniently, the Metamath formalized math system has accumulated about 40,000 theorems across mathematics, all with hand-constructed proofs based ultimately on ZFC set theory. And within this system we can find the theorem dependency graph for the Pythagorean theorem:

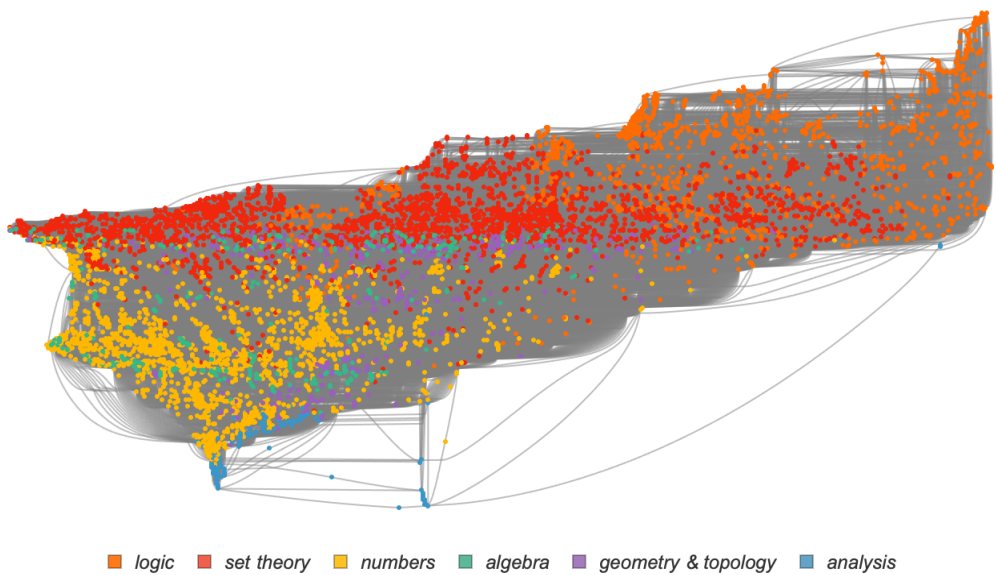

■ *logic*   ■ *set theory*   ■ *numbers*   ■ *algebra*   ■ *geometry & topology*   ■ *analysis*

Altogether it involves 6970 intermediate theorems, or about 18% of all theorems in Metamath—including ones from many different areas of mathematics. But how does it ultimately depend on the axioms? First, we need to talk about what the axioms actually are. In addition to "pure ZFC set theory", we need axioms for (predicate) logic, as well as ones that define real and complex numbers. And the way things are set up in Metamath's "set.mm" there are (essentially) 49 basic axioms (9 for pure set theory, 15 for logic and 25 related to numbers).

And much as in Euclid's *Elements* we found that the Pythagorean theorem depended on all the axioms, so now here we find that the Pythagorean theorem depends on 48 of the 49 axioms—with the only missing axiom being the Axiom of Choice.

Just like in the Euclid's *Elements* case, we can imagine "unrolling" things to see how many copies of each axiom are used. Here are the results—together with the "depth" to reach each axiom:

| Axiom | Copies | Depth |
|---|---|---|
| Modus ponens | $6.0 \times 10^{36}$ | 2 |
| Axiom of Simplification | $1.8 \times 10^{36}$ | 3 |
| Frege's Axiom | $1.9 \times 10^{36}$ | 4 |
| Principle of Transposition | $2.7 \times 10^{35}$ | 7 |
| Rule of Generalization | $1.3 \times 10^{32}$ | 6 |
| Axiom of Quantified Implication | $7.1 \times 10^{31}$ | 6 |
| Axiom of Distinctness | $4.2 \times 10^{31}$ | 4 |
| Axiom of Existence | $3.8 \times 10^{31}$ | 8 |
| Axiom of Equality | $5.2 \times 10^{31}$ | 8 |
| Binary predicate left equality | $4.4 \times 10^{24}$ | 12 |
| Binary predicate right equality | $1.1 \times 10^{30}$ | 10 |
| Axiom of Quantified Negation | $1.9 \times 10^{30}$ | 10 |
| Axiom of Quantifier Commutation | $1.8 \times 10^{29}$ | 9 |
| Axiom of Substitution | $2.2 \times 10^{30}$ | 9 |
| Axiom of Quantified Equality | $5.2 \times 10^{29}$ | 9 |
| Axiom of Extensionality | $5.5 \times 10^{29}$ | 6 |
| Axiom of Replacement | $2.8 \times 10^{10}$ | 10 |
| Axiom of Separation | $4.2 \times 10^{26}$ | 7 |
| Null Set Axiom | $2.7 \times 10^{26}$ | 6 |
| Axiom of Power Sets | $5.8 \times 10^{23}$ | 7 |
| Axiom of Pairing | $4.1 \times 10^{26}$ | 10 |
| Axiom of Union | $2.1 \times 10^{24}$ | 8 |
| Axiom of Infinity | $2.4 \times 10^{9}$ | 7 |
| Complex numbers form set | $1.8 \times 10^{16}$ | 6 |
| Real numbers subset of complex | $3.4 \times 10^{22}$ | 3 |
| 1 is a complex number | $2.1 \times 10^{22}$ | 3 |
| i is a complex number | $1.0 \times 10^{22}$ | 3 |
| Closure of complex addition | $1.0 \times 10^{22}$ | 3 |
| Closure of real addition | $2.1 \times 10^{22}$ | 5 |
| Closure of complex multiplication | $1.1 \times 10^{22}$ | 3 |
| Closure of real multiplication | $5.7 \times 10^{21}$ | 5 |
| Commutativity of complex multiplication | $6.3 \times 10^{21}$ | 4 |
| Associativity of complex addition | $2.1 \times 10^{21}$ | 5 |
| Associativity of complex multiplication | $3.7 \times 10^{21}$ | 6 |
| Distributivity for complex numbers | $3.7 \times 10^{21}$ | 5 |
| $i^2 = -1$ | $7.8 \times 10^{21}$ | 4 |
| 1 and 0 are distinct | $5.0 \times 10^{21}$ | 4 |
| 1 is multiplicative identity | $1.6 \times 10^{21}$ | 6 |
| Negative reals exist | $5.9 \times 10^{21}$ | 4 |
| Reciprocals of nonzero reals exist | $6.2 \times 10^{21}$ | 4 |
| Complex numbers expressible by 2 reals | $6.1 \times 10^{21}$ | 5 |
| Ordering on reals has strict trichotomy | $7.6 \times 10^{21}$ | 6 |
| Ordering on reals is transitive | $1.9 \times 10^{21}$ | 7 |
| Ordering property of addition on reals | $5.7 \times 10^{21}$ | 7 |
| Product of 2 positive reals is positive | $2.8 \times 10^{17}$ | 5 |
| Supremum property for set of reals | $3.0 \times 10^{13}$ | 10 |
| Addition applies to complex numbers | $1.4 \times 10^{7}$ | 10 |
| Multiplication applies to complex numbers | $1.4 \times 10^{7}$ | 9 |

And, yes, the numbers of copies of most of the axioms required to establish the Pythagorean theorem are extremely large.

There are several additional wrinkles that we should discuss. First, we've so far only considered overall theorem dependency—or in effect "coarse-grained entailment". But the Metamath system ultimately gives complete proofs in terms of explicit substitutions (or, effectively, bisubstitutions) on symbolic expressions. So, for example, while the first-level "whole-theorem-dependency" graph for the Pythagorean theorem is

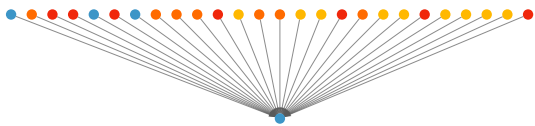

the full first-level entailment structure based on the detailed proof is (where the black vertices indicate "internal structural elements" in the proof—such as variables, class specifications and "inputs"):

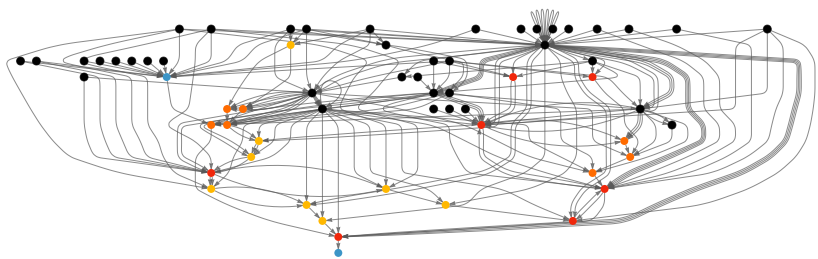

Another important wrinkle has to do with the concept of definitions. The Pythagorean theorem, for example, refers to squaring numbers. But what is squaring? What are numbers? Ultimately all these things have to be defined in terms of the "raw data structures" we're using.

In the case of Boolean algebra, for example, we could set things up just using **Nand** (say denoted ∘), but then we could define **And** and **Or** in terms of **Nand** (say as $(p \circ q) \circ (p \circ q)$ and $(p \circ p) \circ (q \circ q)$ respectively). We could still write expressions using **And** and **Or**—but with our definitions we'd immediately be able to convert these to pure **Nand**s. Axioms—say about **Nand**—give us transformations we can use repeatedly to make derivations. But definitions are transformations we use "just once" (like macro expansion in programming) to reduce things to the point where they involve only constructs that appear in the axioms.

In Metamath's "set.mm" there are about 1700 definitions that effectively build up from "pure set theory" (as well as logic, structural elements and various axioms about numbers) to give the mathematical constructs one needs. So, for example, here is the definition dependency graph for addition ("+" or **Plus**):

At the bottom are the basic constructs of logic and set theory—in terms of which things like order relations, complex numbers and finally addition are defined. The definition dependency graph for GCD, for example, is somewhat larger, though has considerable overlap at lower levels:

Different constructs have definition dependency graphs of different sizes—in effect reflecting their "definitional distance" from set theory and the underlying axioms being used:

| construct | distance |
|---:|---:|
| First | 12 |
| ∞ | 58 |
| Plus | 60 |
| 1 | 61 |
| 2 | 67 |
| EvenQ | 83 |
| GCD | 91 |
| Power | 92 |
| PrimePi | 102 |

| construct | distance |
|---:|---:|
| Sinh | 102 |
| SymmetricGroup | 104 |
| π | 105 |
| Log | 112 |
| Abs | 114 |
| Eigenvalues | 119 |
| Zeta | 123 |
| Area | 133 |
| Det | 191 |

In our physicalized approach to metamathematics, though, something like set theory is not our ultimate foundation. Instead, we imagine that everything is eventually built up from the raw ruliad, and that all the constructs we're considering are formed from what amount to configurations of emes in the ruliad. We discussed above how constructs like numbers and logic can be obtained from a combinator representation of the ruliad.

We can view the definition dependency graph above as being an empirical example of how somewhat higher-level definitions can be built up. From a computer science perspective, we can think of it as being like a type hierarchy. From a physics perspective, it's as if we're starting from atoms, then building up to molecules and beyond.

It's worth pointing out, however, that even the top of the definition hierarchy in something like Metamath is still operating very much at an axiomatic kind of level. In the analogy we've been using, it's still for the most part "formulating math at the molecular dynamics level" not at the more human "fluid dynamics" level.

We've been talking about "the Pythagorean theorem". But even on the basis of set theory there are many different possible formulations one can give. In Metamath, for example, there is the pythag version (which is what we've been using), and there is also a (somewhat more general) pythi version. So how are these related? Here's their combined theorem dependency graph (or at least the first two levels in it)—with red indicating theorems used only in deriving pythag, blue indicating ones used only in deriving pythi, and purple indicating ones used in both:

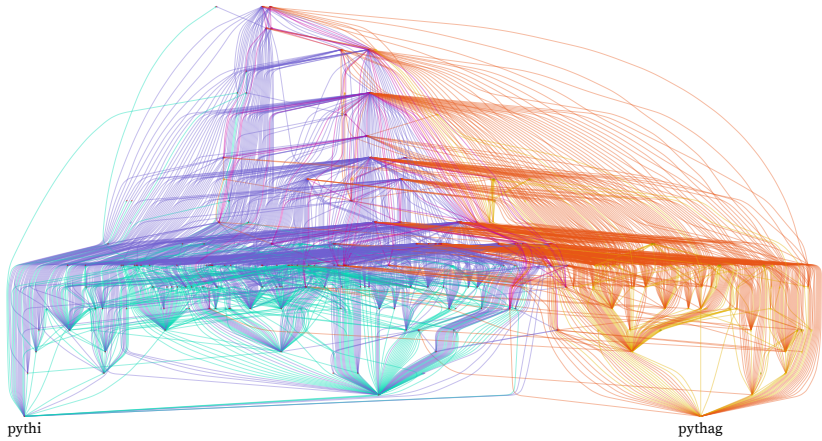

And what we see is there's a certain amount of "lower-level overlap" between the derivations of these variants of the Pythagorean theorem, but also some discrepancy—indicating a certain separation between these variants in metamathematical space.

So what about other theorems? Here's a table of some famous theorems from all over mathematics, sorted by the total number of theorems on which proofs of them formulated in Metamath depend—giving also the number of axioms and definitions used in each case:

| | total | axioms | definitions |
|---:|---:|---:|---:|
| Schröder–Bernstein Theorem | 1474 | 21 | 48 |
| Principle of Mathematical Induction | 1532 | 20 | 47 |
| Cantor's Theorem | 1719 | 21 | 55 |
| Number of subsets of a set | 1737 | 21 | 59 |
| Infinitude of primes | 3394 | 45 | 98 |
| Non–denumerability of continuum | 3508 | 45 | 101 |
| Triangle Inequality | 3626 | 44 | 106 |
| Intermediate Value Theorem | 3782 | 44 | 111 |
| Irrationality of $\sqrt{2}$ | 3816 | 44 | 108 |
| Bezout's Theorem | 3856 | 44 | 111 |
| Denumerability of rationals | 3948 | 45 | 108 |
| Binomial as number of subsets | 3961 | 44 | 110 |
| Euclid's GCD algorithm | 3966 | 44 | 113 |
| Formula for Pythagorean triples | 4281 | 44 | 118 |
| Sum of a geometric series | 4459 | 46 | 121 |
| Fermat's Little Theorem | 4477 | 45 | 123 |
| Fundamental Theorem of Arithmetic | 4497 | 45 | 121 |
| Binomial Theorem | 4522 | 46 | 123 |
| Lebesgue Integration Theorem | 4543 | 46 | 127 |
| Sum of an arithmetic series | 4545 | 46 | 123 |
| Divisibility by 3 Rule | 4561 | 46 | 129 |
| Principle of Inclusion/Exclusion | 4651 | 46 | 124 |
| Sum of triangular–number reciprocals | 4662 | 46 | 126 |
| Sum of kth powers | 4698 | 46 | 124 |
| Four–Squares Theorem | 4710 | 45 | 125 |
| Lagrange's Theorem | 4738 | 46 | 141 |
| Divergence of harmonic series | 4868 | 46 | 127 |
| Ramsey's Theorem | 4916 | 46 | 130 |
| Königsberg bridge theorem | 4996 | 45 | 145 |
| Polynomial Factor Theorem | 5015 | 47 | 133 |
| De Moivre's Theorem | 5217 | 48 | 133 |
| Divergence of inverse prime series | 5368 | 46 | 135 |
| Sylow's Theorem | 5445 | 46 | 153 |
| Wilson's Theorem | 5544 | 48 | 175 |
| Cramer's Rule | 6192 | 47 | 218 |
| Mean Value Theorem | 6438 | 48 | 216 |
| Law of Quadratic Reciprocity | 6448 | 48 | 214 |
| L'Hôpital's Rule | 6498 | 48 | 217 |
| Cayley–Hamilton Theorem | 6551 | 47 | 230 |
| Friendship Graph Theorem | 6853 | 47 | 181 |
| Ptolemy's Theorem | 7103 | 48 | 229 |
| Liouville's Theorem | 7179 | 48 | 238 |
| Primes 4k+1 are sums of 2 squares | 7246 | 48 | 243 |
| Taylor's Theorem | 7253 | 48 | 249 |
| Law of Cosines | 7269 | 48 | 231 |
| Pythagorean Theorem | 7271 | 48 | 231 |
| Heron's Formula | 7285 | 48 | 231 |
| Sum of angles of a triangle | 7295 | 48 | 231 |
| Isosceles Triangle Theorem | 7309 | 48 | 231 |
| Fundamental Theorem of Calculus | 7338 | 49 | 223 |
| Product of segments of chords | 7361 | 48 | 231 |
| Solution of quartic equations | 7461 | 48 | 233 |
| Solution of cubic equations | 7490 | 48 | 233 |
| Value of $\zeta(2)$ | 7522 | 48 | 238 |
| Fundamental Theorem of Algebra | 7542 | 48 | 238 |
| Perfect Number Theorem | 7579 | 48 | 237 |
| Arithmetic–geometric mean inequality | 7750 | 48 | 253 |
| Leibniz' series for $\pi$ | 7870 | 48 | 238 |
| Bertrand's Postulate | 7945 | 48 | 240 |
| Birthday Problem probability | 8001 | 48 | 237 |
| Prime Number Theorem | 8477 | 48 | 248 |
| Dirichlet's Theorem | 10 006 | 48 | 312 |

The Pythagorean theorem (here the pythi formulation) occurs solidly in the second half. Some of the theorems with the fewest dependencies are in a sense very structural theorems. But it's interesting to see that theorems from all sorts of different areas soon start appearing, and then are very much mixed together in the remainder of the list. One might have thought that theorems involving "more sophisticated concepts" (like Ramsey's theorem) would appear later than "more elementary" ones (like the sum of angles of a triangle). But this doesn't seem to be true.

There's a distribution of what amount to "proof sizes" (or, more strictly, theorem dependency sizes)—from the Schröder–Bernstein theorem which relies on less than 4% of all theorems, to Dirichlet's theorem that relies on 25%:

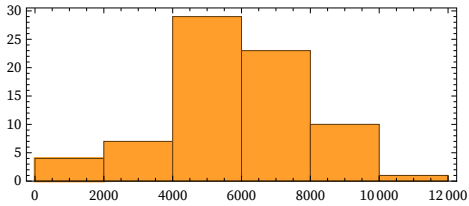

If we look not at "famous" theorems, but at all theorems covered by Metamath, the distribution becomes broader, with many short-to-prove "glue" or essentially "definitional" lemmas appearing:

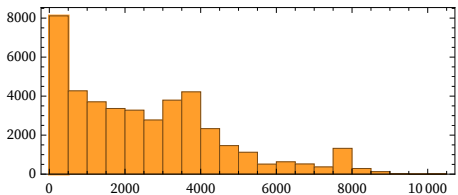

But using the list of famous theorems as an indication of the "math that mathematicians care about" we can conclude that there is a kind of "metamathematical floor" of results that one needs to reach before "things that we care about" start appearing. It's a bit like the situation in our Physics Project—where the vast majority of microscopic events that happen in the universe seem to be devoted merely to knitting together the structure of space, and only "on top of that" can events which can be identified with things like particles and motion appear.

And indeed if we look at the "prerequisites" for different famous theorems, we indeed find that there is a large overlap (indicated by lighter colors)—supporting the impression that in a sense one first has "knit together metamathematical space" and only then can one start generating "interesting theorems":

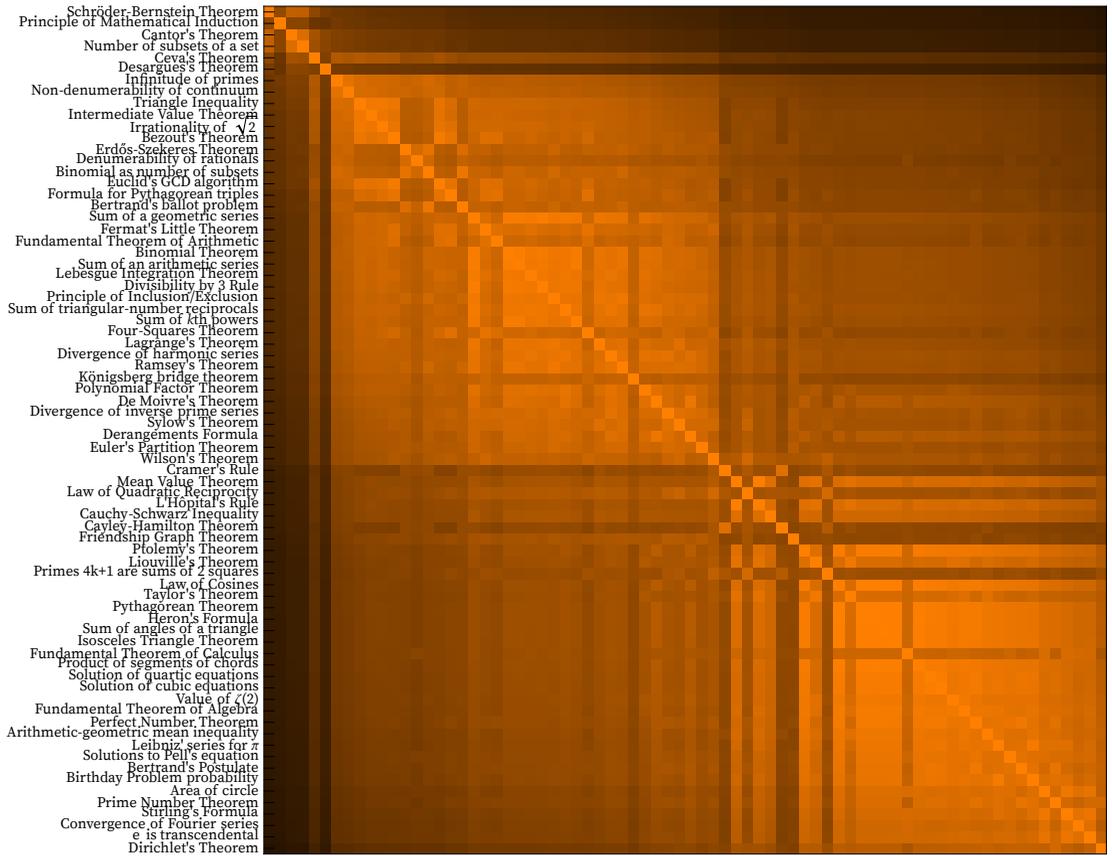

Another way to see "underlying overlap" is to look at what axioms different theorems ultimately depend on (the colors indicate the "depth" at which the axioms are reached):

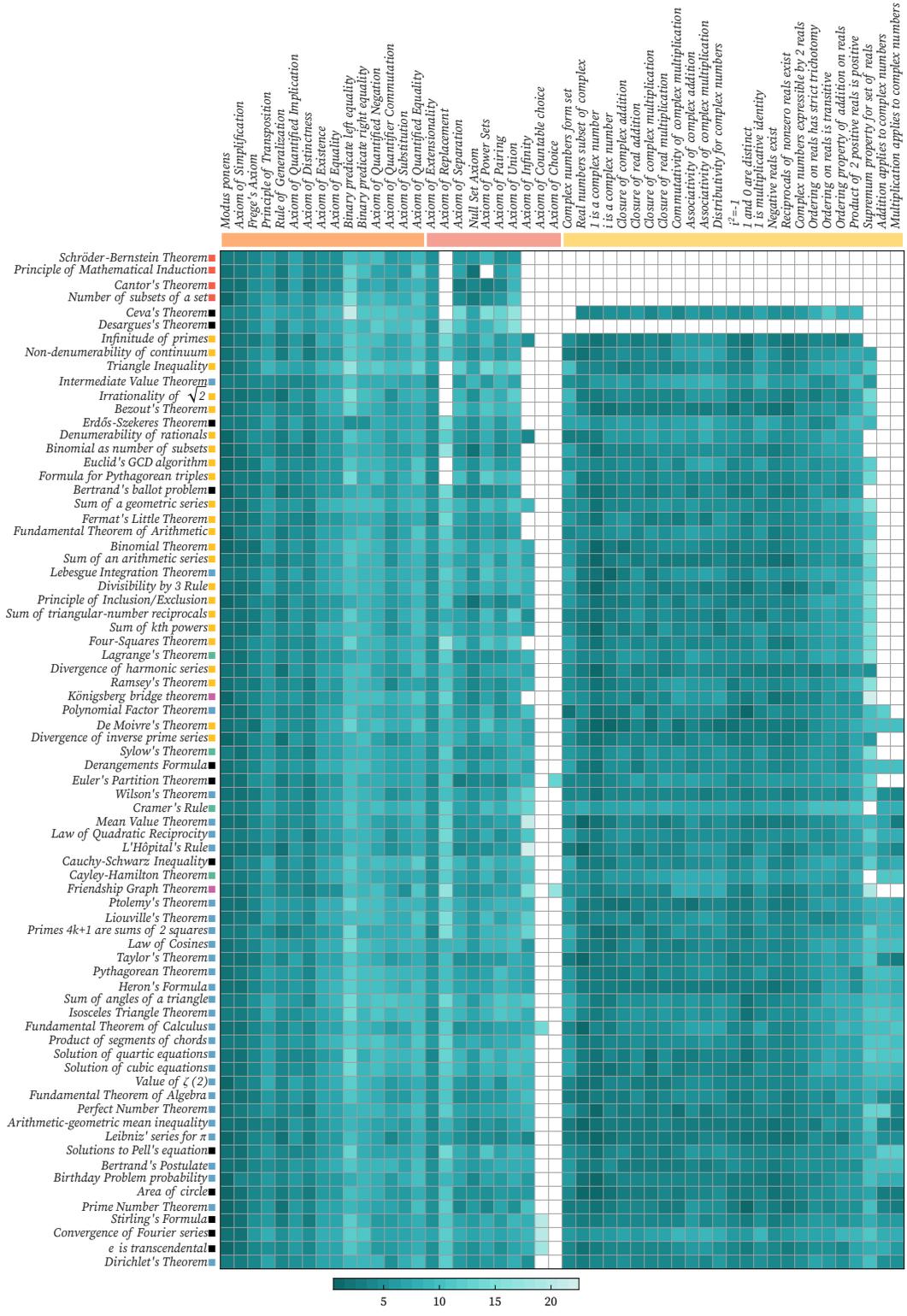

The theorems here are again sorted in order of "dependency size". The "very-set-theoretic" ones at the top don't depend on any of the various number-related axioms. And quite a few "integer-related theorems" don't depend on complex number axioms. But otherwise, we see that (at least according to the proofs in set.mm) most of the "famous theorems" depend on almost all the axioms. The only axiom that's rarely used is the Axiom of Choice—on which only things like "analysis-related theorems" such as the Fundamental Theorem of Calculus depend.

If we look at the "depth of proof" at which axioms are reached, there's a definite distribution:

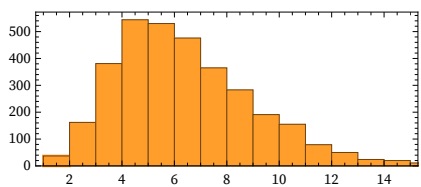

And this may be about as robust as any a "statistical characteristic" of the sampling of metamathematical space corresponding to mathematics that is "important to humans". If we were, for example, to consider all possible theorems in the entailment cone we'd get a very different picture. But potentially what we see here may be a characteristic signature of what's important to a "mathematical observer like us".

Going beyond "famous theorems" we can ask, for example, about all the 42,000 or so identified theorems in the Metamath set.mm collection. Here's a rough rendering of their theorem dependency graph, with different colors indicating theorems in different fields of math (and with explicit edges removed):

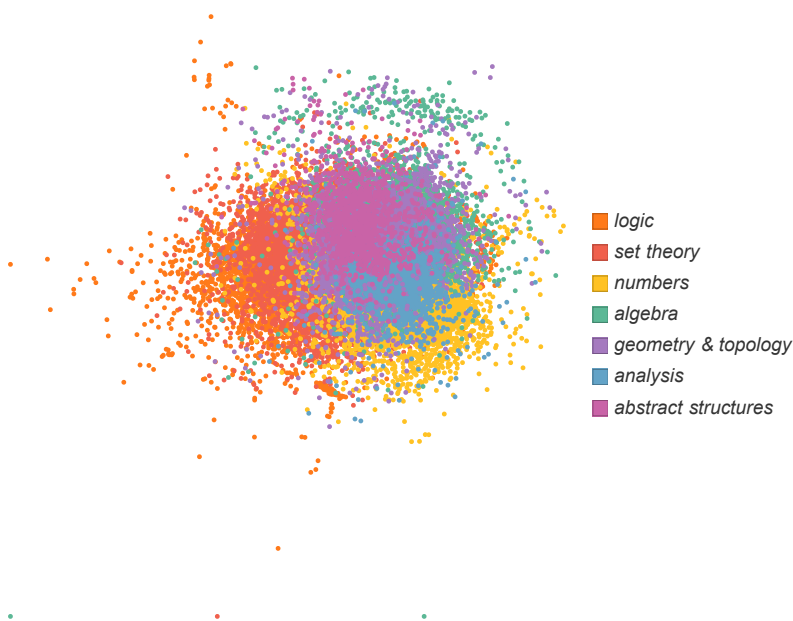

There's some evidence of a certain overall uniformity, but we can see definite "patches of metamathematical space" dominated by different areas of mathematics. And here's what happens if we zoom in on the central region, and show where famous theorems lie:

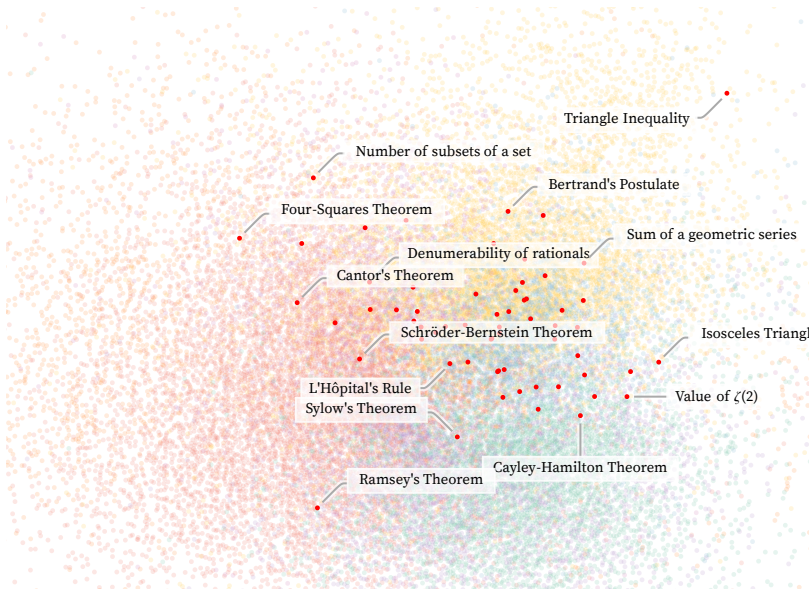

A bit like we saw for the named theorems of Boolean algebra clumps of famous theorems appear to somehow "stake out their own separate metamathematical territory". But notably the famous theorems seem to show some tendency to congregate near "borders" between different areas of mathematics.

To get more of a sense of the relation between these different areas, we can make what amounts to a highly coarsened branchial graph, effectively laying out whole areas of mathematics in metamathematical space, and indicating their cross-connections:

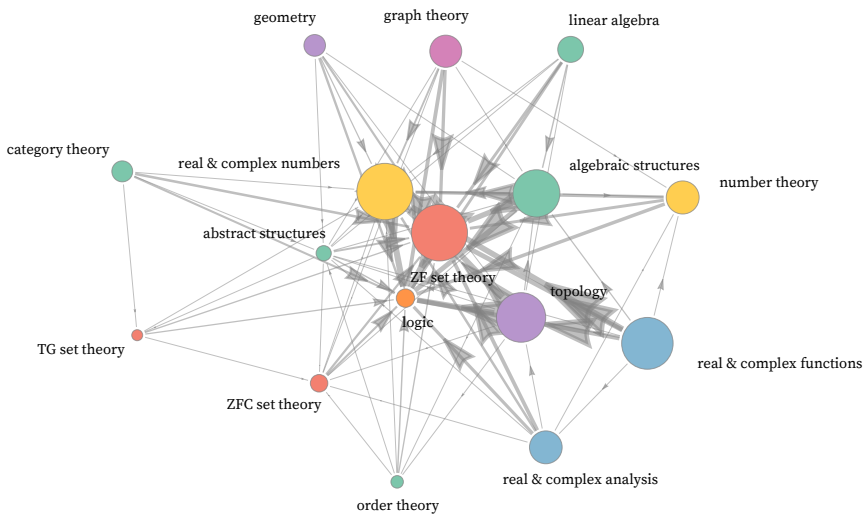

We can see "highways" between certain areas. But there's also a definite "background entanglement" between areas, reflecting at least a certain background uniformity in meta-mathematical space, as sampled with the theorems identified in Metamath.

It's not the case that all these areas of math "look the same"—and for example there are differences in their distributions of theorem dependency sizes:

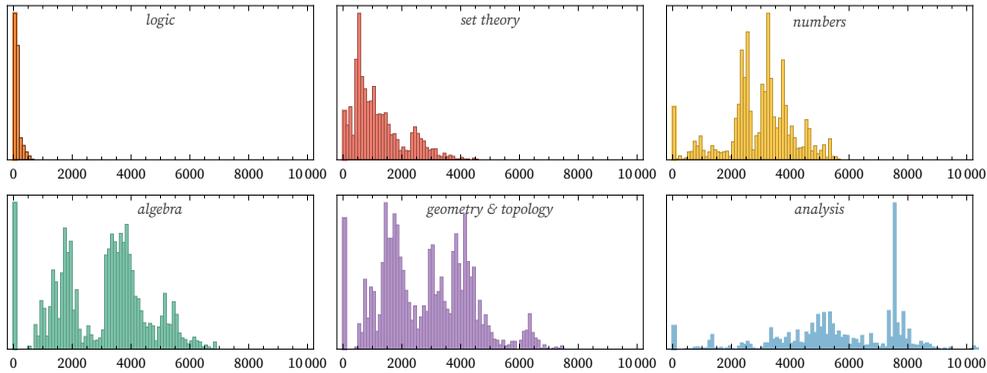

In areas like algebra and number theory, most proofs are fairly long, as revealed by the fact that they have many dependencies. But in set theory there are plenty of short proofs, and in logic all the proofs of theorems that have been included in Metamath are short.

What if we look at the overall dependency graph for all theorems in Metamath? Here's the adjacency matrix we get:

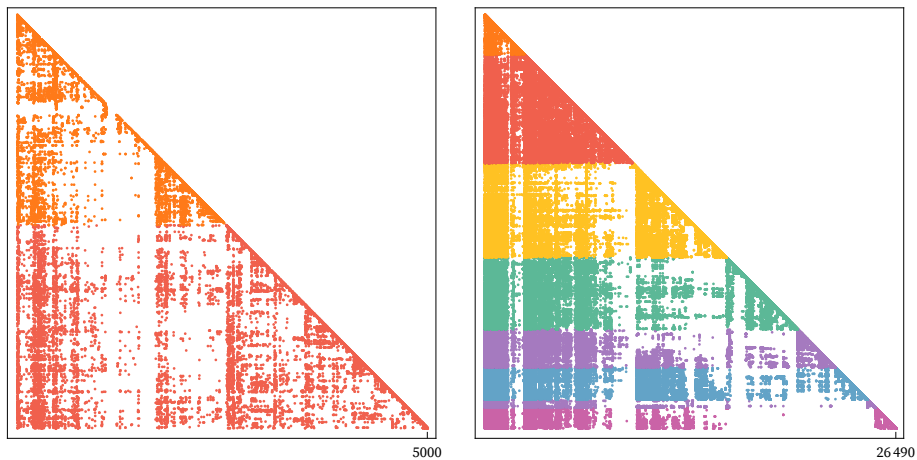

■ logic ■ set theory ■ numbers ■ algebra ■ geometry & topology ■ analysis ■ abstract structures

The results are triangular because theorems in the Metamath database are arranged so that later ones only depend on earlier ones. And while there's considerable patchiness visible, there still seems to be a certain overall background level of uniformity.

In doing this empirical metamathematics we're sampling metamathematical space just through particular "human mathematical settlements" in it. But even from the distribution of these "settlements" we potentially begin to see evidence of a certain background uniformity in metamathematical space.

Perhaps in time as more connections between different areas of mathematics are found human mathematics will gradually become more "uniformly settled" in metamathematical space—and closer to what we might expect from entailment cones and ultimately from the raw ruliad. But it's interesting to see that even with fairly basic empirical metamathematics—operating on a current corpus of human mathematical knowledge—it may already be possible to see signs of some features of physicalized metamathematics.

One day, no doubt, we'll be able do experiments in physics that take our "parsing" of the physical universe in terms of things like space and time and quantum mechanics—and reveal "slices" of the raw ruliad underneath. But perhaps something similar will also be possible in empirical metamathematics: to construct what amounts to a metamathematical microscope (or telescope) through which we can see aspects of the ruliad.

# 27 | Invented or Discovered? How Mathematics Relates to Humans

It's an old and oft-asked question: is mathematics ultimately something that is invented, or something that is discovered? Or, put another way: is mathematics something arbitrarily set up by us humans, or something inevitable and fundamental and in a sense "preexisting", that we merely get to explore? In the past it's seemed as if these were two fundamentally incompatible possibilities. But the framework we've built here in a sense blends them both into a rather unexpected synthesis.

The starting point is the idea that mathematics—like physics—is rooted in the ruliad, which is a representation of formal necessity. Actual mathematics as we "experience" it is—like physics—based on the particular sampling we make of the ruliad. But then the crucial point is that very basic characteristics of us as "observers" are sufficient to constrain that experience to be our general mathematics—or our physics.

At some level we can say that "mathematics is always there"—because every aspect of it is ultimately encoded in the ruliad. But in another sense we can say that the mathematics we have is all "up to us"—because it's based on how we sample the ruliad. But the point is that that sampling is not somehow "arbitrary": if we're talking about mathematics for us humans then it's us ultimately doing the sampling, and the sampling is inevitably constrained by general features of our nature.

A major discovery from our Physics Project is that it doesn't take much in the way of constraints on the observer to deeply constrain the laws of physics they will perceive. And similarly we posit here that for "observers like us" there will inevitably be general ("physicalized") laws of mathematics, that make mathematics inevitably have the general kinds of characteristics we perceive it to have (such as the possibility of doing mathematics at a high level, without always having to drop down to an "atomic" level).

Particularly over the past century there's been the idea that mathematics can be specified in terms of axiom systems, and that these axiom systems can somehow be "invented at will". But our framework does two things. First, it says that "far below" axiom systems is the raw ruliad, which in a sense represents all possible axiom systems. And second, it says that whatever axiom systems we perceive to be "operating" will be ones that we as observers can pick out from the underlying structure of the ruliad.

At a formal level we can "invent" an arbitrary axiom system (and it'll be somewhere in the ruliad), but only certain axiom systems will be ones that describe what we as "mathematical observers" can perceive. In a physics setting we might construct some formal physical theory that talks about detailed patterns in the atoms of space (or molecules in a gas), but the kind of "coarse-grained" observations that we can make won't capture these. Put another way, observers like us can perceive certain kinds of things, and can describe things in terms of these perceptions. But with the wrong kind of theory—or "axioms"—these descriptions won't be sufficient—and only an observer who's "shredded" down to a more "atomic" level will be able to track what's going on.

There's lots of different possible math—and physics—in the ruliad. But observers like us can only "access" a certain type. Some putative alien not like us might access a different type— and might end up with both a different math and a different physics. Deep underneath they—like us—would be talking about the ruliad. But they'd be taking different samples of it, and describing different aspects of it.

For much of the history of mathematics there was a close alignment between the mathematics that was done and what we perceive in the world. For example, Euclidean geometry— with its whole axiomatic structure—was originally conceived just as an idealization of geometrical things that we observe about the world. But by the late 1800s the idea had emerged that one could create "disembodied" axiomatic systems with no particular grounding in our experience in the world.

And, yes, there are many possible disembodied axiom systems that one can set up. And in doing ruliology and generally exploring the computational universe it's interesting to investigate what they do. But the point is that this is something quite different from mathematics as mathematics is normally conceived. Because in a sense mathematics—like physics—is a "more human" activity that's based on what "observers like us" make of the raw formal structure that is ultimately embodied in the ruliad.

When it comes to physics there are, it seems, two crucial features of "observers like us". First, that we're computationally bounded. And second, that we have the perception that we're persistent—and have a definite and continuous thread of experience. At the level of atoms of space, we're in a sense constantly being "remade". But we nevertheless perceive it as always being the "same us".

This single seemingly simple assumption has far-reaching consequences. For example, it leads us to experience a single thread of time. And from the notion that we maintain a continuity of experience from every successive moment to the next we are inexorably led to the idea of a perceived continuum—not only in time, but also for motion and in space. And when combined with intrinsic features of the ruliad and of multicomputation in general, what comes out in the end is a surprisingly precise description of how we'll perceive our universe to operate—that seems to correspond exactly with known core laws of physics.

What does that kind of thinking tell us about mathematics? The basic point is that—since in the end both relate to humans—there is necessarily a close correspondence between physical and mathematical observers. Both are computationally bounded. And the assumption of persistence in time for physical observers becomes for mathematical observers the concept of maintaining coherence as more statements are accumulated. And when combined with intrinsic features of the ruliad and multicomputation this then turns out to imply the kind of physicalized laws of mathematics that we've discussed.

In a formal axiomatic view of mathematics one just imagines that one invents axioms and sees their consequences. But what we're describing here is a view of mathematics that is ultimately just about the ways that we as mathematical observers sample and experience the ruliad. And if we use axiom systems it has to be as a kind of "intermediate language" that helps us make a slightly higher-level description of some corner of the raw ruliad. But actual "human-level" mathematics—like human-level physics—operates at a higher level.

Our everyday experience of the physical world gives us the impression that we have a kind of "direct access" to many foundational features of physics, like the existence of space and the phenomenon of motion. But our Physics Project implies that these are not concepts that are in any sense "already there"; they are just things that emerge from the raw ruliad when you "parse" it in the kinds of ways observers like us do.

In mathematics it's less obvious (at least to all but perhaps experienced pure mathematicians) that there's "direct access" to anything. But in our view of mathematics here, it's ultimately just like physics—and ultimately also rooted in the ruliad, but sampled not by physical observers but by mathematical ones.

So from this point view there's just as much that's "real" underneath mathematics as there is underneath physics. The mathematics is sampled slightly differently (though very similarly)—but we should not in any sense consider it "fundamentally more abstract".

When we think of ourselves as entities within the ruliad, we can build up what we might consider a "fully abstract" description of how we get our "experience" of physics. And we can basically do the same thing for mathematics. So if we take the commonsense point of view that physics fundamentally exists "for real", we're forced into the same point of view for mathematics. In other words, if we say that the physical universe exists, so must we also say that in some fundamental sense, mathematics also exists.

It's not something we as humans "just make", but it is something that is made through our particular way of observing the ruliad, that is ultimately defined by our particular characteristics as observers, with our particular core assumptions about the world, our particular kinds of sensory experience, and so on.

So what can we say in the end about whether mathematics is "invented" or "discovered"? It is neither. Its underpinnings are the ruliad, whose structure is a matter of formal necessity. But its perceived form for us is determined by our intrinsic characteristics as observers. We neither get to "arbitrarily invent" what's underneath, nor do we get to "arbitrarily discover" what's already there. The mathematics we see is the result of a combination of formal necessity in the underlying ruliad, and the particular forms of perception that we—as entities like us—have. Putative aliens could have quite different mathematics, but not because the underlying ruliad is any different for them, but because their forms of perception might be different. And it's the same with physics: even though they "live in the same physical universe" their perception of the laws of physics could be quite different.

## 28 | What Axioms Can There Be for Human Mathematics?

When they were first developed in antiquity the axioms of Euclidean geometry were presumably intended basically as a kind of "tightening" of our everyday impressions of geometry—that would aid in being able to deduce what was true in geometry. But by the mid-1800s—between non-Euclidean geometry, group theory, Boolean algebra and quaternions—it had become clear that there was a range of abstract axiom systems one could in principle consider. And by the time of Hilbert's program around 1900 the pure process of deduction was in effect being viewed as an end in itself—and indeed the core of mathematics—with axiom systems being seen as "starter material" pretty much just "determined by convention".

In practice even today very few different axiom systems are ever commonly used—and indeed in *A New Kind of Science* I was able to list essentially all of them comfortably on a couple of pages. But why these axiom systems and not others? Despite the idea that axiom systems could ultimately be arbitrary, the concept was still that in studying some particular area of mathematics one should basically have an axiom system that would provide a "tight specification" of whatever mathematical object or structure one was trying to talk about. And so, for example, the Peano axioms are what became used for talking about arithmetic-style operations on integers.

In 1931, however, Gödel's theorem showed that actually these axioms weren't strong enough to constrain one to be talking only about integers: there were also other possible models of the axiom system, involving all sorts of exotic "non-standard arithmetic". (And moreover, there was no finite way to "patch" this issue.) In other words, even though the Peano axioms had been invented—like Euclid's axioms for geometry—as a way to describe a definite "intuitive" mathematical thing (in this case, integers) their formal axiomatic structure "had a life of its own" that extended (in some sense, infinitely) beyond its original intended purpose.

Both geometry and arithmetic in a sense had foundations in everyday experience. But for set theory dealing with infinite sets there was never an obvious intuitive base rooted in everyday experience. Some extrapolations from finite sets were clear. But in covering infinite sets various axioms (like the Axiom of Choice) were gradually added to capture what seemed like "reasonable" mathematical assertions.

But one example whose status for a long time wasn't clear was the Continuum Hypothesis—which asserts that the "next distinct possible cardinality" $\aleph_1$ after the cardinality $\aleph_0$ of the integers is $2^{\aleph_0}$: the cardinality of real numbers (i.e. of "the continuum"). Was this something that followed from previously accepted axioms of set theory? And if it was added, would it even be consistent with them? In the early 1960s it was established that actually the Continuum Hypothesis is independent of the other axioms.

With the axiomatic view of the foundations of mathematics that's been popular for the past century or so it seems as if one could, for example, just choose at will whether to include the Continuum Hypothesis (or its negation) as an axiom in set theory. But with the approach to the foundations of mathematics that we've developed here, this is no longer so clear.

Recall that in our approach, everything is ultimately rooted in the ruliad—with whatever mathematics observers like us "experience" just being the result of the particular sampling we do of the ruliad. And in this picture, axiom systems are a particular representation of fairly low-level features of the sampling we do of the raw ruliad.

If we could do any kind of sampling we want of the ruliad, then we'd presumably be able to get all possible axiom systems—as intermediate-level "waypoints" representing different kinds of slices of the ruliad. But in fact by our nature we are observers capable of only certain kinds of sampling of the ruliad.

We could imagine "alien observers" not like us who could for example make whatever choice they want about the Continuum Hypothesis. But given our general characteristics as observers, we may be forced into a particular choice. Operationally, as we've discussed above, the wrong choice could, for example, be incompatible with an observer who "maintains coherence" in metamathematical space.

Let's say we have a particular axiom stated in standard symbolic form. "Underneath" this axiom there will typically be at the level of the raw ruliad a huge cloud of possible configurations of emes that can represent the axiom. But an "observer like us" can only deal with a coarse-grained version in which all these different configurations are somehow considered equivalent. And if the entailments from "nearby configurations" remain nearby, then everything will work out, and the observer can maintain a coherent view of what's going, for example just in terms of symbolic statements about axioms.

But if instead different entailments of raw configurations of emes lead to very different places, the observer will in effect be "shredded"—and instead of having definite coherent "single-minded" things to say about what happens, they'll have to separate everything into all the different cases for different configurations of emes. Or, as we've said it before, the observer will inevitably end up getting "shredded"—and not be able to come up with definite mathematical conclusions.

So what specifically can we say about the Continuum Hypothesis? It's not clear. But conceivably we can start by thinking of $\aleph_0$ as characterizing the "base cardinality" of the ruliad, while $2^{\aleph_0}$ characterizes the base cardinality of a first-level hyperruliad that could for example be based on Turing machines with oracles for their halting problems. And it could be that for us to conclude that the Continuum Hypothesis is false, we'd have to somehow be straddling the ruliad and the hyperruliad, which would be inconsistent with us maintaining a coherent view of mathematics. In other words, the Continuum Hypothesis might somehow be equivalent to what we've argued before is in a sense the most fundamental "contingent fact"—that just as we live in a particular location in physical space—so also we live in the ruliad and not the hyperruliad.

We might have thought that whatever we might see—or construct—in mathematics would in effect be "entirely abstract" and independent of anything about physics, or our experience in the physical world. But particularly insofar as we're thinking about mathematics as done by humans we're dealing with "mathematical observers" that are "made of the same stuff" as physical observers. And this means that whatever general constraints or features exist for physical observers we can expect these to carry over to mathematical observers—so it's no coincidence that both physical and mathematical observers have the same core characteristics, of computational boundedness and "assumption of coherence".

And what this means is that there'll be a fundamental correlation between things familiar from our experience in the physical world and what shows up in our mathematics. We might have thought that the fact that Euclid's original axioms were based on our human perceptions of physical space would be a sign that in some "overall picture" of mathematics they should be considered arbitrary and not in any way central. But the point is that in fact our notions of space are central to our characteristics as observers. And so it's inevitable that "physical-experience-informed" axioms like those for Euclidean geometry will be what appear in mathematics for "observers like us".

# 29 | Counting the Emes of Mathematics and Physics

How does the "size of mathematics" compare to the size of our physical universe? In the past this might have seemed like an absurd question, that tries to compare something abstract and arbitrary with something real and physical. But with the idea that both mathematics and physics as we experience them emerge from our sampling of the ruliad, it begins to seem less absurd.

At the lowest level the ruliad can be thought of as being made up of atoms of existence that we call emes. As physical observers we interpret these emes as atoms of space, or in effect the ultimate raw material of the physical universe. And as mathematical observers we interpret them as the ultimate elements from which the constructs of mathematics are built.

As the entangled limit of all possible computations, the whole ruliad is infinite. But we as physical or mathematical observers sample only limited parts of it. And that means we can meaningfully ask questions like how the number of emes in these parts compare—or, in effect, how big is physics as we experience it compared to mathematics.

In some ways an eme is like a bit. But the concept of emes is that they're "actual atoms of existence"—from which "actual stuff" like the physical universe and its history are made—rather than just "static informational representations" of it. As soon as we imagine that everything is ultimately computational we are immediately led to start thinking of representing it in terms of bits. But the ruliad is not just a representation. It's in some way something lower level. It's the "actual stuff" that everything is made of. And what defines our particular experience of physics or of mathematics is the particular samples we as observers take of what's in the ruliad.

So the question is now how many emes there are in those samples. Or, more specifically, how many emes "matter to us" in building up our experience.

Let's return to an analogy we've used several times before: a gas made of molecules. In the volume of a room there might be $10^{27}$ individual molecules, each on average colliding every $10^{-10}$ seconds. So that means that our "experience of the room" over the course of a minute or so might sample $10^{39}$ collisions. Or, in terms closer to our Physics Project, we might say that there are perhaps $10^{39}$ "collision events" in the causal graph that defines what we experience.

But these "collision events" aren't something fundamental; they have what amounts to "internal structure" with many associated parameters about location, time, molecular configuration, etc.

Our Physics Project, however, suggests that—far below for example our usual notions of space and time—we can in fact have a truly fundamental definition of what's happening in the universe, ultimately in terms of emes. We don't yet know the "physical scale" for this—and in the end we presumably need experiments to determine that. But rather rickety estimates based on a variety of assumptions suggest that the elementary length might be around $10^{-90}$ meters, with the elementary time being around $10^{-100}$ seconds.

And with these estimates we might conclude that our "experience of a room for a minute" would involve sampling perhaps $10^{370}$ update events, that create about this number of atoms of space.

But it's immediately clear that this is in a sense a gross underestimate of the total number of emes that we're sampling. And the reason is that we're not accounting for quantum mechanics, and for the multiway nature of the evolution of the universe. We've so far only considered one "thread of time" at one "position in branchial space". But in fact there are many threads of time, constantly branching and merging. So how many of these do we experience?

In effect that depends on our size in branchial space. In physical space "human scale" is of order a meter—or perhaps $10^{90}$ elementary lengths. But how big is it in branchial space?

The fact that we're so large compared to the elementary length is the reason that we consistently experience space as something continuous. And the analog in branchial space is that if we're big compared to the "elementary branchial distance between branches" then we won't experience the different individual histories of these branches, but only an aggregate "objective reality" in which we conflate together what happens on all the branches. Or, put another way, being large in branchial space is what makes us experience classical physics rather than quantum mechanics.

Our estimates for branchial space are even more rickety than for physical space. But conceivably there are on the order of $10^{120}$ "instantaneous parallel threads of time" in the universe, and $10^{20}$ encompassed by our instantaneous experience—implying that in our minute-long experience we might sample a total of on the order of close to $10^{500}$ emes.

But even this is a vast underestimate. Yes, it tries to account for our extent in physical space and in branchial space. But then there's also rulial space—which in effect is what "fills out" the whole ruliad. So how big are we in that space? In essence that's like asking how many different possible sequences of rules there are that are consistent with our experience.

The total conceivable number of sequences associated with $10^{500}$ emes is roughly the number of possible hypergraphs with $10^{500}$ nodes—or around $(10^{500})^{10^{500}}$. But the actual number consistent with our experience is smaller, in particular as reflected by the fact that we attribute specific laws to our universe. But when we say "specific laws" we have to recognize that there is a finiteness to our efforts at inductive inference which inevitably makes these laws at least somewhat uncertain to us. And in a sense that uncertainty is what represents our "extent in rulial space".

But if we want to count the emes that we "absorb" as physical observers, it's still going to be a huge number. Perhaps the base may be lower—say $10^{10}$—but there's still a vast exponent, suggesting that if we include our extent in rulial space, we as physical observers may experience numbers of emes like $(10^{10})^{10^{500}}$.

But let's say we go beyond our "everyday human-scale experience". For example, let's ask about "experiencing" our whole universe. In physical space, the volume of our current universe is about $10^{78}$ times larger than "human scale" (while human scale is perhaps $10^{270}$ times larger than the "scale of the atoms of space"). In branchial space, conceivably our current universe is $10^{100}$ times larger than "human scale". But these differences absolutely pale in comparison to the sizes associated with rulial space.

We might try to go beyond "ordinary human experience" and for example measure things using tools from science and technology. And, yes, we could then think about "experiencing" lengths down to $10^{-22}$ meters, or something close to "single threads" of quantum histories. But in the end, it's still the rulial size that dominates, and that's where we can expect most of the vast number of emes that form of our experience of the physical universe to come from.

OK, so what about mathematics? When we think about what we might call human-scale mathematics, and talk about things like the Pythagorean theorem, how many emes are there "underneath"? "Compiling" our theorem down to typical traditional mathematical axioms, we've seen that we'll routinely end up with expressions containing, say, $10^{20}$ symbolic elements. But what happens if we go "below that", compiling these symbolic elements—which might include things like variables and operators—into "pure computational elements" that we can think of as emes? We've seen a few examples, say with combinators, that suggest that for the traditional axiomatic structures of mathematics, we might need another factor of maybe roughly $10^{10}$.

These are incredibly rough estimates, but perhaps there's a hint that there's "further to go" to get from human-scale for a physical observer down to atoms of space that correspond to emes, than there is to get from human-scale for a mathematical observer down to emes.

Just like in physics, however, this kind of "static drill-down" isn't the whole story for mathematics. When we talk about something like the Pythagorean theorem, we're really referring to a whole cloud of "human-equivalent" points in metamathematical space. The total number of "possible points" is basically the size of the entailment cone that contains something like the Pythagorean theorem. The "height" of the entailment cone is related to typical lengths of proofs—which for current human mathematics might be perhaps hundreds of steps.

And this would lead to overall sizes of [entailment cones](#) of very roughly $10^{100}$ theorems. But within this "how big" is the cloud of variants corresponding to particular "human-recognized" theorems? Empirical metamathematics could provide additional data on this question. But if we very roughly imagine that half of every proof is "flexible", we'd end up with things like $10^{50}$ variants. So if we asked how many emes correspond to the "experience" of the Pythagorean theorem, it might be, say, $10^{80}$.

To give an analogy of "everyday physical experience" we might consider a mathematician thinking about mathematical concepts, and maybe in effect pondering a few tens of theorems per minute—implying according to our extremely rough and speculative estimates that while typical "specific human-scale physics experience" might involve $10^{500}$ emes, specific human-scale mathematics experience might involve $10^{80}$ emes (a number comparable, for example, to the number of physical atoms in our universe).

What if instead of considering "everyday mathematical experience" we consider [all humanly explored mathematics](#)? On the scales we're describing, the factors are not large. In the history of human mathematics, only a few million theorems have been published. If we think about all the computations that have been done in the service of mathematics, it's a somewhat larger factor. I suspect Mathematica is the dominant contributor here—and we can estimate that the total number of Wolfram Language operations corresponding to "human-level mathematics" done so far is perhaps $10^{20}$.

But just like for physics, all these numbers pale in comparison with those introduced by rulial sizes. We've talked essentially about a particular path from emes through specific axioms to theorems. But the ruliad in effect contains all possible axiom systems. And if we start thinking about enumerating these—and effectively "populating all of rulial space"—we'll end up with exponentially more emes.

But as with the perceived laws of physics, in mathematics as done by humans it's actually just a narrow slice of rulial space that we're sampling. It's like a generalization of the idea that something like arithmetic as we imagine it can be derived from a whole cloud of possible axiom systems. It's not just one axiom system; but it's also not all possible axiom systems.

One can imagine doing some combination of ruliology and empirical metamathematics to get an estimate of "how broad" human-equivalent axiom systems (and their construction from emes) might be. But the answer seems likely to be much smaller than the kinds of sizes we have been estimating for physics.

It's important to emphasize that what we've discussed here is extremely rough—and speculative. And indeed I view its main value as being to provide an example of how to imagine thinking through things in the context of the ruliad and the framework around it. But on the basis of what we've discussed, we might make the very tentative conclusion that "human-experienced physics" is bigger than "human-experienced mathematics". Both involve vast numbers of emes. But physics seems to involve a lot more. In a sense—even with all its abstraction—the suspicion is that there's "less ultimately in mathematics" as far as we're concerned than there is in physics. Though by any ordinary human standards, mathematics still involves absolutely vast numbers of emes.

# 30 | Some Historical (and Philosophical) Background

The human activity that we now call "mathematics" can presumably trace its origins into prehistory. What might have started as "a single goat", "a pair of goats", etc. became a story of abstract numbers that could be indicated purely by things like tally marks. In Babylonian times the practicalities of a city-based society led to all sorts of calculations involving arithmetic and geometry—and basically everything we now call "mathematics" can ultimately be thought of as a generalization of these ideas.

The tradition of philosophy that emerged in Greek times saw mathematics as a kind of reasoning. But while much of arithmetic (apart from issues of infinity and infinitesimals) could be thought of in explicit calculational ways, precise geometry immediately required an idealization—specifically the concept of a point having no extent, or equivalently, the continuity of space. And in an effort to reason on top of this idealization, there emerged the idea of defining axioms and making abstract deductions from them.

But what kind of a thing actually was mathematics? Plato talked about things we sense in the external world, and things we conceptualize in our internal thoughts. But he considered mathematics to be at its core an example of a third kind of thing: something from an abstract world of ideal forms. And with our current thinking, there is an immediate resonance between this concept of ideal forms and the concept of the ruliad.

But for most of the past two millennia of the actual development of mathematics, questions about what it ultimately was lay in the background. An important step was taken in the late 1600s when Newton and others "mathematicized" mechanics, at first presenting what they did in the form of axioms similar to Euclid's. Through the 1700s mathematics as a practical field was viewed as some kind of precise idealization of features of the world—though with an increasingly elaborate tower of formal derivations constructed in it. Philosophy, meanwhile, typically viewed mathematics—like logic—mostly as an example of a system in which there was a formal process of derivation with a "necessary" structure not requiring reference to the real world.

But in the first half of the 1800s there arose several examples of systems where axioms—while inspired by features of the world—ultimately seemed to be "just invented" (e.g. group theory, curved space, quaternions, Boolean algebra, ...). A push towards increasing rigor (especially for calculus and the nature of real numbers) led to more focus on axiomatization and formalization—which was still further emphasized by the appearance of a few non-constructive "purely formal" proofs.

But if mathematics was to be formalized, what should its underlying primitives be? One obvious choice seemed to be logic, which had originally been developed by Aristotle as a kind of catalog of human arguments, but two thousand years later felt basic and inevitable. And so it was that Frege, followed by Whitehead and Russell, tried to start "constructing mathematics" from "pure logic" (along with set theory). Logic was in a sense a rather low-level "machine code", and it took hundreds of pages of unreadable (if impressive-looking) "code" for Whitehead and Russell, in their 1910 *Principia Mathematica*, to get to 1+1=2.

*Courtesy of archive.org*

Meanwhile, starting around 1900, Hilbert took a slightly different path, essentially representing everything with what we would now call symbolic expressions, and setting up axioms as relations between these. But what axioms should be used? Hilbert seemed to feel that the core of mathematics lay not in any "external meaning" but in the pure formal structure built up from whatever axioms were used. And he imagined that somehow all the truths of mathematics could be "mechanically derived" from axioms, a bit, as he said in a certain resonance with our current views, like the "great calculating machine, Nature" does it for physics.

Not all mathematicians, however, bought into this "formalist" view of what mathematics is. And in 1931 Gödel managed to prove from inside the formal axiom system traditionally used for arithmetic that this system had a fundamental incompleteness that prevented it from ever having anything to say about certain mathematical statements. But Gödel's theorem seems to have maintained a more Platonic belief about mathematics: that even though the axiomatic method falls short, the truths of mathematics are in some sense still "all there", and it's potentially possible for the human mind to have "direct access" to them. And while this is not quite the same as our picture of the mathematical observer accessing the ruliad, there's again some definite resonance here.

But, OK, so how has mathematics actually conducted itself over the past century? Typically there's at least lip service paid to the idea that there are "axioms underneath"—usually assumed to be those from set theory. There's been significant emphasis placed on the idea of formal deduction and proof—but not so much in terms of formally building up from axioms as in terms of giving narrative expositions that help humans understand why some theorem might follow from other things they know.

There's been a field of "mathematical logic" concerned with using mathematics-like methods to explore mathematics-like aspects of formal axiomatic systems. But (at least until very recently) there's been rather little interaction between this and the "mainstream" study of mathematics. And for example phenomena like undecidability that are central to mathematical logic have seemed rather remote from typical pure mathematics—even though many actual long-unsolved problems in mathematics do seem likely to run into it.

But even if formal axiomatization may have been something of a sideshow for mathematics, its ideas have brought us what is without much doubt the single most important intellectual breakthrough of the twentieth century: the abstract concept of computation. And what's now become clear is that computation is in some fundamental sense much more general than mathematics.

At a philosophical level one can view the ruliad as containing all computation. But mathematics (at least as it's done by humans) is defined by what a "mathematical observer like us" samples and perceives in the ruliad.

The most common "core workflow" for mathematicians doing pure mathematics is first to imagine what might be true (usually through a process of intuition that feels a bit like making "direct access to the truths of mathematics")—and then to "work backwards" to try to construct a proof. As a practical matter, though, the vast majority of "mathematics done in the world" doesn't follow this workflow, and instead just "runs forward"—doing computation. And there's no reason for at least the innards of that computation to have any "humanized character" to it; it can just involve the raw processes of computation.

But the traditional pure mathematics workflow in effect depends on using "human-level" steps. Or if, as we described earlier, we think of low-level axiomatic operations as being like molecular dynamics, then it involves operating at a "fluid dynamics" level.

A century ago efforts to "globally understand mathematics" centered on trying to find common axiomatic foundations for everything. But as different areas of mathematics were explored (and particularly ones like algebraic topology that cut across existing disciplines) it began to seem as if there might also be "top-down" commonalities in mathematics, in effect directly at the "fluid dynamics" level. And within the last few decades, it's become increasingly common to use ideas from category theory as a general framework for thinking about mathematics at a high level.

But there's also been an effort to progressively build up—as an abstract matter—formal "higher category theory". A notable feature of this has been the appearance of connections to both geometry and mathematical logic—and for us a connection to the ruliad and its features.

The success of category theory has led in the past decade or so to interest in other high-level structural approaches to mathematics. A notable example is homotopy type theory. The basic concept is to characterize mathematical objects not by using axioms to describe properties they should have, but instead to use "types" to say "what the objects are" (for example, "mapping from reals to integers"). Such type theory has the feature that it tends to look much more "immediately computational" than traditional mathematical structures and notation—as well as making explicit proofs and other metamathematical concepts. And in fact questions about types and their equivalences wind up being very much like the questions we've discussed for the multiway systems we're using as metamodels for mathematics.

Homotopy type theory can itself be set up as a formal axiomatic system—but with axioms that include what amount to metamathematical statements. A key example is the univalence axiom which essentially states that things that are equivalent can be treated as the same. And now from our point of view here we can see this being essentially a statement of metamathematical coarse graining—and a piece of defining what should be considered "mathematics" on the basis of properties assumed for a mathematical observer.

When Plato introduced ideal forms and their distinction from the external and internal world the understanding of even the fundamental concept of computation—let alone multicomputation and the ruliad—was still more than two millennia in the future. But now our picture is that everything can in a sense be viewed as part of the world of ideal forms that is the ruliad—and that not only mathematics but also physical reality are in effect just manifestations of these ideal forms.

But a crucial aspect is how we sample the "ideal forms" of the ruliad. And this is where the "contingent facts" about us as human "observers" enter. The formal axiomatic view of mathematics can be viewed as providing one kind of low-level description of the ruliad. But the point is that this description isn't aligned with what observers like us perceive—or with what we will successfully be able to view as human-level mathematics.

A century ago there was a movement to take mathematics (as well, as it happens, as other fields) beyond its origins in what amount to human perceptions of the world. But what we now see is that while there is an underlying "world of ideal forms" embodied in the ruliad that has nothing to do with us humans, mathematics as we humans do it must be associated with the particular sampling we make of that underlying structure.

And it's not as if we get to pick that sampling "at will"; the sampling we do is the result of fundamental features of us as humans. And an important point is that those fundamental features determine our characteristics both as mathematical observers and as physical observers. And this fact leads to a deep connection between our experience of physics and our definition of mathematics.

Mathematics historically began as a formal idealization of our human perception of the physical world. Along the way, though, it began to think of itself as a more purely abstract pursuit, separated from both human perception and the physical world. But now, with the general idea of computation, and more specifically with the concept of the ruliad, we can in a sense see what the limit of such abstraction would be. And interesting though it is, what we're now discovering is that it's not the thing we call mathematics. And instead, what we call mathematics is something that is subtly but deeply determined by general features of human perception—in fact, essentially the same features that also determine our perception of the physical world.

The intellectual foundations and justification are different now. But in a sense our view of mathematics has come full circle. And we can now see that mathematics is in fact deeply connected to the physical world and our particular perception of it. And we as humans can do what we call mathematics for basically the same reason that we as humans manage to parse the physical world to the point where we can do science about it.

# 31 | Implications for the Future of Mathematics

Having talked a bit about historical context let's now talk about what the things we've discussed here mean for the future of mathematics—both in theory and in practice.

At a theoretical level we've characterized the story of mathematics as being the story of a particular way of exploring the [ruliad](). And from this we might think that in some sense the ultimate limit of mathematics would be to just deal with the ruliad as a whole. But observers like us—at least doing mathematics the way we normally do it—simply can't do that. And in fact, with the limitations we have as mathematical observers we can inevitably sample only tiny slices of the ruliad.

But as we've discussed, it is exactly this that leads us to experience the kinds of "general laws of mathematics" that we've talked about. And it is from these laws that we get a picture of the "large-scale structure of mathematics"—that turns out to be in many ways similar to the picture of the large-scale structure of our physical universe that we get from physics.

As we've discussed, what corresponds to the coherent structure of physical space is the possibility of doing mathematics in terms of high-level concepts—without always having to drop down to the "atomic" level. Effective uniformity of metamathematical space then leads to the idea of "pure metamathematical motion", and in effect the possibility of translating at a high level between different areas of mathematics. And what this suggests is that in some sense "all high-level areas of mathematics" should ultimately be connected by "high-level dualities"—some of which have already been seen, but many of which remain to be discovered.

Thinking about metamathematics in physicalized terms also suggests another phenomenon: essentially an [analog of gravity for metamathematics](). As we discussed earlier, in direct analogy to the way that "larger densities of activity" in the spatial hypergraph for physics lead to a deflection in geodesic paths in physical space, so also larger "entailment density" in metamathematical space will lead to deflection in geodesic paths in metamathematical space. And when the entailment density gets sufficiently high, it presumably becomes inevitable that these paths will all converge, leading to what one might think of as a "metamathematical singularity".

In the spacetime case, a typical analog would be a place where all geodesics have finite length, or in effect "time stops". In our view of metamathematics, it corresponds to a situation where "all proofs are finite"—or, in other words, where everything is decidable, and there is no more "fundamental difficulty" left.

Absent other effects we might imagine that in the physical universe the effects of gravity would eventually lead everything to collapse into black holes. And the analog in metamathematics would be that everything in mathematics would "collapse" into decidable theories. But among the effects not accounted for is continued expansion—or in effect the creation of new physical or metamathematical space, formed in a sense by underlying raw computational processes.

What will observers like us make of this, though? In statistical mechanics an observer who does coarse graining might perceive the "heat death of the universe". But at a molecular level there is all sorts of detailed motion that reflects a continued irreducible process of computation. And inevitably there will be an infinite collection of possible "slices of reducibility" to be found in this—just not necessarily ones that align with any of our current capabilities as observers.

What does this mean for mathematics? Conceivably it might suggest that there's only so much that can fundamentally be discovered in "high-level mathematics" without in effect "expanding our scope as observers"—or in essence changing our definition of what it is we humans mean by doing mathematics.

But underneath all this is still raw computation—and the ruliad. And this we know goes on forever, in effect continually generating "irreducible surprises". But how should we study "raw computation"?

In essence we want to do unfettered exploration of the computational universe, of the kind I did in *A New Kind of Science*, and that we now call the science of ruliology. It's something we can view as more abstract and more fundamental than mathematics—and indeed, as we've argued, it's for example what's underneath not only mathematics but also physics.

Ruliology is a rich intellectual activity, important for example as a source of models for many processes in nature and elsewhere. But it's one where computational irreducibility and undecidability are seen at almost every turn—and it's not one where we can readily expect "general laws" accessible to observers like us, of the kind we've seen in physics, and now see in mathematics.

We've argued that with its foundation in the ruliad mathematics is ultimately based on structures lower level than axiom systems. But given their familiarity from the history of mathematics, it's convenient to use axiom systems—as we have done here—as a kind of "intermediate-scale metamodel" for mathematics.

But what is the "workflow" for using axiom systems? One possibility in effect inspired by ruliology is just to systematically construct the entailment cone for an axiom system,

*The Physicalization of Metamathematics and Its Implications for the Foundations of Mathematics* | 187progressively generating all possible theorems that the axiom system implies. But while doing this is of great theoretical interest, it typically isn't something that will in practice reach much in the way of (currently) familiar mathematical results.

But let's say one's thinking about a particular result. A proof of this would correspond to a path within the entailment cone. And the idea of automated theorem proving is to systematically find such a path—which, with a variety of tricks, can usually be done vastly more efficiently than just by enumerating everything in the entailment cone. In practice, though, despite half a century of history, automated theorem proving has seen very little use in mainstream mathematics. Of course it doesn't help that in typical mathematical work a proof is seen as part of the high-level exposition of ideas—but automated proofs tend to operate at the level of "axiomatic machine code" without any connection to human-level narrative.

But if one doesn't already know the result one's trying to prove? Part of the intuition that comes from *A New Kind of Science* is that there can be "interesting results" that are still simple enough that they can conceivably be found by some kind of explicit search—and then verified by automated theorem proving. But so far as I know, only one significant unexpected result has so far ever been found in this way with automated theorem proving: my 2000 result on the simplest axiom system for Boolean algebra.

And the fact is that when it comes to using computers for mathematics, the overwhelming fraction of the time they're used not to construct proofs, but instead to do "forward computations" and "get results" (yes, often with Mathematica). Of course, within those forward computations, there are many operations—like Reduce, SatisfiableQ, PrimeQ, etc.—that essentially work by internally finding proofs, but their output is "just results" not "why-it's-true explanations". (FindEquationalProof—as its name suggests—is a case where an actual proof is generated.)

Whether one's thinking in terms of axioms and proofs, or just in terms of "getting results", one's ultimately always dealing with computation. But the key question is how that computation is "packaged". Is one dealing with arbitrary, raw, low-level constructs, or with something higher level and more "humanized"?

As we've discussed, at the lowest level, everything can be represented in terms of the ruliad. But when we do both mathematics and physics what we're perceiving is not the raw ruliad, but rather just certain high-level features of it. But how should these be represented? Ultimately we need a language that we humans understand, that captures the particular features of the underlying raw computation that we're interested in.

From our computational point of view, mathematical notation can be thought of as a rough attempt at this. But the most complete and systematic effort in this direction is the one I've worked towards for the past several decades: what's now the full-scale computational language that is the Wolfram Language (and Mathematica).

Ultimately the Wolfram Language can represent any computation. But the point is to make it easy to represent the computations that people care about: to capture the high-level constructs (whether they're polynomials, geometrical objects or chemicals) that are part of modern human thinking.

The process of language design (on which, yes, I've spent immense amounts of time) is a curious mixture of art and science, that requires both drilling down to the essence of things, and creatively devising ways to make those things accessible and cognitively convenient for humans. At some level it's a bit like deciding on words as they might appear in a human language—but it's something more structured and demanding.

And it's our best way of representing "high-level" mathematics: mathematics not at the axiomatic (or below) "machine code" level, but instead at the level human mathematicians typically think about it.

We've definitely not "finished the job", though. Wolfram Language currently has around 7000 built-in primitive constructs, of which at least a couple of thousand can be considered "primarily mathematical". But while the language has long contained constructs for algebraic numbers, random walks and finite groups, it doesn't (yet) have built-in constructs for algebraic topology or K-theory. In recent years we've been slowly adding more kinds of pure-mathematical constructs—but to reach the frontiers of modern human mathematics might require perhaps a thousand more. And to make them useful all of them will have to be carefully and coherently designed.

The great power of the Wolfram Language comes not only from being able to represent things computationally, but also being able to compute with things, and get results. And it's one thing to be able to represent some pure mathematical construct—but quite another to be able to broadly compute with it.

The Wolfram Language in a sense emphasizes the "forward computation" workflow. Another workflow that's achieved some popularity in recent years is the proof assistant one—in which one defines a result and then as a human one tries to fill in the steps to create a proof of it, with the computer verifying that the steps correctly fit together. If the steps are low level then what one has is something like typical automated theorem proving—though now being attempted with human effort rather than being done automatically.

In principle one can build up to much higher-level "steps" in a modular way. But now the problem is essentially the same as in computational language design: to create primitives that are both precise enough to be immediately handled computationally, and "cognitively convenient" enough to be usefully understood by humans. And realistically once one's done the design (which, after decades of working on such things, I can say is hard), there's likely to be much more "leverage" to be had by letting the computer just do computations than by expending human effort (even with computer assistance) to put together proofs.

One might think that a proof would be important in being sure one's got the right answer. But as we've discussed, that's a complicated concept when one's dealing with human-level mathematics. If we go to a full axiomatic level it's very typical that there will be all sorts of pedantic conditions involved. Do we have the "right answer" if underneath we assume that 1/0=0? Or does this not matter at the "fluid dynamics" level of human mathematics?

One of the great things about computational language is that—at least if it's written well—it provides a clear and succinct specification of things, just like a good "human proof" is supposed to. But computational language has the great advantage that it can be run to create new results—rather than just being used to check something.

It's worth mentioning that there's another potential workflow beyond "compute a result" and "find a proof". It's "here's an object or a set of constraints for creating one; now find interesting facts about this". Type into Wolfram|Alpha something like sin^4(x) (and, yes, there's "natural math understanding" needed to translate something like this to precise Wolfram Language). There's nothing obvious to "compute" here. But instead what Wolfram|Alpha does is to "say interesting things" about this—like what its maximum or its integral over a period is.

In principle this is a bit like exploring the entailment cone—but with the crucial additional piece of picking out which entailments will be "interesting to humans". (And implementationally it's a very deeply constrained exploration.)

It's interesting to compare these various workflows with what one can call experimental mathematics. Sometimes this term is basically just applied to studying explicit examples of known mathematical results. But the much more powerful concept is to imagine discovering new mathematical results by "doing experiments".

Usually these experiments are not done at the level of axioms, but rather at a considerably higher level (e.g. with things specified using the primitives of Wolfram Language). But the typical pattern is to enumerate a large number of cases and to see what happens—with the most exciting result being the discovery of some unexpected phenomenon, regularity or irregularity.

This type of approach is in a sense much more general than mathematics: it can be applied to anything computational, or anything described by rules. And indeed it is the core methodology of ruliology, and what it does to explore the computational universe—and the ruliad.

One can think of the typical approach in pure mathematics as representing a gradual expansion of the entailment fabric, with humans checking (perhaps with a computer) statements they consider adding. Experimental mathematics effectively strikes out in some "direction" in metamathematical space, potentially jumping far away from the entailment fabric currently within the purview of some mathematical observer.

And one feature of this—very common in ruliology—is that one may run into undecidability. The "nearby" entailment fabric of the mathematical observer is in a sense "filled in enough" that it doesn't typically have infinite proof paths of the kind associated with undecidability. But something reached by experimental mathematics has no such guarantee.

What's good of course is that experimental mathematics can discover phenomena that are "far away" from existing mathematics. But (like in automated theorem proving) there isn't necessarily any human-accessible "narrative explanation" (and if there's undecidability there may be no "finite explanation" at all).

So how does this all relate to our whole discussion of new ideas about the foundations of mathematics? In the past we might have thought that mathematics must ultimately progress just by working out more and more consequences of particular axioms. But what we've argued is that there's a fundamental infrastructure even far below axiom systems—whose low-level exploration is the subject of ruliology. But the thing we call mathematics is really something higher level.

Axiom systems are some kind of intermediate modeling layer—a kind of "assembly language" that can be used as a wrapper above the "raw ruliad". In the end, we've argued, the details of this language won't matter for typical things we call mathematics. But in a sense the situation is very much like in practical computing: we want an "assembly language" that makes it easiest to do the typical high-level things we want. In practical computing that's often achieved with RISC instruction sets. In mathematics we typically imagine using axiom systems like ZFC. But—as reverse mathematics has tended to indicate—there are probably much more accessible axiom systems that could be used to reach the mathematics we want. (And ultimately even ZFC is limited in what it can reach.)

But if we could find such a "RISC" axiom system for mathematics it has the potential to make practical more extensive exploration of the entailment cone. It's also conceivable—though not guaranteed—that it could be "designed" to be more readily understood by humans. But in the end actual human-level mathematics will typically operate at a level far above it.

And now the question is whether the "physicalized general laws of mathematics" that we've discussed can be used to make conclusions directly about human-level mathematics. We've identified a few features—like the very possibility of high-level mathematics, and the expectation of extensive dualities between mathematical fields. And we know that basic commonalities in structural features can be captured by things like category theory. But the question is what kinds of deeper general features can be found, and used.

In physics our everyday experience immediately makes us think about "large-scale features" far above the level of atoms of space. In mathematics our typical experience so far has been at a lower level. So now the challenge is to think more globally, more metamathematically and, in effect, more like in physics.

In the end, though, what we call mathematics is what mathematical observers perceive. So if we ask about the future of mathematics we must also ask about the future of mathematical observers.

If one looks at the history of physics there was already much to understand just on the basis of what we humans could "observe" with our unaided senses. But gradually as more kinds of detectors became available—from microscopes to telescopes to amplifiers and so on—the domain of the physical observer was expanded, and the perceived laws of physics with it.

And today, as the practical computational capability of observers increases, we can expect that we'll gradually see new kinds of physical laws (say associated with hitherto "it's just random" molecular motion or other features of systems).

As we've discussed above, we can see our characteristics as physical observers as being associated with "experiencing" the ruliad from one particular "vantage point" in rulial space (just as we "experience" physical space from one particular vantage point in physical space). Putative "aliens" might experience the ruliad from a different vantage point in rulial space—leading them to have laws of physics utterly incoherent with our own. But as our technology and ways of thinking progress, we can expect that we'll gradually be able to expand our "presence" in rulial space (just as we do with spacecraft and telescopes in physical space). And so we'll be able to "experience" different laws of physics.

We can expect the story to be very similar for mathematics. We have "experienced" mathematics from a certain vantage point in the ruliad. Putative aliens might experience it from another point, and build their own "paramathematics" utterly incoherent with our mathematics. The "natural evolution" of our mathematics corresponds to a gradual expansion in the entailment fabric, and in a sense a gradual spreading in rulial space. Experimental mathematics has the potential to launch a kind of "metamathematical space probe" which can discover quite different mathematics. At first, though, this will tend to be a piece of "raw ruliology". But, if pursued, it potentially points the way to a kind of "colonization of rulial space" that will gradually expand the domain of the mathematical observer.

The physicalized general laws of mathematics we've discussed here are based on features of current mathematical observers (which in turn are highly based on current physical observers). What these laws would be like with "enhanced" mathematical observers we don't yet know.

Mathematics as it is today is a great example of the "humanization of raw computation". Two other examples are theoretical physics and computational language. And in all cases there is the potential to gradually expand our scope as observers. It'll no doubt be a mixture of technology and methods along with expanded cognitive frameworks and understanding. We can use ruliology—or experimental mathematics—to "jump out" into the raw ruliad. But most of what we'll see is "non-humanized" computational irreducibility.

But perhaps somewhere there'll be another slice of computational reducibility: a different "island" on which "alien" general laws can be built. But for now we exist on our current "island" of reducibility. And on this island we see the particular kinds of general laws that we've discussed. We saw them first in physics. But there we discovered that they could emerge quite generically from a lower-level computational structure—and ultimately from the very general structure that we call the ruliad. And now, as we've discussed here, we realize that the thing we call mathematics is actually based on exactly the same foundations—with the result that it should show the same kinds of general laws.

It's a rather different view of mathematics—and its foundations—than we've been able to form before. But the deep connection with physics that we've discussed allows us to now have a physicalized view of metamathematics, which informs both what mathematics really is now, and what the future can hold for the remarkable pursuit that we call mathematics.

## *Some Personal History: The Evolution of These Ideas*

It's been a long personal journey to get to the ideas described here—stretching back nearly 45 years. Parts have been quite direct, steadily building over the course of time. But other parts have been surprising—even shocking. And to get to where we are now has required me to rethink some very long-held assumptions, and adopt what I had believed was a rather different way of thinking—even though, ironically, I've realized in the end that many aspects of this way of thinking pretty much mirror what I've done all along at a practical and technological level.

Back in the late 1970s as a young theoretical physicist I had discovered the "secret weapon" of using computers to do mathematical calculations. By 1979 I had outgrown existing systems and decided to build my own. But what should its foundations be? A key goal was to represent the processes of mathematics in a computational way. I thought about the methods I'd found effective in practice. I studied the history of mathematical logic. And in the end I came up with what seemed to me at the time the most obvious and direct approach: that everything should be based on transformations for symbolic expressions.

I was pretty sure this was actually a good general approach to computation of all kinds—and the system we released in 1981 was named SMP ("Symbolic Manipulation Program") to reflect this generality. History has indeed borne out the strength of the symbolic expression paradigm—and it's from that we've been able to build the huge tower of technology that is the modern Wolfram Language. But all along mathematics has been an important use case—and in effect we've now seen four decades of validation that the core idea of transformations on symbolic expressions is a good metamodel of mathematics.

When Mathematica was first released in 1988 we called it "A System for Doing Mathematics by Computer", where by "doing mathematics" we meant doing computations in mathematics and getting results. People soon did all sorts of experiments on using Mathematica to create and present proofs. But the overwhelming majority of actual usage was for directly computing results—and almost nobody seemed interested in seeing the inner workings, presented as a proof or otherwise.

But in the 1980s I had started my work on exploring the computational universe of simple programs like cellular automata. And doing this was all about looking at the ongoing behavior of systems—or in effect the (often computationally irreducible) history of computations. And even though I sometimes talked about using my computational methods to do "experimental mathematics", I don't think I particularly thought about the actual progress of the computations I was studying as being like mathematical processes or proofs.

In 1991 I started working on what became *A New Kind of Science*, and in doing so I tried to systematically study possible forms of computational processes—and I was soon led to substitution systems and symbolic systems which I viewed in their different ways as being minimal idealizations of what would become Wolfram Language, as well as to multiway systems. There were some areas to which I was pretty sure the methods of *A New Kind of Science* would apply. Three that I wasn't sure about were biology, physics and mathematics.

But by the late 1990s I had worked out quite a bit about the first two, and started looking at mathematics. I knew that Mathematica and what would become Wolfram Language were good representations of "practical mathematics". But I assumed that to understand the foundations of mathematics I should look at the traditional low-level representation of mathematics: axiom systems.

And in doing this I was soon able to simplify to multiway systems—with proofs being paths:

I had long wondered what the detailed relationships between things like my idea of computational irreducibility and earlier results in mathematical logic were. And I was pleased at how well many things could be clarified—and explicitly illustrated—by thinking in terms of multiway systems.

My experience in exploring simple programs in general had led to the conclusion that computational irreducibility and therefore undecidability were quite ubiquitous. So I considered it quite a mystery why undecidability seemed so rare in the mathematics that mathematicians typically did. I suspected that in fact undecidability was lurking close at hand—and I got some evidence of that by doing experimental mathematics. But why weren't mathematicians running into this more? I came to suspect that it had something to do with the history of mathematics, and with the idea that mathematics had tended to expand its subject matter by asking "How can this be generalized while still having such-and-such a theorem be true?"

But I also wondered about the particular axiom systems that had historically been used for mathematics. They all fit easily on a couple of pages. But why these and not others? Following my general "ruliological" approach of exploring all possible systems I started just enumerating possible axiom systems—and soon found out that many of them had rich and complicated implications.

But where among these possible systems did the axiom systems historically used in mathematics lie? I did searches, and at about the 50,000th axiom was able to find the simplest axiom system for Boolean algebra. Proving that it was correct gave me my first serious experience with automated theorem proving.

But what kind of a thing was the proof? I made some attempt to understand it, but it was clear that it wasn't something a human could readily understand—and reading it felt a bit like trying to read machine code. I recognized that the problem was in a sense a lack of "human connection points"—for example of intermediate lemmas that (like words in a human language) had a contextualized significance. I wondered about how one could find lemmas that "humans would care about"? And I was surprised to discover that at least for the "named theorems" of Boolean algebra a simple criterion could reproduce them.

Quite a few years went by. Off and on I thought about two ultimately related issues. One was how to represent the execution histories of Wolfram Language programs. And the other was how to represent proofs. In both cases there seemed to be all sorts of detail, and it seemed difficult to have a structure that would capture what would be needed for further computation—or any kind of general understanding.

Meanwhile, in 2009, we released Wolfram|Alpha. One of its features was that it had "step-by-step" math computations. But these weren't "general proofs": rather they were narratives synthesized in very specific ways for human readers. Still, a core concept in Wolfram|Alpha—and the Wolfram Language—is the idea of integrating in knowledge about as many things as possible in the world. We'd done this for cities and movies and lattices and animals and much more. And I thought about doing it for mathematical theorems as well.

We did a pilot project—on theorems about continued fractions. We trawled through the mathematical literature assessing the difficulty of extending the "natural math understanding" we'd built for Wolfram|Alpha. I imagined a workflow which would mix automated theorem generation with theorem search—in which one would define a mathematical scenario, then say "tell me interesting facts about this". And in 2014 we set about engaging the mathematical community in a large-scale curation effort to formalize the theorems of mathematics. But try as we might, only people already involved in math formalization seemed to care; with few exceptions working mathematicians just didn't seem to consider it relevant to what they did.

We continued, however, to push slowly forward. We worked with proof assistant developers. We curated various kinds of mathematical structures (like function spaces). I had estimated that we'd need more than a thousand new Wolfram Language functions to cover "modern pure mathematics", but without a clear market we couldn't motivate the huge design (let alone implementation) effort that would be needed—though, partly in a nod to the intellectual origins of mathematics, we did for example do a project that has succeeded in finally making Euclid-style geometry computable.

Then in the latter part of the 2010s a couple more "proof-related" things happened. Back in 2002 we'd started using equational logic automated theorem proving to get results in functions like FullSimplify. But we hadn't figured out how to present the proofs that were generated. In 2018 we finally introduced FindEquationalProof—allowing programmatic access to proofs, and making it feasible for me to explore collections of proofs in bulk.

I had for decades been interested in what I've called "symbolic discourse language": the extension of the idea of computational language to "everyday discourse"—and to the kind of thing one might want for example to express in legal contracts. And between this and our involvement in the idea of computational contracts, and things like blockchain technology, I started exploring questions of AI ethics and "constitutions". At this point we'd also started to introduce machine-learning-based functions into the Wolfram Language. And—with my "human incomprehensible" Boolean algebra proof as "empirical data"—I started exploring general questions of explainability, and in effect proof.

And not long after that came the surprise breakthrough of our Physics Project. Extending my ideas from the 1990s about computational foundations for fundamental physics it suddenly became possible finally to understand the underlying origins of the main known laws of physics. And core to this effort—and particularly to the understanding of quantum mechanics—were multiway systems.

At first we just used the knowledge that multiway systems could also represent axiomatic mathematics and proofs to provide analogies for our thinking about physics ("quantum observers might in effect be doing critical-pair completions", "causal graphs are like higher categories", etc.) But then we started wondering whether the phenomenon of the emergence that we'd seen for the familiar laws of physics might also affect mathematics—and whether it could give us something like a "bulk" version of metamathematics.

I had long studied the transition from discrete "computational" elements to "bulk" behavior, first following my interest in the Second Law of thermodynamics, which stretched all the way back to age 12 in 1972, then following my work on cellular automaton fluids in the mid-1980s, and now with the emergence of physical space from underlying hypergraphs in our Physics Project. But what might "bulk" metamathematics be like?

One feature of our Physics Project—in fact shared with thermodynamics—is that certain aspects of its observed behavior depend very little on the details of its components. But what did they depend on? We realized that it all had to do with the observer—and their interaction (according to what I've described as the 4th paradigm for science) with the general

"multicomputational" processes going on underneath. For physics we had some idea what characteristics an "observer like us" might have (and actually they seemed to be closely related to our notion of consciousness). But what might a "mathematical observer" be like?

In its original framing we talked about our Physics Project as being about "finding the rule for the universe". But right around the time we launched the project we realized that that wasn't really the right characterization. And we started talking about rulial multiway systems that instead "run every rule"—but in which an observer perceives only some small slice, that in particular can show emergent laws of physics.

But what is this "run every rule" structure? In the end it's something very fundamental: the entangled limit of all possible computations—that I call the ruliad. The ruliad basically depends on nothing: it's unique and its structure is a matter of formal necessity. So in a sense the ruliad "necessarily exists"—and, I argued, so must our universe.

But we can think of the ruliad not only as the foundation for physics, but also as the foundation for mathematics. And so, I concluded, if we believe that the physical universe exists, then we must conclude—a bit like Plato—that mathematics exists too.

But how did all this relate to axiom systems and ideas about metamathematics? I had two additional pieces of input from the latter half of 2020. First, following up on a note in *A New Kind of Science*, I had done an extensive study of the "empirical metamathematics" of the network of the theorems in Euclid, and in a couple of math formalization systems. And second, in celebration of the 100th anniversary of their invention essentially as "primitives for mathematics", I had done an extensive ruliological and other study of combinators.

I began to work on this current piece in the fall of 2020, but felt there was something I was missing. Yes, I could study axiom systems using the formalism of our Physics Project. But was this really getting at the essence of mathematics? I had long assumed that axiom systems really were the "raw material" of mathematics—even though I'd long gotten signals they weren't really a good representation of how serious, aesthetically oriented pure mathematicians thought about things.

In our Physics Project we'd always had as a target to reproduce the known laws of physics. But what should the target be in understanding the foundations of mathematics? It always seemed like it had to revolve around axiom systems and processes of proof. And it felt like validation when it became clear that the same concepts of "substitution rules applied to expressions" seemed to span my earliest efforts to make math computational, the underlying structure of our Physics Project, and "metamodels" of axiom systems.

But somehow the ruliad—and the idea that if physics exists so must math—made me realize that this wasn't ultimately the right level of description. And that axioms were some kind of intermediate level, between the "raw ruliad", and the "humanized" level at which pure mathematics is normally done. At first I found this hard to accept; not only had axiom systems dominated thinking about the foundations of mathematics for more than a century, but they also seemed to fit so perfectly into my personal "symbolic rules" paradigm.

But gradually I got convinced that, yes, I had been wrong all this time—and that axiom systems were in many respects missing the point. The true foundation is the ruliad, and axiom systems are a rather-hard-to-work-with "machine-code-like" description below the inevitable general "physicalized laws of metamathematics" that emerge—and that imply that for observers like us there's a fundamentally higher-level approach to mathematics.

At first I thought this was incompatible with my general computational view of things. But then I realized: "No, quite the opposite!" All these years I've been building the Wolfram Language precisely to connect "at a human level" with computational processes—and with mathematics. Yes, it can represent and deal with axiom systems. But it's never felt particularly natural. And it's because they're at an awkward level—neither at the level of the raw ruliad and raw computation, nor at the level where we as humans define mathematics.

But now, I think, we begin to get some clarity on just what this thing we call mathematics really is. What I've done here is just a beginning. But between its explicit computational examples and its conceptual arguments I feel it's pointing the way to a broad and incredibly fertile new understanding that—even though I didn't see it coming—I'm very excited is now here.

## *Notes & Thanks*

For more than 25 years Elise Cawley has been telling me her thematic (and rather Platonic) view of the foundations of mathematics—and that basing everything on constructed axiom systems is a piece of modernism that misses the point. From what's described here, I now finally realize that, yes, despite my repeated insistence to the contrary, what she's been telling me has been on the right track all along!

I'm grateful for extensive help on this project from James Boyd and Nik Murzin, with additional contributions by Brad Klee and Mano Namuduri. Some of the early core technical ideas here arose from discussions with Jonathan Gorard, with additional input from Xerxes Arsiwalla and Hatem Elshatlawy. (Xerxes and Jonathan have now also been developing connections with homotopy type theory.)

I've had helpful background discussions (some recently and some longer ago) with many people, including Richard Assar, Jeremy Avigad, Andrej Bauer, Kevin Buzzard, Mario Carneiro, Greg Chaitin, Harvey Friedman, Tim Gowers, Tom Hales, Lou Kauffman, Maryanthe Malliaris, Norm Megill, Assaf Peretz, Dana Scott, Matthew Szudzik, Michael Trott and Vladimir Voevodsky.

I'd like to recognize Norm Megill, creator of the Metamath system used for some of the empirical metamathematics here, who died in December 2021. (Shortly before his death he was also working on simplifying the proof of my axiom for Boolean algebra.)

Most of the specific development of this report has been livestreamed or otherwise recorded, and is available—along with archives of working notebooks—at the Wolfram Physics Project website.

The Wolfram Language code to produce all the images here is directly available by clicking each image. And I should add that this project would have been impossible without the Wolfram Language, both its practical manifestation, and the ideas that it has inspired and clarified. So thanks to everyone involved in the 40+ years of its development and gestation!

## Graphical Key

| | |
|---|---|
| state/expression | 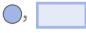 |
| axiom | 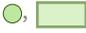 |
| statement/theorem | 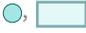 |
| notable theorem | 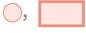 |
| hypothesis | 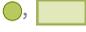 |
| substitution event |  |
| cosubstitution event |  |
| bisubstitution event |  |
| multiway/entailment graph | 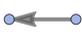 |
| accumulative evolution graph | 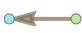 |
| token–event graph | 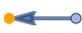 |
| branchial/metamathematical graph | 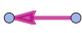 |

## Glossary

*A glossary of terms that are either new, or used in unfamiliar ways*

### accumulative system

A system in which states are rules and rules update rules. Successive steps in the evolution of such a system are collections of rules that can be applied to each other.

### axiomatic level

The traditional foundational way to represent mathematics using axioms, viewed here as being intermediate between the raw ruliad and human-scale mathematics.

### bisubstitution

The combination of substitution and cosubstitution that corresponds to the complete set of possible transformations to make on expressions containing patterns.

*branchial space*

Space corresponding to the limit of a branchial graph that provides a map of common ancestry (or entanglement) in a multiway graph.

*cosubstitution*

The dual operation to substitution, in which a pattern expression that is to be transformed is specialized to allow a given rule to match it.

*eme*

The smallest element of existence according to our framework. In physics it can be identified as an "atom of space", but in general it is an entity whose only internal attribute is that it is distinct from others.

*entailment cone*

The expanding region of a multiway graph or token-event graph affected by a particular node. The entailment cone is the analog in metamathematical space of a light cone in physical space.

*entailment fabric*

A piece of metamathematical space constructed by knitting together many small entailment cones. An entailment fabric is a rough model for what a mathematical observer might effectively perceive.

*entailment graph*

A combination of entailment cones starting from a collection of initial nodes.

*expression rewriting*

The process of rewriting (tree-structured) symbolic expressions according to rules for symbolic patterns. (Called "operator systems" in *A New Kind of Science*. Combinators are a special case.)

*mathematical observer*

An entity sampling the ruliad as a mathematician might effectively do it. Mathematical observers are expected to have certain core human-derived characteristics in common with physical observers.

*metamathematical space*

The space in which mathematical expressions or mathematical statements can be considered to lie. The space can potentially acquire a geometry as a limit of its construction through a branchial graph.

200 | Stephen Wolfram

*multiway graph*

A graph that represents an evolution process in which there are multiple outcomes from a given state at each step. Multiway graphs are central to our Physics Project and to the multicomputational paradigm in general.

*paramathematics*

Parallel analogs of mathematics corresponding to different samplings of the ruliad by putative aliens or others.

*pattern expression*

A symbolic expression that involves pattern variables (x_ etc. in Wolfram Language, or ∀ quantifiers in mathematical logic).

*physicalization of metamathematics*

The concept of treating metamathematical constructs like elements of the physical universe.

*proof cone*

Another term for the entailment cone.

*proof graph*

The subgraph in a token-event graph that leads from axioms to a given statement.

*proof path*

The path in a multiway graph that shows equivalence between expressions, or the subgraph in a token-event graph that shows the constructibility of a given statement.

*ruliad*

The entangled limit of all possible computational processes, that is posited to be the ultimate foundation of both physics and mathematics.

*rulial space*

The limit of rulelike slices taken from a foliation of the ruliad in time. The analog in the rulelike "direction" of branchial space or physical space.

*shredding of observers*

The process by which an observer who has aggregated statements in a localized region of metamathematical space is effectively pulled apart by trying to cover consequences of these statements.

*statement*

A symbolic expression, often containing a two-way rule, and often derivable from axioms, and thus representing a lemma or theorem.

*substitution event*

An update event in which a symbolic expression (which may be a rule) is transformed by substitution according to a given rule.

*token-event graph*

A graph indicating the transformation of expressions or statements ("tokens") through updating events.

*two-way rule*

A transformation rule for pattern expressions that can be applied in both directions (indicated with ↔).

*uniquification*

The process of giving different names to variables generated through different events.

## Annotated Bibliography

The text above includes direct links to specific documents and references. Here we'll give a slightly more general bibliography, though most of it should be considered background, since the approach taken here represents a significant departure from traditional directions, and builds more or less directly on extremely low-level concepts.

*The earliest known large-scale axiomatic presentation of mathematics was:*

Euclid (300 BC), Στοιχεῖα (in Ancient Greek) [*Elements*].

*Empirical metamathematics from this was given in:*

S. Wolfram (2020), "The Empirical Metamathematics of Euclid and Beyond". arXiv: 2107.07337.

*The concept that there is underlying reality in mathematics was discussed in:*

Plato (375 BC), πολιτεία (in Ancient Greek) [*The Republic*].

Plato (360 BC), Τίμαιος (in Ancient Greek) [*Timaeus*].

*Modern explorations of these ideas include:*

M. Balaguer (1998), *Platonism and Anti-Platonism in Mathematics*, Oxford University Press.

J. Gray (2008), *Plato's Ghost: The Modernist Transformation of Mathematics,* Princeton University Press.

R. Tieszen (2011), *After Gödel: Platonism and Rationalism in Mathematics and Logic*, Oxford University Press.

*The contemporary axiomatic formulation of mathematics was developed in:*

F. L. G. Frege (1879), *Begriffsschrift: eine der arithmetischen nachgebildete Formelsprache des reinen Denkens* (in German), Verlag von Louis Neber. (Translated in J. v. Heijenoort (1967), as "Begriffsschrift: A Formal Language, Modeled upon That of Arithmetic, for Pure Thought" in *From Frege to Gödel: A Source Book in Mathematical Logic, 1879–1931*, Harvard University Press, 1–82.)

R. Dedekind (1888), *Was sind und was sollen die Zahlen?* (in German), F. Vieweg und Sohn. (Translated in H. Pogorzelski, et al. (1995), as *What Are Numbers and What Should They Be?*, Research Institute for Mathematics.)

G. Peano (1889), *Arithmetices principia, nova methodo exposita* (in Italian), Fratres Bocca. (Translated by H. C. Kennedy (1973), as "The Principles of Arithmetic, Presented by a New Method", in *Selected Works of Giuseppe Peano*, University of Toronto Press, 101–134.)

D. Hilbert (1903), *Grundlagen der geometrie* (in German), B. G. Teubner. (Translated by E. G. Townsend (1902), as *The Foundations of Geometry,* Open Court.)